\documentclass[a4paper]{amsart}
\usepackage{fullpage}
\usepackage{amsmath,amssymb,mathtools,booktabs,tikz}
\usetikzlibrary{arrows.meta,calc,decorations.markings,positioning}
\usepackage{mathrsfs}
\usepackage[T1]{fontenc}
\usepackage[utf8]{inputenc}
\usepackage{lmodern}
\usepackage{microtype}
\usepackage{needspace}
\usepackage{tabularx,array,longtable}
\usepackage[hypertexnames=false]{hyperref}
\usepackage{doi}
\usepackage{float}
\hypersetup{colorlinks=true,linkcolor=blue,citecolor=blue,urlcolor=blue,
pdftitle={Kushel--Tyaglov's problem, sharp Schoenberg-type inequalities and majorization relations between zeros and critical points},pdfauthor={Teng Zhang}}
\newtheorem{thm}{Theorem}[section]
\newtheorem{lem}[thm]{Lemma}
\newtheorem{prop}[thm]{Proposition}
\newtheorem{cor}[thm]{Corollary}
\newtheorem{rem}[thm]{Remark}
\newtheorem{prob}[thm]{Problem}
\newtheorem{defin}[thm]{Definition}

\newtheorem*{schthm}{Schmeisser's theorem}
\newtheorem*{pmthm}{Pereira--Malamud's theorem}
\newcommand{\C}{\mathbb C}
\newcommand{\R}{\mathbb R}

\newcommand{\dd}{\,\mathrm d}

\DeclareMathOperator{\tr}{Tr}\DeclareMathOperator{\diag}{diag}
\DeclareMathOperator{\rank}{rank}

\renewcommand{\Re}{\operatorname{Re}}\renewcommand{\Im}{\operatorname{Im}}
\numberwithin{equation}{section}
\allowdisplaybreaks[2]
\title{Kushel--Tyaglov's problem, sharp Schoenberg-type inequalities and majorization relations between zeros and critical points}
\author{Teng Zhang}
\address{School of Mathematics and Statistics, Xi'an Jiaotong University,
Xi'an 710049, P. R. China}
\email{teng.zhang@stu.xjtu.edu.cn}
\date{November 23, 2026}
\begin{document}
\begin{abstract}
In this paper, we give a purely complex-analytic proof of
Schoenberg's conjecture. We also determine the sharp constants
in Schoenberg-type inequalities for every real exponent
\(p>0\) and every degree \(n\ge2\), and characterize all
corresponding extremal configurations. In particular, we confirm two conjectures posed in \cite{TZ25} and provide a complete solution to a problem of Kushel and Tyaglov. Furthermore, we introduce the tool called a weighted weak majorization and establish, under the zero-centroid condition, new weak-majorization relations between the moduli of the critical points and those of the zeros. As consequences of these relations, we obtain refinements of several known results.
\end{abstract}
\maketitle

\section{Introduction}\label{sec:introduction}

\subsection{Kushel--Tyaglov's problem and sharp Schoenberg-type inequalities}
Let 
$
p(z) = z^n + a_1 z^{n-1} + \cdots + a_{n-1} z + a_n
$
be a complex polynomial of degree $n \ge 2$, with zeros $z_1, \dots, z_n$ and \emph{critical points} $w_1, \dots, w_{n-1}$, i.e., the zeros of $p'(z)$.

In 1986, Schoenberg \cite{Sch86} conjectured that if the centroid of the zeros lies at the origin, i.e.\ $\sum_{j=1}^n z_j = 0$, then the quadratic inequality
\begin{equation}\label{eq:schoenberg}
	\sum_{k=1}^{n-1} |w_k|^2 \le \frac{n-2}{n} \sum_{j=1}^n |z_j|^2
\end{equation}
holds, with equality if and only if all zeros are collinear in the complex plane. Following the equivalent reformulation in~\cite{DIS99},
Schoenberg's conjecture was independently proved by
Pereira~\cite{Per03}, using the theory of differentiators,
and Malamud~\cite{Mal05}, using inverse spectral theory
and majorization for normal matrices.
Alternative proofs were subsequently given by
Cheung and Ng~\cite{CN06}, using companion matrices,
and Kushel and Tyaglov~\cite{KT16}, using circulant matrices.

Although Schoenberg's conjecture arises naturally in complex analysis, as discussed above, all existing proofs rely on methods from matrix theory. 
This naturally leads to the following question.
\begin{prob}\label{prob:1}
Does Schoenberg's conjecture admit a purely complex-analytic proof?
\end{prob}

In this paper, we give a purely complex-analytic proof of
Schoenberg's conjecture in Section~\ref{sec:complex-schoenberg}, thereby answering
Problem~\ref{prob:1}. The main steps of the proof are as follows.

\begin{itemize}
    \item Using Newton's identities, we reduce Schoenberg's
    inequality to an inequality involving the imaginary parts
    of the critical points.

    \item We introduce a holomorphic family of polynomials
    and construct a subharmonic function on the upper half-plane.

    \item We apply the maximum principle to establish the
    desired inequality. The equality case follows from
    the strong maximum principle.
\end{itemize}

Schoenberg's inequality~\eqref{eq:schoenberg} has since motivated a broad
family of inequalities relating higher moments of
the zeros and critical points, commonly known as
\emph{Schoenberg-type inequalities}.
The fourth-order case was established by Cheung
and Ng~\cite{CN06} and subsequently sharpened
by Kushel and Tyaglov~\cite{KT16}. 
Tang~\cite{Tan25} obtained inequalities
of orders one and six. More recently, Tang and the author \cite{TZ25} confirmed two conjectures of
Georgiev, G\'{o}mez-Serrano, Tao, and Wagner \cite{GGTW26} on fourth-order Schoenberg type inequalities, giving a
complete solution to the
\href{https://google-deepmind.github.io/alphaevolve_repository_of_problems/problems/24.html}{de~Bruin--Sharma problem}, which originated in~\cite{DS99}. 
For further developments on higher-order Schoenberg-type
inequalities, we refer the reader to
\cite{LXZ21,Zha26b}.
Beyond their intrinsic interest, these inequalities
also have applications to the geometry of polynomial
critical points. As an application of \eqref{eq:schoenberg},
the author~\cite{Zha26a,Zha27} obtained a substantial
improvement of a bound due to
Paw{\l}owski~\cite{Paw98} on the distance
from the centroid of the zeros to the
nearest critical point, a problem closely related to Schmeisser's conjecture~\cite{Sch77},
which was subsequently resolved by the author \cite{Zha26d}.

Kushel and Tyaglov~\cite[Section~5]{KT16} emphasized the difficulty of obtaining
higher-order trace expressions in terms of the zeros and raised the following
problem.
  \begin{prob}[\cite{KT16}] \label{prob:kt} Let \(p(z)\) be a complex polynomial of degree \(n\ge2\) with zeros
\(z_1,\ldots,z_n\) and critical points \(w_1,\ldots,w_{n-1}\), counted
with multiplicity. Assume  that \(\sum_{j=1}^n z_j=0\).      Is it possible to estimate $\sum_{k=1}^{n-1}|w_k|^m$ for $m\in\mathbb N$?
   \end{prob} 
More generally, we are interested in the following question.
   \begin{prob}\label{prob:main}
Let \(p(z)\) be a complex polynomial of degree \(n\ge2\), with zeros
\(z_1,\ldots,z_n\) and critical points
\(w_1,\ldots,w_{n-1}\), counted with multiplicity. Assume  that \(\sum_{j=1}^n z_j=0\).   Determine the optimal constants
\(C_{n,p}\)  such that
\[
\sum_{k=1}^{n-1}|w_k|^p
\le C_{n,p}\sum_{j=1}^n|z_j|^p.
\]
What is the optimal constant \(C_{n,p}\) if we drop the zero-centroid condition?
\end{prob}

In this paper, we give a complete solution to
Problem~\ref{prob:main}, thereby answering
Problem~\ref{prob:kt}. The optimal constants in both
the centered and uncentered cases are summarized in
Table~\ref{tab:sharp-constants}.
In the table \(m\) is a positive integer. For \(1\le k\le\lfloor n/2\rfloor\),
the ratio for a centered list
with multiplicities \(k\) and \(n-k\) is
\[
  R_{n,p}(k)=
  \frac{(k-1)(n-k)^p+(n-k-1)k^p+(n-2k)^p}
       {k(n-k)^p+(n-k)k^p}.
\]

\begin{table}[htbp]
\centering
\small
\renewcommand{\arraystretch}{1.35}
\begin{tabularx}{\textwidth}{
|>{\centering\arraybackslash}p{5.3cm}
|>{\centering\arraybackslash}p{3.7cm}
|>{\centering\arraybackslash}p{1.5cm}
|>{\centering\arraybackslash}X|}
\hline
$C(n,p)$ & \(n,p\) & Centroid & Reference\\
\hline

\(\displaystyle
\left(\frac{n-2}{n}\right)^{\frac{p}{2}}
\)
&
\(n\ge2,\quad 0<p\le1\)
&
\(0\)
&
Theorem~\ref{thm:centered-subunit}
\\
\hline

\(\displaystyle
\max\left\{
\left(\frac{n-2}{n}\right)^{\frac{p}{2}},
\max_k R_{n,p}(k)
\right\}
\)
&
\(n\ge2,\quad 1<p<2\)
&
\(0\)
&
Theorem~\ref{thm:centered-subquadratic}
\\
\hline

\(\displaystyle\frac{n-2}{n}\)
&
\(n\ge2,\quad p=2\)
&
\(0\)
&
Pereira~\cite{Per03},
Malamud~\cite{Mal05}
\\
\hline

\(\displaystyle\frac{n-2}{n}\)
&
\(n=2m,\quad p>2\)
&
\(0\)
&
Theorem~\ref{thm:centered-superquadratic}
\\
\hline

\(\displaystyle R_{3,p}(1)\)
&
\(n=3,\quad 2<p\le4\)
&
\(0\)
&
Proposition~\ref{prop:cubic-sharp}
\\
\hline

\(\displaystyle 3^{-\frac{p}{2}}\)
&
\(n=3,\quad p>4\)
&
\(0\)
&
Proposition~\ref{prop:cubic-sharp}
\\
\hline

\(\displaystyle R_{n,p}\bigl((n-1)/2\bigr)\)
&
\(n=2m+3,\quad p>2\)
&
\(0\)
&
Theorem~\ref{thm:centered-odd}
\\
\hline

\(\displaystyle
\left(\frac{n-1}{n}\right)^p
\)
&
\(n\ge2,\quad 0<p<1\)
& arbitrary
&
Theorem~\ref{thm:general-subunit}
\\
\hline

\(\displaystyle\frac{n-1}{n}\)
&
\(n\ge2,\quad p\ge1\)
& arbitrary
&
de~Bruijn--Springer~\cite{DS47}
\\
\hline

\end{tabularx}
\caption{Optimal constants in power-sum inequalities for the zeros
and critical points of complex polynomials.}
\label{tab:sharp-constants}
\end{table}
    \begin{thm}\label{thm:centered-subunit}
Let \(p(z)\) be a complex polynomial of degree \(n\ge2\) with zeros
\(z_1,\ldots,z_n\) and critical points \(w_1,\ldots,w_{n-1}\), counted
with multiplicity. Assume that \(\sum_{j=1}^n z_j=0\) and \(0<p\le1\). Then
	\[
	\sum_{k=1}^{n-1} |w_k|^p \le \left(\frac{n-2}{n}\right)^{\frac{p}{2}} \sum_{j=1}^n |z_j|^p .
    \]
Equality holds exactly for \(p(z)=c z^{n-2}(z^2-a^2)\),
where \(c\in\C\setminus\{0\}\) and \(a\in\C\).
    \end{thm}
For \(n\ge4\), Tang and Zhang~\cite[Theorem~4.7]{TZ25}
obtained a bound for \(1< p<2\) by complex interpolation.
The next theorem gives the sharp constant for \(1<p<2\).
\begin{thm}\label{thm:centered-subquadratic}
Let \(p(z)\) be a complex polynomial of degree \(n\ge2\) with zeros
\(z_1,\ldots,z_n\) and critical points \(w_1,\ldots,w_{n-1}\), counted
with multiplicity. Assume that \(\sum_{j=1}^n z_j=0\) and \(1<p<2\). Then
\begin{equation*}
\sum_{k=1}^{n-1}|w_k|^p
\le
\max\Bigg\{
\left(\frac{n-2}{n}\right)^{\frac{p}{2}},\max_{1\le k\le\lfloor n/2\rfloor}
\frac{
(k-1)(n-k)^p+(n-k-1)k^p+(n-2k)^p
}{
k(n-k)^p+(n-k)k^p
}
\Bigg\}
\sum_{j=1}^n|z_j|^p.
\end{equation*}
The constant is optimal. Apart from \(p(z)=c z^n\), equality holds
for precisely the following polynomials, with \(a,c\in\C\setminus\{0\}\):
\(c z^{n-2}(z^2-a^2)\) if the first term realizes the maximum,
and \(c(z-(n-k)a)^k(z+ka)^{n-k}\) for each maximizing index \(k\)
in the second term. All maximizing families are included in a tie.
\end{thm}

Our Theorems~\ref{thm:centered-superquadratic}
and~\ref{thm:centered-odd} confirm two conjectures in
\cite[Conjectures~1.6 and~1.7]{TZ25}, respectively. 
\begin{thm}\label{thm:centered-superquadratic}
Let \(p(z)\) be a complex polynomial of degree \(n\ge2\) with zeros
\(z_1,\ldots,z_n\) and critical points \(w_1,\ldots,w_{n-1}\), counted
with multiplicity. Assume that \(\sum_{j=1}^n z_j=0\) and \(p\ge2\). Then
	\[
	\sum_{k=1}^{n-1} |w_k|^p \le \frac{n-2}{n} \sum_{j=1}^n |z_j|^p .
    \]
For \(p=2\), equality holds precisely when the zeros are collinear.
For \(p>2\) and a nonzero root list, equality holds precisely when
\(n\) is even and \(p(z)=c(z^2-a^2)^{\frac{n}{2}}\) for nonzero \(a,c\in\C\).
\end{thm}

\begin{thm}\label{thm:centered-odd}
Let \(p(z)\) be a complex polynomial of odd degree \(n\ge5\) with zeros
\(z_1,\ldots,z_n\) and critical points \(w_1,\ldots,w_{n-1}\), counted
with multiplicity. Assume that \(\sum_{j=1}^n z_j=0\) and \(p\ge2\). Then
	\[
	\sum_{k=1}^{n-1} |w_k|^p \le \frac{(n-3)(n+1)^{p}+(n-1)^{p+1}+2^{p+1}}
	{(n-1)(n+1)^{p}+(n+1)(n-1)^{p}} \sum_{j=1}^n |z_j|^p .
    \]
For \(p=2\), equality holds precisely when the zeros are collinear.
For \(p>2\), all equality cases with a nonzero root list are
\(p(z)=c(z-(n+1)a)^{\frac{n-1}{2}}(z+(n-1)a)^{\frac{n+1}{2}}\),
where \(a,c\in\C\setminus\{0\}\).
\end{thm}

\begin{rem}
In \cite{TZ25}, the authors established the \(n=4\) cases of
Theorems~\ref{thm:centered-superquadratic} and
\ref{thm:centered-odd}. They also confirmed two conjectures
of Georgiev, G\'omez-Serrano, Tao, and Wagner~\cite{GGTW26}
concerning fourth-order Schoenberg-type inequalities.
These results resolved the de Bruin--Sharma problem.
\end{rem}

For \(p\ge1\), de Bruijn and Springer~\cite[Theorem~4]{DS47}
proved \eqref{eq:Pereira1} without any assumption on the centroid.
The same inequality follows from Pereira's convex-function
inequality~\cite[Corollary~5.5]{Per03}. An alternative proof based
on complex interpolation is given in~\cite{TZ25}.
\begin{equation}\label{eq:Pereira1}
  \sum_{k=1}^{n-1}|w_k|^p
  \le
  \frac{n-1}{n}\sum_{j=1}^n|z_j|^p.
\end{equation}
 For \(p=1\), equality holds precisely
when all zeros lie on one closed ray from the origin. For \(p>1\),
it holds precisely when all zeros coincide. Taking \(p(z)=(z-a)^n\)
with \(a\ne0\) gives the ratio \((n-1)/n\) for every \(p\ge1\).

The following theorem extends
\eqref{eq:Pereira1} to the range $0<p<1$,
with a sharp constant.
\begin{thm}\label{thm:general-subunit}
Let \(p(z)\) be a complex polynomial of degree \(n\ge2\), with zeros
\(z_1,\ldots,z_n\) and critical points \(w_1,\ldots,w_{n-1}\),
counted with multiplicity. If \(0<p<1\), then
\[
  \sum_{k=1}^{n-1} |w_k|^p \le \left(\frac{n-1}{n}\right)^p \sum_{j=1}^n |z_j|^p.    
\]
For each fixed pair \(n,p\), the constant is sharp.
If $z_1,\ldots,z_n$ are not all zero, equality holds
if and only if $p(z)$ is a nonzero scalar multiple of
$z^{n-1}(z-a)$ for some $a\in\mathbb C\setminus\{0\}$.
\end{thm}

\subsection{Majorization relations between zeros and critical points}\label{sub:major}
For \(x=(x_1,\ldots,x_n)\in\mathbb R^n\), let
\(x^\downarrow=(x_1^\downarrow,\ldots,x_n^\downarrow)\)
denote the non-increasing rearrangement of its components.

Given \(x,y\in\mathbb R^n\), we say that \(x\) is
\emph{weakly majorized} by \(y\), denoted by \(x\prec_w y\), if
\[
\sum_{j=1}^k x_j^\downarrow
\le
\sum_{j=1}^k y_j^\downarrow,
\qquad 1\le k\le n.
\]
If, in addition,
$
\sum_{j=1}^n x_j=\sum_{j=1}^n y_j,
$
then \(x\) is \emph{majorized} by \(y\), written \(x\prec y\).

Similarly, for \(x,y\in \mathbb R_+^n\), we say that
\(x\) is \emph{weakly log-majorized} by \(y\),
denoted by \(x\prec_{w\log}y\), if
\[
\prod_{j=1}^k x_j^\downarrow
\le
\prod_{j=1}^k y_j^\downarrow,
\qquad 1\le k\le n.
\]
If, in addition,
$
\prod_{j=1}^n x_j=\prod_{j=1}^n y_j,
$
then \(x\) is \emph{log-majorized} by \(y\),
written \(x\prec_{\log}y\). 
In particular, weak log-majorization implies weak majorization; that is,
\begin{equation}\label{eq:log-majorization}
x\prec_{w\log}y \quad\Longrightarrow\quad x\prec_w y.
\end{equation}
For a comprehensive account of majorization theory,
we refer to Marshall--Olkin--Arnold
\cite{MOA11}.

Schmeisser~\cite[Corollary~4]{Sch04}
established the following majorization relations
between the moduli of the zeros and critical points of a polynomial.
\begin{schthm}
   Let \(p(z)\) be a complex polynomial of degree \(n\ge2\) with zeros
\(z_1,\ldots,z_n\) and critical points \(w_1,\ldots,w_{n-1}\), counted
with multiplicity. 
Then
\[
(|w_1|,\ldots,|w_{n-1}|,0)
\prec_{w\log}
(|z_1|,\ldots,|z_n|).
\]
\end{schthm}

By using \eqref{eq:log-majorization}, Schmeisser's theorem implies that
\((|w_1|,\ldots,|w_{n-1}|,0) \prec_w (|z_1|,\ldots,|z_n|)\).

A stronger geometric majorization relation was obtained independently by
Malamud~\cite[Theorem~4.7]{Mal05} and
Pereira~\cite[Theorem~5.4]{Per03}.

\begin{pmthm}
    Let \(p(z)\) be a complex polynomial of degree \(n\ge2\) with zeros
\(z_1,\ldots,z_n\) and critical points \(w_1,\ldots,w_{n-1}\), counted
with multiplicity. 
Then
\[
        (w_1,\ldots,w_{n-1}, n^{-1}\sum_{j=1}^n z_j)
        \prec
        (z_1,\ldots,z_n),
\]
where the majorization is understood in the complex-vector sense, equivalently
through the existence of a doubly stochastic matrix.
\end{pmthm}
By Pereira's convexity corollary
\cite[Corollary~5.5]{Per03}, every convex function
\(\varphi:\C\to\R\) satisfies
$$\sum_{k=1}^{n-1}\varphi(w_k) \le \frac{n-1}{n}\sum_{j=1}^{n}\varphi(z_j).$$
Taking \(\varphi(\zeta)=|\zeta|^p\),
\(\varphi(\zeta)=|\Re\zeta|^p\), and
\(\varphi(\zeta)=|\Im\zeta|^p\), respectively, gives, for \(p\ge1\),
\begin{equation}\label{eq:unrestricted-projections}
\sum_{k=1}^{n-1}|w_k|^p
\le \frac{n-1}{n}\sum_{j=1}^{n}|z_j|^p,\quad
\sum_{k=1}^{n-1}|\Re w_k|^p
\le \frac{n-1}{n}\sum_{j=1}^{n}|\Re z_j|^p,\quad
\sum_{k=1}^{n-1}|\Im w_k|^p
\le \frac{n-1}{n}\sum_{j=1}^{n}|\Im z_j|^p.
\end{equation}

As discussed above, a central case in the study of
Schoenberg-type inequalities is when the zeros have
centroid at the origin. Under this additional assumption,
the constants in the preceding estimates can be improved.
To this end, we develop a weighted weak-majorization
framework that yields sharper constants for polynomials
with zero centroid. In particular, when
\(\sum_{j=1}^n z_j=0\), we improve the coefficients
in~\eqref{eq:unrestricted-projections}.

\begin{defin}[Weighted weak majorization]
\label{def:weighted-majorization}
Let \(x=(x_1,\ldots,x_m)\in\R_+^m\),
\(y=(y_1,\ldots,y_n)\in\R_+^n\), and
\(\gamma\ge0\). We say that \(x\) is
\emph{weakly majorized by \(y\) with weight \(\gamma\)},
denoted by \(x\prec_{w,\gamma}y\), if
\begin{equation}\label{eq:wmaj-def}
\sum_{i=1}^{m}(x_i-t)_+
\le
\gamma\sum_{j=1}^{n}(y_j-t)_+,
\qquad t\ge0,
\end{equation}
where \((s)_+=\max\{s,0\}\).
\end{defin}

For \(x,y\in\R_+^n\), the usual weak majorization
\(x\prec_w y\) is equivalent to
\[
\sum_{i=1}^{n}(x_i-t)_+
\le
\sum_{j=1}^{n}(y_j-t)_+,
\qquad t\ge0.
\]
Thus, when \(\gamma=1\) and \(m=n\),
Definition~\ref{def:weighted-majorization} reduces to
the usual weak majorization. In general,
\(x\prec_{w,\gamma}y\) is not equivalent to
\(x\prec_w\gamma y\).

The following lemma characterizes weighted weak majorization
in terms of increasing convex functions; its proof is given
in Section~\ref{sec:matrix-inequalities}.

\begin{lem}\label{lem:weighted-convex}
Let \(x\in\R_+^m\), \(y\in\R_+^n\), and \(\gamma\ge0\).
Then \(x\prec_{w,\gamma}y\) if and only if, for every
increasing convex function
\(\Phi:\R_+\to\R\) satisfying \(\Phi(0)=0\),
\[
\sum_{i=1}^{m}\Phi(x_i)
\le
\gamma\sum_{j=1}^{n}\Phi(y_j).
\]
\end{lem}

We first establish a weighted weak-majorization relation
between the moduli of the critical points and those of the zeros.

\begin{thm}\label{thm:main-modulus-wmaj}
 Let \(p(z)\) be a complex polynomial of degree \(n\ge2\) with zeros
\(z_1,\ldots,z_n\) and critical points \(w_1,\ldots,w_{n-1}\), counted
with multiplicity. 
Assume that \(\sum_{j=1}^n z_j=0\). Then
\[
(|w_1|,\ldots,|w_{n-1}|)
\prec_{w,\sqrt{\frac{n-2}{n}}}
(|z_1|,\ldots,|z_n|).
\]
\end{thm}
\begin{rem}
Under the hypotheses of Theorem~\ref{thm:main-modulus-wmaj}, Lemma~\ref{lem:weighted-convex} gives
\begin{equation}\label{eq:centered-radial-convex}
\sum_{k=1}^{n-1}\Phi(|w_k|)
\le
\sqrt{\frac{n-2}{n}}
\sum_{j=1}^{n}\Phi(|z_j|)
\end{equation}
for every increasing convex function
\(\Phi:\R_+\to\R\) with \(\Phi(0)=0\).
For \(n=2\) both sides vanish, since the sole critical point is zero.

Taking \(\Phi(t)=t\) in \eqref{eq:centered-radial-convex} recovers
Tang's order-one inequality~\cite[Theorem~3.4]{Tan25}.
Tang~\cite[Remark~3.5]{Tan25} asks about a convex-function
extension of his order-one inequality. Theorem~\ref{thm:main-modulus-wmaj}
gives such an extension for radial convex functions vanishing at the
origin. The normalization matters, since the constant function \(1\)
would require \(n-1\le n\sqrt{(n-2)/n}\).
For the normalized radial functions, the centered coefficient is
smaller than \((n-1)/n\) in Pereira's
inequality~\cite[Corollary~5.5]{Per03}.
\end{rem}

Theorem~\ref{thm:main-modulus-wmaj} compares the moduli.
For their squares the weight improves to \((n-2)/n\).

\begin{thm}\label{thm:main-square-wmaj}
Let \(p(z)\) be a complex polynomial of degree \(n\ge2\) with zeros
\(z_1,\ldots,z_n\) and critical points \(w_1,\ldots,w_{n-1}\), counted
with multiplicity. Assume that \(\sum_{j=1}^n z_j=0\). Then
\[
(|w_1|^2,\ldots,|w_{n-1}|^2)
\prec_{w,\frac{n-2}{n}}
(|z_1|^2,\ldots,|z_n|^2).
\]
\end{thm}
\begin{rem}
Under the hypotheses of Theorem~\ref{thm:main-square-wmaj}, Lemma~\ref{lem:weighted-convex} gives
\(\sum_{k=1}^{n-1}\Phi(|w_k|^2) \le \frac{n-2}{n} \sum_{j=1}^{n}\Phi(|z_j|^2)\)
for every increasing convex function
\(\Phi:\R_+\to\R\) satisfying \(\Phi(0)=0\).
Thus Pereira's inequality~\cite[Corollary~5.5]{Per03} has a smaller
coefficient for functions of the form \(\zeta\mapsto\Phi(|\zeta|^2)\)
under the zero-centroid assumption. For \(p\ge2\), taking
\(\Phi(t)=t^{\frac{p}{2}}\) gives the inequality in
Theorem~\ref{thm:centered-superquadratic}.
\end{rem}

We next consider projections onto lines through the origin. For
\(\theta\in\R\) and \(\zeta\in\C\), write \(L_\theta(\zeta)=\Re(e^{-i\theta}\zeta)\).
The next theorem compares the absolute values of the projected zeros and critical points.

\begin{thm}\label{thm:main-directional-first}
 Let \(p(z)\) be a complex polynomial of degree \(n\ge2\) with zeros
\(z_1,\ldots,z_n\) and critical points \(w_1,\ldots,w_{n-1}\), counted
with multiplicity. 
Assume that \(\sum_{j=1}^n z_j=0\).
Then, for every \(\theta\in\R\),
\[
\bigl(|L_\theta(w_1)|,\ldots,|L_\theta(w_{n-1})|\bigr)
\prec_{w,\sqrt{\frac{n-2}{n}}}
\bigl(|L_\theta(z_1)|,\ldots,|L_\theta(z_n)|\bigr).
\]
\end{thm}
\begin{rem}
Fix \(\theta\in\R\). Theorem~\ref{thm:main-directional-first}
and Lemma~\ref{lem:weighted-convex} give
\begin{equation}\label{eq:centered-directional-convex}
\sum_{k=1}^{n-1}\Phi(|L_\theta(w_k)|)
\le
\sqrt{\frac{n-2}{n}}
\sum_{j=1}^{n}\Phi(|L_\theta(z_j)|)
\end{equation}
for every increasing convex function
\(\Phi:\R_+\to\R\) satisfying \(\Phi(0)=0\).
Since \(\zeta\mapsto\Phi(|L_\theta(\zeta)|)\) is convex on \(\C\),
\eqref{eq:centered-directional-convex} improves Pereira's coefficient
for functions of one projection~\cite[Corollary~5.5]{Per03}.
Together, \eqref{eq:centered-radial-convex} and \eqref{eq:centered-directional-convex}
extend Tang's order-one bound~\cite[Theorem~3.4]{Tan25} to radial functions
and functions of a single projection.
\end{rem}

As in the modulus comparison, a stronger estimate
is available for the squares of the projected
absolute values.

\begin{thm}\label{thm:main-directional-square}
 Let \(p(z)\) be a complex polynomial of degree \(n\ge2\) with zeros
\(z_1,\ldots,z_n\) and critical points \(w_1,\ldots,w_{n-1}\), counted
with multiplicity. 
Assume that \(\sum_{j=1}^n z_j=0\).
Then, for every \(\theta\in\R\),
\[
\bigl(|L_\theta(w_1)|^2,\ldots,|L_\theta(w_{n-1})|^2\bigr)
\prec_{w,\frac{n-2}{n}}
\bigl(|L_\theta(z_1)|^2,\ldots,|L_\theta(z_n)|^2\bigr).
\]
\end{thm}

\begin{rem}
For every \(\theta\in\R\) and every increasing convex
function \(\Phi:\R_+\to\R\) satisfying \(\Phi(0)=0\),
Theorem~\ref{thm:main-directional-square} yields
\(\sum_{k=1}^{n-1}\Phi(|L_\theta(w_k)|^2) \le \frac{n-2}{n} \sum_{j=1}^{n}\Phi(|L_\theta(z_j)|^2)\).
For the convex functions \(\zeta\mapsto\Phi(|L_\theta(\zeta)|^2)\),
the coefficient \((n-2)/n\) is smaller than \((n-1)/n\) in
Pereira's inequality~\cite[Corollary~5.5]{Per03}.
\end{rem}

\subsection{Organization of this paper}\label{subsec:organization}
Section~\ref{sec:complex-schoenberg} answers Problem~\ref{prob:1}
by applying the maximum principle to prove Schoenberg's conjecture~\eqref{eq:schoenberg}.
Section~\ref{sec:subunit-proof} proves Theorem~\ref{thm:centered-subunit}
for centered zeros and \(0<p\le1\), and Theorem~\ref{thm:general-subunit}
for zeros with arbitrary centroid and \(0<p<1\).
For \(1<p<2\), Section~\ref{sec:subquadratic-proof} proves
Theorem~\ref{thm:centered-subquadratic} by reducing to polynomials
with real zeros and applying the moment estimates of Appendix~\ref{app:sub12-details}.

Section~\ref{sec:matrix-inequalities} proves Lemma~\ref{lem:weighted-convex}
and establishes singular-value bounds for compressions of diagonal matrices.
Section~\ref{sec:modulus-comparisons} proves Theorem~\ref{thm:main-modulus-wmaj}
on the moduli of the zeros and critical points.
Section~\ref{sec:square-directional-proofs} proves
Theorems~\ref{thm:centered-superquadratic}, \ref{thm:main-square-wmaj},
\ref{thm:main-directional-first}, and~\ref{thm:main-directional-square}
on power sums, squared moduli, and projections onto lines through the origin.
Proposition~\ref{prop:cubic-sharp} gives the sharp constant for centered cubics
when \(p>2\).

Section~\ref{sec:odd-all-exponents} proves Theorem~\ref{thm:centered-odd}
for odd degrees \(n\ge5\). For \(2<p<3\), the proof uses the logarithmic
comparison in Appendix~\ref{app:odd23-details}. For \(p\ge3\),
Appendix~\ref{app:odd-ge-three-details} reduces the real extremal problem
to polynomials with at most three distinct zeros and proves the cubic tail
inequalities. Appendix~\ref{app:odd-exact-certificates} gives the block
allocations and the exact-arithmetic certificates. The accompanying
\href{https://github.com/zhangteng2000/schoenberg-validation/tree/e7435c5521dcd6a10131ba04a6da500b1404800b/certificates}{Python scripts}
are available on GitHub~\cite{Zha26c}.

\subsection{Acknowledgments and AI tools disclosure} 

This work was supported by the China Scholarship Council,
by the Young Elite Scientists Sponsorship Program for PhD Students
of the China Association for Science and Technology,
and by the Fundamental Research Funds for the Central Universities
at Xi'an Jiaotong University (Grant No.~xzy022024045).

While seeking a proof of Theorem~\ref{thm:centered-superquadratic},
the author developed the weighted weak-majorization framework
used in this paper and established the modulus and projection
comparisons presented in Subsection~\ref{sub:major};
see the \href{https://zhangteng2000.github.io/files/weighted_weak_majorization.pdf}{author's earlier note}
on weighted weak majorization.
Before consulting ChatGPT~6.0 (OpenAI) on this project,
the author had also independently obtained proofs of
Theorems~\ref{thm:centered-subunit} and
\ref{thm:general-subunit}.

The author had conjectured the bound established in
Theorem~\ref{thm:centered-subquadratic} and, jointly with
Tang, proposed the odd-degree bound established in
Theorem~\ref{thm:centered-odd}.
By applying Pereira's projection theorem~\cite[Theorem~4.6]{Per03},
the author reduced both conjectures to inequalities for
polynomials with real zeros. However, the author had not
found proofs of either of these real-rooted inequalities.
ChatGPT~6.0 suggested the approaches used in this paper
to establish both inequalities. It also suggested the
complex-analytic proof of Schoenberg's conjecture presented
in Section~\ref{sec:complex-schoenberg}, which is based on
the maximum principle for subharmonic functions.

ChatGPT 6.0 was used for English-language editing, proofreading, grammatical corrections and exploratory discussions of possible approaches to selected parts of this paper.  The mathematical arguments and proofs in the original manuscript were independently developed, checked, and written by the author.

\section{A complex-analytic proof of Schoenberg's conjecture}
\label{sec:complex-schoenberg}

We prove \eqref{eq:schoenberg} by applying the maximum
principle to a subharmonic function on the upper half-plane.
The argument relies on the following identity.

\begin{lem}\label{lem:complex-laplacian}
Let $f$ be holomorphic in a domain contained in the upper
half-plane. Then
\begin{equation}\label{eq:complex-laplacian}
\Delta_t\left(\frac{(\Im f(t))^2}{\Im t}\right)
=
\frac{2}{\Im t}
\left|f'(t)-\frac{\Im f(t)}{\Im t}\right|^2,
\end{equation}
where $\Delta_t$ denotes the planar Laplacian.
\end{lem}

\begin{proof}
Write $t=u+iv$ with $v>0$, and set
$h(u,v)=\Im f(u+iv)$. Since $h$ is harmonic,
direct differentiation gives
\(\Delta_t\left(\frac{h^2}{v}\right) = \frac{2}{v} \left(h_u^2+\left(h_v-\frac{h}{v}\right)^2\right)\).
The Cauchy--Riemann equations give
$h_u=\Im f'(t)$ and $h_v=\Re f'(t)$,
proving \eqref{eq:complex-laplacian}.
\end{proof}

\begin{proof}[Proof of Schoenberg's conjecture
\eqref{eq:schoenberg}]
If all zeros coincide, the centroid condition gives
$p(z)=z^n$, and equality holds. Otherwise, a rotation
makes the imaginary parts of the distinct zeros
pairwise different. Such a rotation preserves both
sides of \eqref{eq:schoenberg} and the equality
condition.

\medskip
\noindent\emph{Step 1: Reduction to a projection inequality.}
Let $\eta_1,\ldots,\eta_{n-1}$ be the critical
points of the real polynomial
$\prod_{j=1}^{n}(x-\Im z_j).$
The first two Newton identities
\cite[p.~5, Exercise~E.2]{BE95}, applied to $p(z)$,
$p'(z)/n$, and the projected polynomial, give
\begin{equation}\label{eq:complex-newton}
\sum_{k=1}^{n-1}w_k^2
=
\frac{n-2}{n}\sum_{j=1}^{n}z_j^2,
\qquad
\sum_{k=1}^{n-1}\eta_k^2
=
\frac{n-2}{n}\sum_{j=1}^{n}(\Im z_j)^2.
\end{equation}
Since
$
|z|^2=\Re(z^2)+2(\Im z)^2,
$
equation \eqref{eq:complex-newton} implies
\begin{equation}\label{eq:complex-deficit}
\frac{n-2}{n}\sum_{j=1}^{n}|z_j|^2
-\sum_{k=1}^{n-1}|w_k|^2=
2\left(
\sum_{k=1}^{n-1}\eta_k^2
-\sum_{k=1}^{n-1}(\Im w_k)^2
\right).
\end{equation}
Thus \eqref{eq:schoenberg} is equivalent to
\begin{equation}\label{eq:complex-projection}
\sum_{k=1}^{n-1}(\Im w_k)^2
\le
\sum_{k=1}^{n-1}\eta_k^2.
\end{equation}

\medskip
\noindent\emph{Step 2: A holomorphic family and its asymptotic behavior.}
For $t\in\mathbb C$, consider the polynomial
\begin{equation}\label{eq:complex-family}
p_t(w)=
\prod_{j=1}^{n}
\bigl(w-\Re z_j-t\Im z_j\bigr),
\end{equation}
and denote its critical points by
$w_1(t),\ldots,w_{n-1}(t)$.
At \(t=i\), \eqref{eq:complex-family} gives \(p_i(z)=p(z)\); for real \(t\)
all zeros and critical points of $p_t$ are real.

For $\Im t>0$, the Gauss--Lucas theorem
\cite[p.~18, Theorem~1.3.1]{BE95} gives
\begin{equation}\label{eq:complex-strip}
\min_{1\le j\le n}\Im z_j
\le
\frac{\Im w_k(t)}{\Im t}
\le
\max_{1\le j\le n}\Im z_j,
\qquad 1\le k\le n-1,
\end{equation}
because each critical point belongs to the convex
hull of the zeros of $p_t$.

We next determine the behavior of the critical
points as $|t|\to\infty$.
Let $a_1,\ldots,a_s$ be the distinct zeros of $p(z)$,
with multiplicities $m_1,\ldots,m_s$, indexed so
that
$
\Im a_1<\cdots<\Im a_s.
$
For every \(t\), the derivative of \(p_t\) factors as
\begin{equation}\label{eq:complex-free-factor}
 p_t'(w)=\prod_{j=1}^s(w-\Re a_j-t\Im a_j)^{m_j-1}
 \sum_{j=1}^s m_j\prod_{\ell\ne j}(w-\Re a_\ell-t\Im a_\ell).
\end{equation}
Two distinct values \(\Re a_j+t\Im a_j\) can coincide only for real \(t\),
because the imaginary parts \(\Im a_j\) are pairwise distinct. Away from
these finitely many parameters, the first factor supplies the inherited
critical points, with multiplicities \(m_j-1\). The other \(s-1\) critical
points are the solutions of
\begin{equation}\label{eq:complex-free-critical}
\sum_{j=1}^{s}
\frac{m_j}
{w-\Re a_j-t\Im a_j}=0.
\end{equation}

Setting $w=t\xi$ and letting $|t|\to\infty$
gives the limiting equation
\begin{equation}\label{eq:complex-limit-critical}
\sum_{j=1}^{s}
\frac{m_j}{\xi-\Im a_j}=0.
\end{equation}
The left-hand side of
\eqref{eq:complex-limit-critical} is strictly
decreasing on each interval
$(\Im a_j,\Im a_{j+1})$, with limits
$+\infty$ and $-\infty$ at the respective
endpoints. It therefore has exactly $s-1$
simple real zeros, one in each interval.
Together with \(m_j-1\) copies of each \(\Im a_j\), the simple zeros
of \eqref{eq:complex-limit-critical} form the list \(\eta_1,\ldots,\eta_{n-1}\).

The holomorphic implicit function theorem
\cite[p.~34, Theorem~7.6]{FG02}, applied with
parameter $1/t$ near each simple root of
\eqref{eq:complex-limit-critical}, shows that
the moving critical points can be labelled
for sufficiently large $|t|$ so that
\begin{equation}\label{eq:complex-series}
w_k(t)
=
\eta_k t+b_k+
\sum_{\ell=1}^{\infty}c_{k,\ell}t^{-\ell},
\qquad
b_k,c_{k,\ell}\in\mathbb R.
\end{equation}
The coefficients are real because the equations
defining the branches have real coefficients
for real $t$. The series in
\eqref{eq:complex-series} converge for all
sufficiently large $|t|$; for the inherited
branches the coefficients $c_{k,\ell}$ vanish.

For every integer $\ell\ge1$ and $\Im t>0$,
\(|\Im(t^{-\ell})| \le \frac{\ell\,\Im t}{|t|^{\ell+1}}\).
Taking imaginary parts in
\eqref{eq:complex-series} therefore yields
\begin{equation}\label{eq:complex-asymptotics}
\frac{\Im w_k(t)}{\Im t}
=
\eta_k+O(|t|^{-2}),
\qquad
|t|\to\infty,\quad\Im t>0.
\end{equation}
The error in \eqref{eq:complex-asymptotics}
is uniform even when $t$ approaches the real
axis.

The sum in \eqref{eq:complex-free-factor} has degree \(s-1\) in \(w\),
leading coefficient \(n\), and polynomial coefficients in \(t\).
Its discriminant is therefore a polynomial in \(t\).
It is nonzero for sufficiently large $|t|$,
by the simple-root property of
\eqref{eq:complex-limit-critical}.
Hence the free critical points can collide
only at finitely many parameter values.

\medskip
\noindent\emph{Step 3: Application of the maximum principle.}
For $\Im t>0$, define
\begin{equation}\label{eq:complex-subharmonic}
F(t)
=
\sum_{k=1}^{n-1}
\frac{(\Im w_k(t))^2}{\Im t}
-(\Im t)\sum_{k=1}^{n-1}\eta_k^2.
\end{equation}
Away from the collision parameters, the
critical points admit local holomorphic
branches, so Lemma~\ref{lem:complex-laplacian}
gives
\begin{equation}\label{eq:complex-positive-laplacian}
\Delta_tF(t)
=
\frac{2}{\Im t}
\sum_{k=1}^{n-1}
\left|
w_k'(t)-\frac{\Im w_k(t)}{\Im t}
\right|^2
\ge0.
\end{equation}

The multiset of critical points depends
continuously on $t$, and hence $F$ is
continuous at every collision parameter.
To check subharmonicity at a collision \(t_0\), take a disk
\(|t-t_0|\le r\) containing no other collision and contained in the
upper half-plane. Let \(h\) be harmonic in the disk and continuous on
its closure, with \(F\le h\) on its boundary. For \(\varepsilon>0\),
\(F-h+\varepsilon\log(|t-t_0|/r)\) is nonpositive on the outer circle
and on every sufficiently small inner circle, since \(F-h\) is bounded.
The comparison principle~\cite[p.~265, Corollary~3.5]{Con78} gives the
same bound on the intervening annulus. Shrinking the inner circle and
then letting \(\varepsilon\downarrow0\) gives \(F\le h\) throughout
the disk by continuity. Thus \(F\) is subharmonic at \(t_0\).

By \eqref{eq:complex-strip}, the function $F$
extends continuously to the real axis with
boundary value zero. In addition,
\eqref{eq:complex-asymptotics} and
\eqref{eq:complex-subharmonic} give
\begin{equation}\label{eq:complex-decay}
F(t)
=
O\left(\frac{\Im t}{|t|^2}\right),
\qquad
|t|\to\infty,\quad\Im t>0.
\end{equation}

Apply the maximum principle for subharmonic
functions
\cite[p.~264, Theorem~3.3]{Con78}
to the upper half-disk
\(\{t\in\mathbb C:\Im t>0,\ |t|<R\}\).
The boundary value on the diameter is zero,
and \eqref{eq:complex-decay} bounds $F$
by $O(R^{-1})$ on the semicircular boundary.
Letting $R\to\infty$, we obtain
\begin{equation}\label{eq:complex-max-conclusion}
F(t)\le0,\qquad\Im t>0.
\end{equation}
At \(t=i\), this gives \eqref{eq:complex-projection}.
Equation~\eqref{eq:complex-deficit} now proves \eqref{eq:schoenberg}.

\medskip
\noindent\emph{Step 4: The equality case.}
Suppose now that equality holds in
\eqref{eq:schoenberg}. Equations
\eqref{eq:complex-deficit} and
\eqref{eq:complex-subharmonic} give $F(i)=0$.
The strong maximum principle
\cite[p.~264, Theorem~3.2]{Con78},
together with \eqref{eq:complex-max-conclusion},
implies that $F$ vanishes identically.

By \eqref{eq:complex-positive-laplacian},
every local holomorphic branch of a critical
point satisfies
\begin{equation}\label{eq:complex-affine}
w_k'(t)
=
\frac{\Im w_k(t)}{\Im t}.
\end{equation}
The derivative $w_k'$ is holomorphic and
real-valued, so it is constant. Equation
\eqref{eq:complex-affine} then shows that
each branch is affine in $t$, with real
coefficients.

Choose the free branch with limiting slope
$
\eta\in(\Im a_{s-1},\Im a_s).
$
By \eqref{eq:complex-series} and
\eqref{eq:complex-affine}, it has the form
\[
w(t)=\eta t+b,\qquad b\in\mathbb R.
\]
Substitution into
\eqref{eq:complex-free-critical} gives the
rational identity
\begin{equation}\label{eq:complex-rational}
\sum_{j=1}^{s}
\frac{m_j}
{(\eta-\Im a_j)t+b-\Re a_j}=0.
\end{equation}

The $j$th summand has a real pole at
$
\frac{\Re a_j-b}{\eta-\Im a_j}
$
with residue
$
\frac{m_j}{\eta-\Im a_j}.
$
By the choice of $\eta$, the first $s-1$
residues are positive and the last is negative.
Since the sum in \eqref{eq:complex-rational}
vanishes identically, every pole carrying
a positive residue must coincide with the
pole carrying the negative residue.
All poles therefore coincide.

Thus there exists $\tau\in\mathbb R$ such that
\begin{equation}\label{eq:complex-collinearity}
\Re a_j+\tau\Im a_j=b+\tau\eta,
\qquad 1\le j\le s.
\end{equation}
Equation \eqref{eq:complex-collinearity}
shows that all zeros of $p(z)$ lie on a
straight line.

Conversely, if the centered zeros are
collinear, a rotation makes all zeros real.
Rolle's theorem then gives real critical
points, and \eqref{eq:complex-newton}
yields equality in \eqref{eq:schoenberg}.
This completes the proof.
\end{proof}

\section{Proof of Theorems~\ref{thm:centered-subunit} and \ref{thm:general-subunit}}\label{sec:subunit-proof}

For \(k=1\), the elementary symmetric function bound in the next
lemma is due to Erd\H{o}s and Niven~\cite{EN48}. The general case
follows from Malamud's convex-function inequality.

\begin{lem}[{\cite[Theorem~4.10]{Mal05}}]
\label{lem:elementary-symmetric}
Let \(z_1,z_2,\ldots,z_n\) be the zeros of a complex
polynomial \(p(z)\) of degree \(n\ge2\), and let
\(w_1,w_2,\ldots,w_{n-1}\) be its critical points.
Then
\begin{equation*}
e_k(|w_1|,\ldots,|w_{n-1}|)
\le
\frac{n-k}{n}
e_k(|z_1|,\ldots,|z_n|),
\qquad 1\le k\le n-1,
\end{equation*}
where \(e_k\) denotes the \(k\)-th elementary
symmetric polynomial.
\end{lem}
\begin{lem}\label{lem:logpower}
Let $0<p<1$. Then, for every $a\ge0$,
\begin{equation}\label{eq:log-integral-power}
a^p=\frac{p\sin(\pi p)}{\pi}
\int_0^\infty \log(1+ta)t^{-p-1}\dd t.
\end{equation}
\end{lem}
\begin{proof}
The identity is immediate when $a=0$. Assume that $a>0$.
Making the change of variables $u=ta$, we obtain
\begin{equation}\label{eq:logpower-substitution}
\int_0^\infty \log(1+ta)t^{-p-1}\dd t
=a^p\int_0^\infty \log(1+u)u^{-p-1}\dd u.
\end{equation}
The integral on the right converges, since
\[
\log(1+u)u^{-p-1}
=
\begin{cases}
O(u^{-p}), & u\to0^+,\\[2mm]
O(u^{-p-1}\log u), & u\to\infty,
\end{cases}
\]
and $0<p<1$.

To evaluate this integral, we integrate by parts on
$[\varepsilon,R]$, where $0<\varepsilon<R$:
\[
\int_\varepsilon^R \log(1+u)u^{-p-1}\dd u =-\frac1p \left[u^{-p}\log(1+u)\right]_\varepsilon^R+\frac1p \int_\varepsilon^R\frac{u^{-p}}{1+u}\dd u.
\]
The boundary terms vanish because \(u^{-p}\log(1+u)\sim u^{1-p}\)
as \(u\downarrow0\) and \(u^{-p}\log(1+u)\to0\) as \(u\to\infty\).
Hence
\begin{equation}\label{eq:logpower-parts}
\int_0^\infty \log(1+u)u^{-p-1}\dd u
=\frac1p\int_0^\infty\frac{u^{-p}}{1+u}\dd u.
\end{equation}

Finally, the substitution $s=u/(1+u)$ gives
\begin{equation}\label{eq:logpower-beta}
\int_0^\infty\frac{u^{-p}}{1+u}\dd u
=\int_0^1 s^{-p}(1-s)^{p-1}\dd s=B(1-p,p)=\Gamma(1-p)\Gamma(p)=\frac{\pi}{\sin(\pi p)},
\end{equation}
where the last equality follows from Euler's reflection
formula. Combining \eqref{eq:logpower-substitution},
\eqref{eq:logpower-parts}, and \eqref{eq:logpower-beta}
proves the integral formula in Lemma~\ref{lem:logpower}.
\end{proof}

\begin{proof}[Proof of Theorem~\ref{thm:centered-subunit}]
When \(n=2\), the centroid condition implies that
the unique critical point is zero. The inequality
and the equality characterization are therefore
immediate. Henceforth, assume \(n\ge3\).

\medskip
\noindent\emph{Case 1: The inequality for \(p=1\).}
By Lemma~\ref{lem:elementary-symmetric}, we have
\[
e_k(|w_1|,\ldots,|w_{n-1}|)
\le
\frac{n-k}{n}e_k(|z_1|,\ldots,|z_n|),
\qquad 1\le k\le n-1.
\]
In particular,
\(e_2(|w_1|,\ldots,|w_{n-1}|) \le \frac{n-2}{n}e_2(|z_1|,\ldots,|z_n|)\).
Combining this with Schoenberg's inequality~\eqref{eq:schoenberg},
we obtain
\[
\begin{aligned}
\left(\sum_{k=1}^{n-1}|w_k|\right)^2
&=\sum_{k=1}^{n-1}|w_k|^2
+2e_2(|w_1|,\ldots,|w_{n-1}|)\\
&\le
\frac{n-2}{n}
\left(
\sum_{j=1}^n|z_j|^2
+2e_2(|z_1|,\ldots,|z_n|)
\right)\\
&=
\frac{n-2}{n}
\left(\sum_{j=1}^n|z_j|\right)^2.
\end{aligned}
\]
Taking square roots gives
\begin{equation}\label{eq:centered-first-power}
\sum_{k=1}^{n-1}|w_k|
\le
\left(\frac{n-2}{n}\right)^{\frac{1}{2}}
\sum_{j=1}^n|z_j|.
\end{equation}
Equation~\eqref{eq:centered-first-power} proves the case \(p=1\).

For \(2\le k\le n\), Bernoulli's inequality yields
\(\left(\frac{n-2}{n}\right)^{\frac{k}{2}} = \left(1-\frac2n\right)^{\frac{k}{2}} \ge 1-\frac{k}{n}=\frac{n-k}{n}\).
It follows from Lemma~\ref{lem:elementary-symmetric}
and \eqref{eq:centered-first-power} that
\begin{equation}\label{eq:elementary-scaled}
e_k(|w_1|,\ldots,|w_{n-1}|,0)
\le
\left(\frac{n-2}{n}\right)^{\frac{k}{2}}
e_k(|z_1|,\ldots,|z_n|),
\qquad 1\le k\le n,
\end{equation}
where the case \(k=n\) is trivial.

Expanding the products in terms of elementary
symmetric polynomials and comparing their
coefficients in \eqref{eq:elementary-scaled},
we conclude that
\begin{equation}\label{eq:product-comparison}
\prod_{k=1}^{n-1}(1+t|w_k|)
\le
\prod_{j=1}^n
\left(
1+t\sqrt{\frac{n-2}{n}}\,|z_j|
\right),
\qquad t\ge0.
\end{equation}

\medskip
\noindent\emph{Case 2: The inequality for \(0<p<1\).}
Now assume that \(0<p<1\).
Taking logarithms in
\eqref{eq:product-comparison} gives
\begin{equation}\label{eq:log-product-comparison}
\sum_{k=1}^{n-1}\log(1+t|w_k|)
\le
\sum_{j=1}^n
\log\left(
1+t\sqrt{\frac{n-2}{n}}\,|z_j|
\right).
\end{equation}
Apply the integral representation \eqref{eq:log-integral-power} of
Lemma~\ref{lem:logpower} term by term.
Since \(p\sin(\pi p)>0\), multiplying
\eqref{eq:log-product-comparison} by
\(\frac{p\sin(\pi p)}{\pi}t^{-p-1}\)
and integrating over \((0,\infty)\),
we obtain from \eqref{eq:log-integral-power}
\(\sum_{k=1}^{n-1}|w_k|^p \le \left(\frac{n-2}{n}\right)^{\frac{p}{2}} \sum_{j=1}^n|z_j|^p\).
Together with \eqref{eq:centered-first-power},
this establishes the inequality for
every \(0<p\le1\).

It remains to determine the equality cases.

\medskip
\noindent\emph{Case 1: Equality for \(p=1\).}
Assume that equality holds in
\eqref{eq:centered-first-power}.
Since both Schoenberg's inequality and the
elementary symmetric function inequality
used in its derivation have nonnegative
deficits, equality must hold in each.
In particular,
\begin{equation}\label{eq:e2-equality}
e_2(|w_1|,\ldots,|w_{n-1}|)
=
\frac{n-2}{n}
e_2(|z_1|,\ldots,|z_n|).
\end{equation}

We show that \eqref{eq:e2-equality} is
impossible if at least three zeros are
nonzero. Set $A=QDQ$, where $Q,D$ are defined in Lemma~\ref{lem:zpprimez}.
Let \(e_1,\ldots,e_n\) denote the standard
orthonormal basis of \(\C^n\), and define
\[
v_{ij}=
(\bigwedge\nolimits^2Q)(e_i\wedge e_j),
\qquad 1\le i<j\le n.
\]
By the multiplicativity of exterior powers,
\begin{equation}\label{eq:exterior-expansion}
\bigwedge\nolimits^2 A
=
\sum_{1\le i<j\le n}
z_iz_jv_{ij}v_{ij}^*.
\end{equation}
The entries of \(Q\) are \(Q_{ij}=\delta_{ij}-1/n\).
The determinant formula for inner products in exterior powers gives
\begin{equation}\label{eq:exterior-inner-products}
\|v_{ij}\|^2=\frac{n-2}{n},
\qquad
|\langle v_{ij},v_{ik}\rangle|
=\frac1n
\end{equation}
whenever \(i,j,k\) are pairwise distinct.

The eigenvalues of \(A\) are \(w_1,\ldots,w_{n-1},0\), counted
with algebraic multiplicity. The eigenvalues of \(\bigwedge^2 A\)
are therefore their pairwise products. Using the trace-norm bound
for the sum of the absolute values
of the eigenvalues, followed by
\eqref{eq:exterior-expansion} and
\eqref{eq:exterior-inner-products}, we get
\begin{equation}\label{eq:exterior-trace-bound}
e_2(|w_1|,\ldots,|w_{n-1}|)\le
\left\|\bigwedge\nolimits^2 A\right\|_1\le
\sum_{i<j}|z_iz_j|\|v_{ij}v_{ij}^*\|_1=
\frac{n-2}{n}
e_2(|z_1|,\ldots,|z_n|).
\end{equation}

By \eqref{eq:e2-equality}, both inequalities
in \eqref{eq:exterior-trace-bound} must
be equalities. In particular, equality
holds in the trace-norm triangle inequality
applied to \eqref{eq:exterior-expansion}.

By the dual characterization of the
trace norm, there exists a contraction
\(U\) such that
\(\Re\tr\left( U^*\bigwedge\nolimits^2 A \right) = \left\|\bigwedge\nolimits^2 A\right\|_1\).
Equality in the triangle inequality then
forces every nonzero term in
\eqref{eq:exterior-expansion} to attain
equality in the same dual estimate.
Consequently,
\begin{equation}\label{eq:exterior-equality-condition}
Uv_{ij}
=
\frac{z_iz_j}{|z_iz_j|}v_{ij}
\qquad\text{whenever }z_iz_j\ne0.
\end{equation}
Indeed, this follows from equality in
the Cauchy--Schwarz inequality and
\(\|U\|\le1\).

Recall that eigenvectors of a contraction
corresponding to distinct eigenvalues of
modulus one are orthogonal. If
\(z_i,z_j,z_k\) are all nonzero,
then \eqref{eq:exterior-inner-products}
shows that \(v_{ij}\) and \(v_{ik}\)
are not orthogonal. Thus
\eqref{eq:exterior-equality-condition} implies
\(\frac{z_iz_j}{|z_iz_j|} = \frac{z_iz_k}{|z_iz_k|}\),
and therefore
\(\frac{z_j}{|z_j|} = \frac{z_k}{|z_k|}\).
If at least three zeros are nonzero,
we may choose \(i\) distinct from any
given pair \(j,k\) of nonzero zeros.
It follows that all nonzero zeros
have the same argument. This contradicts
the centroid condition
\(\sum_{j=1}^n z_j=0\).

Therefore, equality for \(p=1\) requires
at most two nonzero zeros.

\medskip
\noindent\emph{Case 2: Equality for \(0<p<1\).}
Assume that at least three of the zeros
are nonzero. Then
\(e_3(|z_1|,\ldots,|z_n|)>0\).
Since \(n\ge3\),
\(\left(\frac{n-2}{n}\right)^{\frac{3}{2}} > \frac{n-3}{n}\).
For \(n\ge4\), apply Lemma~\ref{lem:elementary-symmetric} with \(k=3\).
For \(n=3\), \(e_3(|w_1|,|w_2|,0)=0\). In either case,
\[
e_3(|w_1|,\ldots,|w_{n-1}|,0)
\le
\frac{n-3}{n}e_3(|z_1|,\ldots,|z_n|)<
\left(\frac{n-2}{n}\right)^{\frac{3}{2}}
e_3(|z_1|,\ldots,|z_n|).
\]
Thus the coefficient of \(t^3\) in
\eqref{eq:product-comparison} is strictly
smaller on the left, while all other
coefficient inequalities remain valid.
Consequently, \eqref{eq:product-comparison},
and hence \eqref{eq:log-product-comparison},
is strict for every \(t>0\).
The positivity of the integration weight
in \eqref{eq:log-integral-power} then shows
that the inequality in Theorem~\ref{thm:centered-subunit} is strict.

Therefore, equality for \(0<p<1\) also
requires at most two nonzero zeros.

\medskip\noindent\emph{Sharpness.}
In either case, the centroid condition
implies that if there are exactly two
nonzero zeros, they must be \(a\) and
\(-a\) for some \(a\ne0\).
If all zeros vanish, the polynomial is
a nonzero constant multiple of \(z^n\).
Thus all equality cases are of the form
\[
p(z)=c\,z^{n-2}(z^2-a^2),
\qquad c\in\C\setminus\{0\},\quad a\in\C.
\]

Conversely, differentiating this polynomial
gives
\(p'(z)=c\,z^{n-3} \bigl(nz^2-(n-2)a^2\bigr)\).
Its only possible nonzero critical points
are
$
\pm\sqrt{\frac{n-2}{n}}\,a.
$
Consequently, for every \(0<p\le1\),
\[
\sum_{k=1}^{n-1}|w_k|^p=
2\left(\frac{n-2}{n}\right)^{\frac{p}{2}}|a|^p=
\left(\frac{n-2}{n}\right)^{\frac{p}{2}}
\sum_{j=1}^n|z_j|^p.
\]
Hence the constant is sharp, and the
equality characterization is complete.

\end{proof}

As in the proof of Theorem~\ref{thm:centered-subunit}, we next prove Theorem~\ref{thm:general-subunit}.

\begin{proof}[Proof of Theorem~\ref{thm:general-subunit}]
By Lemma~\ref{lem:elementary-symmetric}, we have
\[
e_k(|w_1|,\ldots,|w_{n-1}|)
\le
\frac{n-k}{n}e_k(|z_1|,\ldots,|z_n|),
\qquad 1\le k\le n-1.
\]
On the other hand, Bernoulli's inequality gives
\[
\left(\frac{n-1}{n}\right)^k
=
\left(1-\frac1n\right)^k
\ge 1-\frac{k}{n}
=\frac{n-k}{n},
\qquad 1\le k\le n.
\]
Consequently,
\begin{equation}\label{eq:general-elementary-scaled}
e_k(|w_1|,\ldots,|w_{n-1}|,0)
\le
\left(\frac{n-1}{n}\right)^k
e_k(|z_1|,\ldots,|z_n|),
\qquad 1\le k\le n,
\end{equation}
where the case \(k=n\) is immediate.

Expanding the products in terms of elementary
symmetric polynomials and comparing their
coefficients in \eqref{eq:general-elementary-scaled},
we obtain
\begin{equation}\label{eq:general-product-comparison}
\prod_{k=1}^{n-1}(1+t|w_k|)
\le
\prod_{j=1}^n
\left(1+\frac{n-1}{n}t|z_j|\right),
\qquad t\ge0.
\end{equation}

Taking logarithms in
\eqref{eq:general-product-comparison}, we get
\begin{equation}\label{eq:general-log-comparison}
\sum_{k=1}^{n-1}\log(1+t|w_k|)
\le
\sum_{j=1}^n
\log\left(1+\frac{n-1}{n}t|z_j|\right).
\end{equation}
Since \(0<p<1\), the weight
$
\frac{p\sin(\pi p)}{\pi}t^{-p-1}
$
is strictly positive for every \(t>0\).
Multiplying \eqref{eq:general-log-comparison}
by this weight and integrating over
\((0,\infty)\), we obtain, by
Lemma~\ref{lem:logpower},
\[
\sum_{k=1}^{n-1}|w_k|^p
\le
\sum_{j=1}^n
\left(\frac{n-1}{n}|z_j|\right)^p=
\left(\frac{n-1}{n}\right)^p
\sum_{j=1}^n|z_j|^p.
\]

\medskip
\noindent\emph{Equality and sharpness.}
Assume that at least two zeros are nonzero.
Then
\(e_2(|z_1|,\ldots,|z_n|)>0\).
Since
\(\left(\frac{n-1}{n}\right)^2 -\frac{n-2}{n} =\frac1{n^2}>0\),
Lemma~\ref{lem:elementary-symmetric}
implies, for \(n\ge3\),
\[
e_2(|w_1|,\ldots,|w_{n-1}|,0)
\le
\frac{n-2}{n}e_2(|z_1|,\ldots,|z_n|)<
\left(\frac{n-1}{n}\right)^2
e_2(|z_1|,\ldots,|z_n|).
\]
The same strict inequality holds when
\(n=2\), since its left-hand side vanishes.

Thus the coefficient of \(t^2\) in
\eqref{eq:general-product-comparison}
is strictly smaller on the left, while
all the remaining coefficient inequalities
are valid. It follows that
\eqref{eq:general-product-comparison},
and hence \eqref{eq:general-log-comparison},
is strict for every \(t>0\).
Since the integration weight in
Lemma~\ref{lem:logpower} is strictly positive,
the resulting power inequality is also strict.

Therefore, equality is possible only if
at most one zero is nonzero. If the zeros
are not all zero, this means that, up to
a nonzero constant factor,
\[
p(z)=z^{n-1}(z-a),
\qquad a\in\C\setminus\{0\}.
\]

Conversely, for
$
p(z)=c\,z^{n-1}(z-a),
$ with $c\ne0$ and $a\ne0,$
we have
\(p'(z)=c\,z^{n-2} \bigl(nz-(n-1)a\bigr)\).
Hence the only nonzero critical point is
$
\frac{n-1}{n}a,
$
and therefore
\(\sum_{k=1}^{n-1}|w_k|^p = \left(\frac{n-1}{n}\right)^p|a|^p = \left(\frac{n-1}{n}\right)^p \sum_{j=1}^n|z_j|^p\).
Thus the constant is sharp, and the
equality characterization is complete.
\end{proof}

\section{Proof of Theorem~\ref{thm:centered-subquadratic}}
\label{sec:subquadratic-proof}

Let $n\ge2$ be an integer, and let $z_1,\ldots,z_n$ and
$w_1,\ldots,w_{n-1}$ be the zeros and critical points, respectively,
of a complex polynomial \(p(z)\) of degree \(n\). For a real exponent
$p\ge1$, put
\begin{equation}\label{eq:notation}
Z_p:=\sum_{j=1}^{n}|z_j|^p,\qquad W_p:=\sum_{k=1}^{n-1}|w_k|^p.
\end{equation}

We begin by showing that the optimal constant in
Theorem~\ref{thm:centered-subquadratic} can be determined from
polynomials with real zeros. The reduction uses Pereira's
majorization theorem.

\begin{lem}[{\cite[Theorem~4.6]{Per03}}]
\label{lem:pereira-projection}
Let $w_1,\ldots,w_{n-1}$ and $\eta_1,\ldots,\eta_{n-1}$ be the
critical points of $\prod_{j=1}^{n}(z-z_j)$ and
$\prod_{j=1}^{n}(x-\Re z_j)$, respectively. Then
\[
(\Re w_1,\ldots,\Re w_{n-1})\prec(\eta_1,\ldots,\eta_{n-1}).
\]
 Consequently, for every convex
function $\varphi:\mathbb R\to\mathbb R$,
\begin{equation}\label{eq:pereira-projection}
\sum_{k=1}^{n-1}\varphi(\Re w_k)\le\sum_{k=1}^{n-1}\varphi(\eta_k).
\end{equation}
\end{lem}

Recall $C_{n,p}$ is the least constant for which $W_p\le C_{n,p}Z_p$
holds for every real-rooted polynomial of degree $n$ satisfying
$\sum_{j=1}^{n}z_j=0$. Consider now a complex polynomial with
$\sum_{j=1}^{n}z_j=0$. For $\theta\in[0,2\pi]$, let
$\eta_1(\theta),\ldots,\eta_{n-1}(\theta)$ be the critical points of
$\prod_{j=1}^{n}(x-\Re(e^{-i\theta}z_j))$. Applying
inequality~\eqref{eq:pereira-projection} to the polynomial with zeros
\(e^{-i\theta}z_j\) and the convex function \(\varphi(x)=|x|^p\) gives
the first inequality in \eqref{eq:complex-projection-transfer}.
The second follows from the real-rooted bound, since
\(\sum_j\Re(e^{-i\theta}z_j)=0\).
\begin{equation}\label{eq:complex-projection-transfer}
\sum_{k=1}^{n-1}|\Re(e^{-i\theta}w_k)|^p
\le\sum_{k=1}^{n-1}|\eta_k(\theta)|^p
\le C_{n,p}\sum_{j=1}^{n}|\Re(e^{-i\theta}z_j)|^p.
\end{equation}
Integrating \eqref{eq:complex-projection-transfer} over
$0\le\theta\le2\pi$ and using
\begin{equation}\label{eq:angular}
\frac1{2\pi}\int_0^{2\pi}|\Re(e^{-i\theta}z)|^p\,d\theta
=\frac{|z|^p}{2\pi}\int_0^{2\pi}|\cos\theta|^p\,d\theta,
\end{equation}
we obtain $W_p\le C_{n,p}Z_p$ for centered complex polynomials:
the positive factor
$\frac1{2\pi}\int_0^{2\pi}|\cos\theta|^p\,d\theta$
cancels from both sides. Conversely, every centered real-rooted
polynomial is a centered complex polynomial. Thus the sharp
constants in the real and complex cases coincide. The equality
cases for complex zeros will be considered after the sharp
constant has been determined.

For the remainder of the section, fix $1<p<2$. The cases $n=2,3$,
including their equality cases, are treated in
Proposition~\ref{sub12:prop:turn3-degrees-two-three}. We henceforth
assume $n\ge4$. For real zeros, the quotient $W_p/Z_p$ is invariant
under nonzero real scaling and attains its maximum on the compact
set $\sum_{j=1}^{n}z_j=0$, $\sum_{j=1}^{n}z_j^2=1$. Choose a
maximizing configuration. Dividing the zeros by half the distance
between the smallest and largest zero, and reflecting them if
necessary, we may assume that the two extreme zeros are $t-1$ and
$t+1$ with $t\ge0$. Because the remaining $n-2$ zeros lie between
$t-1$ and $t+1$, the centroid condition implies
$0\le t\le(n-2)/n$.

We use a compression representation for critical points. Related
matrix representations appear in Komarova and Rivin~\cite{KR03}.
For diagonal \(D\), \(e=n^{-\frac{1}{2}}\mathbf1\) satisfies
\(e^*D^j e=n^{-1}\tr(D^j)\) for all integers \(j\ge0\).
Thus the compression also follows from Pereira's characterization of
differentiators~\cite[Theorem~2.5]{Per03}.

\begin{lem}\label{lem:zpprimez}
Let $p(z)=\prod_{j=1}^{n}(z-z_j)$, and write
$\mathbf1=(1,\ldots,1)^T$,
$D=\operatorname{diag}(z_1,\ldots,z_n)$, and
$Q=I-\mathbf1\mathbf1^T/n$. Then
\[
\det(zI-DQ)=\frac{z}{n}p'(z).
\]
Consequently, the compression $A=QDQ|_{\mathbf1^\perp}$ has the
critical points of $p(z)$ as its eigenvalues, counted with
multiplicity.
\end{lem}

The next lemma converts the first-order optimality condition for
$W_p/Z_p$ into an integral identity with three prescribed moments.

\begin{lem}\label{lem:real-spectral-reduction}
Let a centered real-rooted polynomial maximize $W_p/Z_p$, and
normalize its smallest and largest zeros to $t-1$ and $t+1$, where
$0\le t\le(n-2)/n$. Then there exists a probability measure $\mu$
on $[-1,1]$ such that
\begin{equation}\label{eq:section4-gradient}
\frac{W_p}{Z_p}
=\frac{2\displaystyle\int_{-1}^{1}y\,\operatorname{sgn}(t+y)|t+y|^{p-1}\,d\mu(y)}
{(1-t)^{p-1}+(1+t)^{p-1}},
\end{equation}
and
\begin{equation}\label{eq:section4-moments}
\int_{-1}^{1}y\,d\mu(y)=0,\qquad
\int_{-1}^{1}y^2\,d\mu(y)=\frac{n-2}{n},\qquad
\int_{-1}^{1}y^3\,d\mu(y)=-\frac{2t}{n}.
\end{equation}
\end{lem}

\begin{proof}
Choose indices $i,j$ with $z_i=t-1$ and $z_j=t+1$, and retain
$D,Q,A$ from Lemma~\ref{lem:zpprimez}. As the zeros are real, $A$
is self-adjoint. Let $X=A-tI$ on $\mathbf1^\perp$ and
$f=(e_j-e_i)/\sqrt2$, where $e_1,\ldots,e_n$ are the standard basis
vectors. A calculation using $z_i=t-1$ and $z_j=t+1$ gives
$Xf=Q(e_i+e_j)/\sqrt2$.

Vary the two extreme zeros by $z_i(\varepsilon)=z_i-\varepsilon$
and $z_j(\varepsilon)=z_j+\varepsilon$, keeping the other zeros
fixed. This perturbation preserves the centroid and gives
$A'(0)=f(Xf)^T+(Xf)f^T$. Put
$h(s)=\operatorname{sgn}(s)|s|^{p-1}$. Since $p>1$, the function
$|s|^p$ is continuously differentiable, and the trace derivative
formula for self-adjoint matrices yields
$W_p'(0)=p\operatorname{Tr}(h(A)A'(0))$. Direct differentiation also
gives $Z_p'(0)=p((1-t)^{p-1}+(1+t)^{p-1})>0$.
The derivative of $W_p(\varepsilon)/Z_p(\varepsilon)$ vanishes at
$\varepsilon=0$ by maximality. Therefore
\begin{equation}\label{eq:real-variational}
\frac{W_p}{Z_p}
=\frac{W_p'(0)}{Z_p'(0)}
=\frac{2\langle h(tI+X)f,Xf\rangle}
{(1-t)^{p-1}+(1+t)^{p-1}}.
\end{equation}

The eigenvalues of the self-adjoint compression $A$ lie in
$[t-1,t+1]$, so $\operatorname{spec}(X)\subseteq[-1,1]$. Let $\mu$
be the spectral measure of $X$ associated with the unit vector $f$;
in particular,
$\int_{-1}^{1}\varphi(y)\,d\mu(y)=\langle\varphi(X)f,f\rangle$
for every continuous function $\varphi$ on $[-1,1]$. Applying the
spectral identity to \eqref{eq:real-variational} gives
\eqref{eq:section4-gradient}.

For the first two moments, $\langle Xf,f\rangle=0$ and
$\|Xf\|^2=(n-2)/n$. For the third moment, the explicit expression
for $Xf$ and $\sum_{r=1}^{n}z_r=0$ give
\begin{equation}\label{eq:real-third-moment}
\int_{-1}^{1}y^3\,d\mu(y)
=\langle X(Xf),Xf\rangle
=\frac12\sum_{r=1}^{n}(z_r-t)
\left(\delta_{ri}+\delta_{rj}-\frac2n\right)^2
=-\frac{2t}{n}.
\end{equation}
In \eqref{eq:real-third-moment} we used \(z_i-t=-1\), \(z_j-t=1\), and
$\sum_{r\ne i,j}(z_r-t)=-nt$. The three moment identities in
\eqref{eq:section4-moments} follow.
\end{proof}

The conditions in \eqref{eq:section4-moments} hold for the spectral
measure of an extremal polynomial. They do not, by themselves,
retain the restriction that a polynomial has an integer number of
zeros at each location. That restriction enters through a fourth
moment.

\begin{lem}\label{lem:section4-fourth}
Under the assumptions of Lemma~\ref{lem:real-spectral-reduction},
write the remaining $n-2$ zeros as $t+y_1,\ldots,t+y_{n-2}$, where
$-1\le y_\nu\le1$. Then
\begin{equation}\label{eq:section4-root-defect}
\int_{-1}^{1}(1-y^2)(y-t)^2\,d\mu(y)
=\frac{2}{n^2}\sum_{\nu=1}^{n-2}(1-y_\nu^2).
\end{equation}
If $\ell\le t\le u$, where
\begin{equation}\label{eq:section4-cell}
\ell=1-\frac{2(k+1)}n,\qquad
u=1-\frac{2k}n,\qquad 1\le k<\frac n2,
\end{equation}
then
\begin{equation}\label{eq:section4-fourth}
\int_{-1}^{1}(1-y^2)(y-t)^2\,d\mu(y)
\ge2(t-\ell)(u-t).
\end{equation}
\end{lem}

\begin{proof}
The centroid condition gives $\sum_{\nu=1}^{n-2}y_\nu=-nt$.
For the unit vector $f$ in Lemma~\ref{lem:real-spectral-reduction},
we have $Xf=Q(e_i+e_j)/\sqrt2$. Since
$X^2f=Q(D-tI)Xf$, evaluating its Euclidean norm gives
\begin{equation}\label{eq:section4-fourth-identity}
\int_{-1}^{1}y^4\,d\mu(y)
=\|X^2f\|^2
=\left(\frac{n-2}{n}\right)^2
+\frac2{n^2}
\left(\sum_{\nu=1}^{n-2}y_\nu^2-nt^2\right).
\end{equation}
Expanding $(1-y^2)(y-t)^2$ and substituting
\eqref{eq:section4-moments} and \eqref{eq:section4-fourth-identity}
shows that its integral is
$2(n-2-\sum_{\nu=1}^{n-2}y_\nu^2)/n^2$, which is
\eqref{eq:section4-root-defect}.

To prove \eqref{eq:section4-fourth}, maximize
$\sum_{\nu=1}^{n-2}y_\nu^2$ subject to $-1\le y_\nu\le1$ and
$\sum_{\nu=1}^{n-2}y_\nu=-nt$. At a maximizing point, at most one
coordinate lies strictly between $-1$ and $1$: if two coordinates
were interior, changing them by equal and opposite amounts would
preserve their sum, while their squared sum would be strictly
convex in the perturbation parameter. For $\ell\le t\le u$, a
maximizing list has $k-1$ entries equal to $1$, $n-k-2$ entries
equal to $-1$, and one entry equal to $v=n(u-t)-1$. Since
$-1\le v\le1$ and $1-v^2=n^2(t-\ell)(u-t)$, we obtain
$\sum_{\nu=1}^{n-2}(1-y_\nu^2)\ge n^2(t-\ell)(u-t)$.
Equation~\eqref{eq:section4-root-defect} now yields
\eqref{eq:section4-fourth}.
\end{proof}

For $0<r\le1$, define
\begin{equation}\label{eq:section4-H}
H_p(r)=\frac{(1+r)\{1+r^p-(1-r)^p\}}{r+r^p}.
\end{equation}
The centered polynomial
$p(z)=(z-(n-k)a)^k(z+ka)^{n-k}$, with $a\ne0$ and
$1\le k\le\lfloor n/2\rfloor$, has
\begin{equation}\label{eq:section4-binary}
R_{n,p}(k)
=\frac{(k-1)(n-k)^p+(n-k-1)k^p+(n-2k)^p}
{k(n-k)^p+(n-k)k^p}
=1-\frac1nH_p\left(\frac{k}{n-k}\right).
\end{equation}
Theorem~\ref{sub12:thm:subquad-H-unimodal} shows that $H_p$ is
strictly decreasing before its unique minimizer $r_p\in(0,1)$
and strictly increasing afterward. Put $\theta_p=r_p/(1+r_p)$.
Because $k/(n-k)$ increases with $k$, equation
\eqref{eq:section4-binary} shows that the largest two-value ratio
occurs at an admissible integer adjacent to $n\theta_p$.

When $n\theta_p$ is not an integer, the continuous minimum of
$H_p$ need not correspond to a polynomial with integer root
multiplicities. The fourth-moment bound \eqref{eq:section4-fourth}
provides the correction needed to compare a noninteger value with
its neighboring integer values.

\begin{lem}\label{lem:section4-penalty}
Assume that $n\theta_p>1$ is not an integer, and put
$k=\lfloor n\theta_p\rfloor$. Let $\ell,u$ be as in
\eqref{eq:section4-cell}. If the maximizing configuration in
Lemma~\ref{lem:real-spectral-reduction} has $\ell\le t\le u$, then
\begin{equation}\label{eq:section4-penalty}
\frac{W_p}{Z_p}
\le1-\frac1nH_p\left(\frac{1-t}{1+t}\right)
-\frac{4E_*(t)}{(1-t)^{p-1}+(1+t)^{p-1}}
(t-\ell)(u-t),
\end{equation}
where $E_*(t)$ is the largest admissible nonnegative quartic
coefficient from Theorem~\ref{sub12:thm:quartic-anchor}.
Theorem~\ref{sub12:thm:turn3-peak-anchor-positive}
gives $E_*(t)>0$ on $[\ell,u]$.
\end{lem}

\begin{proof}
Write $\phi_t(y)=y\operatorname{sgn}(t+y)|t+y|^{p-1}$, and let
$Q_0$ be the cubic Hermite interpolant determined by
$Q_0(-1)=\phi_t(-1)$, $Q_0(1)=\phi_t(1)$,
$Q_0(t)=\phi_t(t)$, and $Q_0'(t)=\phi_t'(t)$.
By Theorem~\ref{sub12:thm:quartic-anchor} and the positivity in
Theorem~\ref{sub12:thm:turn3-peak-anchor-positive},
$Q_0(y)-E_*(t)(1-y^2)(y-t)^2\ge\phi_t(y)$ for $-1\le y\le1$.

The three moment identities in \eqref{eq:section4-moments} imply,
for every polynomial $q$ of degree at most three,
\begin{equation}\label{eq:section4-three-point-evaluation}
\int_{-1}^{1}q(y)\,d\mu(y)
=\frac{n-2+nt}{2n(1+t)}q(-1)
+\frac{2}{n(1-t^2)}q(t)
+\frac{n-2-nt}{2n(1-t)}q(1).
\end{equation}
Indeed, the coefficients in \eqref{eq:section4-three-point-evaluation}
reproduce the moments of $1,y,y^2,y^3$ in
\eqref{eq:section4-moments}. Applying
\eqref{eq:section4-three-point-evaluation} to $q=Q_0$ and using
the interpolation values at $-1,t,1$, we obtain
\begin{equation}\label{eq:section4-cubic-expectation}
\frac{2\displaystyle\int_{-1}^{1}Q_0(y)\,d\mu(y)}
{(1-t)^{p-1}+(1+t)^{p-1}}
=1-\frac1nH_p\left(\frac{1-t}{1+t}\right).
\end{equation}
Integrating the pointwise quartic majorant against $\mu$ and
applying \eqref{eq:section4-gradient},
\eqref{eq:section4-cubic-expectation}, and
\eqref{eq:section4-fourth} yields \eqref{eq:section4-penalty}.
\end{proof}

The reduction outside the integer cell is proved in the argument
for Theorem~\ref{sub12:thm:turn3-peak-cell-sufficient}. Inside the
cell, Theorem~\ref{sub12:thm:turn3-integer-anchor-branch} and the
chord inequality \eqref{sub12:eq:turn3-right-chord-theta} control
the correction term in \eqref{eq:section4-penalty}. Together they
give the sharp estimate for real zeros.

\begin{lem}\label{lem:section4-integer-comparison}
Let \(n\ge4\) and \(1<p<2\). For every centered real-rooted polynomial
of degree \(n\) maximizing \(W_p/Z_p\) among root lists with \(Z_p>0\),
\begin{equation}\label{eq:section4-final-bound}
\frac{W_p}{Z_p}
\le\max\left\{
\left(\frac{n-2}{n}\right)^{\frac{p}{2}},
\max_{1\le k\le\lfloor n/2\rfloor}R_{n,p}(k)
\right\}.
\end{equation}
\end{lem}

\begin{proof}
If $n\theta_p\le1$, the exact three-moment optimization in
Theorem~\ref{sub12:thm:turn3-uncapped-lp3-exact} gives
$W_p/Z_p\le\max\{((n-2)/n)^{\frac{p}{2}},R_{n,p}(1)\}$.
If $n\theta_p>1$ is an integer,
Theorem~\ref{sub12:thm:turn3-uncapped-lp3-exact} gives
$W_p/Z_p\le1-H_p(r_p)/n=R_{n,p}(n\theta_p)$.
Both bounds imply \eqref{eq:section4-final-bound}.

Assume that $n\theta_p>1$ is not an integer. Put
$k=\lfloor n\theta_p\rfloor$ and define $\ell,u$ by
\eqref{eq:section4-cell}. The three-moment classification and
Radau boundary analysis in the proof of
Theorem~\ref{sub12:thm:turn3-peak-cell-sufficient} show that
\eqref{eq:section4-final-bound} already holds when
$t\notin[\ell,u]$.

For $\ell\le t\le u$, apply \eqref{eq:section4-penalty}. The
quartic coefficient is $E_*(t)=\min\{E_-(t),E_+(t)\}$, where
\(E_-(t)\) and \(E_+(t)\) are given by \eqref{eq:quartic-endpoint-coefficients}. If
$E_-(t)\le E_+(t)$, then
Theorem~\ref{sub12:thm:turn3-integer-anchor-branch} gives
$H_p((1-t)/(1+t))>H_p(k/(n-k))$. Together with
$E_*(t)>0$, equation \eqref{eq:section4-penalty} implies
$W_p/Z_p\le R_{n,p}(k)$.

If $E_+(t)<E_-(t)$, put $\theta=(1-t)/2$,
$a=k/n$, and $b=(k+1)/n$. Then
$t-\ell=2(b-\theta)$ and $u-t=2(\theta-a)$.
The strict backward-remainder theorem
\ref{sub12:thm:turn3-backward-remainder-proved} yields the
chord inequality \eqref{sub12:eq:turn3-right-chord-theta}.
Combining \eqref{sub12:eq:turn3-right-chord-theta} with
\eqref{eq:section4-penalty} gives
$W_p/Z_p\le\max\{R_{n,p}(k),R_{n,p}(k+1)\}$.
This proves \eqref{eq:section4-final-bound} in the last case.
\end{proof}

\begin{proof}[Proof of Theorem~\ref{thm:centered-subquadratic}]
For $n\ge4$, Lemma~\ref{lem:section4-integer-comparison} gives
the upper bound
\begin{equation}\label{eq:section4-sharp}
C_{n,p}\le\max\left\{
\left(\frac{n-2}{n}\right)^{\frac{p}{2}},
\max_{1\le k\le\lfloor n/2\rfloor}R_{n,p}(k)
\right\}.
\end{equation}
The polynomial $p(z)=z^{n-2}(z^2-a^2)$ attains
$((n-2)/n)^{\frac{p}{2}}$, while
$p(z)=(z-(n-k)a)^k(z+ka)^{n-k}$ attains $R_{n,p}(k)$,
for every $a\in\mathbb C\setminus\{0\}$ and
$1\le k\le\lfloor n/2\rfloor$. Both families are centered and
have real zeros when $a$ is real. Hence equality holds in
\eqref{eq:section4-sharp}, and $C_{n,p}$ equals its right-hand side.
The Pereira reduction, proved through
\eqref{eq:complex-projection-transfer} and \eqref{eq:angular},
transfers the same sharp constant to centered complex polynomials.

Theorem~\ref{sub12:thm:turn3-full-finite-equality} gives the complete
equality classification for $n\ge4$: besides the zero root list,
equality is possible only for nonzero scalar multiples of
$p(z)=z^{n-2}(z^2-a^2)$ when the first term in
\eqref{eq:section4-sharp} is maximal, or of
$p(z)=(z-(n-k)a)^k(z+ka)^{n-k}$ when $R_{n,p}(k)$ is maximal.
Here $a\in\mathbb C\setminus\{0\}$, and all maximizing families
are included when several terms in \eqref{eq:section4-sharp} tie.
The cases $n=2,3$ were established in
Proposition~\ref{sub12:prop:turn3-degrees-two-three}.
This completes the proof.
\end{proof}

\section{Matrix inequalities and weighted majorization}
\label{sec:matrix-inequalities}
\begin{proof}[Proof of Lemma~\ref{lem:weighted-convex}]
When \(\gamma=0\), either condition forces every \(x_i=0\), by taking
\(t=0\) or \(\Phi(u)=u\), respectively. Conversely, \(x=0\) satisfies
both conditions. Suppose now that \(\gamma>0\) and \(x\prec_{w,\gamma}y\), and put
$
        M=\max\{x_i,y_j:1\le i\le m,\ 1\le j\le n\}.
$
On \([0,M]\), every increasing convex \(\Phi\) with \(\Phi(0)=0\) has the
Stieltjes representation
\[
        \Phi(u)=au+\int_{(0,M]}(u-t)_+\,d\mu(t),
        \qquad 0\le u\le M,
\]
where \(a=\Phi'_+(0)\ge0\), and \(\mu\) is the positive measure induced by the
nondecreasing right derivative \(\Phi'_+\). The defining tail inequalities,
including the case \(t=0\), give
\[
\sum_{i=1}^{m}\Phi(x_i) =a\sum_{i=1}^{m}x_i+ \int_{(0,M]}\sum_{i=1}^{m}(x_i-t)_+\,d\mu(t) \le \gamma a\sum_{j=1}^{n}y_j+ \gamma\int_{(0,M]}\sum_{j=1}^{n}(y_j-t)_+\,d\mu(t) =\gamma\sum_{j=1}^{n}\Phi(y_j).
\]
Conversely, taking \(\Phi(u)=(u-t)_+\) for each fixed \(t\ge0\) recovers
\eqref{eq:wmaj-def}. This proves the equivalence.
\end{proof}

\begin{lem}\label{lem:singular-interlacing}
Let \(A=QDQ\), where $Q,D$ are defined by Lemma~\ref{lem:zpprimez}.
Then
\begin{equation*}
s_{j+2}(D)\le s_j(A)
\qquad 1\le j\le n-2,
\end{equation*}
and
\[
s_j(A)\le s_j(D),
\qquad 1\le j\le n.
\]
\end{lem}

\begin{proof}
The first inequality is vacuous for \(n=2\). For \(n\ge3\), observe that
\(D-A=D-QDQ=(I-Q)D+QD(I-Q)\).
Since \(I-Q=\frac1nJ\) has rank one,
both terms on the right have rank at most
one. Consequently,
\begin{equation}\label{eq:rank-two-difference}
\rank(D-A)\le2.
\end{equation}
Applying Weyl's inequality for singular values~\cite[p.~75, Problem~III.6.5]{Bha97} to
$
D=A+(D-A)$
with \(k=3\), we obtain
\(s_{j+2}(D) \le s_j(A)+s_3(D-A)\).
By \eqref{eq:rank-two-difference},
\(s_3(D-A)=0\). Hence
\[
s_{j+2}(D)\le s_j(A),
\qquad 1\le j\le n-2.
\]

For the second inequality, \(Q\) is an orthogonal projection, so
we have \(\|Q\|_\infty=1\).
By the standard singular-value inequalities
for matrix products (see \cite[p.~75, Problem~III.6.2]{Bha97}),
we obtain
\(s_j(A)=s_j(QDQ) \le \|Q\|_\infty s_j(D)\|Q\|_\infty=s_j(D)\),
for every \(1\le j\le n\).
\end{proof}

Let \(\mathbb M_n\) denote the set of all \(n\times n\) complex matrices,
and let \(\mathbb H_n\) denote the set of all \(n\times n\) Hermitian matrices.
For \(A\in\mathbb M_n\), let \(\tr A\) denote its trace.
For \(A,B\in\mathbb H_n\), we write \(A\le B\) if
\(B-A\) is positive semidefinite. This defines the
Loewner order on \(\mathbb H_n\).

For \(A\in\mathbb H_n\), define its positive and negative parts by
\[
A_+=\frac{|A|+A}{2},
\qquad
A_-=\frac{|A|-A}{2},
\]
where \(|A|=(A^2)^{\frac{1}{2}}\). Then
\(A=A_+-A_-\) and \(A_+,A_-\ge0\).

By Weyl's monotonicity principle, we immediately obtain the following lemma.

\begin{lem}\label{lem:positive-trace}
Let \(X, Y\in \mathbb H_n\).  If \(X\le Y\), then
$
        \tr X_+\le \tr Y_+.
$
\end{lem}

\begin{lem}[Weyl's Majorant Theorem
{\cite[p.~42, Theorem~II.3.6]{Bha97}}]
\label{lem:weyl-horn}
Let \(A\in \mathbb M_n\), and let
\(\lambda_1(A),\ldots,\lambda_n(A)\) and
\(s_1(A),\ldots,s_n(A)\) denote its eigenvalues
and singular values, respectively. Then
\[
(|\lambda_1(A)|,\ldots,|\lambda_n(A)|)
\prec_w
(s_1(A),\ldots,s_n(A)).
\]
\end{lem}

\begin{lem}\label{lem:hermitian-majorization}
Let \(A\in \mathbb{M}_n\), and fix \(\theta\in\R\).  Put
$
        R_\theta(A)={(e^{-i\theta}A+e^{i\theta}A^*)}/{2}.
$
If \(\lambda_1,\ldots,\lambda_n\) are the eigenvalues of \(A\), and
\(\nu_1,\ldots,\nu_n\) are the eigenvalues of the Hermitian matrix
\(R_\theta(A)\), then
\[
        (\Re(e^{-i\theta}\lambda_1),\ldots,\Re(e^{-i\theta}\lambda_n))
        \prec
        (\nu_1,\ldots,\nu_n).
\]
\end{lem}

\begin{proof}
Choose a unitary matrix \(U\in\mathbb M_n\) for which \(U^*AU\) is
upper triangular with diagonal entries \(\lambda_1,\ldots,\lambda_n\),
as in \cite[p.~5, Schur's triangularization theorem]{Bha97}.
Consequently,
\(U^*R_\theta(A)U =\frac{e^{-i\theta}U^*AU+e^{i\theta}U^*A^*U}{2}\)
is Hermitian with diagonal entries
$
\Re(e^{-i\theta}\lambda_1),\ldots,
\Re(e^{-i\theta}\lambda_n).
$
The diagonal of a Hermitian matrix is majorized by its eigenvalue
list, by Schur's inequality
\cite[p.~35, Exercise~II.1.12]{Bha97}.
\end{proof}

\section{Proof of Theorem~\ref{thm:main-modulus-wmaj}}\label{sec:modulus-comparisons}

When \(n=2\), the centroid condition in Theorem~\ref{thm:main-modulus-wmaj} implies that
the unique critical point is zero, so the result
is immediate. Henceforth, assume \(n\ge3\).

Write \(\mathbf z=(z_1,\ldots,z_n)^T\),
\(e=n^{-\frac{1}{2}}(1,\ldots,1)^T\), \(D=\diag(z_1,\ldots,z_n)\), and
\(Q=I-ee^*\), and \(A=QDQ\). Since \(\sum_{j=1}^n z_j=0\), the vector
\(\mathbf z\) lies in \(e^\perp\). By Lemma~\ref{lem:zpprimez}, the eigenvalues
of \(A\) are \(w_1,\ldots,w_{n-1},0\),
counted with algebraic multiplicity.
Weyl's Majorant Theorem (Lemma~\ref{lem:weyl-horn}) therefore gives
\((|w_1|,\ldots,|w_{n-1}|,0) \prec_w (s_1(A),\ldots,s_n(A))\). 
Since \(u\mapsto(u-t)_+\) is increasing
and convex on \([0,\infty)\), we obtain
\begin{equation}\label{eq:modulus-eigen-singular}
\sum_{k=1}^{n-1}(|w_k|-t)_+
\le \sum_{i=1}^n(s_i(A)-t)_+,
\qquad t\ge0.
\end{equation}

For \(t\ge0\), define
\[
F(\mathbf z,t)
=\sum_{i=1}^n(s_i(A)-t)_+,
\qquad
H(\mathbf z,t)
=\sum_{j=1}^n(|z_j|-t)_+.
\]
We first prove the bound at \(t=0\) for \(n\ge3\). For \(n=3\) it
will give the desired critical-point bound directly. For \(n\ge4\),
\eqref{eq:modulus-eigen-singular} reduces the problem to showing
\begin{equation}\label{eq:modulus-F-bound}
F(\mathbf z,t)
\le
\sqrt{\frac{n-2}{n}}\,H(\mathbf z,t),
\qquad t\ge0.
\end{equation}

\subsection{Step 1. The case \texorpdfstring{\(t=0\)}{t = 0}}
At \(t=0\) the required bound is the trace-norm estimate of
Tang and Zhang~\cite[Theorem~4.1]{TZ25}. In our notation it reads
\begin{equation}\label{eq:modulus-trace-endpoint}
\sum_{i=1}^ns_i(A)
\le
\sqrt{\frac{n-2}{n}}
\sum_{j=1}^n|z_j|.
\end{equation}

Since \(\sum_{j=1}^n z_j=0\), we have
\(QD^*e=D^*e\). Consequently,
\(A^*A =Q|D|^2Q-D^*ee^*D\).
Taking traces gives
\begin{equation}\label{eq:modulus-HS-identity}
\sum_{i=1}^n s_i(A)^2=\tr(A^*A)=\frac{n-2}{n}\sum_{j=1}^n|z_j|^2.
\end{equation}

Let \(e_1,\ldots,e_n\) be the standard
orthonormal basis of \(\C^n\), and put
\[
v_{ij}
=(\bigwedge\nolimits^2Q)(e_i\wedge e_j),
\qquad 1\le i<j\le n.
\]
By the multiplicativity of exterior powers,
\begin{equation}\label{eq:modulus-exterior-expansion}
\bigwedge\nolimits^2A
=
\sum_{i<j}z_iz_jv_{ij}v_{ij}^*.
\end{equation}
Moreover,
\[
\|v_{ij}\|^2=
\det
\begin{pmatrix}
1-\frac1n&-\frac1n\\
-\frac1n&1-\frac1n
\end{pmatrix}=\frac{n-2}{n}.
\]

The singular values of \(\bigwedge^2A\)
are precisely the pairwise products
\(s_i(A)s_j(A)\), \(i<j\). Hence, by
\eqref{eq:modulus-exterior-expansion}
and the trace-norm triangle inequality,
\[
e_2(s_1(A),\ldots,s_n(A))=\left\|\bigwedge\nolimits^2A\right\|_1\le
\sum_{i<j}|z_iz_j|\|v_{ij}v_{ij}^*\|_1=\frac{n-2}{n}
e_2(|z_1|,\ldots,|z_n|).
\]
Together with \eqref{eq:modulus-HS-identity},
this yields
\[
\begin{aligned}
(\sum_{i=1}^ns_i(A))^2
&=\sum_{i=1}^n s_i(A)^2
+2e_2(s_1(A),\ldots,s_n(A))\\
&\le
\frac{n-2}{n}
\left(
\sum_{j=1}^n|z_j|^2
+2e_2(|z_1|,\ldots,|z_n|)
\right)\\
&=
\frac{n-2}{n}
\left(\sum_{j=1}^n|z_j|\right)^2.
\end{aligned}
\]
Taking square roots proves
\eqref{eq:modulus-trace-endpoint}.

\subsection{Step 2. Cubics}
For \(n=3\) the critical points have the form \(w,-w\).
By \eqref{eq:modulus-eigen-singular}
and \eqref{eq:modulus-trace-endpoint},
\(2|w| \le\frac1{\sqrt3}\sum_{j=1}^3|z_j|\).

If \(t\ge |w|\), then \(2(|w|-t)_+=0\). For \(0\le t<|w|\),
using \((u-t)_+\ge u-t\) for \(u\ge0\),
we obtain
\[
\frac1{\sqrt3}
\sum_{j=1}^3(|z_j|-t)_+
\ge
\frac1{\sqrt3}\sum_{j=1}^3|z_j|
-\sqrt3\,t\ge 2|w|-\sqrt3\,t\ge 2(|w|-t),
\]
since \(t\ge0\) and \(\sqrt3\le2\).
Thus \(2(|w|-t)_+\le H(\mathbf z,t)/\sqrt3\) for all \(t\ge0\),
which proves Theorem~\ref{thm:main-modulus-wmaj} for cubics.
For the rest of this section assume \(n\ge4\).

\subsection{Step 3. Reduction to extreme points}
Lemma~\ref{lem:vari} gives a variational proof of the subadditivity of \(F\).
\begin{lem}\label{lem:vari}
For \(t\ge0\) and \(A=QDQ\),
\[
 F(\mathbf z,t)=\sup_{\|X\|_\infty\le1}
 \{\Re\tr(X^*A)-t\|X\|_1\}.
\]
\end{lem}
\begin{proof}
    For any \(X\in\mathbb M_n\) with
\(\|X\|_\infty\le1\), von Neumann's trace
inequality~\cite[p.~77, Problem~III.6.12]{Bha97} gives
\[
\Re\tr(X^*A)-t\|X\|_1
\le
\sum_{i=1}^n s_i(X)(s_i(A)-t)\le
\sum_{i=1}^n(s_i(A)-t)_+.
\]
Conversely, let
$
A=U\diag(s_1(A),\ldots,s_n(A))V^*
$
be a singular-value decomposition, and set
\[
X=U\diag(\varepsilon_1,\ldots,\varepsilon_n)V^*,
\qquad
\varepsilon_i=
\begin{cases}
1,&s_i(A)>t,\\
0,&s_i(A)\le t.
\end{cases}
\]
Then \(\|X\|_\infty\le1\) and
$
\Re\tr(X^*A)-t\|X\|_1
=
\sum_{i=1}^n(s_i(A)-t)_+.
$
This proves the variational lemma.
\end{proof}

\begin{rem}\label{rem:Fsub}
    It easily follows from Lemma~\ref{lem:vari} that \(F\) is positively
homogeneous and subadditive in the
joint variables \((\mathbf z,t)\).
In particular,
\begin{equation}\label{eq:modulus-subadditivity}
F(\mathbf z+\mathbf z',t+t')
\le
F(\mathbf z,t)+F(\mathbf z',t')
\end{equation}
for centered vectors \(\mathbf z,\mathbf z'\in \mathbb C^n\)
and \(t,t'\ge0\).
\end{rem}

\begin{defin}[Recession cone]
\label{def:recession-cone}
Let \(\mathscr C\) be a nonempty convex subset of
a finite-dimensional real vector space.
A vector \(v\) is called a \emph{recession direction}
of \(\mathscr C\) if
$
x+\lambda v\in\mathscr C
$
for every \(x\in\mathscr C\) and every \(\lambda\ge0\).
The set of all recession directions of
\(\mathscr C\) is called its \emph{recession cone},
denoted by
\[
\operatorname{rec}(\mathscr C)
=
\{v:x+\lambda v\in\mathscr C
\text{ for all }x\in\mathscr C,\ \lambda\ge0\}.
\]
\end{defin}

\begin{lem}[{\cite[p.~166, Theorem~18.5]{Roc70}}]
\label{lem:extreme-recession-representation}
Let \(\mathscr C\) be a nonempty closed convex
subset of a finite-dimensional real vector space.
Assume that \(\mathscr C\) contains no lines.
Then every \(x\in\mathscr C\) can be expressed as
\begin{equation}\label{eq:extreme-recession-representation}
x=\sum_{\nu=1}^N\alpha_\nu x^{(\nu)}+v,
\end{equation}
where \(x^{(1)},\ldots,x^{(N)}\) are extreme
points of \(\mathscr C\),
$
\alpha_\nu\ge0$, $
\sum_{\nu=1}^N\alpha_\nu=1,
$
and \(v\in\operatorname{rec}(\mathscr C)\).
\end{lem}
Consider the closed convex set
\begin{equation}\label{eq:mathcalC}
    \mathcal C=
\left\{
(\mathbf z,t)\in\C^n\times[0,\infty):
\sum_{j=1}^n z_j=0,\quad
H(\mathbf z,t)\le1
\right\}.
\end{equation}
The set \(\mathcal C\) contains no lines.
Indeed, assume that
\[
(\mathbf z+\lambda\mathbf v,t+\lambda\tau)
\in\mathcal C,
\qquad \lambda\in\R,
\]
for some \(\mathbf v\in\C^n\) and \(\tau\in\R\).
Since \(t+\lambda\tau\ge0\) for every
\(\lambda\in\R\), we must have \(\tau=0\).
Moreover, \(H(\mathbf z+\lambda\mathbf v,t)\le1\)
implies
\[
|z_j+\lambda v_j|\le t+1,
\qquad 1\le j\le n,\quad \lambda\in\R.
\]
Letting \(|\lambda|\to\infty\), we obtain
\(v_j=0\) for every \(j\).
Thus \(\mathbf v=0\), proving that
\(\mathcal C\) contains no lines.

\begin{lem}[Recession cone and vanishing of \(F\)]
\label{lem:recession-cone}
Let \(\mathcal C\) be the convex set defined by \eqref{eq:mathcalC}. Then
\begin{equation}\label{eq:recession-cone}
\operatorname{rec}(\mathcal C)
=
\left\{
(\mathbf v,\tau)\in\C^n\times\R:
\sum_{j=1}^n v_j=0,\quad
\tau\ge0,\quad
|v_j|\le\tau\ (1\le j\le n)
\right\}.
\end{equation}
Moreover,
$
F(\mathbf v,\tau)=0
$
for every $
(\mathbf v,\tau)\in\operatorname{rec}(\mathcal C).
$
\end{lem}

\begin{proof}
We first establish \eqref{eq:recession-cone}.
Assume that
\((\mathbf v,\tau)\in\operatorname{rec}(\mathcal C)\).
Since \((\mathbf 0,0)\in\mathcal C\), the definition
of the recession cone implies that
$
(\lambda\mathbf v,\lambda\tau)\in\mathcal C, \lambda\ge0.
$
In particular,
$
\sum_{j=1}^n v_j=0,
\tau\ge0.
$
Moreover, by the positive homogeneity of \(H\),
\(\lambda\sum_{j=1}^n(|v_j|-\tau)_+ =H(\lambda\mathbf v,\lambda\tau) \le1\)
for every \(\lambda>0\).
Letting \(\lambda\to\infty\), we obtain
\(\sum_{j=1}^n(|v_j|-\tau)_+=0\).
Since each summand is nonnegative, it follows that
\(|v_j|\le\tau\) for every \(1\le j\le n\).
This proves one inclusion in
\eqref{eq:recession-cone}.

Conversely, assume that
\[
\sum_{j=1}^n v_j=0,\qquad
\tau\ge0,\qquad
|v_j|\le\tau,\quad 1\le j\le n.
\]
Let \((\mathbf z,t)\in\mathcal C\) and
\(\lambda\ge0\). Clearly,
\[
\sum_{j=1}^n(z_j+\lambda v_j)=0,
\qquad
t+\lambda\tau\ge0.
\]
Furthermore, the triangle inequality gives
\[
|z_j+\lambda v_j|-(t+\lambda\tau)
\le |z_j|+\lambda|v_j|-t-\lambda\tau\le |z_j|-t.
\]
Taking positive parts and summing over \(j\),
we obtain
\[
H(\mathbf z+\lambda\mathbf v,t+\lambda\tau) =\sum_{j=1}^n \bigl(|z_j+\lambda v_j|-(t+\lambda\tau)\bigr)_+ \le\sum_{j=1}^n(|z_j|-t)_+ =H(\mathbf z,t)\le1.
\]
Thus
$
(\mathbf z+\lambda\mathbf v,t+\lambda\tau)
\in\mathcal C
$
for every \((\mathbf z,t)\in\mathcal C\)
and every \(\lambda\ge0\).
By the definition of the recession cone,
\((\mathbf v,\tau)\in\operatorname{rec}(\mathcal C)\).
This completes the proof of
\eqref{eq:recession-cone}.

Finally, let
\((\mathbf v,\tau)\in\operatorname{rec}(\mathcal C)\).
Since \(Q\) is an orthogonal projection,
by \eqref{eq:recession-cone},
\[
\|Q\diag(v_1,\ldots,v_n)Q\|_\infty
\le\|\diag(v_1,\ldots,v_n)\|_\infty=\max_{1\le j\le n}|v_j|\le\tau.
\]
Consequently, every singular value of
\(Q\diag(v_1,\ldots,v_n)Q\) is at most
\(\tau\). It follows from the definition
of \(F\) that
\(F(\mathbf v,\tau) =\sum_{i=1}^n \bigl(s_i(Q\diag(v_1,\ldots,v_n)Q)-\tau\bigr)_+=0\).
This completes the proof.
\end{proof}

By Lemmas~\ref{lem:extreme-recession-representation}
and \ref{lem:recession-cone}, every
\((\mathbf z,t)\in\mathcal C\) admits a representation
\((\mathbf z,t) = \sum_{\nu=1}^N \alpha_\nu(\mathbf z^{(\nu)},t_\nu) +(\mathbf v,\tau)\),
where \((\mathbf z^{(\nu)},t_\nu)\) are extreme
points of \(\mathcal C\), with
\(\alpha_\nu\ge0\) and \(\sum_{\nu=1}^N\alpha_\nu=1\),
and \((\mathbf v,\tau)\in\operatorname{rec}(\mathcal C)\).
By the subadditivity and positive homogeneity
of \(F\) (Remark~\ref{rem:Fsub}), together with
Lemma~\ref{lem:recession-cone}, we have
\[
F(\mathbf z,t)
\le
\sum_{\nu=1}^N
\alpha_\nu F(\mathbf z^{(\nu)},t_\nu)
+F(\mathbf v,\tau)=
\sum_{\nu=1}^N
\alpha_\nu F(\mathbf z^{(\nu)},t_\nu).
\]
Consequently, it suffices to prove
\begin{equation}\label{eq:modulus-extreme-bound}
F(\mathbf z,t)
\le\sqrt{\frac{n-2}{n}}
\end{equation}
at every extreme point of \(\mathcal C\).

Let \((\mathbf z,t)\) be such an extreme point.
If \(t=0\), then \eqref{eq:modulus-trace-endpoint}
and \(H(\mathbf z,0)\le1\) give
\(F(\mathbf z,0) \le \sqrt{\frac{n-2}{n}}H(\mathbf z,0) \le \sqrt{\frac{n-2}{n}}\).
We may therefore assume that \(t>0\).

\begin{lem}\label{lem:extreme-H-equality}
Let \(\mathcal C\) be the convex set defined by
\eqref{eq:mathcalC}. If \((\mathbf z,t)\) is an
extreme point of \(\mathcal C\) with \(t>0\), then
$
H(\mathbf z,t)=1.
$
\end{lem}

\begin{proof}
Suppose, to the contrary, that
\(H(\mathbf z,t)<1\). Since \(H\) is continuous
and \(t>0\), there exists \(0<\varepsilon<t\)
such that
\[
H(\mathbf z,t-\varepsilon)<1,
\qquad
H(\mathbf z,t+\varepsilon)<1.
\]
Both \((\mathbf z,t-\varepsilon)\) and
\((\mathbf z,t+\varepsilon)\) therefore belong
to \(\mathcal C\). However,
\((\mathbf z,t) = \frac12(\mathbf z,t-\varepsilon) +\frac12(\mathbf z,t+\varepsilon)\).
The endpoints of this segment are distinct, contradicting the extremality of
\((\mathbf z,t)\) in \(\mathcal C\).
\end{proof}

\begin{lem}\label{lem:extreme-dimension}
Let \(\mathcal C\) be the convex set defined by
\eqref{eq:mathcalC}, and let \((\mathbf z,t)\)
be an extreme point of \(\mathcal C\) with \(t>0\).
Let \(r\) and \(s\) denote the numbers of components
of \(\mathbf z\) whose moduli are greater than
and less than \(t\), respectively. That is,
\(r=\#\{1\le j\le n:|z_j|>t\}\) and
\(s=\#\{1\le j\le n:|z_j|<t\}\).
 Then
$
(r,s)=(1,0)
$ or $(2,0)$
\end{lem}

\begin{proof}
By Lemma~\ref{lem:extreme-H-equality},
we have \(H(\mathbf z,t)=1\), so \(r\ge1\).
We shall prove that \(r+2s<3\).

Rearranging the components of \(\mathbf z\),
we may assume that
\[
\begin{cases}
|z_j|>t, &1\le j\le r,\\
|z_j|=t, &r+1\le j\le n-s,\\
|z_j|<t, &n-s+1\le j\le n.
\end{cases}
\]
For \(1\le j\le n-s\), write
\(z_j=|z_j|e^{i\theta_j}\).

Assume that \(r+2s\ge3\).
We shall show that there exists a nonzero vector
\((\mathbf v,\tau)\in\C^n\times\R\)
such that both
\[
(\mathbf z+\varepsilon\mathbf v,t+\varepsilon\tau)
\quad\text{and}\quad
(\mathbf z-\varepsilon\mathbf v,t-\varepsilon\tau)
\]
belong to \(\mathcal C\) for all sufficiently
small \(\varepsilon>0\).
This would contradict the extremality of
\((\mathbf z,t)\).

By the definition of \(\mathcal C\), we only
need to verify that
\[
\sum_{j=1}^n(z_j\pm\varepsilon v_j)=0,
\qquad
t\pm\varepsilon\tau\ge0,\qquad\text{and}\qquad
H(\mathbf z\pm\varepsilon\mathbf v,
t\pm\varepsilon\tau)\le1.
\]

Since \(\sum_{j=1}^n z_j=0\), the first
condition is equivalent to
\begin{equation}\label{eq:perturb-centroid}
\sum_{j=1}^n v_j=0.
\end{equation}
The second condition holds automatically
for sufficiently small \(\varepsilon>0\),
because \(t>0\).

To deal with the last condition, choose
\[
v_j=
\begin{cases}
b_je^{i\theta_j},
&1\le j\le r,\\
\tau e^{i\theta_j},
&r+1\le j\le n-s,\\
x_j+iy_j,
&n-s+1\le j\le n,
\end{cases}
\]
where \(b_j,\tau,x_j,y_j\) are real numbers
to be determined.

For sufficiently small \(\varepsilon>0\),
the components remain on the same side
of the perturbed threshold. Therefore,
\[
H(\mathbf z+\varepsilon\mathbf v, t+\varepsilon\tau)= \sum_{j=1}^r \bigl(|z_j|-t+\varepsilon(b_j-\tau)\bigr)= 1+\varepsilon \left(\sum_{j=1}^r b_j-r\tau\right),
\]
and similarly,
\(H(\mathbf z-\varepsilon\mathbf v, t-\varepsilon\tau)= 1-\varepsilon \left(\sum_{j=1}^r b_j-r\tau\right)\).
Thus both quantities are at most \(1\)
if and only if
\begin{equation}\label{eq:perturb-tail}
\sum_{j=1}^r b_j-r\tau=0.
\end{equation}

Substituting the chosen components of
\(\mathbf v\) into
\eqref{eq:perturb-centroid} and separating
real and imaginary parts, we see that
\eqref{eq:perturb-centroid} and
\eqref{eq:perturb-tail} are equivalent
to the following system:
\begin{equation}\label{eq:extreme-linear-system}
\begin{cases}
\displaystyle
\sum_{j=1}^r b_j\cos\theta_j
+\tau\sum_{j=r+1}^{n-s}\cos\theta_j
+\sum_{j=n-s+1}^n x_j=0,\\[8pt]
\displaystyle
\sum_{j=1}^r b_j\sin\theta_j
+\tau\sum_{j=r+1}^{n-s}\sin\theta_j
+\sum_{j=n-s+1}^n y_j=0,\\[8pt]
\displaystyle
\sum_{j=1}^r b_j-r\tau=0.
\end{cases}
\end{equation}

This is a homogeneous system of three real
linear equations in \(1+r+2s\) real unknowns,
namely
\[
b_1,\ldots,b_r,\quad
\tau,\quad
x_{n-s+1},y_{n-s+1},\ldots,x_n,y_n.
\]
Since \(r+2s\ge3\), the number of unknowns
is at least four. By the rank-nullity theorem,
\eqref{eq:extreme-linear-system} admits
a nonzero solution.

Fix such a solution and let
\((\mathbf v,\tau)\) be the corresponding
perturbation direction. Since \(t>0\),
the strict inequalities defining the
components with \(|z_j|>t\) and \(|z_j|<t\)
remain valid for sufficiently small
\(\varepsilon>0\). The components with
\(|z_j|=t\) continue to have modulus
\(t\pm\varepsilon\tau\).

By \eqref{eq:perturb-centroid}, both
perturbed configurations have centroid zero.
Moreover, \eqref{eq:perturb-tail} gives
\(H(\mathbf z\pm\varepsilon\mathbf v, t\pm\varepsilon\tau)=1\).
Both perturbed points belong to \(\mathcal C\), since
\(t\pm\varepsilon\tau>0\). The perturbation is nonzero and the midpoint
is \((\mathbf z,t)\), contradicting its extremality.

Hence \(r+2s<3\). Since \(r\ge1\)
and \(r,s\) are nonnegative integers,
the only possibilities are
\[
(r,s)=(1,0)
\qquad\text{or}\qquad
(r,s)=(2,0).
\]
This completes the proof.
\end{proof}

\subsection{Step 4. One modulus above the threshold}
\label{subsec:modulus-one-excess}

Assume that $|z_1|=t+d$ and $|z_2|=\cdots=|z_n|=t$, where
$d>0$. By Lemma~\ref{lem:singular-interlacing}, the singular values
of $A=QDQ$ are $t+L$, $t$ with multiplicity $n-3$, $a$, and $0$,
where $L\ge0$ and $0\le a\le t$. Thus $F(\mathbf z,t)=L$ and
$H(\mathbf z,t)=d$. The Hilbert--Schmidt identity
\eqref{eq:modulus-HS-identity} gives
$2tL+L^2+a^2=\frac{n-2}{n}(2td+d^2)$, and hence
$2tL+L^2\le 2t\sqrt{\frac{n-2}{n}}\,d+\frac{n-2}{n}d^2$.
Since $s\mapsto2ts+s^2$ is increasing on $[0,\infty)$, we obtain
\begin{equation}\label{eq:modulus-one-excess}
F(\mathbf z,t)\le\sqrt{\frac{n-2}{n}}\,H(\mathbf z,t).
\end{equation}

\subsection{Step 5. Two moduli above the threshold}
\label{subsec:modulus-two-excess}

Assume that $|z_1|=t+d_1$, $|z_2|=t+d_2$, and
$|z_3|=\cdots=|z_n|=t$, where $d_1,d_2>0$.
Then $H(\mathbf z,t)=d_1+d_2$.

If at most one singular value of $A$ exceeds $t$,
Lemma~\ref{lem:singular-interlacing} shows that its singular values
are $t+L$, $t$ with multiplicity $n-3$, $a$, and $0$, where
$L\ge0$ and $0\le a\le t$. By \eqref{eq:modulus-HS-identity},
$2tL+L^2+a^2=\frac{n-2}{n}\bigl(2t(d_1+d_2)+d_1^2+d_2^2\bigr)$.
Since $d_1^2+d_2^2\le(d_1+d_2)^2$, the monotonicity argument in
Subsection~\ref{subsec:modulus-one-excess} gives
$F(\mathbf z,t)=L\le\sqrt{\frac{n-2}{n}}(d_1+d_2)$.

It remains to treat the case in which exactly two singular values
of $A$ exceed $t$. Write $z_j=(t+d_j)u_j$ with $|u_j|=1$ for
$j=1,2$, and put $\Delta_j=d_j(2t+d_j)>0$. Since
$\sum_{j=1}^nz_j=0$, we have $QD^*e=D^*e$. Consequently,
\begin{equation}\label{eq:modulus-rank-three}
A^*A=t^2Q+\sum_{j=1}^{2}\Delta_j(Qe_j)(Qe_j)^*-D^*ee^*D.
\end{equation}
Here $e=n^{-\frac{1}{2}}(1,\ldots,1)^T$, as in the definition of $A$.
The difference $A^*A-t^2Q$ has rank at most three. On
$\mathbf1^\perp$, the eigenvalue $t^2$ therefore has multiplicity
at least $n-4$. Denote the remaining three eigenvalues of $A^*A$,
counted with multiplicity, by $y_1\ge y_2\ge y_3\ge0$.
Lemma~\ref{lem:singular-interlacing} gives $y_3\le t^2$, while
our assumption implies $y_2>t^2$. Hence
\begin{equation}\label{eq:modulus-F-three}
F(\mathbf z,t)=\sqrt{y_1}+\sqrt{y_2}-2t.
\end{equation}

We next compare the singular values for different phases of
$z_1$ and $z_2$, keeping all the moduli fixed. The Gram matrix
of $Qe_1,Qe_2,D^*e$ is
\[
\begin{pmatrix}
1-\frac1n&-\frac1n&\frac{\overline{z_1}}{\sqrt n}\\
-\frac1n&1-\frac1n&\frac{\overline{z_2}}{\sqrt n}\\
\frac{z_1}{\sqrt n}&\frac{z_2}{\sqrt n}&\frac1n\sum_{j=1}^n|z_j|^2
\end{pmatrix}.
\]
The three eigenvalues of the rank-three perturbation in
\eqref{eq:modulus-rank-three}, supplemented with zeros if necessary,
coincide with the eigenvalues of $\operatorname{diag}(\Delta_1,\Delta_2,-1)$
times this Gram matrix. Its trace and the sum of its principal
minors of order two depend only on the moduli. The only
phase-dependent term in its determinant is
$2\Delta_1\Delta_2\Re(z_1\overline{z_2})/n^2$.
It follows, upon adding $t^2$ to each eigenvalue, that
$S:=y_1+y_2+y_3$ and $E:=y_1y_2+y_1y_3+y_2y_3$ depend only
on the moduli, and that
\begin{equation}\label{eq:modulus-phase-product}
y_1y_2y_3=C_0+\frac{2\Delta_1\Delta_2}{n^2}\Re(z_1\overline{z_2}),
\end{equation}
where $C_0$ also depends only on the moduli.
Since $\Delta_1\Delta_2>0$, the product in
\eqref{eq:modulus-phase-product} is smallest when $u_1=-u_2$.

We show that a smaller value of $y_1y_2y_3$, with $S$ and $E$
fixed, cannot decrease $\sqrt{y_1}+\sqrt{y_2}$. Write $x=y_3$.
Then $y_1+y_2=S-x$ and $y_1y_2=E-Sx+x^2$, so
\[
y_1y_2y_3=x(E-Sx+x^2),\qquad
\frac{d}{dx}\bigl(x(E-Sx+x^2)\bigr)=(y_1-x)(y_2-x)\ge0.
\]
Thus decreasing the product cannot increase the smallest eigenvalue
$x$. Furthermore,
\begin{equation}\label{eq:modulus-two-root-sum}
(\sqrt{y_1}+\sqrt{y_2})^2=S-x+2\sqrt{E-Sx+x^2}.
\end{equation}
The right-hand side of \eqref{eq:modulus-two-root-sum} is
nonincreasing in $x$, because its derivative, where defined, is
$-1+(2x-S)/\sqrt{E-Sx+x^2}\le0$.
Also,
$y_2=\bigl(S-x-\sqrt{S^2+2Sx-3x^2-4E}\bigr)/2$
is nonincreasing in $x$: its derivative is nonpositive since
$S-3x=y_1+y_2-2x\ge0$. The corresponding assertions at repeated
or zero eigenvalues follow by continuity. Therefore, decreasing
$y_1y_2y_3$ cannot decrease $\sqrt{y_1}+\sqrt{y_2}$ or cause
$y_2$ to fall below $t^2$. By \eqref{eq:modulus-F-three}, choosing
opposite phases for $z_1,z_2$ cannot decrease $F(\mathbf z,t)$.

We now check that opposite phases are compatible with the centroid
condition. Since $\sum_{j=1}^nz_j=0$ and $|z_j|=t$ for $j\ge3$,
$|z_1+z_2|\le(n-2)t$. The reverse triangle inequality gives
$|d_1-d_2|\le|z_1+z_2|\le(n-2)t$.
For \(r\ge2\) the sums of \(r\) complex numbers of modulus \(t\)
fill the disk of radius \(rt\). For \(r=2\) this follows by varying
the angle between the summands. Adding one summand to a disk of
radius \(rt\) fills the disk of radius \((r+1)t\), since each point
in the disk of radius \((r+1)t\) has distance at most \(rt\) from the circle of radius
\(t\). Apply this with \(r=n-2\).
Consequently, the phases of $z_3,\ldots,z_n$ can be chosen so that
the root list remains centered after setting $u_1=-u_2$.
We may therefore assume
$z_1=(t+d_1)u$ and $z_2=-(t+d_2)u$ for some $|u|=1$.

Interchanging the first two zeros if necessary, assume $d_1\ge d_2$.
We then have the decomposition
\begin{equation}\label{eq:modulus-two-split}
\mathbf z=(d_2u,-d_2u,0,\ldots,0)+((t+d_1-d_2)u,-tu,z_3,\ldots,z_n).
\end{equation}
Both vectors on the right of \eqref{eq:modulus-two-split} are centered.
The first satisfies the trace-norm bound
\eqref{eq:modulus-trace-endpoint} at threshold zero. If \(d_1>d_2\),
the second vector satisfies \eqref{eq:modulus-one-excess}. If
\(d_1=d_2\), its largest modulus is \(t\), so its compression has
operator norm at most \(t\) and its contribution is zero.
Subadditivity~\eqref{eq:modulus-subadditivity} gives
\[
\begin{aligned}
F(\mathbf z,t)
&\le F((d_2u,-d_2u,0,\ldots,0),0)
 +F(((t+d_1-d_2)u,-tu,z_3,\ldots,z_n),t)\\
&\le 2\sqrt{\frac{n-2}{n}}\,d_2
 +\sqrt{\frac{n-2}{n}}(d_1-d_2)
 =\sqrt{\frac{n-2}{n}}H(\mathbf z,t).
\end{aligned}
\]
Thus \eqref{eq:modulus-extreme-bound} holds for the case
$(r,s)=(2,0)$ in Lemma~\ref{lem:extreme-dimension}.

\subsection{Step 6. Completion of the proof}
\label{subsec:modulus-completion}

The extreme-point classification in Lemma~\ref{lem:extreme-dimension}
and the bounds proved in Subsections~\ref{subsec:modulus-one-excess}
and~\ref{subsec:modulus-two-excess} establish
\eqref{eq:modulus-extreme-bound} for every extreme point of
$\mathcal C$ with $t>0$. For $t=0$, the same conclusion follows
from \eqref{eq:modulus-trace-endpoint} and the definition
\eqref{eq:mathcalC}. The decomposition
\eqref{eq:extreme-recession-representation}, the vanishing of \(F\)
on the recession cone from Lemma~\ref{lem:recession-cone},
and subadditivity \eqref{eq:modulus-subadditivity} therefore give
$F(\mathbf z,t)\le\sqrt{(n-2)/n}$ for every
$(\mathbf z,t)\in\mathcal C$.

Now let $\mathbf z$ be any centered vector and $t\ge0$. If
$H(\mathbf z,t)>0$, then
$(\mathbf z/H(\mathbf z,t),t/H(\mathbf z,t))\in\mathcal C$.
Positive homogeneity yields \eqref{eq:modulus-F-bound}. If
$H(\mathbf z,t)=0$, then $|z_j|\le t$ for all $j$, whence
$\|A\|_\infty\le\|D\|_\infty\le t$ and $F(\mathbf z,t)=0$.
Thus \eqref{eq:modulus-F-bound} holds for every centered root list
and every $t\ge0$.

Combining \eqref{eq:modulus-F-bound} with
\eqref{eq:modulus-eigen-singular} gives
\[
\sum_{k=1}^{n-1}(|w_k|-t)_+
\le\sqrt{\frac{n-2}{n}}\sum_{j=1}^n(|z_j|-t)_+,
\qquad t\ge0.
\]
By Definition~\ref{def:weighted-majorization}, the assertion of
Theorem~\ref{thm:main-modulus-wmaj} follows for $n\ge4$.
The cases $n=2$ and $n=3$, treated before the extreme-point reduction,
complete the proof.

\section{Proof of Theorems~\ref{thm:centered-superquadratic},
  \ref{thm:main-square-wmaj}, \ref{thm:main-directional-first},
  and~\ref{thm:main-directional-square}}
\label{sec:square-directional-proofs}

\subsection{Squared moduli and power sums}
\begin{proof}[{Proof of Theorem~\ref{thm:main-square-wmaj}}]
For $n=2$, the only critical point of a centered quadratic is zero.
Assume $n\ge3$, and put $e=n^{-\frac{1}{2}}(1,\ldots,1)^T$,
$Q=I-ee^*$, $D=\diag(z_1,\ldots,z_n)$, and
$A=(QDQ)|_{e^\perp}$. By Lemma~\ref{lem:zpprimez}, the
spectrum of $A$ consists of $w_1,\ldots,w_{n-1}$, counted with
multiplicity. Weyl's majorant theorem (Lemma~\ref{lem:weyl-horn})
and the increasing convex function $s\mapsto(s^2-t)_+$ on
$[0,\infty)$ reduce the proof to
\begin{equation}\label{eq:square-singular-tail}
\sum_{i=1}^{n-1}(s_i(A)^2-t)_+\le\frac{n-2}{n}\sum_{j=1}^n(|z_j|^2-t)_+,\qquad t\ge0.
\end{equation}

The centroid condition implies $D^*e\in e^\perp$. Thus, on
$e^\perp$,
\begin{equation}\label{eq:square-compression-identity}
A^*A=(Q|D|^2Q)|_{e^\perp}-(D^*e)(D^*e)^*,\qquad
\tr(A^*A)=\frac{n-2}{n}\sum_{j=1}^n|z_j|^2.
\end{equation}
The trace identity proves \eqref{eq:square-singular-tail} when
$t=0$. Fix $t>0$. The case $z_1=\cdots=z_n=0$ is immediate.
Otherwise, let $b_j=(|z_j|^2-t)_+$,
$r_j=\min\{|z_j|^2,t\}$, and
$\rho=n^{-1}\sum_{j=1}^n r_j>0$. In particular,
$|z_j|^2=b_j+r_j$ and $\rho\le t$.
Since $|z_j|^2-t\le b_j$, equation
\eqref{eq:square-compression-identity} gives
\begin{equation}\label{eq:square-rank-one-bound}
A^*A-tI_{e^\perp}\le T:=(Q\diag(b_1,\ldots,b_n)Q)|_{e^\perp}-(D^*e)(D^*e)^*.
\end{equation}
A direct trace calculation yields
$\tr T=\frac{n-2}{n}\sum_{j=1}^n b_j-\rho$.

We claim that $T\ge-\rho I_{e^\perp}$. If $b_j=0$, then
$|z_j|^2/(b_j+\rho)=r_j/\rho$. If $b_j>0$, then
$r_j=t\ge\rho$, so
$|z_j|^2/(b_j+\rho)=1+(r_j-\rho)/(b_j+\rho)\le r_j/\rho$.
Consequently, $\sum_j|z_j|^2/(b_j+\rho)\le n$.
For $h\in e^\perp$, the weighted Cauchy--Schwarz inequality gives
\[
|\langle D^*e,h\rangle|^2\le\frac1n\left(\sum_{j=1}^n\frac{|z_j|^2}{b_j+\rho}\right)
\left(\sum_{j=1}^n(b_j+\rho)|h_j|^2\right)
\le\sum_{j=1}^n(b_j+\rho)|h_j|^2.
\]
The definition of $T$ in \eqref{eq:square-rank-one-bound}
therefore implies $\langle(T+\rho I_{e^\perp})h,h\rangle\ge0$.

The operator $T$ has at most one negative eigenvalue, because it
is a rank-one negative perturbation of a positive semidefinite
operator. Since $T\ge-\rho I_{e^\perp}$, we have
$\tr(T_-)\le\rho$. Applying Lemma~\ref{lem:positive-trace}
to \eqref{eq:square-rank-one-bound}, we obtain
\[
\sum_{i=1}^{n-1}(s_i(A)^2-t)_+
\le\tr(T_+)=\tr T+\tr(T_-)
\le\frac{n-2}{n}\sum_{j=1}^n b_j.
\]
Thus \eqref{eq:square-singular-tail} holds for every $t\ge0$.
Finally, Lemma~\ref{lem:weyl-horn} gives the first inequality in
\begin{equation}\label{eq:square-critical-tail}
\sum_{k=1}^{n-1}(|w_k|^2-t)_+
\le\sum_{i=1}^{n-1}(s_i(A)^2-t)_+
\le\frac{n-2}{n}\sum_{j=1}^n(|z_j|^2-t)_+.
\end{equation}
Since \eqref{eq:square-critical-tail} holds for all $t\ge0$,
Definition~\ref{def:weighted-majorization} proves the theorem.
\end{proof}
\begin{proof}[Proof of Theorem~\ref{thm:centered-superquadratic}]
For \(n=2\), the zeros are opposite and the sole critical point is zero. For \(n\ge3\), apply
Lemma~\ref{lem:weighted-convex} to Theorem~\ref{thm:main-square-wmaj}
with \(\Phi(u)=u^{\frac{p}{2}}\). This proves the asserted bound. For \(p=2\),
the equality condition is that of \eqref{eq:schoenberg}.

Suppose \(p>2\) and the roots are not all zero. Put \(q=p/2\).
Since \(u^q=q(q-1)\int_0^\infty(u-t)_+t^{q-2}\,dt\), equality forces
the continuous nonnegative difference between the two outer terms of
\eqref{eq:square-critical-tail} to vanish for every \(t\ge0\).
At \(t=0\) this is equality in \eqref{eq:schoenberg}, so we may rotate
the roots onto the real axis. Write \(M=\max_j|z_j|>0\).
Let \(r\) count roots of modulus \(M\), with multiplicity, and let \(d\)
be the number of distinct such roots. By strict interlacing between
distinct real root values, exactly \(r-d\) critical points have modulus
\(M\). Taking \(t<M^2\) sufficiently close to \(M^2\) gives
\(r-d=(n-2)r/n\), so \(r=nd/2\). If \(d=1\), half the roots
equal one endpoint and the other half have strictly smaller modulus,
contrary to centering. Thus \(d=2\), \(r=n\), and centering gives
equal multiplicities at \(M\) and \(-M\). Conversely, for the roots
of \(c(z^2-a^2)^{\frac{n}{2}}\), the critical points are \(0\) and
\(\pm a\), the latter each with multiplicity \(n/2-1\).
This gives equality for every \(p\ge2\).
\end{proof}

\begin{prop}\label{prop:cubic-sharp}
For centered cubics and \(p>2\), the least constant is
\[
 C_{3,p}=
 \begin{cases}
  2/(2^p+2),&2<p<4,\\
  1/9,&p=4,\\
3^{-\frac{p}{2}},&p>4.
 \end{cases}
\]
Up to permutation and a common nonzero complex factor, the nonzero
extremal root lists are \((1,1,-2)\) for \(2<p<4\) and
\((1,-1,0)\) for \(p>4\). At \(p=4\) they are all centered
collinear lists.
\end{prop}
\begin{proof}
The critical points are \(\pm w\), where \(6w^2=\sum_jz_j^2\).
Rotate so that this sum is real and nonnegative and write
\(z_j=x_j+i y_j\). Then \(6|w|^2=\sum_j x_j^2-\sum_j y_j^2
\le\sum_j x_j^2\) and \(\sum_j|z_j|^p\ge\sum_j|x_j|^p\).
Since \(\sum_jx_j=0\), it suffices to consider real roots.
For a nonzero extremizer the first comparison is strict unless
all \(y_j=0\).

Fix \(S=\sum_j x_j^2>0\) and put \(u_j=x_j^2\) and
\(R=(x_1x_2x_3)^2\). Since \(\sum_jx_j=0\), we have
\(\sum_{i<j}u_i u_j=(\sum_{i<j}x_ix_j)^2=S^2/4\).
Thus the \(u_j\) solve \(t^3-St^2+(S^2/4)t-R=0\).
The discriminant of the cubic with roots \(x_j\) is
\(S^3/2-27R\). Hence \(0\le R\le S^3/54\), with extremal
endpoints given respectively by a zero coordinate and by a repeated
coordinate. Every intermediate \(R\) occurs by continuity on the
circle \(\sum x_j=0,\ \sum x_j^2=S\).

Put \(q=p/2\). In the interior the \(u_j\) are positive and
distinct, and implicit differentiation gives
\(\frac{d}{dR}\sum_{j=1}^3u_j^q =q\sum_{j=1}^3\frac{u_j^{q-1}}{\prod_{i\ne j}(u_j-u_i)}\).
The sum on the right is the second divided difference of \(t^{q-1}\), so
it is negative for \(1<q<2\), zero for \(q=2\), and positive for
\(q>2\), by the sign of the second derivative. The numerator
\(2(S/6)^q\) is fixed. The minimizing denominator is therefore
attained at the repeated-root endpoint for \(2<p<4\), at every list
for \(p=4\), and at the zero-coordinate endpoint for \(p>4\).
Substitution of \((1,1,-2)\) and \((1,-1,0)\) gives the stated constants.
\end{proof}

\subsection{Directional comparisons}
\begin{proof}[{Proof of Theorem~\ref{thm:main-directional-first}}]
The assertion is immediate for $n=2$. Assume $n\ge3$. We first
prove a tail inequality for the critical points of a centered
real-rooted polynomial. Let $x_1,\ldots,x_n\in\mathbb R$ satisfy
$\sum_jx_j=0$, put $e=n^{-\frac{1}{2}}(1,\ldots,1)^T$ and $Q=I-ee^*$,
and write $q_j=Qe_j$, where $e_1,\ldots,e_n$ are the standard
coordinate vectors. The self-adjoint compression
$B=(Q\diag(x_1,\ldots,x_n)Q)|_{e^\perp}$ has the critical
points of $\prod_j(z-x_j)$ as its eigenvalues, by
Lemma~\ref{lem:zpprimez}.

Fix a positive contraction $S$ on $e^\perp$, and put
$h_j=\langle Sq_j,q_j\rangle$. For distinct $j,k$, the Hermitian
operator $q_jq_j^*-q_kq_k^*$ has trace zero and two nonzero
eigenvalues $\pm\sqrt{(n-2)/n}$, because
$\|q_j\|^2=(n-1)/n$ and $\langle q_j,q_k\rangle=-1/n$.
The variational characterization of the positive and negative
parts therefore gives
$|h_j-h_k|\le\sqrt{(n-2)/n}$.

Set $a=\min_j h_j$ and $g_j=h_j-a$. Then
$0\le g_j\le\sqrt{(n-2)/n}$, and
$\sum_jg_j=\tr S-na\le\tr S$, since
$\sum_jq_jq_j^*=I_{e^\perp}$. For $t\ge0$, centering yields
\[
\tr\bigl(S(B-tI_{e^\perp})\bigr)
=\sum_j(x_j-t)g_j-t\bigl(\tr S-\sum_jg_j\bigr)
\le\sqrt{\frac{n-2}{n}}\sum_j(x_j-t)_+.
\]
Maximizing over positive contractions $S$ gives
$\tr\bigl((B-tI_{e^\perp})_+\bigr)\le\sqrt{(n-2)/n}\sum_j(x_j-t)_+$.
Applying the same argument to $-B$ and adding yields
\begin{equation}\label{eq:directional-real-tail}
\sum_{k=1}^{n-1}(|\eta_k|-t)_+
\le\sqrt{\frac{n-2}{n}}\sum_{j=1}^n(|x_j|-t)_+,\qquad t\ge0,
\end{equation}
where $\eta_1,\ldots,\eta_{n-1}$ are the critical points of
$\prod_j(z-x_j)$.

Fix $\theta\in\mathbb R$, and apply
\eqref{eq:directional-real-tail} with
$x_j=L_\theta(z_j)=\Re(e^{-i\theta}z_j)$. The projected zeros
satisfy $\sum_jx_j=0$. By
Lemma~\ref{lem:pereira-projection}, the real parts of the
rotated critical points are majorized by
$\eta_1,\ldots,\eta_{n-1}$. Applying the convex function
$s\mapsto(|s|-t)_+$ to that majorization and then using
\eqref{eq:directional-real-tail}, we obtain
\begin{equation}\label{eq:directional-first-tail}
\sum_{k=1}^{n-1}(|L_\theta(w_k)|-t)_+
\le\sum_{k=1}^{n-1}(|\eta_k|-t)_+
\le\sqrt{\frac{n-2}{n}}\sum_{j=1}^n(|L_\theta(z_j)|-t)_+.
\end{equation}
The parameter $t\ge0$ and direction $\theta\in\mathbb R$ are
arbitrary. Equation~\eqref{eq:directional-first-tail} is exactly
the assertion of Theorem~\ref{thm:main-directional-first} under
Definition~\ref{def:weighted-majorization}.
\end{proof}

\begin{proof}[{Proof of Theorem~\ref{thm:main-directional-square}}]
Fix $\theta\in\mathbb R$, and let
$x_j=L_\theta(z_j)=\Re(e^{-i\theta}z_j)$.
The centroid condition gives $\sum_{j=1}^n x_j=0$. Let
$\eta_1,\ldots,\eta_{n-1}$ be the critical points of the real
polynomial $\prod_{j=1}^n(z-x_j)$.

Lemma~\ref{lem:pereira-projection} shows that the vector of
projected critical points $(L_\theta(w_k))_{k=1}^{n-1}$ is
majorized by $(\eta_k)_{k=1}^{n-1}$. Since
$s\mapsto(s^2-t)_+$ is convex on $\mathbb R$ for every $t\ge0$,
Pereira's majorization gives an upper bound in terms of the
$\eta_k$. Theorem~\ref{thm:main-square-wmaj}, applied to the
centered real polynomial $\prod_{j=1}^n(z-x_j)$, bounds the
resulting sum by the squared-modulus tails of the $x_j$.
Consequently,
\begin{equation}\label{eq:directional-square-tail}
\sum_{k=1}^{n-1}(|L_\theta(w_k)|^2-t)_+
\le\frac{n-2}{n}\sum_{j=1}^n(|L_\theta(z_j)|^2-t)_+,\qquad t\ge0.
\end{equation}
Since $\theta$ is arbitrary, Definition~\ref{def:weighted-majorization}
shows that \eqref{eq:directional-square-tail} is equivalent to the
claimed weighted weak majorization for every direction.
\end{proof}

\section{Proof of Theorem~\ref{thm:centered-odd}}
\label{sec:odd-all-exponents}

Let \(n=2m+1\ge5\), and let \(p(z)\) be a polynomial whose zeros
have centroid zero. We use the power sums $Z_p$ and $W_p$ defined
in \eqref{eq:notation}, and set
\begin{equation}\label{eq:coefficient}
N_p=(m-1)(m+1)^p+m^{p+1}+1,\qquad D_p=m(m+1)^p+(m+1)m^p,\qquad B_{n,p}=\frac{N_p}{D_p}.
\end{equation}
At $p=2$, the coefficient in \eqref{eq:coefficient} satisfies
$B_{n,2}=(n-2)/n$, so \eqref{eq:schoenberg} settles \(p=2\).
It remains to consider \(p>2\).

We first prove the inequality for real zeros. For any function
$f:[0,\infty)\to\mathbb R$ with $f(0)=0$, let
\begin{equation}\label{eq:odd-all-functional}
\begin{aligned}
\mathcal L[f]={}&(m-1)\sum_{j=1}^n f((m+1)|z_j|)
+m\sum_{j=1}^n f(m|z_j|)+\sum_{j=1}^n f(|z_j|)\\
&-m\sum_{k=1}^{n-1}f((m+1)|w_k|)
-(m+1)\sum_{k=1}^{n-1}f(m|w_k|).
\end{aligned}
\end{equation}
Substitution of $f(u)=u^p$ into \eqref{eq:odd-all-functional}
gives
\begin{equation}\label{eq:odd-all-power-identity}
\mathcal L[u^p]=N_pZ_p-D_pW_p.
\end{equation}
For centered real zeros, Newton's identities give
$W_2=(n-2)Z_2/n$. Since $N_2/D_2=(n-2)/n$,
\begin{equation}\label{eq:odd-all-quadratic}
\mathcal L[u^2]=0.
\end{equation}
Since $D_p>0$, equations \eqref{eq:coefficient} and
\eqref{eq:odd-all-power-identity} show that the real-rooted
inequality $W_p\le B_{n,p}Z_p$ is equivalent to
$\mathcal L[u^p]\ge0$. The logarithmic comparison in
Theorem~\ref{odd23:thm:global-log} handles $2<p<3$;
Proposition~\ref{prop:cubes} handles $p\ge3$ after the
extremal reduction in Lemma~\ref{lem:reduction}.

\subsection{The range \texorpdfstring{\(2<p<3\)}{2 < p < 3}}

For $s>0$, put $\ell_s(u)=u/s-\log(1+u/s)$ for $u\ge0$.
For every centered real-rooted polynomial of degree $n$,
Theorem~\ref{odd23:thm:global-log}, proved in
Appendix~\ref{app:odd23-details}, gives
$\mathcal J_s\ge0$. The functional $\mathcal J_s$, defined
in Subsection~\ref{odd23:sec:log-framework}, is formed from
the reference root list $((m+1)^{[m]},(-m)^{[m+1]})$ and its
critical-point list $((m+1)^{[m-1]},(-m)^{[m]},1)$.
Expanding the associated product multisets yields
$\mathcal J_s=-\mathcal L[\ell_s]$. Consequently,
\begin{equation}\label{eq:odd-all-log-comparison}
\mathcal L[\ell_s]\le0,\qquad s>0.
\end{equation}
The bound \eqref{eq:odd-all-log-comparison} holds for every
centered real configuration, regardless of the number of
distinct zero values.

For $u\ge0$ and $2<p<3$, the Mellin identity is
\begin{equation}\label{eq:odd-all-mellin}
u^p=\frac{p\sin(\pi p)}{\pi}\int_0^\infty
s^{p-1}\left[\frac{u^2}{2s^2}-\ell_s(u)\right] \,ds.
\end{equation}
For \(u>0\), substitute \(t=u/s\) and put
\(g(t)=\log(1+t)-t+t^2/2\). The integral on the right before
multiplication by its prefactor becomes
\(u^p\int_0^\infty g(t)t^{-p-1}\,dt\).
Since \(g(t)=O(t^3)\) at zero and \(g(t)=O(t^2)\) at infinity,
the integral converges and both boundary terms vanish. Integration
by parts, using \(g'(t)=t^2/(1+t)\), and Euler's beta integral give
\(\int_0^\infty g(t)t^{-p-1}\,dt =\frac1p\int_0^\infty\frac{t^{2-p}}{1+t}\,dt =\frac{\pi}{p\sin(\pi p)}\).
This proves \eqref{eq:odd-all-mellin}, including \(u=0\) by direct
substitution. Its prefactor is positive for \(2<p<3\).
Applying the linear functional \eqref{eq:odd-all-functional}
to \eqref{eq:odd-all-mellin}, and then using
\eqref{eq:odd-all-quadratic} and
\eqref{eq:odd-all-log-comparison}, yields
\begin{equation}\label{eq:odd-all-subcubic-bound}
\mathcal L[u^p]=-\frac{p\sin(\pi p)}{\pi}
\int_0^\infty s^{p-1}\mathcal L[\ell_s] \,ds\ge0.
\end{equation}
The identity \eqref{eq:odd-all-power-identity} therefore
implies $W_p\le B_{n,p}Z_p$ for every centered real-rooted
polynomial when $2<p<3$.

\subsection{The range \texorpdfstring{\(p\ge3\)}{p at least 3}}

Fix $p\ge3$. For nonzero real root lists the quotient
$W_p/Z_p$ is invariant under common nonzero scaling.
It therefore attains a maximum on the compact set
$\sum_{j=1}^n z_j=0$, $\sum_{j=1}^n z_j^2=1$.
Choose a maximizing polynomial. By
Lemma~\ref{lem:reduction}, proved in Appendix~\ref{app:odd-ge-three-details}, its zeros take at most three
distinct values; if there are three, the middle value has
multiplicity one. After reflection and scaling, the root
values are $-\ell-t$, $(n-1)t$, and $k-t$, with respective
multiplicities $k,1,\ell$, where $k+\ell=n-1$,
$1\le k\le m$, and $-\ell/n\le t\le k/n$.
This is the parameterization in \eqref{eq:triple}. At
\(t=-\ell/n\) and \(t=k/n\), two of the root values coincide.

For $s\ge0$, define $f_s(u)=(u^3-s)_+$. The functional $\mathcal D_f$ in \eqref{eq:cross} is
constructed from the normalized three-value list
\eqref{eq:triple} after multiplying all root and critical-point
values by $n$. Since $f_{n^3s}(nu)=n^3f_s(u)$,
the formulas \eqref{eq:odd-all-functional} and \eqref{eq:cross}
give
\begin{equation}\label{eq:odd-all-cubic-tail}
\mathcal L[f_s]=n^{-3}\mathcal D_{f_{n^3s}}\ge0,
\qquad s\ge0.
\end{equation}
The inequality in \eqref{eq:odd-all-cubic-tail} is
Proposition~\ref{prop:cubes}, equation~\eqref{eq:cubic-tail},
proved in Appendix~\ref{app:odd-ge-three-details}.

When $p=3$, set $s=0$ in
\eqref{eq:odd-all-cubic-tail} to obtain
\(\mathcal L[u^3]\ge0\). For \(p>3\) we use the integral representation
\begin{equation}\label{eq:odd-all-cubic-mellin}
u^p=\frac p3\left(\frac p3-1\right)
\int_0^\infty (u^3-s)_+s^{\frac{p}{3}-2}\,ds,
\qquad u\ge0.
\end{equation}
Applying \(\mathcal L\) to \eqref{eq:odd-all-cubic-mellin} and using \eqref{eq:odd-all-cubic-tail} gives
\begin{equation}\label{eq:odd-all-supercubic-bound}
\mathcal L[u^p]=\frac p3\left(\frac p3-1\right)
\int_0^\infty\mathcal L[f_s]s^{\frac{p}{3}-2}\,ds\ge0.
\end{equation}
For the chosen real maximizer, \eqref{eq:odd-all-power-identity}
and the bound at \(s=0\) for \(p=3\), or \eqref{eq:odd-all-supercubic-bound} for \(p>3\), imply
$W_p/Z_p\le B_{n,p}$. Since the configuration maximizes
the quotient, $W_p\le B_{n,p}Z_p$ holds for every centered
real-rooted polynomial of degree $n$.

We have now proved \(W_p\le B_{n,p}Z_p\) for real zeros in both ranges:
\eqref{eq:odd-all-subcubic-bound} treats \(2<p<3\),
\eqref{eq:odd-all-cubic-tail} at \(s=0\) treats \(p=3\), and
\eqref{eq:odd-all-supercubic-bound} treats \(p>3\).
To pass to complex zeros, fix \(\theta\in[0,2\pi]\).
Lemma~\ref{lem:pereira-projection}, applied to the rotated
polynomial and the convex function $x\mapsto|x|^p$, followed
by the real-rooted inequality, gives
\begin{equation}\label{eq:odd-all-projection-transfer}
\sum_{k=1}^{n-1}|\Re(e^{-i\theta}w_k)|^p
\le B_{n,p}\sum_{j=1}^n|\Re(e^{-i\theta}z_j)|^p.
\end{equation}
Integrating \eqref{eq:odd-all-projection-transfer} over
$[0,2\pi]$ and using \eqref{eq:angular} proves
\begin{equation}\label{eq:odd-all-conclusion}
W_p\le B_{n,p}Z_p,\qquad n=2m+1\ge5,\quad p>2.
\end{equation}

To prove sharpness, choose $a,\lambda\in\mathbb C\setminus\{0\}$
and consider
\begin{equation}\label{eq:equality}
p(z)=\lambda(z-(m+1)a)^m(z+ma)^{m+1}.
\end{equation}
The critical points are $(m+1)a$ with multiplicity $m-1$,
$-ma$ with multiplicity $m$, and $a$ once. Hence
$Z_p=D_p|a|^p$ and $W_p=N_p|a|^p$, so equality holds in
\eqref{eq:odd-all-conclusion} for every $p>2$.

For $p\ge3$, the strict assertion in
Proposition~\ref{prop:cubes} excludes a maximizing real-rooted
polynomial with three distinct zeros. Lemma~\ref{lem:reduction}
excludes more than three distinct zero values. The two-value
comparison in Subsection~\ref{sec:finish} shows that equality
among nonzero real-rooted polynomials is attained only for the
multiplicities $m$ and $m+1$. If a complex root list has at least
three distinct values, its real projections have at least three
distinct values for almost every direction. The real inequality
is therefore strict for almost every projection, and integration
of \eqref{eq:odd-all-projection-transfer} gives strictness in
\eqref{eq:odd-all-conclusion}. A centered two-value complex
list lies on a line through the origin, so the real equality
classification applies after rotation. The nonzero equality
polynomials for $p\ge3$ are consequently exactly those in
\eqref{eq:equality}, as asserted in
Theorem~\ref{thm:oddcomplete}.
For \(2<p<3\), equality for a real root list forces the integral
in \eqref{eq:odd-all-subcubic-bound} to vanish. Its weight is
positive, and \(-\mathcal L[\ell_s]\) is continuous and nonnegative
for \(s>0\), so \(\mathcal L[\ell_s]=0\) for every \(s>0\).
For fixed \(u\), as \(s\to\infty\),
\(\ell_s(u)=\frac{u^2}{2s^2}-\frac{u^3}{3s^3}+O(s^{-4})\).
The functional in \eqref{eq:odd-all-functional} is a finite sum.
By \eqref{eq:odd-all-quadratic}, multiplying by \(s^3\) and passing
to the limit gives \(\mathcal L[u^3]=0\).
The cubic equality case of Theorem~\ref{thm:oddcomplete} now gives
\eqref{eq:equality} for nonzero real root lists. A complex list
with three or more distinct values has three distinct real projections
in almost every direction. Integrating the strict real inequality
in those directions rules out equality for that list. The remaining
two-value lists are collinear. The identities \(Z_p=D_p|a|^p\) and
\(W_p=N_p|a|^p\) for \eqref{eq:equality} give sharpness.

\appendix
\section{Integer-sensitive moment estimates for \texorpdfstring{\(1<p<2\)}{1 < p < 2}}\label{app:sub12-details}
\subsection{The unique minimizing two-value proportion}

\begin{thm}\label{sub12:thm:subquad-H-unimodal}
For every $1<p<2$, the function
\[
H_p(r)=\frac{(1+r)(1+r^p-(1-r)^p)}{r+r^p}, \qquad 0<r\le1,
\]
has a unique critical point $r_p\in(0,1)$. It is strictly decreasing
on $(0,r_p)$ and strictly increasing on $(r_p,1)$.
In particular, $r_p$ is the unique global minimum.
\end{thm}

\begin{proof}
Put $q=p-1\in(0,1)$ and
$J(r)=(r+r^p)^2H_p'(r)$ and $K(r)=\frac{J(r)}{r^{q+1}}$.
Direct differentiation gives
\[
J(r)=(1-r)^q\bigl[1+qr+pr^2+(p+qr+r^2)r^q\bigr] -1-p(1-r^2)r^q+r^{2q+2}.
\]
We show that $K$ is strictly concave on $(0,1)$.

Write the binomial expansion
$(1-r)^q=1-\sum_{j=1}^{\infty}a_jr^j$ with $a_j>0$ and $a_1=q$,
and set
$c_2=p-q^2-a_2=\frac{2+q-q^2}{2}>1$ and $b_2=1+p-q^2-pa_2$.
For $j\ge3$, put
$A_j=a_j+qa_{j-1}+pa_{j-2}>0$ and $B_j=pa_j+qa_{j-1}+a_{j-2}>0$.
Expanding the two factors and collecting equal powers of $r$ gives
\begin{equation}
\label{sub12:eq:subquad-K-series} K(r)=-q^2+c_2r^{1-q}+b_2r+r^{q+1} -\sum_{j=3}^{\infty}A_jr^{j-q-1} -\sum_{j=3}^{\infty}B_jr^{j-1}.
\end{equation}
Both series may be differentiated termwise on every compact
subinterval of $(0,1)$. For $j\ge3$, the exponents $j-q-1$ and
$j-1$ exceed $1$, so their second-derivative contributions are
strictly negative. Keeping only the $j=3$ term of the second sum,
$K''(r)\le -c_2q(1-q)r^{-q-1} +q(q+1)r^{q-1}-2B_3$.
Concavity of $z^{\frac{1-q}{1+q}}$ gives
$q(q+1)r^{q-1} \le q(1-q)r^{-q-1}+2q^2$.
Since $B_3=pa_3+qa_2+a_1>q>q^2$ and $c_2>1$, it follows that
$K''(r) \le-(c_2-1)q(1-q)r^{-q-1}-2(B_3-q^2)<0$.
This proves strict concavity.

Equation~\eqref{sub12:eq:subquad-K-series} gives
$K(0+)=-q^2<0$. The explicit expression for $J$ gives $K(1-)=0$,
and, more precisely,
$\lim_{r\uparrow1}\frac{J(r)}{(1-r)^q}=4p>0$.
Indeed the bracket multiplying $(1-r)^q$ tends to $4p$, while
$1+p(1-r^2)r^q-r^{2q+2}=O(1-r)$ and $q<1$.
Hence $K(r)>0$ for all $r<1$ sufficiently close to $1$.

A strictly concave function with $K(0+)<0$, $K(1-)=0$, and
$K(r)>0$ for $r<1$ sufficiently close to $1$ has exactly one
zero in $(0,1)$: two interior zeros together with
the endpoint zero would contradict strict concavity at the middle
zero. Its sign is negative before the interior zero and positive
after it. Since $H_p'(r)$ and $K(r)$ have the same sign for $0<r<1$,
the unique zero of $K$ is the asserted minimizer $r_p$.
\end{proof}

\begin{cor}\label{sub12:cor:subquad-two-value-integers}
For centered two-value lists of degree $n$, the largest ratio is
\[
\mathcal B_{n,p}=\max_{1\le k\le\lfloor n/2\rfloor} \frac{(k-1)(n-k)^p+(n-k-1)k^p+(n-2k)^p} {k(n-k)^p+(n-k)k^p}.
\]
By Theorem~\ref{sub12:thm:subquad-H-unimodal}, it suffices to test the
integers \(\lfloor n\theta_p\rfloor\) and \(\lceil n\theta_p\rceil\),
where \(\theta_p=r_p/(1+r_p)\). Replace either integer outside
\([1,\lfloor n/2\rfloor]\) by the nearest endpoint of the interval.
\end{cor}
\begin{proof}
The ratio is exactly $1-H_p(k/(n-k))/n$, and $k/(n-k)$ increases
strictly with $k$. The asserted monotonicity on either side of
$n\theta_p$ follows from Theorem~\ref{sub12:thm:subquad-H-unimodal}.
\end{proof}

\paragraph{Range of the minimum.}
For $0<r<1$, $H_p(r)>1$, since
$(1+r)(1-r)^p<(1+r)(1-r)<1<1+r^{p+1}$.
At zero,
$H_p(r)=p-(p-1)r^{p-1}+o(r^{p-1})$.
Together with the attained interior minimum this proves
$1<H_p(r_p)<p$.

\subsection{An upper bound for the minimizing two-value proportion}
\label{sub12:sec:subquad-upper-location}

The next estimate controls the location of the minimum in
Theorem~\ref{sub12:thm:subquad-H-unimodal}. It is also used to
locate the integer cell requiring the quartic correction.

\begin{lem}\label{sub12:lem:subquad-upper-location}
Let $1<p<2$, set $q=p-1$, and let $r_p$ be the unique minimizer of $H_p$.
Then
\begin{equation}\label{sub12:eq:subquad-explicit-b-location}
0<\frac{r_p}{1-r_p}<\min\left\{\frac{q^2}{4},\frac1{16}\right\}.
\end{equation}
\end{lem}

\begin{proof}
Put $b_p=r_p/(1-r_p)$. We first prove $b_p<q^2/4$ for
$0<q\le1/2$, using the function $K$ and the series
\eqref{sub12:eq:subquad-K-series}. For $q\ge1/2$, we prove
$b_p<1/16$ by a fixed-radius derivative estimate.

\medskip\noindent\emph{The range $0<q\le1/2$.}
Retain $K$ from the proof of Theorem~\ref{sub12:thm:subquad-H-unimodal}.
By its definition, $K(r)$ and $H_p'(r)$ have the same sign, and
\[
 K(r)=-q^2+c_2r^{1-q}+b_2r+r^{q+1}
 -\sum_{j=3}^{\infty}A_jr^{j-q-1}-\sum_{j=3}^{\infty}B_jr^{j-1},
\]
where
\[
 c_2=1+\frac{q(1-q)}2,\qquad
 b_2=2+\frac{q(1-q)^2}2,\qquad 0<A_j,B_j\le2(1+q)q.
\]
The bounds on \(A_j\) and \(B_j\) follow from $0<a_j\le q$ in the binomial expansion
$(1-r)^q=1-\sum_{j=1}^{\infty}a_jr^j$.
Set $r=q^2/(4+q^2)$, corresponding to $b=q^2/4$.
If $T$ denotes the two positive tails divided by $r$, then
\[
 T\le2(1+q)q\frac{r^{1-q}+r}{1-r},\qquad
 \frac{K(r)}r=-(4+q^2)+c_2r^{-q}+b_2+r^q-T.
\]
For $0<q\le1/4$, one has $r\le1/65$,
$r^{1-q}<1/20$, and therefore $T<q/6$.
Using $r^{-q}+r^q\ge2$ gives
\[
 \frac{K(r)}r\ge
 q\left(1-\frac52q+\frac12q^2\right)-T
 \ge\left(\frac{13}{32}-\frac16\right)q>0.
\]
For $1/4\le q\le1/2$, one has $r\le1/17$ and
$r^{-q}>2$. Thus $r^{-q}+r^q>5/2$, while
$r^{1-q}<1/4$ gives $T<63q/64\le63/128$.
The remaining main term is at least
\[
 \frac12-q^2+q(1-q)+\frac12q(1-q)^2
 \ge\frac14+\frac3{16}+\frac1{16}=\frac12.
\]
Hence $K(r)/r>1/128$. The unique zero of $K$ lies to the left,
which proves $b_p<q^2/4$ for this range of $q$.

\medskip\noindent\emph{The range $1/2\le q<1$.}
It remains to prove $H_p'(1/17)>0$. Put $t=1-q\in(0,1/2]$,
$\alpha=\log(17/16)$, and $\beta=\log17$. The exact numerator
$J=(r+r^p)^2H_p'(r)$ at $r=1/17$ becomes
\[
J(t)=\left(\frac{4928}{4913}-\frac{288}{4913}t\right)e^{\alpha t} +\left(\frac{9536}{83521}-\frac{288}{4913}t\right)e^{(\alpha+\beta)t} -1+\left(-\frac{576}{4913}+\frac{288}{4913}t\right)e^{\beta t} +\frac1{83521}e^{2\beta t}.
\]
The constant Taylor coefficient vanishes. Writing
$J(t)=\sum_{j=1}^{\infty}c_jt^j$, it suffices to verify
\[
 c_1>\frac{49}{100000},\quad
 c_j>0\ (2\le j\le5),\quad
 c_j>-\frac1{1000}\ (6\le j\le11),
\]
and, for $0<t\le1/2$,
\(\left|\sum_{j=12}^{\infty}c_jt^j\right|<\frac{t}{1000000}\).
The rational-arithmetic certificate
\href{https://github.com/zhangteng2000/schoenberg-validation/blob/e7435c5521dcd6a10131ba04a6da500b1404800b/certificates/subquad_04_the_upper_location_estimate.py}{\nolinkurl{subquad_04_the_upper_location_estimate.py}}
verifies the coefficient and tail bounds by enclosing the logarithms using
\[
 \log\frac{1+z}{1-z}
 =2\sum_{j=0}^{23}\frac{z^{2j+1}}{2j+1}+R,
 \qquad 0<R<\frac{2z^{49}}{49(1-z^2)},
\]
using $z=1/33$ and $z=1/3$, and uses
$\log17=4\log2+\log(17/16)$.
The exponential tails are bounded geometrically with
$\alpha<1/16$, $\alpha+\beta<3$, $\beta<3$, $2\beta<6$.
The checker uses exact rational arithmetic only.
It follows that
\(J(t)>\left(\frac{49}{100000}-\frac1{16000} -\frac1{1000000}\right)t>0\).
Theorem~\ref{sub12:thm:subquad-H-unimodal} now yields $b_p<1/16$
for $1/2\le q<1$. Since $q^2/4\ge1/16$
on this range, and since $q^2/4\le1/16$ for $q\le1/2$,
the two estimates prove \eqref{sub12:eq:subquad-explicit-b-location}.

\end{proof}

\subsection{A uniform two-value location bound through degree nineteen}
\label{sub12:sec:turn3-nineteen-location}

\begin{thm}\label{sub12:thm:turn3-nineteen-location}
Let $1<p<2$ and let $r_p$ be the unique minimizer of $H_p$ in
Theorem~\ref{sub12:thm:subquad-H-unimodal}. Put
$\theta_p=r_p/(1+r_p)$. Then
\[
r_p<\frac1{18},\qquad \theta_p<\frac1{19}.
\]
\end{thm}

\begin{proof}
Write $q=p-1$.  For $0<q\le12/25$, Lemma~\ref{sub12:lem:subquad-upper-location}
gives $r_p<q^2/(4+q^2)\le36/661<1/18$.
For $12/25\le q<1$, we prove $H_p'(1/18)>0$. By
Theorem~\ref{sub12:thm:subquad-H-unimodal}, this implies
$r_p<1/18$.

Set $r=1/18$, $\varepsilon=1-q\in(0,13/25]$,
$\alpha=\log(18/17)$ and $\beta=\log18$.  The exact derivative
numerator $J=(r+r^p)^2H_p'(r)$ is
\begin{align*}
 J(\varepsilon)={}&
 (1-r)\{1+r+2r^2-\varepsilon(r+r^2)\}e^{\alpha\varepsilon}\\
 &+(1-r)r\{2+r+r^2-\varepsilon(1+r)\}
                       e^{(\alpha+\beta)\varepsilon}\\
 &-1+r(1-r^2)(-2+\varepsilon)e^{\beta\varepsilon}
       +r^4e^{2\beta\varepsilon}.
\end{align*}
Its constant coefficient vanishes.  Write
$J(\varepsilon)=\sum_{j=1}^{\infty}c_j\varepsilon^j$.
The following finite rational bounds hold:
$c_1>\frac1{5000}$, $c_j>0\quad(2\le j\le5)$, $c_j>-\frac1{1000}\quad(6\le j\le11)$,
and
\[
\left|\sum_{j=12}^{\infty}c_j\varepsilon^j\right| <\frac{\varepsilon}{500000} \qquad\left(0<\varepsilon\le\frac{13}{25}\right).
\]
The certificate \href{https://github.com/zhangteng2000/schoenberg-validation/blob/e7435c5521dcd6a10131ba04a6da500b1404800b/certificates/turn3_norm_nineteen_location_certificate.py}{\nolinkurl{turn3_norm_nineteen_location_certificate.py}}
verifies the coefficient and tail bounds using
\[
2\sum_{j=0}^{23}\frac{z^{2j+1}}{2j+1} <\log\frac{1+z}{1-z} <2\sum_{j=0}^{23}\frac{z^{2j+1}}{2j+1} +\frac{2z^{49}}{49(1-z^2)}.
\]
It takes $z=1/35$ for $\alpha$, and uses
$\beta=4\log2+\log(9/8)$ with $z=1/3$ and $z=1/17$.
All coefficient calculations are rational.  The exponential tails
are bounded geometrically using
$\alpha<\frac1{17}$, $\alpha+\beta<3$, $\beta<3$, and $2\beta<6$.

With $T=13/25$, the coefficient and tail bounds for $J(\varepsilon)$ imply
\[
\frac{J(\varepsilon)}{\varepsilon} >\frac1{5000}-\frac1{1000}\frac{T^5}{1-T} -\frac1{500000} =\frac{17401}{146484375}>\frac1{10000}>0.
\]
Thus $H_p'(1/18)>0$ throughout the remaining range.  Strict
unimodality gives $r_p<1/18$, and the stated bound on $\theta_p$
follows by the increasing map $r\mapsto r/(1+r)$.
\end{proof}

In particular, the condition $n\theta_p\le1$ from the uncapped
three-moment theorem holds for every $4\le n\le19$ and every
$1<p<2$.  The implication uses $\theta_p<1/19$ from
Theorem~\ref{sub12:thm:turn3-nineteen-location}; no monotonicity
of $r_p$ in $p$ is required.

For later use write $S_{n,p}=((n-2)/n)^{\frac{p}{2}}$; the maximum
$\mathcal B_{n,p}$ over two-value configurations is defined in
Corollary~\ref{sub12:cor:subquad-two-value-integers}.

\subsection{The fourth-moment cap beyond the integer constraint}\label{sub12:sec:remaining-fourth-moment}

We retain the notation of Lemmas~\ref{lem:real-spectral-reduction}
and~\ref{lem:section4-fourth}. Put $q=p-1$ and $c=(n-2)/n$.
At a maximizing real-rooted configuration, the objective and the
first three moments of $\mu$ are given by
\eqref{eq:section4-gradient} and \eqref{eq:section4-moments}.
The exact fourth-moment identity is \eqref{eq:section4-fourth-identity}.

Write the remaining $n-2$ roots as $t+y_\nu$, where
$y_\nu\in[-1,1]$ and $\sum_{\nu=1}^{n-2}y_\nu=-nt$.
Choose an integer $m\in\{0,\ldots,n-3\}$ such that
$v=-nt-2m+(n-3)$ lies in $[-1,1]$.
The convex maximization on this cube section, together with
Cauchy--Schwarz, yields
\begin{equation}
\label{sub12:eq:finite-gradient-m4-interval} c^2+\frac{4t^2}{n(n-2)}\le\int_{-1}^{1}y^4\,d\mu(y)\le c^2+\frac2{n^2}(n-3+v^2-nt^2).
\end{equation}
Equality in the upper bound requires every $y_\nu$ except possibly
one to belong to $\{-1,1\}$. The fourth-moment defect identity
\eqref{eq:section4-root-defect} consequently gives
\begin{equation}
\label{sub12:eq:finite-gradient-deficiency} \int_{-1}^{1}(1-y^2)(y-t)^2\,d\mu(y)\ge\frac{2(1-v^2)}{n^2}.
\end{equation}

\paragraph{The symmetric endpoint.}

If $t=0$, equations \eqref{eq:section4-gradient} and
\eqref{eq:section4-moments}, together with concavity of $s^{\frac{p}{2}}$,
give $W_p/Z_p=\int_{-1}^{1}|y|^p\,d\mu(y)\le c^{\frac{p}{2}}$.
Equality forces $\mu$ to be supported on $\{\pm\sqrt c\}$,
so $\int y^4\,d\mu=c^2$. By
\eqref{eq:section4-fourth-identity}, this implies
$\sum_{\nu=1}^{n-2}y_\nu^2=0$, since $t=0$.
Thus all internal zeros vanish, and equality holds only for
the sparse pair, up to real scaling and permutation.

\paragraph{The relaxed moment problem.}

Set $\Lambda_{n,p}=\max\{S_{n,p},\mathcal B_{n,p}\}$. It suffices
to prove that every probability measure
$\mu$ on $[-1,1]$ satisfying
\eqref{eq:section4-moments} and the upper bound in
\eqref{sub12:eq:finite-gradient-m4-interval} obeys
\[
2\int y\,\operatorname{sgn}(t+y)|t+y|^q\,d\mu(y) \le \Lambda_{n,p}\bigl((1-t)^q+(1+t)^q\bigr), \qquad 0\le t\le c.
\]
The spectral measure of a maximizing real root list satisfies
\eqref{eq:section4-moments} and \eqref{sub12:eq:finite-gradient-m4-interval}.
Outside the cell containing the continuous maximum we use the three-moment bound. Inside the cell,
\eqref{sub12:eq:finite-gradient-deficiency} bounds the integral of
the quartic correction from below.

For a quartic polynomial $P(y)=\sum_{j=0}^4a_jy^j$ majorizing
$y\operatorname{sgn}(t+y)|t+y|^q$ on the whole interval, its integral
is bounded using the fixed first three moments and, according to
the sign of $a_4$, the appropriate fourth-moment endpoint. A quartic polynomial gives an exact dual certificate only after
its pointwise majorization and the resulting integral inequality
have both been proved on the full interval $[-1,1]$.

\subsubsection{The entire rounding cell lies in the quartic theorem's domain}

Let $r_p$ be the unique minimizer of $H_p$, put
$\theta_p=r_p/(1+r_p)$ and $t_p=1-2\theta_p$.
Equation~\eqref{sub12:eq:subquad-explicit-b-location} implies
$r_p<1/17$ and hence $\theta_p=r_p/(1+r_p)<1/18$.
If $n\theta_p<1$, then $t_p>(n-2)/n$ and the continuous
two-value expression is bounded by the one-outlier endpoint on
the whole admissible interval. If $n\theta_p$ is an integer,
its continuous maximum is already attained at an allowed lattice
point. Otherwise put $k=\lfloor n\theta_p\rfloor\ge1$.
Only the cell
$\left[1-\frac{2(k+1)}n,\ 1-\frac{2k}n\right]$
can contain points where the continuous two-value expression exceeds
all discrete two-value values. Its left endpoint satisfies
$1-\frac{2(k+1)}n >t_p-\frac2n \ge1-4\theta_p>\frac79>\frac34$.
Thus \(t>3/4\) throughout this cell.
Theorem~\ref{sub12:thm:turn3-peak-anchor-positive} gives the sign of \(E_-\).
Theorems~\ref{sub12:thm:turn3-integer-anchor-branch}
and~\ref{sub12:thm:turn3-backward-remainder-proved} give the comparison with the discrete maximum.

\subsection{A global two-contact classification for the mean-only envelope}
\label{sub12:sec:mean-only-high-asymmetry}

Fix $0<q<1$, set $p=q+1$, and let $0<t<1$.  Define
$\phi_t(y)=y\operatorname{sgn}(t+y)|t+y|^q$ and $D=(1-t)^q+(1+t)^q$,
and let $K_t$ be the continuous divided difference
$K_t(y)=\phi_t[-1,y,1] =\frac{D/2+\bigl((1+t)^q-(1-t)^q\bigr)y/2-\phi_t(y)}{1-y^2}$.
The endpoint values are understood by continuity.

Among all probability measures $\nu$ on $[-1,1]$ with mean $t$, the
minimum of $\int K_t\,d\nu$ is attained either by $\delta_t$ or by a
measure supported on $\{-1,b\}$ for some $t<b\le1$.  Thus this mean-only
minimum is exactly
\begin{equation}
\label{sub12:eq:mean-only-envelope-one-variable} \min_{t\le b\le1} \frac{(b-t)K_t(-1)+(1+t)K_t(b)}{1+b}.
\end{equation}
The choice $b=t$ includes the point mass $\delta_t$.
This does not impose the reciprocal or variance constraints in the
endpoint-gradient problem, and therefore does not by itself give the
sharp finite-degree constant.

\paragraph{A global endpoint-chord inequality on the negative branch.}
Let $Q_H$ be the cubic Hermite interpolant of $\phi_t$ at $-1,1$, with
values and first derivatives prescribed at both endpoints.  Put
$A=1-t$, $B=1+t$.  Its coefficients of $y^2,y^3$ are
$d_2=\frac{D+q(A^{q-1}+B^{q-1})}{4}>0$ and $d_3=\frac{q(B^{q-1}-A^{q-1})}{4}<0$.
Write $g(u)=u^{q+1}+tu^q$ and $P(u)=Q_H(-t-u)$.
Since $P$ and $g$ have equal value and derivative at $u=A$,
the cubic Taylor formula gives
$P(u)=g(A)+g'(A)(u-A) +\bigl(d_2-(2+t)d_3-d_3u\bigr)(u-A)^2$.
The useful coefficient identity is
$d_2-(2+t)d_3=qA^{q-1}+\frac{1-q}{4}D$.
Concavity of $u^q$ on $[0,\infty)$ gives
$u^q\le A^q+qA^{q-1}(u-A)$.
Multiplying this inequality by $u\ge0$ yields
$u^{q+1}\le A^{q+1}+(q+1)A^q(u-A) +qA^{q-1}(u-A)^2$.
Adding the tangent bound for $tu^q$ therefore gives
$g(u)\le g(A)+g'(A)(u-A)+qA^{q-1}(u-A)^2$.
Consequently
\begin{equation}
\label{sub12:eq:global-negative-hermite-gap} P(u)-g(u)\ge \left(\frac{1-q}{4}D-d_3u\right)(u-A)^2>0 \qquad(u\ge0,\ u\ne A).
\end{equation}
In particular this applies to $0\le u<A$, which is the entire
negative shifted branch $-1<y\le-t$, including the cusp.
Let $L_H$ be the affine function joining the endpoint values of $K_t$.
The identity
$Q_H(y)-\phi_t(y)=(1-y^2)(K_t(y)-L_H(y))$
and \eqref{sub12:eq:global-negative-hermite-gap} show that
\begin{equation}
\label{sub12:eq:global-negative-K-chord} K_t(y)>\frac{1-y}{2}K_t(-1)+\frac{1+y}{2}K_t(1) \qquad(-1<y\le-t).
\end{equation}
Replacing any mass in this interval by its endpoint mixture preserves
the mean and strictly decreases its contribution to the objective.
Hence no optimizing measure puts positive mass there.

\begin{proof}[Proof of the classification]
Compactness and the elementary description of extreme points of the
probability measures with one prescribed moment show that an optimizer
can be chosen with at most two support points.  Equivalently, the lower
convex envelope of the continuous function $K_t$ at the interior point
$t$ is realized by a chord.  It has a supporting affine function $\ell$,
with $\ell\le K_t$ and equality at every support point of this optimizer.
By \eqref{sub12:eq:global-negative-K-chord}, every interior support point lies
strictly to the right of $-t$.
The cubic polynomial
$Q(y)=D/2+\bigl((1+t)^q-(1-t)^q\bigr)y/2-(1-y^2)\ell(y)$
majorizes $\phi_t$ and automatically has contact at $y=-1,1$.
Write $z=t+y$ and $Q(z-t)=a_0+a_1z+a_2z^2+a_3z^3$.
Majorization forces $a_0>0$: if $a_0=0$, then at $z=-u<0$ small the
residual would equal $-tu^q+O(u)<0$.
On the positive half-axis the residual is
$a_0+t z^q+a_1z-z^{q+1}+a_2z^2+a_3z^3$.
Its coefficient signs in increasing exponent order have at most three
changes.  The generalized Descartes rule therefore bounds its positive
zeros, with multiplicity, by three.  This rule follows from repeated
Rolle arguments after division by the lowest power; zero coefficients
are simply omitted.
Two interior support points would both give double positive-$z$
contacts, contradicting that zero bound.  An optimizer on
$\{a,1\}$ with $-1<a<t$ is impossible as well: its interior atom
satisfies $a>-t$, and contact of $\ell$ with $K_t$ at $1$ makes
$Q-\phi_t=(1-y^2)(K_t-\ell)$ double at the right endpoint.
Together with the double contact at $a$, this again gives four positive
zeros.  The only remaining choices are a single atom at $t$ and a
two-atom measure on $\{-1,b\}$.  Their weights give exactly
\eqref{sub12:eq:mean-only-envelope-one-variable}.
\end{proof}

At $t=0$ the conclusion needed for the envelope follows directly from
$K_0(y)=(1-|y|^p)/(1-y^2)\ge p/2$, with equality at the endpoints.
The minimization in \eqref{sub12:eq:mean-only-envelope-one-variable}
uses only a prescribed first moment. The spectral problem in
\eqref{eq:section4-moments} also fixes the second and third moments;
the mean-only calculation cannot determine its sharp value.
\subsection{An exact global quartic majorant criterion}
\label{sub12:sec:quartic-anchor-majorant}

We construct a quartic majorant for the integrand in
\eqref{eq:section4-gradient}. Its negative fourth-degree correction
can then be controlled by \eqref{eq:section4-fourth}.  Throughout, $0<q<1$, $p=1+q$, and
$0<t<1$.  Put
$\phi_t(y)=y\operatorname{sgn}(t+y)|t+y|^q$ for $-1\le y\le1$.
Let $Q_0$ be the cubic Hermite interpolant specified by
$Q_0(-1)=\phi_t(-1),\quad Q_0(1)=\phi_t(1),\quad Q_0(t)=\phi_t(t),\quad Q_0'(t)=\phi_t'(t)$.
Define
\begin{equation}\label{eq:quartic-endpoint-coefficients}
 E_-=\frac{Q_0'(-1)-\phi_t'(-1)}{2(1+t)^2},\qquad
 E_+=\frac{\phi_t'(1)-Q_0'(1)}{2(1-t)^2}.
\end{equation}

\begin{thm}[Global quartic anchor criterion]
\label{sub12:thm:quartic-anchor}
For $E\ge0$, the polynomial
\[
Q_E(y)=Q_0(y)+E(y^2-1)(y-t)^2
\]
majorizes $\phi_t$ throughout $[-1,1]$ if and only if
\[
E\le\min\{E_-,E_+\}.
\]
In particular, whenever this minimum is nonnegative, it is the exact
largest admissible nonnegative coefficient.
\end{thm}

\begin{proof}
Let $\ell_F$ be the affine interpolant of $\phi_t$ at $-1,1$, and set
$K_t(y)=\frac{\ell_F(y)-\phi_t(y)}{1-y^2}$ and $P_E(y)=K_t(t)+K_t'(t)(y-t)+E(y-t)^2$.
The endpoint values of $K_t$ are defined by continuity.
The interpolation conditions give the exact identity
$Q_E(y)=\ell_F(y)-(1-y^2)P_E(y)$.
Moreover, the endpoint slope conditions are precisely
$E\le E_\pm\quad\Longleftrightarrow\quad P_E(\pm1)\le K_t(\pm1)$.
They are necessary for a majorant because the residual vanishes at
both endpoints.

Assume \(E\le\min\{E_-,E_+\}\).  Since $E\ge0$, the function $P_E$ is
convex and therefore lies below the affine chord joining its endpoint
values.  Hence it lies below the affine chord $L_H$ joining
$K_t(-1),K_t(1)$.  The global negative-branch chord inequality
\eqref{sub12:eq:global-negative-K-chord} gives
$P_E(y)\le L_H(y)<K_t(y)$ for $-1<y\le-t$.
Thus $Q_E-\phi_t$ is strictly positive on this entire interval,
including the cusp $y=-t$.

It remains to consider $z=t+y>0$.  Write
$Q_E(z-t)=a_0+a_1z+a_2z^2+a_3z^3+Ez^4$.  We have $a_0>0$ by the
strict inequality just proved.  The residual on the positive half-axis
is the generalized polynomial
$R(z)=a_0+t z^q+a_1z-z^{q+1}+a_2z^2+a_3z^3+Ez^4$.
For $E>0$, its coefficients, in increasing exponent order, have at
most four sign changes.  The generalized Descartes rule bounds the
number of its positive zeros, counted with multiplicity, by four.
The prescribed contact at $z=2t$ is at least double, and the endpoint
$B=1+t>2t$ is a zero with $R'(B)\le0$.
If $R'(B)=0$, the contacts at \(2t\) and \(B\) already account for at least four
zeros.  If $R'(B)<0$, then $R$ is negative just to the right of $B$,
whereas its positive leading coefficient makes $R(z)>0$ for all
sufficiently large $z$; hence there is an additional zero beyond $B$.
Again four zeros have been accounted for.  Therefore there can be no
additional positive zero in $(0,B)$, and the contact at $2t$ has even
multiplicity.  Since $R(0)=a_0>0$, it follows that $R\ge0$ on $[0,B]$.

When $E=0$, the generalized Descartes rule applied to the cubic residual gives at most three positive
zeros, irrespective of the sign of the cubic coefficient.  The double
contact $2t$ and the endpoint $B$ already account for three.  There is
no further zero and the contact $2t$ cannot have odd multiplicity.
Again $R(0)>0$ proves $R\ge0$ on $[0,B]$.
This proves sufficiency, including zero endpoint slopes and zero
intermediate coefficients.
\end{proof}

\begin{cor}[An arbitrary positive shifted contact]
\label{sub12:cor:quartic-general-contact}
The criterion in Theorem~\ref{sub12:thm:quartic-anchor} remains valid
when the double-contact point $t$ is replaced by
$b\in(-t,1)$.  More precisely, let $Q_{0,b}$ match the values of
$\phi_t$ at $-1,1,b$ and its derivative at $b$, and put
\[
E_{-,b}=\frac{Q_{0,b}'(-1)-\phi_t'(-1)}{2(1+b)^2},\qquad E_{+,b}=\frac{\phi_t'(1)-Q_{0,b}'(1)}{2(1-b)^2}.
\]
For $E\ge0$, the polynomial
$Q_{0,b}(y)+E(y^2-1)(y-b)^2$ majorizes $\phi_t$ on $[-1,1]$
if and only if $E\le\min\{E_{-,b},E_{+,b}\}$.
\end{cor}
\begin{proof}
In the proof of Theorem~\ref{sub12:thm:quartic-anchor}, replace $P_E$ by
$K_t(b)+K_t'(b)(y-b)+E(y-b)^2$.  The endpoint chord argument is
unchanged.  On the positive shifted half-axis, the prescribed double
zero is $z=t+b\in(0,1+t)$, so the generalized Descartes argument
from Theorem~\ref{sub12:thm:quartic-anchor} applies without change.
\end{proof}

\begin{lem}[The right endpoint coefficient is positive]
\label{sub12:lem:anchor-right-positive}
For every $0<q<1$ and $1/2\le t<1$, one has $E_+>0$.
Thus, in this parameter range, Theorem~\ref{sub12:thm:quartic-anchor} says that
the existence of an admissible
nonnegative $E$ is equivalent to $E_-\ge0$.
\end{lem}
\begin{proof}
Put $A=1-t$, $b=2t$, $B=1+t$, and
$g(z)=\frac{z^q+A^q}{z+A}$ and $k(z)=z^q-g(z)$ for $z>0$.
On shifting the nodes by $t$, a first divided difference at the
negative node $-A$ gives
$E_+=k[b,b,B,B]$.
Let $Q=q(q-1)(q-2)=q(1-q)(2-q)>0$ and $r=A/z$.
Direct differentiation yields
$z^{4-q}g'''(z)=\frac{P(r)}{(1+r)^4}$,
where
$P(r)=Q(1+r)^3-3q(q-1)(1+r)^2+6q(1+r)-6(1+r^q)$.
For $0\le r\le1$, the inequality $r^q\ge r$ gives
\begin{align*}
 Q(1+r)^4-P(r)
 &\ge (1-q)(1+r)\bigl(q(2-q)r(1+r)^2-3q(1+r)+6\bigr)>0.
\end{align*}
The last bracket exceeds $6-3q(1+r)>0$.
Consequently, for $z\ge1$ and $A\le1/2$,
$k'''(z)=Qz^{q-3}-g'''(z)>Qz^{q-4}(z-1)\ge0$.
Since $b=2t\ge1$ and $B=1+t>1$, the four nodes $b,b,B,B$
lie in $[1,\infty)$. The Hermite--Genocchi formula and the strict
positivity of $k\prime\prime\prime$ on $(1,\infty)$ give $E_+>0$.
\end{proof}

\paragraph{An exact comparison of the endpoint coefficients.}
The two endpoint coefficients admit explicit expressions in
$r=(1-t)/(1+t)$; no differentiation of the interpolating
polynomial is needed to compare them.
Put $r=(1-t)/(1+t)$, $X=r^q$, and $Z=(1-r)^q$.  Direct substitution gives
\begin{align*}
 \frac{E_-}{(1+t)^{q-2}}
 &=\frac{Xqr^2+Xqr-2Xr^3-Xr^2-Zqr^2-Zqr-2Zr^3-Zr^2-Z+1}
 {4r^2},\\
 \frac{E_+}{(1+t)^{q-2}}
 &=\frac{Xr^3+Zqr^2+Zqr+Zr^3+Zr+2Z+qr^2+qr-r-2}
 {4r^3}.
\end{align*}
In particular,
$\frac{E_--E_+}{(1+t)^{q-2}} =\frac{1+r}{4r^3} \left(2-qr+Xr^2(q-2r)-Z(2+qr+qr^2+2r^3)\right)$.
Thus, when $q>2r$, the condition that the right endpoint determines
the admissible coefficient is a linear lower bound on $X$:
$E_-\ge E_+ \quad\Longleftrightarrow\quad X\ge\frac{Z(2+qr+qr^2+2r^3)-2+qr}{r^2(q-2r)}$.

\paragraph{Moment consequence.}
Suppose a probability measure $\mu$ on $[-1,1]$ satisfies
$m_1=0$, $m_2=c_n=\frac{n-2}{n}$, and $m_3=-\frac{2t}{n}$,
and write
$\Delta=\int(1-y^2)(y-t)^2\,d\mu(y) =c_n-\frac{2t^2}{n}-m_4$.
Then the majorant gives the exact bound
\[
\frac{2\int\phi_t\,d\mu}{(1-t)^q+(1+t)^q} \le 1-\frac{H_p((1-t)/(1+t))}{n} -\frac{2E\Delta}{(1-t)^q+(1+t)^q}.
\]
In the endpoint-gradient application,
$\Delta\ge2(1-u^2)/n^2$ with the integer/parity parameter $u\in[-1,1]$.
The discrete comparison uses Theorem~\ref{sub12:thm:turn3-integer-anchor-branch}
and Theorem~\ref{sub12:thm:turn3-backward-remainder-proved}.

\paragraph{Differentiating the integer-cell objective.}
To differentiate the integer-cell objective, set
\[
D(t)=(1-t)^q+(1+t)^q,\qquad \kappa(t)=H_p\!\left(\frac{1-t}{1+t}\right),\qquad K(t)=D(t)\kappa(t),\qquad L(t)=\frac{D'(t)}{D(t)}.
\]
Then
\begin{equation}\label{eq:integer-cell-derivatives}\begin{aligned}
 K'(t)&=a(t)K(t)-b(t)D(t),\\
 a(t)&=\frac{q+1+(1-q)t^2}{t(1-t^2)},\qquad
 b(t)=\frac{2(q+1)}{t(1-t^2)},\\
 \kappa'(t)&=(a(t)-L(t))\kappa(t)-b(t),\\
 L'(t)&=\frac{q(q-1)-2(q-1)tL(t)}{1-t^2}-L(t)^2.
\end{aligned}\end{equation}
For example, the first identity follows by differentiating
$K(t)=\frac{2(1-t)^q}{1+t}+\frac{2(1+t)^q}{1-t} -\frac{4t(2t)^q}{1-t^2}$.
Moreover
$\frac{E_+(t)}{D(t)}=A(t)\kappa(t)+B(t)+\frac{L(t)}{2(1-t)}$,
where
\[
A(t)=\frac{qt^2-q-4t^2-2t-2}{8t(1-t)^2(1+t)},\qquad B(t)=\frac{2qt^2+qt+q+t^2+t+2}{4t(1-t)^2(1+t)}.
\]
Thus, for fixed lattice endpoints $a_0<b_0$ of distance $2/n$, the
function
$\mathcal F(t)=\kappa(t)+4n\frac{E_+(t)}{D(t)}(t-a_0)(b_0-t)$
and each of its first two derivatives are affine in $\kappa(t)$
by \eqref{eq:integer-cell-derivatives}. In particular,
if $\alpha=E_+/D$ and $u=n(t-(a_0+b_0)/2)$, then
$\mathcal F''(t)=\kappa''(t)+\frac4n(1-u^2)\alpha''(t) -16u\alpha'(t)-8n\alpha(t)$.

\subsection{A necessary condition for joint stationary points of the three-moment relaxation}
Fix $n\ge3$, $q\in(0,1)$, $c=(n-2)/n$, and write
$F_t(y)=y\operatorname{sgn}(t+y)|t+y|^q,\quad D_t=(1-t)^q+(1+t)^q,\quad d(t)=-2t/n$.
The three-moment relaxed objective is $2I/D_t$, where
$I=\int F_t\,d\mu$ and $\mu$ is a probability measure on $[-1,1]$ with
moments $m_1=0,m_2=c,m_3=d(t)$.
Suppose $0<t<c$ and a finitely supported measure has all its spectral
atoms in $(-1,1)$.  Assume explicitly that there is a cubic dual polynomial
$Q(y)=\lambda_0+\lambda_1y+\lambda_2y^2+\lambda_3y^3\ge F_t(y) \quad(-1\le y\le1)$
that touches the integrand at every atom and that the joint $t$ multiplier equation holds.
The dual contact and stationarity conditions hold at an interior global maximizer over all
feasible probability measures and $t$; they are not asserted for a local
maximum within a prescribed support-size family.
Then
\begin{equation}
\label{sub12:eq:lp3-lambda0-identity} \lambda_0=\frac I2\left(1-q+\frac{tD_t'}{D_t}\right), \qquad 1-q+\frac{tD_t'}{D_t}\ge0.
\end{equation}
For the asserted application to a finite global maximizer, cusp atoms can
also be excluded before using smooth multipliers.  If at least three atoms
are present, two noncusp atoms supply the rank-four confluent moment
Jacobian; moving an atom at $y=-t$ and correcting the moments with those
other atoms produces an odd leading term of order $|\varepsilon|^q$ in
the objective. One sign strictly increases it. If a feasible two-atom
measure had an atom at \(-t\), its zero mean and second moment \(c\)
would force the other atom to be \(c/t>1\), outside the support interval.

The multiplier condition is regular for distinct interior atoms: a cubic
whose value and derivative vanish at two distinct atoms is zero, so the
confluent moment Jacobian has rank four.  At a global measure maximizer,
this rank also allows insertion of a new atom with small positive mass
while adjusting old weights and positions to preserve all four moments.
The resulting one-sided first-variation inequality gives the global
majorant $Q\ge F_t$.  The possible cusp $y=-t$ cannot
be a contact, because for $t>0$ the leading local term of $F_t$ is
$-t\operatorname{sgn}(y+t)|y+t|^q$, which exceeds any smooth contact
polynomial on one side.

For the proof, interior contact implies $F_t'(y)=Q'(y)$ at every atom.
Euler homogeneity gives
$t\partial_tF_t+y\partial_yF_t=(q+1)F_t$.
Integrating and using the moments therefore gives
$t\int\partial_tF_t\,d\mu =(q+1)I-2c\lambda_2-3d\lambda_3$.
Joint stationarity of $I/D_t$ gives
$\int\partial_tF_t\,d\mu-\frac{ID_t'}{D_t} +\lambda_3d'(t)=0$.
Using \(td'=d\) and \(I=\int Q\,d\mu=\lambda_0+c\lambda_2+d\lambda_3\),
we eliminate \(\int\partial_tF_t\,d\mu\) and \(\lambda_2\) from the homogeneity
and stationarity identities to obtain \eqref{sub12:eq:lp3-lambda0-identity}.
Finally, $Q(0)\ge F_t(0)=0$. The strict positivity $I>0$ follows from
monotonicity of $h(x)=\operatorname{sgn}(x)|x|^q$, since
$I=\operatorname{Cov}(Y,h(t+Y))$ and $\operatorname{Var}Y=c>0$.

The function
$1-q+tD_t'/D_t$ strictly decreases from $1-q$ to $-\infty$ on $(0,1)$:
$D_t'<0$ and $D_t''<0$, so $\log D_t$ is strictly concave there.
It has a unique zero $\tau_q$.  Thus a joint maximizer with all spectral
atoms interior must have $t\le\tau_q$.
For $q=1/2$, direct simplification gives $\tau_{1/2}=\sqrt3/2$.
This criterion excludes only high-$t$ interior-contact maximizers; it does
not by itself classify the remaining low-$t$ or endpoint-contact cases.

If the dual polynomial also contacts an endpoint, the Euler and
stationarity calculation includes boundary derivatives and yields
$\lambda_0=\frac I2\left(1-q+\frac{tD_t'}{D_t}\right)+\frac e2$,
where
$e=\mu(\{1\})\{F_t'(1)-Q'(1)\} -\mu(\{-1\})\{F_t'(-1)-Q'(-1)\}\ge0$.
The inequality follows from the one-sided derivative conditions for a
majorant at the two endpoints. The term \(e/2\) is the boundary
contribution in the Radau and Lobatto cases of
\eqref{sub12:eq:lp3-lambda0-identity}.

\paragraph{The continuous two-value optimum.}
Let
\[
\mathcal H_p(t)= \frac{2\{(1-t)^p+(1+t)^p-(2t)^p\}} {(1-t^2)\{(1-t)^{p-1}+(1+t)^{p-1}\}},\qquad 0<t<1.
\]
This is the two-value deficit written in the parameter
$r=(1-t)/(1+t)$.  Its unique minimum is $H_p(r_p)<p$ at $t=t_p$.
Then $t_p>\tau_q$.
To prove this directly, put
$J(t)=(1-t)^p+(1+t)^p-(2t)^p$, so
$\mathcal H_p=2J/[(1-t^2)D_t]$.
The elementary identity
$tJ'(t)=p\{J(t)-D_t\}$
and $\mathcal H_p'(t_p)=0$ give
$1-q+\frac{t_pD_{t_p}'}{D_{t_p}} =\frac{2(1-p/H_p(r_p))}{1-t_p^2}<0$.
Strict decrease of the left side proves $t_p>\tau_q$.
Thus the continuous two-value maximum lies in a region where a
configuration with only interior contact points cannot be jointly stationary.
\subsection{Three contact points in the high-asymmetry regime}
Let $0<q<1$, $p=q+1$, and $p/3\le t<1$.
Suppose a cubic polynomial $Q(y)$ majorizes
$F_t(y)=y\operatorname{sgn}(t+y)|t+y|^q$
on $[-1,1]$.  Suppose a probability measure supported on the contact set
has moments $m_1=0$, $m_2=c<1$, and
$m_3=-(1-c)t$.  Then this contact set has at most three points that can
support such a measure; in fact the full contact set has at most three
points once the necessarily present interior positive contact is used.
Consequently, in this regime a cubic dual certificate for the
three-moment problem reduces every optimizer to at most three atoms.
No claim for $t<p/3$ is made here.

Write $x=t+y$ and
$P(x)=Q(x-t)=a_0+a_1x+a_2x^2+a_3x^3$.
Majorization at the cusp forces $a_0>0$: if $a_0=0$, then at $x=-u<0$
small the residual would be $-tu^q+O(u)<0$.
For $x>0$ the residual is the generalized polynomial
$E_+(x)=a_0+t x^q+a_1x-x^{q+1}+a_2x^2+a_3x^3$.
The generalized Descartes rule (obtained from Rolle's theorem by dividing
by the first power and differentiating successively) gives at most three
positive zeros, counting multiplicities, because the coefficient signs
in increasing exponent order are
$(+, +,\operatorname{sgn}a_1,-,\operatorname{sgn}a_2,
\operatorname{sgn}a_3)$.
If $a_3\ge0$, this count improves to two.
If $a_3<0$ and $a_2\le0$, it improves to one.
Zero coefficients are simply omitted.

The moment identity
$\int y(1-y^2)\,d\mu=-(m_3)=(1-c)t>0$
forces an interior contact with $y>0$, hence with $x>0$.
Every such interior contact has multiplicity at least two in $E_+$.
It follows that there is at most one interior positive-$x$ contact.
When $a_3\ge0$, there are no other positive-$x$ contacts.
When $a_3<0$, one must have $a_2>0$, and there can be at most one further
positive-$x$ contact, necessarily the upper endpoint $x=t+1$.

It remains to count contacts with $x=-u<0$, $0<u\le1-t$.
Put $g(u)=u^{q+1}+tu^q$.  Since $t\ge p/3$ and $u\le1-t$,
$g'''(u)=q(q-1)u^{q-3}\{(q+1)u-t(2-q)\}>0$
in the open interval.  Hence $g''$ is increasing and
$g''(u)\le q(p-2t)(1-t)^{q-2}$.
If $a_3<0$, choose the already established positive interior contact
$x=v$, where $t<v<t+1$.  Contact of both value and derivative gives
$a_2=\frac{a_0}{v^2}-2a_3v+qv^{q-1}+(1-q)t v^{q-2}$.
The last two terms strictly decrease as a function of $v>0$.
Using $a_0>0$, $a_3<0$, and $v<t+1$, we obtain
$\{P(-u)\}''=2a_2-6a_3u> \frac{2(q+t)}{(1+t)^{2-q}}$.
The elementary inequality
\begin{equation}
\label{sub12:eq:high-t-scalar-curvature} \frac{2(q+t)}{(1+t)^{2-q}}> \frac{q(q+1-2t)}{(1-t)^{2-q}},\qquad 0<q<1,\quad 0<t<1,
\end{equation}
therefore proves $\{P(-u)-g(u)\}''>0$.
A strictly convex nonnegative function has at most one zero on this
interval.  Combined with the at most two positive-$x$ contacts, this
proves the claim in this case.

\eqref{sub12:eq:high-t-scalar-curvature} holds even without the
restriction $t\ge p/3$.  It is immediate if $q+1-2t\le0$.
Otherwise let $\mathcal R(q,t)$ be the ratio of its two positive sides.
For fixed $t$ and $q$ with $q+1-2t>0$,
\[
\partial_q\log\mathcal R =-\frac{t}{q(q+t)}-\frac1{q+1-2t} +\log\frac{1+t}{1-t} \le-\frac{t}{1+t}-\frac1{2(1-t)} +\log\frac{1+t}{1-t}<0.
\]
To justify the last strict sign, put $U=(1+t)/(1-t)>1$ and $v=\sqrt U$.
The elementary inequality $\log U\le\sqrt U-1/\sqrt U$ follows by
integrating $(v-1)^2/v^2\ge0$.  Moreover,
$\frac{U+3-2/U}{4}-\left(\sqrt U-\frac1{\sqrt U}\right) =\frac{2+(v-1)(v-2)^2(v+1)}{4v^2}>0$.
Since $t/(1+t)+1/[2(1-t)]=(U+3-2/U)/4$, the derivative is negative.
Finally $\mathcal R(1,t)=1$, so $\mathcal R(q,t)>1$ for $q<1$.

If $a_3\ge0$, then
$\{P(-u)-g(u)\}'''=-6a_3-g'''(u)<0$.
A nonnegative function with strictly negative third derivative cannot
have two distinct interior zeros: two interior minima and a maximum
between them would force three zeros of the first derivative and hence
a zero of the third derivative by Rolle's theorem.
There is therefore at most one interior negative-$x$ contact and at most
one additional contact at $u=1-t$; the cusp is not a contact.
Adding the unique positive-$x$ contact again gives at most three.
\subsection{Three contacts for the capped fourth-moment dual}
\label{sub12:sec:quartic-three-contacts}

Let $0<q<1$, $p=q+1$, and $p/4\le t<1$.  Suppose a polynomial $Q$ of
degree at most four majorizes
$\phi_t(y)=y\operatorname{sgn}(t+y)|t+y|^q$
on $[-1,1]$, and suppose the coefficient of $y^4$ is nonnegative.
If its contact set supports a probability measure with
$m_1=0$, $m_2=c<1$, and $m_3=-(1-c)t$,
then the full contact set contains at most three points.  In particular,
whenever an optimal quartic dual certificate exists for the
fourth-moment upper-cap relaxation in this regime, every measure
complementary to that certificate is supported on at most three points.

\paragraph{A scalar curvature bound.}
For $g(u)=u^{q+1}+tu^q$ one has
\begin{equation}
\label{sub12:eq:quartic-contact-global-curvature} \sup_{u>0}g''(u) =q\left(\frac{q+1}{2-q}\right)^{2-q}t^{q-1} <\frac{2(q+t)}{(1+t)^{2-q}} \qquad(t\ge p/4).
\end{equation}
Indeed $g''$ has its unique maximum at $u=(2-q)t/(q+1)$.
The ratio of the right side of \eqref{sub12:eq:quartic-contact-global-curvature}
to its middle expression has logarithmic derivative in $t$ equal to
$\frac{(1-q)(q+2t)}{t(q+t)(1+t)}>0$.
At $t=p/4$ this ratio is
$R(q)=\frac{2(5q+1)}{q(q+1)} \left(\frac{2-q}{q+5}\right)^{2-q}$.
Writing $h(q)=\log R(q)$ gives
\begin{align*}
 h'(q)&=\frac5{5q+1}-\frac1q-\frac1{q+1}
       +\log\frac{q+5}{2-q}-\frac7{q+5},\\
 h''(q)&=\left(\frac1{q^2}-\frac{25}{(5q+1)^2}\right)
       +\frac1{(q+1)^2}+\frac1{q+5}+\frac1{2-q}
       +\frac7{(q+5)^2}>0.
\end{align*}
Since $h(1)=0$ and $h'(1)=\log6-11/6<0$, we have $h(q)>0$ for
$0<q<1$.  The elementary strict bound on $\log6$ follows, for example,
from the first five terms of the exponential series at $11/6$, whose
sum is $186985/31104>6$.

\begin{proof}
Shift $z=t+y$ and write
$Q(z-t)=\sum_{j=0}^4 a_jz^j$.
As in the cubic argument, majorization forces $a_0>0$; the cusp $z=0$
is not a contact.  On the positive half-axis the residual is
$a_0+t z^q+a_1z-z^{q+1}+a_2z^2+a_3z^3+a_4z^4$.
Its coefficient signs have at most four changes if $a_4>0$, and at most
three if $a_4=0$; zero coefficients are omitted.

The identity
$\int y(1-y^2)\,d\mu=(1-c)t>0$
forces a positive-$y$ interior contact.  This is a double positive-$z$
zero of the residual.  There can be at most two distinct positive-$z$
contacts, since three would count at least five zeros: at most one of
them can be the right endpoint.

Suppose first that there are two positive contacts $v<w\le1+t$.
The smaller contact is interior.  The derivative of the residual at
$w$ is zero if $w$ is interior and is nonpositive if $w=1+t$.
Let $H$ be the cubic Hermite interpolant of
$\varphi(z)=z^{q+1}-tz^q$ at $v,w$, and put
$\delta=\varphi'(w)-Q'(w-t)\ge0$.  The exact identity
$Q(z-t)=H(z)-\delta\frac{(z-v)^2(z-w)}{(w-v)^2} +a_4(z-v)^2(z-w)^2$
follows by matching values and first derivatives at both nodes.
As $\varphi'''(z)<0$ for $z>0$, comparison of cubic coefficients gives
$a_3=\varphi[v,v,w,w]-\frac{\delta}{(w-v)^2} -2a_4(v+w)<-2a_4(v+w)\le0$.
The value and derivative equations at $v$ then imply
\[
a_2=\frac{a_0}{v^2}-2a_3v-3a_4v^2 +qv^{q-1}+(1-q)t v^{q-2} >qv^{q-1}+(1-q)t v^{q-2} \ge\frac{q+t}{(1+t)^{2-q}}.
\]
Here $-2a_3v-3a_4v^2>a_4(v^2+4vw)\ge0$, and the last two power
terms decrease strictly with $v$.

For negative $z=-u$, $0<u\le1-t$, set
$R(u)=Q(-t-u)-g(u)\ge0$.  The scalar bound
\eqref{sub12:eq:quartic-contact-global-curvature} now gives
$R''(u)=2a_2-6a_3u+12a_4u^2-g''(u)>0$.
Thus there is at most one negative-$z$ contact, and at most three
contacts in all.

If there is only one positive-$z$ contact, use instead
$g''''(u)=q(q-1)(q-2)u^{q-4} \bigl((q+1)u-t(3-q)\bigr)<0$.
Indeed the bracket is strictly less than $p-4t\le0$ in the open
interval $0<u<1-t$.
It follows that $R''''(u)=24a_4-g''''(u)>0$.
Three distinct negative-$z$ contacts would count at least five zeros,
because at most one can be the left endpoint and all interior contacts
are double.  Repeated Rolle's theorem would then force a zero of
$R''''$, which is impossible.  There are at most two negative contacts,
again giving a total of at most three.
\end{proof}

The nonnegative leading coefficient is precisely the dual sign associated
with an upper bound on the fourth moment.  This contact-count statement
does not by itself evaluate the resulting three-atom optimization.
The boundary classification in Corollary~\ref{sub12:cor:turn3-global-lp4-boundary}
and the integer-cell reduction in Theorem~\ref{sub12:thm:turn3-peak-cell-sufficient}
provide the additional comparison needed to evaluate the capped problem.

\paragraph{Dual attainment and the actual three-atom reduction.}
The quartic contact-count argument also constrains the support of
optimizers for the capped moment problem. Fix $0<c<1$,
$p/4\le t\le c$, $d=-(1-c)t$, and
$e\ge e_{\min}:=c^2+\frac{d^2}{c}$.
Among all probability measures on $[-1,1]$ with moments
$m_1=0$, $m_2=c$, $m_3=d$, and $m_4\le e$, the maximum of
$\int\phi_t\,d\mu$ is attained by a measure with at most three atoms.
In fact, when $t<c$ and $e>e_{\min}$, every maximizing measure is
supported on the at-most-three-point contact set of an attained quartic
dual certificate.

We first prove existence, including the boundary cases.  For every such
measure,
\begin{equation}\label{eq:minimum-fourth-moment}
 \int\left(y^2-\frac dc y-c\right)^2d\mu=m_4-c^2-\frac{d^2}{c}\ge0.
\end{equation}
Both roots of \(y^2-(d/c)y-c\) lie in \([-1,1]\), because
\(|d|\le c(1-c)\). On this pair of roots there is a unique probability
measure of mean zero. Its second, third, and fourth moments are
\(c,d,e_{\min}\), respectively. Thus the feasible set is nonempty and compact, and its
maximum exists.  If $e=e_{\min}$, \eqref{eq:minimum-fourth-moment} forces
this two-atom measure.  If $t=c$, then
$\int(1+y)(y-c)^2\,d\mu=0$,
which forces support in $\{-1,c\}$, again reducing the problem to two
atoms regardless of $e$.

Suppose now that $t<c$ and $e>e_{\min}$.  There is a strictly feasible
measure with an absolutely continuous component whose density is bounded
below by a positive constant on the entire interval.  To see this
explicitly, let $\sigma$ be the uniform measure
on $[-1,1]$ and take a small $\varepsilon>0$.  The residual moments
\[
c_\varepsilon=\frac{c-\varepsilon/3}{1-\varepsilon},\qquad d_\varepsilon=\frac{d}{1-\varepsilon}
\]
still satisfy $0<c_\varepsilon<1$ and
$|d_\varepsilon|<c_\varepsilon(1-c_\varepsilon)$.
Let $\gamma_\varepsilon$ be the corresponding mean-zero two-atom
measure.  Then
$\bar\mu=\varepsilon\sigma+(1-\varepsilon)\gamma_\varepsilon$
has the prescribed first three moments, and its fourth moment tends to
$e_{\min}$ as $\varepsilon\downarrow0$.  Choosing $\varepsilon$ small
therefore gives $\bar m_4<e$ and an absolutely continuous component
with density at least $\varepsilon/2$.

For completeness, finite-dimensional separation now supplies an attained
dual polynomial, rather than only an approximating sequence.  Let
$V(s_1,s_2,s_3,e')$ be the optimal integral for the indicated first three
moments and fourth-moment cap.  The strictly feasible full-support
measure just constructed can be perturbed by small polynomial densities,
so $V$ is finite on a neighborhood of $(0,c,d,e)$.  It is concave there
and on its convex feasible domain.  A supporting hyperplane to its
hypograph at this interior point has coefficients
$a_1,a_2,a_3,a_4$, with $a_4\ge0$ because $V$ is nondecreasing in the
cap.  Applying that supporting inequality to point masses gives
$\phi_t(y)\le Q(y):=a_0+a_1y+a_2y^2+a_3y^3+a_4y^4$,
where
$a_0=V(0,c,d,e)-a_2c-a_3d-a_4e$.
Thus $Q$ is dual feasible and its objective is exactly the primal
optimum.  For a maximizing measure,
$\int(Q-\phi_t)\,d\mu+a_4(e-m_4)=0$.
Both terms are nonnegative, so its support is contained in the contact
set of $Q$.  The contact theorem proves the asserted three-atom
reduction.  All separation used here is in the finite-dimensional
moment space; no compactness assertion about polynomial coefficients
or infinite-dimensional dual attainment is needed.
\subsection{A Gaussian Hermite criterion for the full reciprocal problem}
\label{sub12:sec:gaussian-negative-criterion}

Let $1/2\le c<1$, $0\le t\le c$, $0<q<1$, and put $p=q+1$.
Let $y_-<0<y_+$ be the roots of
\(y^2+\frac{(1-c)t}{c}y-c=0\).
The two nodes define the unique Gaussian measure with
$m_1=0$, $m_2=c$, and $m_3=-(1-c)t$.
Their shifted values have opposite signs; write
\[
 t+y_-=-u,\qquad t+y_+=v,\qquad u,v>0.
\]
Let $P_q$ be the cubic Hermite interpolant at $-u,v$ of
\(\Psi_{t,q}(x)=(x-t)\operatorname{sgn}(x)|x|^q\).
\begin{lem}[Gaussian Hermite criterion]\label{sub12:lem:gaussian-hermite-criterion}
The negative-branch Hermite inequality is
\begin{equation}\label{sub12:eq:gaussian-negative-majorant}
 P_q(x)\ge\Psi_{t,q}(x)\qquad(t-1\le x\le0).
\end{equation}
It is strict except at the prescribed negative interpolation node.
As a consequence, the Hermite polynomial is a majorant on the full
interval $[t-1,t+1]$ if and only if the single endpoint inequality
\begin{equation}\label{sub12:eq:gaussian-right-endpoint-test}
 P_q(1+t)\ge(1+t)^q
\end{equation}
holds.
\end{lem}

\begin{proof}
Put
$\rho=2-\frac1c\in[0,1)$, $S=v-u=(1+\rho)t$, and $P=uv=c-\rho t^2>0$.
For $s>0$, consider the rational kernel
$f_s(x)=\frac{(x-t)x}{s+x^2}$.
Direct Hermite interpolation at $-u,v$ gives the error
\begin{equation}\label{sub12:eq:gaussian-hermite-rational-error}
 H_s(x)-f_s(x)=
 \frac{(x+u)^2(x-v)^2 N(s,x)}
 {(s+x^2)(s+u^2)^2(s+v^2)^2},
\end{equation}
where
$N(s,x)=(s+tx)\bigl((s+P)^2-sS^2\bigr) +2sS(s+P)(t-x)$.
For $x<0$ and $t>0$, this cubic polynomial in $s$ has a negative constant
term and strictly positive coefficients of all positive powers.
Indeed those coefficients, apart from the leading coefficient $1$, are
\begin{align*}
 [s^2]N&=2c+(1-2\rho-\rho^2)t^2-(1+2\rho)tx
          \ge2c-2t^2>0,\\
 [s]N&=P^2+2(1+\rho)Pt^2
           -tx\bigl(2\rho P+(1+\rho)^2t^2\bigr)>0,\\
 [s^0]N&=txP^2<0.
\end{align*}
Thus $N(s,x)$ changes sign exactly once, from negative to positive,
as $s$ increases.

The Stieltjes representation
\(\operatorname{sgn}(x)|x|^q =\frac{\cos(\pi q/2)}\pi \int_0^\infty s^{\frac{q-1}{2}}\frac{x}{s+x^2}\,ds\)
shows that the Hermite error is the corresponding positive multiple of
the integral of \eqref{sub12:eq:gaussian-hermite-rational-error}.
The error integrals are absolutely convergent, so interpolation and this
integral can be interchanged.  The same formula is valid at $q=0$ away
from $x=0$, with $\operatorname{sgn}(x)$ in place of the power.

At $q=0$ the function on the negative half-axis is simply $-x+t$.
Its cubic Hermite interpolant at $-u,v$ satisfies the exact identity
$P_0(x)-(-x+t)= \frac{2(x+u)^2}{(u+v)^3} \left((2t-S)x+2v^2-t(u+3v)\right)$.
Since $2t-S=(1-\rho)t\ge0$, the bracket on $x\ge t-1$ is at least
$2(v-t)^2-(1-\rho)t$.
The positive Gaussian node $v-t=y_+$ decreases with $t$ and equals $c$
at $t=c$.  Hence $v-t\ge c$, and
$2(v-t)^2-(1-\rho)t \ge2c^2-(1-c)=(2c-1)(c+1)\ge0$.
Consequently the error integral at exponent $-1/2$ is nonnegative.

For completeness, the one-sign-change transfer is elementary.  Write
$W(s)=N(s,x)/[(s+x^2)(s+u^2)^2(s+v^2)^2]$ and let $s_0$ be its
unique sign-change point.  For $0<q<1$,
$\int_0^\infty s^{\frac{q-1}{2}}W(s)\,ds >s_0^{\frac{q}{2}}\int_0^\infty s^{-\frac{1}{2}}W(s)\,ds\ge0$,
because $W$ is negative below $s_0$ and positive above it.
When $t=0$, the rational numerator itself is positive and no transfer
is needed.  This proves \eqref{sub12:eq:gaussian-negative-majorant}, with
strictness away from the negative interpolation node.
The value at $x=0$ follows by continuity; it is strictly positive,
since a zero constant term would make the polynomial smaller than
$t|x|^q+|x|^{q+1}$ for negative $x$ sufficiently close to zero.
At $t=0$ the same positivity follows directly from the even Gaussian
Hermite polynomial.

On $x>0$ the residual has the form
$a_0+t x^q+a_1x-x^{q+1}+a_2x^2+a_3x^3$, with $a_0>0$.
It has at most three positive zeros by the generalized Descartes rule.
The positive Gaussian node $v$ is a prescribed double zero and satisfies
$v<1+t$.  If the residual is nonnegative at $1+t$, any negative component
between $0$ and $1+t$ would require at least two additional zeros,
counting endpoint multiplicities, for a total of at least four.
This is impossible.  Thus \eqref{sub12:eq:gaussian-right-endpoint-test} is
sufficient, and it is plainly necessary.
\end{proof}

\paragraph{Consequence when the endpoint criterion holds.}
If \eqref{sub12:eq:gaussian-right-endpoint-test} holds, the Gaussian two-point
measure maximizes $\int\phi_t$ among all measures with the first three
moments \(m_1=0\), \(m_2=c\), and \(m_3=-(1-c)t\).  It therefore also maximizes the capped
fourth-moment problem whenever the cap is at least its minimal fourth
moment.  Its normalized gradient value is explicitly
\(G_{\rm Gauss}(c,t)= \frac{2c}{(1-t)^q+(1+t)^q}\, \frac{u^q+v^q}{u+v}\).
For $c=(n-2)/n$ this is exactly the endpoint-gradient quotient of the
centered root configuration
\[
 t-1,\qquad
 -\frac{2t}{n-2}\quad\text{repeated }n-2\text{ times},\qquad
 t+1.
\]
The derivative estimate \eqref{sub12:eq:finite-outer-derivative-bound}
bounds this quotient by $\max\{S_{n,p},R_{n,p}(1)\}$.
We verify \eqref{sub12:eq:gaussian-right-endpoint-test} for
$0\le t\le1/2$ and every $1<p<2$.

\begin{thm}[Complete Gaussian solution in the low-midpoint region]
\label{sub12:thm:gaussian-low-midpoint-complete}
Let $1/2\le c<1$, $0\le t\le1/2$, and $0<q<1$.
The cubic Gaussian Hermite interpolant $P_q$ from
Lemma~\ref{sub12:lem:gaussian-hermite-criterion} majorizes
$\Psi_{t,q}$ throughout $[t-1,t+1]$.  The right endpoint inequality
\eqref{sub12:eq:gaussian-right-endpoint-test} is strict.
Consequently the Gaussian two-point measure is the unique maximizer of
$\int\phi_t\,d\mu$ under the first three prescribed moment constraints.
It also solves every feasible capped-fourth-moment version of this
problem.  In particular, for $n\ge4$ the full reciprocal-moment
relaxation in the region $0\le t\le1/2$ is bounded by
\[
 \max\{S_{n,p},R_{n,p}(1)\}.
\]
\end{thm}

\begin{proof}
By Lemma~\ref{sub12:lem:gaussian-hermite-criterion}, it suffices
to verify \eqref{sub12:eq:gaussian-right-endpoint-test}.
The case $t=0$ follows directly from the positive rational numerator, so
assume $0<t\le1/2$.  Use \(\rho,P,S\) from Lemma~\ref{sub12:lem:gaussian-hermite-criterion} and put $d=t(1+t)>0$.
At $x=1+t$ the rational Hermite numerator simplifies to
$N(s,1+t)=(s+d)(s+P)^2-sS^2(s+d)-2sS(s+P)$.
We first establish the two elementary estimates
\begin{equation}\label{sub12:eq:gaussian-two-amgm-parameters}
 3S^2<4P,\qquad 9S^2<16Pd.
\end{equation}
For the first, $4P-3S^2=4c-[4\rho+3(1+\rho)^2]t^2$
is decreasing in $t$, and at $t=1/2$ it equals
$\frac{(1-\rho)^2(3\rho+10)}{4(2-\rho)}>0$.
For the second, divide $16Pd-9S^2$ by $t$ and write
$F(t)=16(c-\rho t^2)(1+t)-9(1+\rho)^2t$.
This function is concave on $[0,1/2]$, because
$F''(t)=-32\rho-96\rho t\le0$.  Its endpoint values satisfy
$F(0)=16c>0$ and $F(1/2)=\frac{3(1-\rho)^2(3\rho+10)}{2(2-\rho)}>0$.
Hence $F(t)>0$ throughout the interval, proving
\eqref{sub12:eq:gaussian-two-amgm-parameters}.

For every $s>0$, the elementary arithmetic--geometric mean bounds give
$\frac{sS^2}{(s+P)^2}\le\frac{S^2}{4P}<\frac13$,
and
\(\frac{2sS}{(s+d)(s+P)} \le\frac{2S}{(\sqrt d+\sqrt P)^2} \le\frac{S}{2\sqrt{dP}}<\frac23\).
Therefore
$\frac{N(s,1+t)}{(s+d)(s+P)^2} =1-\frac{sS^2}{(s+P)^2} -\frac{2sS}{(s+d)(s+P)}>0$.
The positive Stieltjes representation now makes the right-endpoint
Hermite error strictly positive.  This proves global majorance.
The negative-branch strictness and the positive zero count show that
the only contact points are the two Gaussian nodes, so the optimizer
is unique.  For $c=(n-2)/n$, the Gaussian gradient is the derivative quotient
in \eqref{sub12:eq:finite-outer-derivative-quotient}; hence
\eqref{sub12:eq:finite-outer-derivative-bound} bounds its value by
$\max\{S_{n,p},R_{n,p}(1)\}$.
\end{proof}

\subsection{An explicit Jacobi reduction of the fourth-moment problem}
\label{sub12:sec:lp4-jacobi}

For $t\ge p/4$, the three-contact argument gives a Jacobi-matrix
parameterization of the capped moment problem.  Fix
$0<c<1$, $0\le t<c$, and put
$d=-(1-c)t$, $\beta=\frac dc$, and $e_{\min}=c^2+\frac{d^2}{c}$.
For a fourth-moment cap $e\ge e_{\min}$, set
$\gamma_{\max}=\frac ec-c-\beta^2$.
For $0\le\gamma\le\gamma_{\max}$ define
\begin{equation}\label{sub12:eq:lp4-jacobi-matrix}
 J_{\gamma,\delta}=
 \begin{pmatrix}
 0&\sqrt c&0\\
 \sqrt c&\beta&\sqrt\gamma\\
 0&\sqrt\gamma&\delta
 \end{pmatrix}.
\end{equation}
Its spectral measure at the first coordinate has moments
$m_1=0,\quad m_2=c,\quad m_3=c\beta=d,\quad m_4=c(c+\beta^2+\gamma)$.
The full matrix is a contraction, $-I\le J\le I$, precisely when
\begin{equation}
\label{sub12:eq:lp4-delta-interval} -1+\frac{\gamma}{1-c+\beta}\le\delta\le 1-\frac{\gamma}{1-c-\beta}.
\end{equation}
The denominators are positive because $t<c$.
The bounds in \eqref{sub12:eq:lp4-delta-interval} follow either from the Schur complements of
$I\pm J$, or by direct principal-minor calculations.  When $\gamma>0$,
irreducibility makes this equivalent to support of the first-coordinate
spectral measure in $[-1,1]$.  When $\gamma=0$, that measure ignores
$\delta$, which can always be chosen in the interval in \eqref{sub12:eq:lp4-delta-interval} without
changing the objective.

Conversely, every measure with at most three atoms and the indicated
first three moments is the first-coordinate spectral measure of
\eqref{sub12:eq:lp4-jacobi-matrix}. To see this, start by
orthogonalizing $1,y,y^2$ in its moment inner product: the first two
off-diagonal Jacobi coefficients are $\sqrt c$ and
$\sqrt{m_4/c-c-\beta^2}$.  If the second is zero, the first two
coordinates already give the unique two-atom measure, and the final
coordinate has zero spectral weight.

Consequently, for $t\ge p/4$ the capped moment optimum is exactly
$\max_{\gamma,\delta} \frac{2\langle e_1,\phi_t(J_{\gamma,\delta})e_1\rangle} {(1-t)^q+(1+t)^q}$,
where the parameters satisfy \(0\le\gamma\le e/c-c-\beta^2\) and \eqref{sub12:eq:lp4-delta-interval}.
The Jacobi parameterization is exact: each admissible pair
represents a feasible spectral measure, and every optimizer can
be represented by an admissible pair.

\paragraph{The parity cap.}
For $c=(n-2)/n$ and
$e=c-\frac{2t^2}{n}-\frac{2(1-u^2)}{n^2}$, where $|u|\le1$,
the active cap corresponds to $\gamma=\gamma_{\max}$, and its allowed
interval becomes
\begin{equation}
\label{sub12:eq:lp4-active-delta} \frac tc-\frac{1-u^2}{n(c-t)}\le\delta\le \frac tc+\frac{1-u^2}{n(c+t)}.
\end{equation}
Thus active-cap optimization requires only the one scalar $\delta$.
At a lattice endpoint $u^2=1$, the interval collapses to
$\delta=t/c$, and the three spectral nodes are $-1,t,1$.

\paragraph{Inactive caps reduce to boundary curves.}
If $t\ge p/3$ and an optimum has $\gamma<\gamma_{\max}$, complementarity
gives a cubic dual polynomial.  Such an optimizer cannot have three
interior spectral nodes.  Indeed two positive shifted contacts would
give four positive zeros of the cubic residual.  With only one positive
contact, if the shifted cubic coefficient is negative the negative-side
residual is strictly convex by the calculation in
\eqref{sub12:eq:high-t-scalar-curvature}. If the
coefficient is nonnegative, the negative-side residual has strictly
negative third derivative because $g'''(u)>0$ for $0<u<1-t$ when
$t\ge p/3$.  In either case two interior negative contacts are impossible.
It follows that the full optimum in this range occurs on one of the
following curves:
\[
\gamma=0,\qquad \gamma=\gamma_{\max},\qquad \delta=-1+\frac\gamma{1-c+\beta},\qquad \delta=1-\frac\gamma{1-c-\beta}.
\]
The first is the Gaussian two-atom measure; the last two are the endpoint
or Radau curves.  This observation does not assert that a curve maximum
occurs at its endpoints.

\paragraph{The inverse source measure.}
In the polynomial application, let $R(z)$ be the first-coordinate
resolvent of $J$.  The Cauchy transform of the remaining $n-2$ roots in
the relaxed source is
\[
T(z)=\frac{2(1-zR(z))}{(z^2-1)R(z)-z} =(n-2)\frac{z-\delta} {z^2-(\delta+\beta/(1-c))z+(\beta\delta-\gamma)/(1-c)}.
\]
For $\gamma>0$ this is a positive two-atom measure on $[-1,1]$; its
two nodes straddle $\delta$.  Its mean and second moment, after
normalization by total mass $n-2$, are
$-\frac tc$ and $\frac{t^2}{c^2}+\frac\gamma{1-c}$.
Under the active parity cap its total endpoint deficiency is exactly
$\int(1-y^2)\,dT(y)=1-u^2$.
The source also contains one unit atom at each of $-1,1$.
In general the two interior source weights are real, rather than
integers.  Thus a three-atom spectral optimizer is a relaxation to four
source values with the two endpoint weights fixed at one.  One must not
identify it with an actual uniform root list without checking those
weights.  If a spectral endpoint is present, one source pole coalesces
with that endpoint, giving a three-value source.  When both spectral
endpoints are present at a lattice value, the source is the corresponding
integer two-value configuration.

\paragraph{An exact interior stationarity equation.}
Let $y_1<y_2<y_3$ be the eigenvalues for $\gamma>0$, and assume none equals
$-t$.  Write
$P_\delta(z)=\det(zI-J) =z^3-(\beta+\delta)z^2+(\beta\delta-c-\gamma)z+c\delta$.
Differentiating the first-coordinate resolvent gives
$\partial_\delta R(z)=\frac{c\gamma}{P_\delta(z)^2}$.
Hence the derivative of its objective numerator is the fifth divided
difference
\begin{equation}
\label{sub12:eq:lp4-delta-stationary} \partial_\delta\langle e_1,\phi_t(J)e_1\rangle =c\gamma\,\phi_t[y_1,y_1,y_2,y_2,y_3,y_3].
\end{equation}
The formula is valid through Hermite interpolation of the scalar
functional calculus; no analyticity across the cusp is assumed.
An interior maximizer in $\delta$ cannot have an eigenvalue equal to
$-t$: that eigenvalue has strictly positive velocity under variation
of $\delta$, and its contribution has a nonzero odd $|s|^q$ cusp, which
allows an increase in one direction.

Equations \eqref{sub12:eq:lp4-active-delta} and
\eqref{sub12:eq:lp4-delta-stationary} make the remaining active-cap scalar
optimization explicit. Its comparison with the sharp discrete constant
is handled by Theorem~\ref{sub12:thm:turn3-peak-cell-sufficient}.
\subsection{A global three-atom upper bound for every finite degree}
\label{sub12:sec:global-three-atom-upper}

Let $n\ge4$, $1<p<2$, $q=p-1$, and $c=(n-2)/n$.
For \(0\le t\le c\), let \(u_n(t)=-nt-2m+n-3\), choosing the integer \(m\)
as in \eqref{sub12:eq:finite-gradient-m4-interval}, and set
$e_n(t)=c-\frac{2t^2}{n}-\frac{2(1-u_n(t)^2)}{n^2}$, $\beta(t)=-\frac{(1-c)t}{c}$,
and $\Gamma_n(t)=\frac{e_n(t)}c-c-\beta(t)^2$.
At a cell boundary either choice of \(m\) gives \(u_n(t)^2=1\), so
\(e_n\) and \(\Gamma_n\) are continuous. The quantity \(\Gamma_n(t)\) is nonnegative: the parity cap
is at least the Gaussian minimum fourth moment.

Define the compact parameter set
\begin{equation}\label{sub12:eq:global-jacobi-compact-domain}
\begin{aligned}
\mathcal D_n=\{(t,\gamma,\delta):\;&0\le t\le c,\quad
 0\le\gamma\le\Gamma_n(t),\quad -1\le\delta\le1,\\
 &\gamma\le(1+\delta)(1-c+\beta(t)),\\
 &\gamma\le(1-\delta)(1-c-\beta(t))\}.
\end{aligned}
\end{equation}
The product inequalities in \eqref{sub12:eq:global-jacobi-compact-domain} remain meaningful at $t=c$.
For each such triple put
\[
J(t,\gamma,\delta)=
 \begin{pmatrix}
 0&\sqrt c&0\\
 \sqrt c&\beta(t)&\sqrt\gamma\\
 0&\sqrt\gamma&\delta
 \end{pmatrix},
\qquad
 \phi_t(y)=y\operatorname{sgn}(t+y)|t+y|^q,
\]
and define the explicit finite-dimensional upper bound
\begin{equation}
\label{sub12:eq:global-three-atom-upper} V_{n,p}=\max_{(t,\gamma,\delta)\in\mathcal D_n} \frac{2\langle e_1,\phi_t(J(t,\gamma,\delta))e_1\rangle} {(1-t)^q+(1+t)^q}.
\end{equation}

\begin{thm}[Exact global reduction of the reciprocal relaxation]
\label{sub12:thm:global-three-atom-reciprocal}
The quantity $V_{n,p}$ is exactly the global optimum of the full
reciprocal-moment relaxation, equivalently of the spectral moment
problem with
\[
m_1=0,\quad m_2=c,\quad m_3=-2t/n,\quad m_4\le e_n(t), \quad\operatorname{supp}\mu\subset[-1,1].
\]
For every $t$ this problem has an optimizer with at most three atoms.
Moreover
\begin{equation}\label{eq:reciprocal-relaxation-bounds}
\Lambda_{n,p}\le C_{n,p}\le V_{n,p}\le U_{n,p} \le1-\frac{H_p(r_p)}{n},
\end{equation}
where \(\Lambda_{n,p}=\max\{S_{n,p},\mathcal B_{n,p}\}\). To define
\(U_{n,p}\), let \(D=(1-t)^q+(1+t)^q\) and let \(K_t\) be the divided
difference in Subsection~\ref{sub12:sec:mean-only-high-asymmetry}. Then
\begin{equation}\label{eq:mean-variance-relaxation}
 U_{n,p}=\max_{t,\nu}\left(1-\frac{4}{nD}\int K_t\,d\nu\right),
 \qquad \int y\,d\nu=t,\quad \int(y-t)^2\,d\nu\ge\frac{1-u_n(t)^2}{n}.
\end{equation}
The maximum is over \(0\le t\le c\) and probability measures on \([-1,1]\).
\end{thm}

\begin{proof}
For $0\le t\le1/2$, Theorem~\ref{sub12:thm:gaussian-low-midpoint-complete}
gives a two-atom optimizer.  For $t>1/2$, one has $t>p/4$, so the quartic contact-count argument in
Subsection~\ref{sub12:sec:quartic-three-contacts} gives an optimizer with at most three
atoms.  This covers every $t\in[0,c]$.

The Jacobi representation of a measure with at most three atoms and
the first three indicated moments is precisely $J(t,\gamma,\delta)$.
Its fourth moment is $c(c+\beta^2+\gamma)$, giving
$\gamma\le\Gamma_n(t)$; the two Schur-complement inequalities for
$-I\le J\le I$ are exactly
\eqref{sub12:eq:global-jacobi-compact-domain}.  When $\gamma=0$, the spectral
measure at $e_1$ ignores the last coordinate, so $\delta$ is redundant.
At $t=c$, the first Schur-complement inequality forces $\gamma=0$;
the unique measure is the two-point measure on $\{-1,c\}$.
Thus all degeneracies are included without division by zero.
Conversely every triple in $\mathcal D_n$ supplies a feasible spectral
measure.  The domain is compact and the scalar functional calculus is
continuous, including at the cusp, so the maximum in \eqref{sub12:eq:global-three-atom-upper} is attained.
This proves the exact reduction of the relaxation.

The original root configurations give feasible spectral measures.
At a real maximizing root list the moment ratio equals its
endpoint-gradient quotient; hence $C_{n,p}\le V_{n,p}$, and the
angular reduction transfers this bound to complex roots.
For a feasible spectral measure \(\mu\), the probability measure
\(d\nu=(n/2)(1-y^2)\,d\mu\) has mean \(t\). Its variance satisfies
\eqref{eq:mean-variance-relaxation} by the fourth-moment cap.
Writing \(\phi_t\) as its affine endpoint interpolant minus \((1-y^2)K_t\)
shows that its normalized integral is \(1-4\int K_t\,d\nu/(nD)\).
This proves \(V_{n,p}\le U_{n,p}\). The sparse and two-value root lists
give \(\Lambda_{n,p}\le C_{n,p}\).

To prove \(U_{n,p}\le1-H_p(r_p)/n\), fix a probability measure \(\nu\) of mean \(t\) supported
on a compact subinterval of \((-1,1)\). Choose an integer \(N\) large enough
that \(N\theta_p>1\), \(t<(N-2)/N\), and both endpoint masses in the
following construction are nonnegative. On \((-1,1)\) take
\(d\mu=2\,d\nu/[N(1-y^2)]\), and give \(\pm1\) masses
\(\frac12-N^{-1}\int(1\mp y)^{-1}\,d\nu\), respectively.
Then \(\mu\) is a probability measure with \(m_1=0\), \(m_2=1-2/N\),
and \(m_3=-2t/N\). The independent three-moment bound of
Theorem~\ref{sub12:thm:turn3-uncapped-lp3-exact}, applied in degree \(N\),
gives \(1-4\int K_t\,d\nu/(ND)\le1-H_p(r_p)/N\). Thus
\(4\int K_t\,d\nu/D\ge H_p(r_p)\). For an arbitrary measure of mean \(t\),
apply the same argument to its image under \(y\mapsto(1-\varepsilon)y+\varepsilon t\)
and let \(\varepsilon\downarrow0\). Continuity of \(K_t\) on \([-1,1]\)
gives the bound for every measure in \eqref{eq:mean-variance-relaxation}.
\end{proof}

\paragraph{The remaining domain is entirely above $t=1/2$.}
Theorem~\ref{sub12:thm:gaussian-low-midpoint-complete}
bounds the region $0\le t\le1/2$ of
\eqref{sub12:eq:global-three-atom-upper}
by $\max\{S_{n,p},R_{n,p}(1)\}$.  The sparse value occurs at $t=0$,
and the one-outlier value is represented at $t=c\ge1/2$.
Consequently one may equivalently write
\[
V_{n,p}=\max\left\{S_{n,p},\
 \max_{\substack{(t,\gamma,\delta)\in\mathcal D_n\\t\ge1/2}}
 \frac{2\langle e_1,\phi_t(J)e_1\rangle}
 {(1-t)^q+(1+t)^q}\right\}.
\]
Corollary~\ref{sub12:cor:turn3-global-lp4-boundary} reduces the
Jacobi parameter domain to its boundary families, and
Theorem~\ref{sub12:thm:turn3-peak-cell-sufficient}
compares the remaining maximum with integer two-value configurations.

\begin{cor}[Rigidity of equality between the two relaxations]
\label{sub12:cor:two-relaxations-equality-rigidity}
Write $\mathcal B_{n,p}=\max_kR_{n,p}(k)$.  Then
\[
V_{n,p}=U_{n,p}\quad\Longleftrightarrow\quad U_{n,p}=\mathcal B_{n,p}.
\]
In that event the original optimal constant is determined exactly:
$C_{n,p}=V_{n,p}=U_{n,p}=\mathcal B_{n,p}$.
If $U_{n,p}>\mathcal B_{n,p}$, then $V_{n,p}<U_{n,p}$ strictly.
\end{cor}

\begin{proof}
Suppose \(V_{n,p}=U_{n,p}\), and choose a maximizing spectral measure
\(\mu\) for \eqref{sub12:eq:global-three-atom-upper}. The associated
\(d\nu=(n/2)(1-y^2)\,d\mu\) has no endpoint atoms and minimizes
\(\int K_t\,d\nu\) at the same \(t\) in \eqref{eq:mean-variance-relaxation}.
If \(t=0\), then \(K_0(y)=(1-|y|^p)/(1-y^2)>p/2\) on \((-1,1)\).
The measure with equal masses at \(-1,1\) is feasible and gives \(p/2\).
Thus \(t>0\).

There is a quadratic minorant of \(K_t\), with nonnegative quadratic
coefficient, that agrees with \(K_t\) on the support of \(\nu\).
Indeed, at fixed mean \(t\) the attainable pairs
\((\int(y-t)^2\,d\nu,\int K_t\,d\nu)\) form a compact convex set.
The endpoint mixture of mean \(t\) has variance \(1-t^2\ge1-c^2>1/n\).
Separation at the constrained minimum therefore gives a finite multiplier
\(\lambda\ge0\) for the variance lower bound. For
\(K_t(y)-\lambda(y-t)^2\), take a supporting line to its lower convex
envelope at the interior point \(t\). Adding \(\lambda(y-t)^2\) back
gives the required minorant, with equality \(\nu\)-almost everywhere.

By convexity, the minorant lies below the chord joining the endpoint
values of \(K_t\). Equation~\eqref{sub12:eq:global-negative-K-chord} excludes
contact on \((-1,-t]\). On \(y>-t\), multiply the gap between \(K_t\)
and the minorant by \(1-y^2\) and set \(z=t+y\). The resulting residual
has the form \(a_0+tz^q+a_1z-z^{q+1}+a_2z^2+a_3z^3+\lambda z^4\),
with \(a_0>0\) by the strict chord inequality at \(y=-t\).
Its coefficients have at most four sign changes. Two interior contacts
would give two double positive zeros in addition to the zero at \(z=1+t\),
contrary to the generalized Descartes rule used in
Theorem~\ref{sub12:thm:quartic-anchor}. Hence \(\nu=\delta_t\).
Its zero variance forces \(u_n(t)^2=1\) in \eqref{eq:mean-variance-relaxation}.
Thus \(t=1-2k/n\) for an integer \(k\), and the objective equals \(R_{n,p}(k)\).
Since \(R_{n,p}(k)\le\mathcal B_{n,p}\), we have \(U_{n,p}\le\mathcal B_{n,p}\).
The reverse inequality follows from \eqref{eq:reciprocal-relaxation-bounds}.
Conversely, if \(U_{n,p}=\mathcal B_{n,p}\), every quantity between
\(\mathcal B_{n,p}\) and \(U_{n,p}\) in \eqref{eq:reciprocal-relaxation-bounds}
equals their common value. Hence \(U_{n,p}>\mathcal B_{n,p}\) forces
\(V_{n,p}<U_{n,p}\).
\end{proof}
\subsection{A Gaussian wedge and a global boundary reduction for the fourth-moment problem}
\label{sub12:sec:turn3-gaussian-wedge}

The following result removes all inactive-cap, three-interior-node
optimizers from the reciprocal relaxation, in every degree and for every
$1<p<2$.  It does not assert the remaining discrete-envelope comparison.

\begin{thm}[Gaussian optimality for $t\le p/3$]
\label{sub12:thm:turn3-gaussian-wedge}
Let $1/2\le c<1$, $0<q<1$, $p=q+1$, and
$0\le t\le\min\{c,p/3\}$.  Let $\mu$ range over probability measures
on $[-1,1]$ satisfying
\[
m_1=0,\qquad m_2=c,\qquad m_3=-(1-c)t.
\]
The Gaussian two-point measure on the roots of
\(y^2+\frac{(1-c)t}{c}y-c=0\)
is the unique maximizer of
$\int y\operatorname{sgn}(t+y)|t+y|^q\,d\mu(y)$.
The Gaussian measure also uniquely maximizes every feasible version
of the problem with an upper bound on the fourth moment.
\end{thm}

\begin{proof}
We use the negative-branch Hermite inequality
\eqref{sub12:eq:global-negative-hermite-gap}: the cubic
Hermite interpolant at the two shifted Gaussian nodes majorizes the
objective on the negative shifted half-axis, and global majorance is
equivalent to the right-endpoint inequality.  We prove that endpoint
inequality strictly under the present assumptions.

Write the shifted Gaussian nodes as $-u,v$, with $u,v>0$, and put
$\rho=2-\frac1c$ and $T=u+v$, so that $v-u=(1+\rho)t$ and $uv=c-\rho t^2$.
Normalize by $v$ and set
$r=\frac uv$, $\tau=\frac tv$, and $\varepsilon=\frac{1+t-v}{v}$.
Let $h_q(z)=\operatorname{sgn}(z)|z|^q$, put
$g_\tau(z)=(z-\tau)h_q(z)$, and let $H_\tau$ be its cubic Hermite
interpolant at $-r,1$.  The sign of the endpoint error is the sign of
$E_\tau(\varepsilon)=H_\tau(1+\varepsilon)-g_\tau(1+\varepsilon)$.
The physical parameters satisfy
\begin{equation}
\label{sub12:eq:turn3-wedge-normalized-domain} \frac15\le r\le1,\qquad \tau\ge\frac{1-r}{2},\qquad 0<\varepsilon\le r,\qquad \frac r{1+r}\ge\frac{2-q}{6}.
\end{equation}
The Gaussian formula
has $T^2=4c+a^2t^2$, where $a=(1-c)/c\in(0,1]$.
The difference
$T^2-(1+t)^2=4c-1-2t-(1-a^2)t^2$
is decreasing in $t$ and equals zero at $t=c$.  Thus $T\ge1+t$,
which is exactly $\varepsilon\le r$.  The positive Gaussian node is
strictly below $1$, so $\varepsilon>0$.
Furthermore
$T\ge2\sqrt c\ge3-\frac1c=1+\rho$.
For the second inequality, multiply by $c$ and note that
$2c^{\frac{3}{2}}-3c+1$ decreases to zero on $[1/2,1]$.
It follows that
$\frac{1-r}{1+r}=\frac{v-u}{T}\le t$ and $\frac r{1+r}\ge\frac{1-t}{2}\ge\frac{2-q}{6}$.
This also gives $r\ge1/5$.
Finally $1-r=(1+\rho)\tau\le2\tau$.

We first show that $E_\tau(\varepsilon)$ is strictly increasing in
$\tau$ for every fixed $r\in[1/5,1]$ and $\varepsilon>0$.
Let $K$ be the cubic Hermite interpolant of $h_q$ at $-r,1$.
By linearity,
$\partial_\tau E_\tau(\varepsilon) =h_q(1+\varepsilon)-K(1+\varepsilon)$.
Put $\alpha=1-q\in(0,1)$, $L=1+r$, and $s=\log(1/r)<2$.
The cubic coefficient $k_3$ of $K$ obeys
$L^3 k_3=(1-r)(e^{\alpha s}-1) -\alpha L(1+e^{\alpha s})<0$,
because $\tanh(\alpha s/2)\le\alpha s/2$ and
$(1-r)s<2L$.
The function $(1-\alpha)e^{\alpha s}$ is concave in $\alpha$,
since its second derivative is
$s e^{\alpha s}\{s(1-\alpha)-2\}\le0$.
Consequently
$q r^{q-1}\le1+\alpha(s-1)$ and $r^q\ge r(1+\alpha s)$.
Using the exact Hermite formula at $1$ gives
\[
L^2\{K''(1)-h_q''(1)\} =2qLr^{q-1}+4qL-6(1+r^q)+q\alpha L^2 \le\alpha\{2(1-2r)s+r^2-4r-5\}<0.
\]
The last strict sign is immediate when $r\ge1/2$.  When $r\le1/2$,
use $s<2$ to bound its braces by $r^2-12r-1<0$.
Since $h_q'''(z)>0$ for $z>0$ and $K'''=6k_3<0$, the residual
$K-h_q$ has negative third derivative on $[1,\infty)$, negative
second derivative at $1$, and zero value and first derivative there.
Thus $K(z)<h_q(z)$ for $z>1$, proving the claimed monotonicity.
It suffices to consider
$\tau_0=\frac{1-r}{2}$.

For this value, extend
$A(\varepsilon)=\frac{E_{\tau_0}(\varepsilon)}{\varepsilon^2}$
continuously to $\varepsilon=0$.
On the positive axis $g_{\tau_0}''''>0$; the integral Taylor formula
therefore shows that $A$ is strictly concave:
\[
A(\varepsilon)=\frac{H_{\tau_0}''(1)}2+a_3\varepsilon -\int_0^1(1-z)g_{\tau_0}''(1+z\varepsilon)\,dz, \qquad A''(\varepsilon)<0.
\]
Here
$a_3=\frac{q(1-r^{q-1})}{2L}$ and $2L A(0)=1+r^q-qLr^{q-1}+\frac{q\alpha L^2}{2}$.
For $r\ge1/3$, the Hermite coefficient estimate and
$r^q\ge r(1+\alpha s)$ give
$2L A(0)\ge\alpha\left(L-s+\frac{qL^2}{2}\right)>0$,
since $L\ge4/3>\log3\ge s$.
For $1/5\le r\le1/3$, one instead has
$A'(0)\le\frac{q\alpha}{12L} \{L(5-r)-6\log(1/r)\}<0$.
Indeed $6\log(1/r)-(1+r)(5-r)$ decreases in $r$, and at $r=1/3$
it is positive because $\log3>13/12>28/27$.
The first strict logarithmic bound follows from the first two terms
of $\log3=2\sum_{j=0}^{\infty}(1/2)^{2j+1}/(2j+1)$.
As $A$ is concave, $A'(0)<0$ implies $A(0)\ge A(r)$ in this range.
Thus, in both ranges, it is enough to show $A(r)>0$.

Put $x=r/(1+r)\in[1/6,1/2]$.  Direct Hermite evaluation gives
\begin{equation}
\label{sub12:eq:turn3-wedge-power-error} E_{\tau_0}(r)=\frac{(1+r)^{q+1}}2 F_q(x),
\end{equation}
where
$F_q(x)=x^{q+1}\{(1-q)x-q\} +(1+x)^2(1-x)^{q-1}\{1-(1-q)x\}-(1+2x)$.
By \eqref{sub12:eq:turn3-wedge-normalized-domain},
$\alpha=1-q\le6x-1$.
We now prove $F_q(x)>0$ on this entire domain.
Write $a=-\log x$, $b=-\log(1-x)$, so that
$F_q(x)=-x^2\{1-(1+x)\alpha\}e^{a\alpha} +(1+x)^2(1-x\alpha)e^{b\alpha}-(1+2x)$.
The function $(1-(1+x)\alpha)e^{a\alpha}$ is strictly concave
in $\alpha\in[0,1]$, since $a\le\log6<2(1+x)$.
The tangent at zero gives
\begin{equation}\label{eq:gaussian-tangent}
 -x^2\{1-(1+x)\alpha\}e^{a\alpha}\ge-x^2+\alpha x^2(1+x-a).
\end{equation}
The second function $h(\alpha)=(1-x\alpha)e^{b\alpha}$ has
$h(0)=h(1)=1$ and
$-h''(\alpha)=b e^{b\alpha}(2x-b+bx\alpha) \ge b(2x-b)>0$,
because $b<2x$ for $0<x\le1/2$.
Thus \(h(\alpha)\ge1+\frac12 b(2x-b)\alpha(1-\alpha)\).
For \(x\ge1/3\), substituting this bound for \(h\) and \eqref{eq:gaussian-tangent}
in \(F_q\) gives \(F_q(x)\ge\alpha x^2(1+x-a)>0\), since
\(1+x-a\ge4/3-\log3>0\).
If $1/6\le x\le1/3$, use $1-\alpha\ge2(1-3x)$ to obtain
$\frac{F_q(x)}\alpha\ge x^2(1+x+\log x) +(1+x)^2b(2x-b)(1-3x)$.
The power series for $b$ gives
$0\le\frac bx-1\le\frac{x}{2(1-x)}$ and $b(2x-b)\ge x^2\left(1-\frac{x^2}{4(1-x)^2}\right)$.
Hence $F_q(x)\ge\alpha x^2 M(x)$, where
$M(x)=1+x+\log x+(1+x)^2(1-3x) \left(1-\frac{x^2}{4(1-x)^2}\right)$.
The function \(M\) is strictly concave on $[1/6,1/3]$.
For a completely rational verification, put $x=(1+z)/6$, $0\le z\le1$.
Direct differentiation gives
$2x^2(1-x)^4 M''(x) =\frac{-z^7+15z^6-57z^5-129z^4+2933z^3-747z^2+3397z-13603}{10368}<0$.
For \(0\le z\le1\), the numerator of \(2x^2(1-x)^4M''(x)\) is at most
$15+2933+3397-13603=-7258$.
Finally,
$M(1/6)=\frac{4417}{2400}-\log6>0$ and $M(1/3)=\frac43-\log3>0$.
For example $\log6<11/6<4417/2400$, as follows from
$\sum_{j=0}^4(11/6)^j/j!=186985/31104>6$.
Also $e^{\frac{4}{3}}>1+4/3+(4/3)^2/2>3$.
Concavity now proves $M>0$ throughout its interval, and hence $F_q>0$.

Equation \eqref{sub12:eq:turn3-wedge-power-error} gives $A(r)>0$.
The strict concavity of $A$, together with the endpoint bounds
$A(0)\ge A(r)>0$ when $r\le1/3$ and $A(0)>0$ when $r\ge1/3$,
yields $A(\varepsilon)>0$ for every
$0\le\varepsilon\le r$.  Monotonicity in $\tau$ proves the strict
physical right-endpoint inequality.  The negative-branch criterion then
gives a global cubic majorant whose only contacts are the two Gaussian
nodes.  Integrating it against measures with the specified first three
moments proves Gaussian optimality and uniqueness.  The Gaussian measure
has the least possible fourth moment, so it is feasible for every
nonempty fourth-cap problem and proves the final assertion as well.
\end{proof}

\begin{cor}[Global boundary classification of the reciprocal relaxation]
\label{sub12:cor:turn3-global-lp4-boundary}
For every $n\ge4$ and $1<p<2$, the complete reciprocal/fourth-moment
relaxation has an optimizer on one of the four boundary families
\[
\gamma=0,\qquad \gamma=\Gamma_n(t),\qquad \delta=-1+\frac{\gamma}{1-c+\beta(t)},\qquad \delta=1-\frac{\gamma}{1-c-\beta(t)}.
\]
The last two families are understood through their cross-multiplied
contraction conditions at degenerate endpoints.
Consequently its global optimum is the maximum over the four boundary
surfaces; an inactive-cap three-node optimizer with all nodes in
$(-1,1)$ never needs to be retained.
\end{cor}

\begin{proof}
For $t\le p/3$, Theorem~\ref{sub12:thm:turn3-gaussian-wedge} supplies the
Gaussian family $\gamma=0$.
For $t>p/3$, the attained quartic dual has zero leading coefficient
whenever the fourth-moment cap is inactive. The cubic contact-count
argument in Subsection~\ref{sub12:sec:lp4-jacobi} excludes three interior
contacts in this range.  A two-node feasible measure is necessarily
the Gaussian measure; a three-node measure must therefore have an
endpoint node and lies on one of the two Radau boundaries.
An active cap gives the second displayed family.
At $t=c$, the moment problem is already the unique two-node Gaussian
measure.  This covers all parameters and all degeneracies.
\end{proof}

Theorem~\ref{sub12:thm:turn3-peak-cell-sufficient} will compare the remaining
boundary maxima with the discrete two-value ratios.
\subsection{Elimination of the left Radau surface}
\label{sub12:sec:turn3-left-radau-elimination}

Let $1/2\le c<1$, $0\le t<c$, $0<q<1$, and set
$\beta=-\frac{(1-c)t}{c},\quad h=1-c+\beta>0,\quad \phi_t(y)=y\operatorname{sgn}(t+y)|t+y|^q$.
Consider the left Radau Jacobi family
\[
J_\gamma=\begin{pmatrix}
 0&\sqrt c&0\\
 \sqrt c&\beta&\sqrt\gamma\\
 0&\sqrt\gamma&-1+\gamma/h
 \end{pmatrix},\qquad
 0\le\gamma\le\overline\gamma,
\]
where $\overline\gamma$ is any nonnegative bound for which all the
matrices are contractions.  Thus the spectral measure has fixed moments
$m_1=0$, $m_2=c$, $m_3=c\beta$ and contains the left endpoint $-1$.

\begin{thm}\label{sub12:thm:turn3-left-radau-endpoints}
The maximum of $\langle e_1,\phi_t(J_\gamma)e_1\rangle$ on this interval
is attained at $\gamma=0$ or $\gamma=\overline\gamma$.
Every nondegenerate interior critical point is a strict local minimum.
In particular, in the fourth-moment relaxation the whole left Radau
surface may be discarded except for its Gaussian endpoint and its
intersection with the active fourth-moment cap.
\end{thm}

\begin{proof}
Write $\delta=-1+\gamma/h$.  The characteristic polynomial factors as
$\det(zI-J_\gamma) =(z+1)\{z^2-(1+\beta+\delta)z+c\delta\}$.
For $\gamma>0$ the other two eigenvalues $a<b$ satisfy
$a<c<b$, where $b=c+\frac{(1-c)(c-t)}{c-a}$.
Both are strictly increasing functions of $\delta$, because
$\frac{da}{d\delta}=\frac{c-a}{b-a}>0$ and $\frac{db}{d\delta}=\frac{b-c}{b-a}>0$.
At $\gamma=0$ they are the two Gaussian nodes, and therefore
$a+b\ge\beta$ and $b\ge c$.
Differentiating the first-coordinate resolvent gives
$\frac d{d\gamma}\langle e_1,\phi_t(J_\gamma)e_1\rangle =c\,\phi_t[-1,a,a,b,b]$
whenever $a\ne-t$.  The removable $\gamma=0$ limit is not needed below.

\paragraph{The interval where $a<-t$.}
Put
$A=1-t,\quad u=-t-a>0,\quad v=t+b>0, \quad S=v-u,\quad P=uv,\quad x=-A$.
Here $0<u<A$, and
$t\le S\le1+t$.
Indeed $S=2t+a+b\ge2t+\beta=(3-1/c)t\ge t$;
the upper bound follows from $u\ge0$ and $b\le1$.

Let $H_q$ be the cubic Hermite interpolant at $-u,v$ of
$f_q(z)=(z-t)\operatorname{sgn}(z)|z|^q$.
We prove $H_q(-A)>f_q(-A)$.  For the rational kernel
$f_s(z)=\frac{(z-t)z}{s+z^2}$,
the Hermite error is
$H_s(x)-f_s(x)= \frac{(x+u)^2(x-v)^2 N(s,x)} {(s+x^2)(s+u^2)^2(s+v^2)^2}$,
where
$N(s,x)=(s+tx)\{(s+P)^2-sS^2\} +2sS(s+P)(t-x)$.
At $x=-A$ its coefficients are
\begin{align*}
 [s^3]N&=1,\\
 [s^2]N&=2P-S^2+2S-t(1-t)>0,\\
 [s]N&=P^2+2SPt+x\{2P(t-S)-tS^2\}>0,\\
 [s^0]N&=txP^2\le0.
\end{align*}
For the quadratic coefficient, concavity of $S(2-S)$ on
$[t,1+t]$ gives
$S(2-S)-t(1-t)\ge\min\{t,1-t\}\ge0$.
For the linear coefficient use $S\ge t$ and $x<0$.
Thus when $t>0$ the rational numerator has exactly one positive zero,
and changes sign there from negative to positive.  When $t=0$ it is
strictly positive for all $s>0$.

The Hermite error at exponent $q=0$ has the exact expression
$H_0(x)-(-x+t)=\frac{2(x+u)^2}{(u+v)^3} \{(2t-S)x+2v^2-t(u+3v)\}$.
At $x=t-1$ the bracket simplifies to
$2b^2+a+b\ge2c^2+\beta \ge2c^2-(1-c)=(2c-1)(c+1)\ge0$.
Consequently the error integral with weight $s^{-\frac{1}{2}}$ is nonnegative.
For \(0<q<1\), use the Stieltjes representation
\(\operatorname{sgn}(z)|z|^q=\frac{\cos(\pi q/2)}{\pi}
\int_0^\infty s^{\frac{q-1}{2}}z/(s+z^2)\,ds\).
If \(t>0\), let \(s_0\) be the unique sign change of the rational error
kernel. Multiplying that kernel by \(s^{\frac{q}{2}}-s_0^{\frac{q}{2}}\)
gives a strictly positive product on both sides of \(s_0\).
Integration against \(s^{-\frac{1}{2}}\,ds\) therefore makes the new error
strictly greater than \(s_0^{\frac{q}{2}}\) times the nonnegative error
integral at \(q=0\). All integrals converge absolutely. For \(t=0\)
the rational kernel itself is positive.

Since the divided difference is the Hermite error with the opposite
sign,
$\phi_t[-1,a,a,b,b] =\frac{f_q(-A)-H_q(-A)}{(A-u)^2(A+v)^2}<0$.
The Radau objective is therefore strictly decreasing wherever $a<-t$.

\paragraph{The interval where $a>-t$.}
Suppose an interior stationary point occurs here.  Let $Q$ be the cubic
Hermite interpolant of $\phi_t$ at $a,b$.  Stationarity says exactly
that $Q(-1)=\phi_t(-1)$ as well.  On $[a,b]$ the shifted variable
$z=t+y$ is positive, and
$\phi_t^{(4)}(y) =q(q-1)(q-2)z^{q-4}\{(q+1)z+(3-q)t\}>0$.
The Hermite remainder, or its repeated-node divided differences,
therefore gives
$\phi_t''(a)-Q''(a)>0$ and $\phi_t''(b)-Q''(b)>0$.
Let $w_a,w_b>0$ be the two movable spectral weights.  Since the
expectation of every cubic polynomial is fixed along the family,
differentiating twice at the stationary point yields
\[
\frac{d^2}{d\gamma^2}\langle e_1,\phi_t(J_\gamma)e_1\rangle =w_a(a')^2\{\phi_t''(a)-Q''(a)\} +w_b(b')^2\{\phi_t''(b)-Q''(b)\}>0.
\]
All terms involving derivatives of the weights vanish because of the
contact conditions, and the node $-1$ is fixed.  Hence every such
stationary point is a strict local minimum.

If $a=-t$ and $t>0$, its contribution has a nonzero odd
$\operatorname{sgn}(\gamma-\gamma_0)|\gamma-\gamma_0|^q$ leading term,
with negative coefficient.  It is not a local maximum.  For $t=0$,
the increasing node $a$ reaches zero only at the full contraction
endpoint.  There are consequently no interior local maxima, proving
the assertion.
\end{proof}

\begin{cor}[Three boundary surfaces for the full fourth-moment envelope]
\label{sub12:cor:turn3-three-boundary-surfaces}
For every $n\ge4$ and $1<p<2$, put $c=(n-2)/n$ and retain the parity
cap $\Gamma_n(t)$ from the fourth-moment Jacobi reduction.  The full
relaxed envelope $V_{n,p}$ is attained on one of
\[
\gamma=0,\qquad \gamma=\Gamma_n(t),\qquad \delta=1-\frac{\gamma}{1-c-\beta}, \quad\beta=-\frac{(1-c)t}{c}.
\]
Thus the left Radau surface need not be optimized separately.
\end{cor}
\begin{proof}
The Gaussian wedge theorem together with the inactive-cap contact
classification gives an optimizer on the four surfaces consisting of
the three displayed surfaces and the left Radau surface.  For $t<c$,
the largest contraction-compatible value of $\gamma$ is
$\gamma_{\rm full}=(1-c)\left(1-\frac{t^2}{c^2}\right)$.
The parity cap satisfies
$0\le\Gamma_n(t) =\gamma_{\rm full}-\frac{2(1-u_n(t)^2)}{n^2c} \le\gamma_{\rm full}$.
Therefore Theorem~\ref{sub12:thm:turn3-left-radau-endpoints}
applies to the whole left Radau interval
$0\le\gamma\le\Gamma_n(t)$.  Its two endpoints already lie on the
Gaussian or active-cap surfaces.  At $t=c$, both
$\gamma_{\rm full}$ and $\Gamma_n(t)$ vanish; the unique effective
spectral measure is the Gaussian two-atom measure.  Equivalently this
case follows by continuity as $t\uparrow c$.  This proves the claimed
global reduction.
\end{proof}
\subsection{Elimination of inactive-cap right Radau optimizers}
\label{sub12:sec:turn3-right-radau-elimination}

Let $1/2\le c<1$, $0<q<1$, $p=q+1$, and $0<t<c$.
Put $\beta=-(1-c)t/c$ and consider the right Radau family
\[
J_\gamma=\begin{pmatrix}
 0&\sqrt c&0\\
 \sqrt c&\beta&\sqrt\gamma\\
 0&\sqrt\gamma&1-\gamma/(1-c-\beta)
 \end{pmatrix}.
\]
Its first-coordinate spectral measure has fixed first three moments
$m_1=0$, $m_2=c$, $m_3=-(1-c)t$, and contains the endpoint $1$.
Let
\[
\phi_t(y)=y\operatorname{sgn}(t+y)|t+y|^q,\quad D=(1-t)^q+(1+t)^q,\quad \Phi=\int\phi_t\,d\mu,\quad G=2\Phi/D.
\]

\begin{thm}\label{sub12:thm:turn3-right-radau-elimination}
Assume that $\mu$ is an optimizer of the moment problem with these
first three moments and an inactive fourth-moment cap, and that it is
a right Radau measure with
$0<\gamma<\gamma_{\mathrm{full}}=(1-c)(1-t^2/c^2)$.
If $t>p/3$, then its derivative
along the right Radau surface, with $c$ and $\gamma$ fixed, satisfies
\[
\partial_tG>0.
\]
Consequently no global optimizer of the finite-degree capped relaxation
can lie on the right Radau surface with an inactive cap and
$\gamma>0$. Together with Gaussian optimality for $t\le p/3$ and the
left Radau elimination, the global relaxation reduces to the Gaussian
surface and the active fourth-moment cap.
\end{thm}

The proof uses the existence of a cubic dual majorant at the
optimizer. Its derivative estimate concerns the feasible variation in
\(t\) along the right Radau surface.

\begin{proof}
Write the two movable nodes as $-x,b$, so that
$-1<-x<-t<0<t<b<1$.
The characteristic polynomial factors at the fixed node $1$. At
$\gamma=0$ the other two nodes are Gaussian, while at the full
contraction endpoint they equal $-1,t$. Both decrease as $\gamma$
increases. In particular,
$x>b$, $x>t$, $b>t$, and $(b+t)-(x-t)\ge3t-1>0$.
The first inequality also follows directly from
$b-x=\beta-\gamma/(1-c-\beta)<0$.
The moment relation is
$c(1-x+b-t)=xb-t$.
Introduce the positive quantities
\begin{gather*}
 u=x-t,\quad v=b+t,\quad A=1-t,\quad B=1+t,\\
 E=1+x,\quad L=1-b,\quad d=x+b,\\
 \varepsilon=1-x,\quad W=1-x+b-t.
\end{gather*}
Thus $u<v$, $A>u$, $B>v$, and
$c=\frac{xb-t}{W}$, $1-c=\frac{\varepsilon(1+b)}W$, and $x=1-\frac{(1-c)(b-t)}{b+c}$.
The spectral weights at $-x,b,1$ are
\[
w_- =\frac{(1+b)(b-t)}{WdE},\qquad w_b =\frac{\varepsilon(x+t)}{WdL},\qquad w_1=\frac{c-xb}{EL}.
\]
All three weights are strictly positive at the optimizer in question.

Because the cap is inactive, an attained dual polynomial $Q$ has degree
at most three and satisfies $Q\ge\phi_t$ on $[-1,1]$. It contacts
$\phi_t$ at $-x,b,1$ and has derivative contacts at $-x,b$.
Equivalently, it is the cubic Hermite interpolant at $-x,b$, and its
right-endpoint contact is the Radau stationarity equation.
Write $Q_3$ for its cubic coefficient. Necessarily
\begin{equation}
\label{sub12:eq:turn3-right-Q3-negative} Q_3<0.
\end{equation}
Indeed, on $[b,1]$ one has
$\phi_t'''(y)=q(q-1)(t+y)^{q-3} \{(q+1)(t+y)+(2-q)t\}<0$.
If $Q_3\ge0$, the residual $Q-\phi_t$ has strictly positive third
 derivative there. Its second derivative at the interior contact $b$
is nonnegative, while its value and first derivative vanish at $b$.
It would therefore be strictly positive at $1$, contradicting contact.

At the optimizer, differentiation of the moment constraints gives
\begin{equation}
\label{sub12:eq:turn3-right-t-derivative} \Phi_t=q\int y|t+y|^{q-1}\,d\mu-(1-c)Q_3, \qquad G_t=\frac2D\left(\Phi_t-\Phi\frac{D'}D\right).
\end{equation}
To justify the first identity, freeze $Q$ at the point being
differentiated. Its value and derivative contacts make the derivative
of the moving measure's expectation of $\phi_t$ equal to the
corresponding derivative for $Q$. Only the third moment changes, with
derivative \(-(1-c)\). Differentiating the integrand in \(t\) gives
the integral in \eqref{sub12:eq:turn3-right-t-derivative}.

\paragraph{A positive coefficient.}
The coefficient of $A^q$ in the numerator
$\Phi_tD-\Phi D'$ is
$\mathcal L=\Phi_t+q\Phi/A>0$.
To bound $A^q$ from below, we first compare the contributions from the two interior nodes.
Their
negative-to-positive ratio is
$\frac{u^q}{v^q}\, \frac{xL}{bE}\, \frac{b^2-t^2}{x^2-t^2} \le \frac{u^q}{v^q}\frac{bL}{xE}<1$,
where $x>b>t$ and $u<v$ were used. The endpoint contribution is
$q w_1 B^q(1/A+1/B)>0$, and the remaining term is
$-(1-c)Q_3>0$ by \eqref{sub12:eq:turn3-right-Q3-negative}.
This proves $\mathcal L>0$.

\paragraph{Reduction to a quadratic.}
Divide all powers by $v^q$ and put
$X=(u/v)^q$, $H=(B/v)^q$, and $Y=(A/v)^q$.
Define
$K=2x^2+x-b,\quad J=\frac{qxEd}{u}-K,\quad J_0=K+\frac{qd(b-x^2)}v$,
$h_1=-\frac{L^2J}{d^3}$ and $h_0=1+\frac{qL}v+\frac{L^2J_0}{d^3}$.
Direct cubic Hermite evaluation at $1$ gives the stationarity equation
$H=h_0+h_1X$.
Concavity of $r^q$ and its weighted harmonic lower bound give
\begin{equation}
\label{sub12:eq:turn3-right-power-bounds} \ell:=\frac{B}{B-qL}\le H\le h:=1+\frac{qL}v, \qquad Y\ge\frac{AX}{M},\qquad M=A-q\varepsilon>0.
\end{equation}
The bound on \(Y\) in \eqref{sub12:eq:turn3-right-power-bounds}
follows from \((A/u)^q\ge A/((1-q)A+qu)\). Applying the same
weighted harmonic mean inequality to \(B/v\) gives the lower bound on \(H\).

One has $J_0>0$. Indeed $K\ge2x^2$ since $x>b$, while
$\frac dv\le\frac{1+b}{b+t}\le2$ for $t>p/3>1/3,\ b>t$.
If $b-x^2\ge0$, then $J_0>0$ by its definition. Otherwise,
$K>2x^2$ and $d/v\le2$ give
$J_0>2x^2-2qx^2\ge0$, since $0<q<1$.
Consequently $h_0>h$. Since $X>0$ and $H\le h$, stationarity forces
$h_1<0$. The bounds on $H$ are therefore equivalent to
$0<X_h:=\frac{h-h_0}{h_1} \le X\le X_\ell:=\frac{\ell-h_0}{h_1}$.

Set
\[
P_0=bw_b,\quad P_1=xw_-,\quad T_0=\frac{qbw_b}{v} -\frac{\varepsilon(1+b)}W \frac{x-b+qdb/v}{d^3},\quad T_1=-\frac{qxw_-}{u} -\frac{\varepsilon(1+b)}W \frac{x-b-qdx/u}{d^3}.
\]
Here $P_0+P_1X$ is the normalized interior part of $\Phi$, and
$T_0+T_1X$ is the normalized $\Phi_t$ with its fixed-endpoint term
$q w_1H/B$ removed.
Define the quadratic polynomial
\begin{align}
 F(X)={}&AX(T_0+T_1X)+q(P_0+P_1X)X\notag\\
 &+(h_0+h_1X)\left\{
 M(T_0+T_1X)-\frac{qM}{B}(P_0+P_1X)
       +q w_1X\left(1+\frac AB\right)\right\}.
 \label{sub12:eq:turn3-right-quadratic}
\end{align}
Expansion of \eqref{sub12:eq:turn3-right-t-derivative}, followed by the lower
bound for $Y$ in \eqref{sub12:eq:turn3-right-power-bounds}, gives
$\frac{\Phi_tD-\Phi D'}{v^{2q}}\ge\frac{F(X)}M$.
The direction of this inequality is justified by $\mathcal L>0$.

The constant coefficient is strictly negative:
\[
F(0)=h_0M\left(T_0-\frac{qP_0}B\right),\qquad T_0-\frac{qP_0}B =\frac{\varepsilon}{Wd^3vB} \left\{(1+b)Bv(b-x)-qb d\varepsilon(d+B)\right\}<0.
\]
The two endpoint inequalities are
\begin{equation}
\label{sub12:eq:turn3-right-rational-endpoints} F(X_h)>0,\qquad F(X_\ell)>0.
\end{equation}
They hold throughout the larger rational domain
$\tfrac12\le c<1,\quad 0<q<1,\quad (1+q)/3<t<1,\quad t<b<1, \quad x=1-\frac{(1-c)(b-t)}{b+c}$,
wherever the rational functions are defined. The certificate for
\eqref{sub12:eq:turn3-right-rational-endpoints} applies to arbitrary
source weights in this domain.

A quadratic $F$ with $F(0)<0$ and positive values at two positive
arguments is positive on the interval joining $X_h$ and $X_\ell$. If its quadratic coefficient is nonpositive, concavity gives the claim. If it is
positive, its two real roots have opposite signs, so positivity at
the smaller positive argument puts the whole interval past its unique
positive root. Thus \eqref{sub12:eq:turn3-right-rational-endpoints} proves
$F(X)>0$, and hence $G_t>0$.

Finally, if the fourth cap is inactive, a sufficiently small increase
in $t$ along this Radau surface remains feasible, including at a parity
cell boundary, since the cap is continuous and the slack is strict.
The only endpoint $t=c$ is Gaussian. This excludes the asserted global
optimizer and proves the theorem.
\end{proof}

\paragraph{Exact rational certificate for the endpoint inequalities.}
The standalone script \href{https://github.com/zhangteng2000/schoenberg-validation/blob/e7435c5521dcd6a10131ba04a6da500b1404800b/certificates/turn3_source_right_radau_certificate.py}{\nolinkurl{turn3_source_right_radau_certificate.py}}
constructs every rational expression in
\eqref{sub12:eq:turn3-right-quadratic} from the Hermite coefficients and endpoint weights.
It substitutes \(c=(1+C)/2\), \(q=Q\),
\(t=(1+Q+(2-Q)T)/3\), and \(b=t+(1-t)Z\), with \(0<C,Z,T,Q<1\).
The case \(C=0\) is checked separately.
It then factors the numerator and denominator of each endpoint value
and verifies the sign of every factor of odd multiplicity by an exact
coefficient calculation. The only even-multiplicity numerator factor is
$Q^2$, strictly positive because $0<q<1$; the checker verifies this
explicitly. Even-multiplicity denominator factors are nonzero wherever
the rational expression is defined and consequently have positive sign.
The only large
odd-multiplicity factors have the following data:
\[
\begin{array}{c|c|c|c}
\text{endpoint}&\text{multidegree in }(C,Z,T,Q)
 &\text{nonzero sign-certificate coefficients}&\text{sign}\\\hline
X_h&(7,11,11,12)&14729&-\\
X_\ell&(7,11,12,13)&17192&+
\end{array}
\]
The numerator and denominator have signs $(-,-)$ and $(+,+)$,
respectively. Hence both rational endpoint values are positive.

For precision, a polynomial $P$ of multidegree $(d_1,\ldots,d_4)$ is
certified by substituting $z_j=s_j/(1+s_j)$ and multiplying by
$\prod_{j=1}^{4}(1+s_j)^{d_j}$. Every coefficient of the resulting polynomial
has the asserted sign, and at least one is nonzero. This is equivalent
to a tensor Bernstein sign certificate on the unit cube. All arithmetic
uses exact rational multivariate polynomials. The checker also verifies
that nonzero coefficients remain after setting $C=0$, so the strict
inequalities include $c=1/2$. It prints a final PASS only after checking
both rational expressions. The certificate is therefore reproducible
from the formulas alone and uses no stored numerical samples.

\begin{cor}[Two remaining surfaces for the complete relaxation]
\label{sub12:cor:turn3-two-surface-lp4}
For every $n\ge4$ and $1<p<2$, put $c=(n-2)/n$ and let
$\Gamma_n(t)$ be the parity fourth-moment cap in Jacobi coordinates.
The global reciprocal/fourth-moment relaxation has exactly the
following reduction:
\[
V_{n,p}=\max_{0\le t\le c}\max\left\{ G(t,0),\quad \max_{\delta_-(t)\le\delta\le\delta_+(t)} G\bigl(t,\Gamma_n(t),\delta\bigr)\right\},
\]
where
\[
\delta_-(t)=\frac tc-\frac{1-u_n(t)^2}{n(c-t)},\qquad \delta_+(t)=\frac tc+\frac{1-u_n(t)^2}{n(c+t)}.
\]
The formulas are understood by continuity at degenerate endpoints.
Thus only the Gaussian family and the active fourth-moment cap remain.
\end{cor}

\begin{proof}
Corollary~\ref{sub12:cor:turn3-global-lp4-boundary}
restricts global optimizers to the Gaussian, active-cap, left Radau,
and right Radau families. The left Radau endpoint theorem
reduces its contribution to Gaussian or active cap. A global maximizer
on the right Radau family with inactive cap is Gaussian when $t\le p/3$
by the Gaussian-wedge theorem, and is excluded for $t>p/3$ by
Theorem~\ref{sub12:thm:turn3-right-radau-elimination}. Strict inactivity gives
$\gamma<\Gamma_n(t)\le\gamma_{\mathrm{full}}$, as required there.
At $t=c$ the moment problem is already Gaussian. This proves the exact
reduction, without asserting an evaluation of its remaining surfaces.
\end{proof}
\subsection{An endpoint bound for the symmetric three-value family}

\begin{thm}\label{sub12:thm:finite-two-simple-outer}
Let $n\ge3$ and $1<p<2$. For every nonzero centered real root list
with three ordered values whose multiplicities are $(1,n-2,1)$,
\[
 \frac{\sum_{j=1}^{n-1}|w_j|^p}{\sum_{i=1}^{n}|x_i|^p}
 <\max\left\{
 S_{n,p},R_{n,p}(1)\right\},
\]
unless its middle value is zero. The two endpoint ratios are
\(S_{n,p}=((n-2)/n)^{\frac{p}{2}}\) for the sparse pair and
\(R_{n,p}(1)=((n-2)^p+n-2)/((n-1)^p+n-1)\) for the two-value list
obtained when two zeros coalesce. On the closed family the larger
endpoint ratio is the maximum. If the ratios are equal, both endpoints
and only those endpoints attain it, up to reflection and scaling.
\end{thm}

\begin{proof}
Write the ordered root values as $x_-<x_0<x_+$ and set
$d=(x_+-x_-)/2>0$. Centering gives
$x_-+x_+=-(n-2)x_0$. Scaling by $d$ and putting
$\tau=nx_0/(2d)$ therefore gives the root list
\begin{equation}\label{sub12:eq:finite-outer-root-list}
-1-c\tau,\qquad (2\tau/n)^{[n-2]},\qquad 1-c\tau,\qquad c=\frac{n-2}{n}.
\end{equation}
The strict ordering is equivalent to $-1<\tau<1$.
Reflection changes $\tau$ to $-\tau$ and preserves the ratio, so
we may take $0\le\tau<1$. Including $\tau=1$ gives the positive
root collision endpoint. Thus \eqref{sub12:eq:finite-outer-root-list}
covers the entire family.

Put \(a=(n-3)/n\) and \(b=1/n\).
Besides the critical point $2\tau/n$ of multiplicity $n-3$,
the remaining two critical points are
$-a\tau\pm\sqrt{c+b^2\tau^2}$.
They have opposite signs for $0\le\tau\le1$, because
$c-(a^2-b^2)\tau^2>0$ on this interval. For $\tau>0$ set
$x=1/\tau\ge1$ and $y=\sqrt{cx^2+b^2}$. Dividing numerator and
denominator of the moment ratio by $\tau^p$ gives $U(x)/V(x)$,
where $U(x)=(y+a)^p+(y-a)^p+(n-3)(2/n)^p$ and
$V(x)=(x+c)^p+(x-c)^p+(n-2)(2/n)^p$.
The four bases $y\pm a$ and $x\pm c$ are positive. With $q=p-1\in(0,1)$,
the ratio of their derivatives is
\begin{equation}\label{sub12:eq:finite-outer-derivative-quotient}
D(x):=\frac{U'(x)}{V'(x)}=\frac{cx}{y}\frac{(y+a)^q+(y-a)^q}{(x+c)^q+(x-c)^q}.
\end{equation}
We prove
\begin{equation}\label{sub12:eq:finite-outer-derivative-bound}
D(x)\le\max\{S_{n,p},R_{n,p}(1)\},\qquad x\ge1.
\end{equation}
The inequality is strict when $1<x<\infty$.

Define
\[
 \mathcal A_q(t)=\frac{(1+\sqrt t)^q+(1-\sqrt t)^q}{2},
 \qquad F_q(t)=-\log\mathcal A_q(t),\quad 0\le t<1.
\]
Since every even binomial coefficient $\binom q{2j}$, $j\ge1$,
is negative, $\mathcal A_q(t)=1-\sum_{j=1}^{\infty}d_jt^j$ with
$d_j>0$. Expanding $-\log(1-z)$ therefore shows that
\begin{equation}\label{sub12:eq:finite-outer-positive-series}
F_q(t)=\sum_{j=1}^{\infty}e_jt^j,\qquad e_j>0.
\end{equation}
The power series for $F_q$ and its derivatives converge on compact subintervals
of $[0,1)$.

Set
\[
 z=1/x^2\in[0,1],\qquad
 \eta=\frac{b^2}{c},\qquad A=\frac{a^2}{c},\qquad B=c^2.
\]
We have $0\le A<B<1$, since
\(n(n-3)^2-(n-2)^3=8-3n<0\).
Equation~\eqref{sub12:eq:finite-outer-derivative-quotient} becomes
\begin{equation}\label{sub12:eq:finite-outer-log-quotient}
\log\frac{D(x)}{c^{\frac{p}{2}}}=-\frac{1-q}{2}\log(1+\eta z)+F_q(Bz)-F_q\left(\frac{Az}{1+\eta z}\right).
\end{equation}
The right side of \eqref{sub12:eq:finite-outer-log-quotient} is strictly convex
in \(z\). The logarithmic term is strictly convex, and
\eqref{sub12:eq:finite-outer-positive-series} expresses the difference of the two \(F_q\) terms as a positive sum of
\(B^jz^j-A^j\left(\frac{z}{1+\eta z}\right)^j\).
Each such function is convex. For $j=1$ this is immediate; for
$j\ge2$, use
\[
 \frac{d^2}{dz^2}\left(\frac{z}{1+\eta z}\right)^j
 =\frac{jz^{j-2}[(j-1)-2\eta z]}{(1+\eta z)^{j+2}}
 \le j(j-1)z^{j-2}
\]
together with $A\le B$. All arguments of $F_q$ in \eqref{sub12:eq:finite-outer-log-quotient} remain in a fixed
compact subinterval of $[0,1)$, so summing the convex terms is
legitimate.

In \eqref{sub12:eq:finite-outer-derivative-quotient},
\(D(x)\to c^{\frac{p}{2}}=S_{n,p}\) as \(x\to\infty\). At \(x=1\),
we have \(y(1)=(n-1)/n\), hence
\(D(1)=\frac{n-2}{n-1}\frac{(n-2)^q+1}{(n-1)^q+1}=R_{n,p}(1)\).
Strict convexity of $\log D$ as a function of $z$ proves \eqref{sub12:eq:finite-outer-derivative-bound}.

Finally let $C_* =\max\{S_{n,p},R_{n,p}(1)\}$. At $x=1$,
the root list has two values and $U(1)=R_{n,p}(1)V(1)$.
Since $V'(x)>0$ and $D(x)\le C_*$,
\((U-C_*V)'=V'(D-C_*)\le0\).
It follows that $U(x)\le C_*V(x)$ for all $x\ge1$.
The derivative is strictly negative for $x>1$, so the inequality
is strict at every finite $x>1$, or equivalently $0<\tau<1$.
At $\tau=0$ the original, unscaled ratio is $S_{n,p}$; at
$\tau=1$ it is $R_{n,p}(1)$. The strict interior derivative bound and the endpoint calculations
prove the stated classification of the maximizers.
\end{proof}

\subsection{The joint three-moment problem in the subquadratic range}
\label{sub12:sec:turn3-uncapped-lp3-exact}

Let $n\ge4$, $1<p<2$, and set $c=(n-2)/n$, $q=p-1$, and
$S_{n,p}=c^{\frac{p}{2}}$ and $R_{n,p}(1)=\frac{(n-2)^p+n-2}{(n-1)^p+n-1}$.
With $H_p$ as in \eqref{eq:section4-H}, put
$R_{\rm cont}(t)=1-H_p((1-t)/(1+t))/n$.
Let $r_p$ be the unique minimizer of $H_p$ and
$\theta_p=r_p/(1+r_p)$, $t_p=1-2\theta_p$.

Define the joint uncapped three-moment optimum by
\[
\mathcal T_{n,p}=\max_{0\le t\le c}\ 
 \max_{\substack{\mu\text{ probability on }[-1,1]\\
 m_1=0,\ m_2=c,\ m_3=-(1-c)t}}
 \frac{2\int y\operatorname{sgn}(t+y)|t+y|^q\,d\mu(y)}
 {(1-t)^q+(1+t)^q}.
\]

\begin{thm}[Exact joint three-moment optimum]
\label{sub12:thm:turn3-uncapped-lp3-exact}
For every $n\ge4$ and $1<p<2$,
\[
 \mathcal T_{n,p} =\max\left\{S_{n,p},\ \max_{0\le t\le c}R_{\rm cont}(t)\right\}.
\]
Equivalently, this equals
\[
\begin{cases}
 \max\{S_{n,p},R_{n,p}(1)\},&n\theta_p\le1,\\[2pt]
 1-H_p(r_p)/n,&n\theta_p>1.
 \end{cases}
\]
In the first case every maximizing moment configuration is Gaussian:
$t=0$ gives the sparse value, and $t=c$ gives the one-outlier value.
Only the endpoint or endpoints whose displayed value is largest occur.
When $n\theta_p>1$, the unique maximizing moment configuration is the
intrinsic-cap measure on $\{-1,t_p,1\}$ displayed in the proof.
\end{thm}

\begin{proof}
The joint feasible set is compact in $[0,c]$ times the weak topology
of probability measures on $[-1,1]$, and the objective is continuous.
Thus a global maximum exists.

The apparently uncapped problem has an intrinsic fourth-moment cap.
Indeed, its first three moments give
$0\le\int(1-y^2)(y-t)^2\,d\mu =c-m_4-(1-c)t^2$.
Thus
$m_4\le e_{\rm full}(t):=c-(1-c)t^2$.
Equality forces support in $\{-1,t,1\}$. Its weights are uniquely
specified by the first two moments:
\[
\alpha_- =\frac{c+t}{2(1+t)},\qquad \alpha_t =\frac{1-c}{1-t^2},\qquad \alpha_+ =\frac{c-t}{2(1-t)}.
\]
They are nonnegative for $0\le t\le c$. Substitution into the objective
gives exactly $R_{\rm cont}(t)$.

For \(t\le p/3\), Theorem~\ref{sub12:thm:turn3-gaussian-wedge} gives
the unique Gaussian maximizer. Suppose \(p/3<t<c\) and the intrinsic cap
\(e=e_{\rm full}(t)\) is inactive. As in the proof of
Corollary~\ref{sub12:cor:turn3-global-lp4-boundary}, complementary slackness
gives a cubic dual majorant. Its contact count excludes three interior
nodes, so the maximizing measure is Gaussian or has an endpoint node.
Theorem~\ref{sub12:thm:turn3-left-radau-endpoints}
reduces the left Radau case to the Gaussian or full-cap endpoint.
At a nondegenerate right Radau optimizer with inactive cap,
Theorem~\ref{sub12:thm:turn3-right-radau-elimination} gives \(G_t>0\),
excluding a joint maximum.

In Jacobi coordinates the intrinsic cap is
$\gamma_{\rm full}(t) =(1-c)\left(1-\frac{t^2}{c^2}\right)$.
At this value the contraction interval for the last diagonal entry
collapses to \(\delta=t/c\), giving the measure on \(\{-1,t,1\}\)
whose weights were determined from the first two moments. Thus a joint
maximizer is Gaussian or is this three-node measure. Gaussian uniqueness
at \(t=0\) follows from Theorem~\ref{sub12:thm:turn3-gaussian-wedge}.
At \(t=c\), the moment problem is
itself Gaussian and unique: the nonnegative polynomial
$(1+y)(y-c)^2$ has zero expectation, forcing support in $\{-1,c\}$.
This also covers both endpoints without invoking interior duality.

For \(0<t<c\), set \(X=c/t>1\). In the normalization of
\eqref{sub12:eq:finite-outer-derivative-quotient}, the middle root is
\(2/n\) and the outer roots are \(-X-c\) and \(X-c\).
Increasing \(X\) moves the outer roots with velocities \(-1\) and \(1\).
The Gaussian objective is therefore the endpoint derivative quotient
\(G_{\rm Gauss}(t)=U'(X)/V'(X)\). Equation~\eqref{sub12:eq:finite-outer-derivative-bound}
and its strictness for \(X>1\) give \(G_{\rm Gauss}(t)<\max\{S_{n,p},R_{n,p}(1)\}\).
At $t=0$ the Gaussian nodes are $\pm\sqrt c$ with equal masses, giving
$G_{\rm Gauss}(0)=S_{n,p}$. At $t=c$ they are $-1,c$, giving
$G_{\rm Gauss}(c)=R_{n,p}(1)=R_{\rm cont}(c)$.
Therefore the global Gaussian contribution is exactly
$\max\{S_{n,p},R_{n,p}(1)\}$, and the claimed first formula follows.

Finally, $R_{\rm cont}$ increases strictly up to $t_p$ and decreases
strictly afterward. If $n\theta_p\le1$, then $t_p\ge c$ and its unique
maximum on the admissible interval is at $t=c$. Otherwise its maximum
is attained at $t=t_p$. In this second case,
$1/(n-1)<r_p$, so strict decrease of $H_p$ from its value $p$ at zero
gives $R_{n,p}(1)>1-p/n$. Concavity of $s^{\frac{p}{2}}$ gives
\(S_{n,p}<1-p/n\), so the sparse value is strictly smaller.
When \(n\theta_p\le1\), \eqref{sub12:eq:finite-outer-derivative-bound}
and the strict increase of \(R_{\rm cont}\) on \(0<t<c\)
leave only the Gaussian endpoints \(t=0,c\). When \(n\theta_p>1\), the continuous
maximum is strictly greater than \(R_{n,p}(1)\) and \(S_{n,p}\), and is unique at $t_p$.
The intrinsic-cap moment measure at $t_p$ is unique, which proves
the additional equality assertion.
\end{proof}

\begin{cor}[Sharp constant whenever the continuous binary optimum lies beyond the degree constraint]
\label{sub12:cor:turn3-no-rounding-cell-exact}
For every $n\ge4$ and $1<p<2$ satisfying $n\theta_p\le1$, the optimal
centered polynomial constant is
\[
C_{n,p}=\max\left\{ \left(\frac{n-2}{n}\right)^{\frac{p}{2}},\quad \frac{(n-2)^p+n-2}{(n-1)^p+n-1} \right\}.
\]
This is sharp for real roots and for complex roots. Apart from the
all-zero list, equality holds precisely, up to permutation and a
common nonzero complex scalar, for
\[
(1,-1,0,\ldots,0) \quad\text{when the first value is maximal},
\]
or
\[
(n-1,(-1)^{[n-1]}) \quad\text{when the second value is maximal}.
\]
At a tie both families give equality.
\end{cor}

\begin{proof}
At a nonzero real maximizer of the critical/root moment ratio,
normalizing the extreme roots to $t-1,t+1$ and reflecting if necessary
gives $0\le t\le c$. The centered endpoint variation makes its ratio
exactly the endpoint gradient for a feasible measure in
$\mathcal T_{n,p}$. Hence $C_{n,p}\le\mathcal T_{n,p}$.
The sparse pair and one-outlier binary list attain the two displayed
values, proving equality.

For equality, Theorem~\ref{sub12:thm:turn3-uncapped-lp3-exact}
forces $t=0$ or
$t=c$, with a Gaussian spectral measure. The inverse-source identity
for the Jacobi reduction gives normalized internal-source variance
$\gamma/(1-c)$. Gaussianity means $\gamma=0$, so all $n-2$ internal
source points equal $-t/c$. At $t=0$ they are all zero. At $t=c$ they
are all $-1$, which also follows directly from their mean and their
lower bound $-1$. This gives exactly the two real configurations.

The angular Schur--Jensen reduction transfers the sharp
bound to complex roots. Its equality step is rigid for $p>1$:
strict convexity of $|x|^p$ forces the Hermitian parts of a Schur
triangular critical-point compression to be diagonal for almost every
rotation, hence forces the compression to be normal. For a centered
root list, writing $D=\operatorname{diag}(z_j)$ and
$e=n^{-\frac{1}{2}}\mathbf1$, the compression commutator is the difference
of the rank-one matrices generated by $De$ and $\overline D e$.
Normality gives \(De=\omega\overline D e\) for some \(|\omega|=1\).
Hence all roots lie on a common line through zero, and rotation gives
the real equality cases.
\end{proof}

\begin{cor}[An explicit exponent range in every degree]
\label{sub12:cor:turn3-square-root-exponent-regime}
For every integer $n\ge4$ and every $1<p<2$ satisfying
\[
(n-2)(p-1)^2\le4, \qquad\text{equivalently}\qquad p-1\le\frac2{\sqrt{n-2}},
\]
one has
\[
C_{n,p}=\max\left\{ \left(\frac{n-2}{n}\right)^{\frac{p}{2}},\quad \frac{(n-2)^p+n-2}{(n-1)^p+n-1} \right\},
\]
with precisely the equality configurations in
Corollary~\ref{sub12:cor:turn3-no-rounding-cell-exact}.
\end{cor}

\begin{proof}
The bound \eqref{sub12:eq:subquad-explicit-b-location},
where $q=p-1$ and $\theta_p=r_p/(1+r_p)$, gives
$\theta_p<\frac{q^2}{4+2q^2}$.
Under $(n-2)q^2\le4$ this implies $n\theta_p<1$.
Apply Corollary~\ref{sub12:cor:turn3-no-rounding-cell-exact}.
In particular every fixed ray $p=1+\tau/n$ enters this explicit
regime for all sufficiently large $n$.
\end{proof}

\begin{cor}[All exponents in degrees four through nineteen]
\label{sub12:cor:turn3-degrees-four-nineteen}
For every integer $4\le n\le19$ and every $1<p<2$,
\[
C_{n,p}=\max\left\{ \left(\frac{n-2}{n}\right)^{\frac{p}{2}},\quad \frac{(n-2)^p+n-2}{(n-1)^p+n-1} \right\},
\]
with exactly the equality configurations in
Corollary~\ref{sub12:cor:turn3-no-rounding-cell-exact}.
\end{cor}

\begin{proof}
The sharpened fixed-point location certificate,
Theorem~\ref{sub12:thm:turn3-nineteen-location}, proves
$r_p<1/18$ for every $1<p<2$. Therefore
$\theta_p=\frac{r_p}{1+r_p}<\frac1{19}$.
Consequently $n\theta_p<1$ throughout the stated degree range, and
Corollary~\ref{sub12:cor:turn3-no-rounding-cell-exact}
applies for every \(1<p<2\).
\end{proof}
\subsection{Strict positivity of the quartic anchor throughout every peak cell}
\label{sub12:sec:turn3-peak-anchor-positive}

\begin{thm}\label{sub12:thm:turn3-peak-anchor-positive}
Let $0<q<1$, $p=1+q$, and let $r_p$ minimize $H_p$.
The left endpoint coefficient of the quartic anchor is strictly
positive whenever
\[
0<r=\frac{1-t}{1+t}\le\frac{2r_p}{1-r_p}.
\]
Consequently $E_-(t)>0$ throughout every integer cell containing
the continuous two-value peak.  The right coefficient is also
strictly positive there, so the quartic anchor is available on
the entire peak cell, including its endpoints.
\end{thm}

\begin{proof}
Put $\alpha=1-q$, $X=r^q$, and $Z=(1-r)^q$.
The exact endpoint formula gives
$\frac{4r^2E_-}{(1+t)^{q-2}}=N_q(r)$,
where
$N_q(r)=1-Z\{1+qr+(q+1)r^2+2r^3\} +r^{q+1}\{q-(1-q)r-2r^2\}$.
We first establish a useful elementary lower bound.  Write
$(1-r)^q=1-\sum_{j=1}^{\infty}a_jr^j$, where $a_j>0$ and $a_1=q$.
Multiplication gives
$Z\{1+qr+(q+1)r^2+2r^3\} =1+c_2r^2+c_3r^3-\sum_{j=4}^{\infty}C_jr^j$,
with
$c_2=\frac{2+q-q^2}{2}$, $c_3=\frac{(1-q)(6+2q-q^2)}3$,
$C_j=a_j+qa_{j-1}+(q+1)a_{j-2}+2a_{j-3}>0$ and $C_4>2q$.
Keeping only the strict bound $C_4>2q$ yields
\begin{equation}
\label{sub12:eq:turn3-anchor-F-lower} N_q(r)>\alpha r^2F_q(r),
\end{equation}
where
$F_q(r)= (q-2r^2)\frac{r^{-\alpha}-1}{\alpha} -r^q-\frac{q+2}{2} -\frac{6+2q-q^2}{3}r-2r^2$.
Whenever $q>2r^2$, differentiation shows that $F_q(r)$ is strictly
decreasing in $r$: every term in
$F_q'(r)=-4r\frac{r^{-\alpha}-1}{\alpha} -(q-2r^2)r^{-\alpha-1}-qr^{q-1} -\frac{6+2q-q^2}{3}-4r$
is negative.  We will prove positivity on three enlarged parameter
domains.

\paragraph{Small exponents: $0<q\le2/5$ and $r\le q^2/2$.}
It suffices to set $R=q^2/2$.  Define
$g(q)=\frac{qR^{q-1}}{1-q}$.
This function is strictly decreasing on $(0,2/5]$.  Indeed
$(\log g)'=\log(q^2/2)+2-\frac1q+\frac1{1-q}$,
and the right side is increasing in $q$ and, at $q=2/5$, is
$\log(2/25)+7/6<-5/6$.
Moreover $g(2/5)>3$, since
$(25/2)^3>(9/2)^5$.
Thus $g(q)>3$ throughout this range.  In the expression for $F_q(R)$,
the first term can be written as
$\left(1-\frac{q^3}{2}\right)g(q) -\frac{q-2R^2}{1-q}$.
For $q\le3/10$, using $R^q\le1$ gives the explicit lower bound
\[
F_q(R)> \frac{1973}{2000}\,3-\frac37-1-\frac{23}{20} -\frac73\frac9{200}-2\left(\frac9{200}\right)^2 =\frac{38063}{140000}>0.
\]
For $3/10\le q\le2/5$, the function $(q^2/2)^q$ is decreasing,
because its logarithmic derivative is $\log(q^2/2)+2<0$.
Also
$(9/200)^{\frac{3}{10}}<\frac25$,
as is seen by raising both sides to the tenth power.  Consequently
\[
F_q(R)> \frac{121}{125}\,3-\frac23-\frac25-\frac65 -\frac73\frac2{25}-2\left(\frac2{25}\right)^2 =\frac{821}{1875}>0.
\]

\paragraph{Two constant-radius domains.}
We next prove
\[
F_q(r)>0\quad
 \begin{cases}
  2/5\le q<1,&0<r\le1/12,\\
  3/5\le q<1,&0<r\le2/17.
 \end{cases}
\]
By monotonicity in $r$, it suffices to use $R=1/12$ and $R=2/17$,
respectively.  Put $L=\log(1/R)$ and write $F_R(\alpha)=F_{1-\alpha}(R)$.
The continuous value at zero is
$F_R(0)=(1-2R^2)L-\frac32-\frac{10R}{3}-2R^2$.
For every Taylor degree $j\ge2$, the coefficient of $\alpha^j$ is
$\frac{L^j}{j!} \left\{\frac{(1-2R^2)L}{j+1}-(1+R)\right\} +\frac R3\,\mathbf 1_{\{j=2\}}$.
The Taylor coefficients of degree \(j\ge2\) are strictly negative.  To check this using
only loose bounds, note that $2<L<5/2$ for both radii.  The braces
are smaller than $-1/6-R$ when $j\ge2$, and for $j=2$ the additional
$R/3$ is dominated since $L^2/2>2$.  Thus $F_R$ is strictly concave
for $\alpha\ge0$.

At zero the elementary bound $L>2$ gives
$F_{1/12}(0)>\frac{13}{72}$ and $F_{2/17}(0)>\frac{43}{1734}$.
At the other endpoints, use
$12^{\frac{3}{5}}>\frac{22}{5}$ and $(17/2)^{\frac{2}{5}}>\frac{47}{20}$.
Raising both sides of each inequality to the fifth power gives a rational comparison.
The coefficients multiplying the exponential terms in $F_R$ are positive,
and direct substitution gives
$F_{1/12}(3/5)>\frac{571}{1350}$ and $F_{2/17}(2/5)>\frac{3421}{57800}$.
Concavity proves positivity on $0\le\alpha\le3/5$ and
$0\le\alpha\le2/5$, respectively.  This establishes both domains.

\paragraph{Location of the peak cell.}
For $q\le2/5$, equation~\eqref{sub12:eq:subquad-explicit-b-location} gives
$\frac{2r_p}{1-r_p}<\frac{q^2}{2}$.
For $q\ge3/5$, Theorem~\ref{sub12:thm:turn3-nineteen-location} gives
$2r_p/(1-r_p)<2/17$.
In the intermediate interval $2/5\le q\le3/5$, the exact location
bound in Lemma~\ref{sub12:lem:turn3-middle-location} is $r_p<1/25$, so
$2r_p/(1-r_p)<1/12$.
The three positivity domains therefore cover every point required
by the theorem.  Equation~\eqref{sub12:eq:turn3-anchor-F-lower} proves
$E_->0$.  Lemma~\ref{sub12:lem:anchor-right-positive} gives $E_+>0$
because every peak cell
has $t>15/19>1/2$.
\end{proof}

\begin{lem}\label{sub12:lem:turn3-middle-location}
For $2/5\le q=p-1\le3/5$, the unique minimizer satisfies
$r_p<1/25$.
\end{lem}

\begin{proof}
Set $r=1/25$, $\varepsilon=1-q\in[2/5,3/5]$,
$\alpha=\log(25/24)$ and $\beta=\log25$.
At $r=1/25$, the derivative numerator $J=(r+r^p)^2H_p'(r)$
is obtained from the exponential formula in the proof of
Theorem~\ref{sub12:thm:turn3-nineteen-location} by substituting
$r=1/25$, $\alpha=\log(25/24)$, and $\beta=\log25$.
Write $J=\sum_{j=1}^{\infty}c_j\varepsilon^j$.
Exact rational logarithmic enclosures yield
$c_1>-\frac3{5000},\quad c_2>0,\quad c_3>\frac{17}{5000},\quad c_4>\frac3{1000},\quad c_5,c_6>0$,
\[
c_j>-\frac3{5000}\quad(7\le j\le11),\qquad \left|\sum_{j=12}^{\infty}c_j\varepsilon^j\right| <\frac{\varepsilon}{100000}.
\]
The checker \href{https://github.com/zhangteng2000/schoenberg-validation/blob/e7435c5521dcd6a10131ba04a6da500b1404800b/certificates/turn3_norm_peak_anchor_certificate.py}{\nolinkurl{turn3_norm_peak_anchor_certificate.py}}
verifies the coefficient and tail bounds by rational arithmetic.  It encloses
$\alpha$ by the atanh series with $z=1/49$, and uses
$\beta=4\log2+2\log(5/4)$ with $z=1/3,1/9$.
For the exponential tails it uses
$\alpha<1/24$, $\alpha+\beta<10/3$, $\beta<13/4$, and $2\beta<13/2$.
It follows that
\[
\frac J\varepsilon> -\frac3{5000}+\frac{17}{5000}\left(\frac25\right)^2 +\frac3{1000}\left(\frac25\right)^3 -\frac3{5000}\frac{(3/5)^6}{1-3/5}-\frac1{100000} =\frac{3501}{62500000}>0.
\]
Strict unimodality of $H_p$ proves the claimed location.
\end{proof}
\subsection{Reduction of the finite-degree theorem to one peak cell}
\label{sub12:sec:turn3-peak-cell-sufficiency}

Theorem~\ref{sub12:thm:turn3-left-radau-endpoints} and
Theorem~\ref{sub12:thm:turn3-right-radau-elimination}
exclude inactive-cap Radau extrema away from the relevant cell.
We now reduce the sharp finite-degree inequality to a scalar comparison
on the integer interval containing the continuous optimum. The required
quartic coefficients are positive by
Theorem~\ref{sub12:thm:turn3-peak-anchor-positive}.

Fix $n\ge4$, $1<p<2$, $q=p-1$, and $c=(n-2)/n$. Retain $R_{\rm cont}(t)=1-H_p((1-t)/(1+t))/n$,
$\mathcal B_{n,p}=\max_jR_{n,p}(j)$.
The continuous two-value ratio attains its unique maximum at \(t_p=1-2\theta_p\).
If \(n\theta_p\le1\), apply Corollary~\ref{sub12:cor:turn3-no-rounding-cell-exact}.
If \(n\theta_p>1\) is an integer, the maximum in
Theorem~\ref{sub12:thm:turn3-uncapped-lp3-exact} is attained by the
two-value list of multiplicities \(n\theta_p\) and \(n(1-\theta_p)\).

We therefore assume that $k=\lfloor n\theta_p\rfloor\ge1$ and
$n\theta_p$ is not an integer. Put
$\ell=1-\frac{2(k+1)}n$ and $u=1-\frac{2k}n$, so that $\ell<t_p<u$.
Theorem~\ref{sub12:thm:turn3-nineteen-location} gives
$\theta_p<1/19$, which yields $\ell>3/4$.
In this case $\mathcal B_{n,p}>S_{n,p}$: indeed $1/(n-1)<r_p$,
so strict decrease of $H_p$ from its endpoint value $H_p(0)=p$
gives $H_p(1/(n-1))<p$. Hence
$\mathcal B_{n,p}\ge R_{n,p}(1)>1-p/n\ge c^{\frac{p}{2}}$.

\begin{thm}[Sufficient scalar comparison on one cell]
\label{sub12:thm:turn3-peak-cell-sufficient}
The endpoint anchor coefficients satisfy
\[
E_-(t)>0,\qquad E_+(t)>0\qquad(\ell\le t\le u).
\]
Set $E_*(t)=\min\{E_-(t),E_+(t)\}$. If
\begin{equation}\label{sub12:eq:turn3-complete-peak-comparison}
 H_p\!\left(\frac{1-t}{1+t}\right)+\frac{4nE_*(t)}{D_t}(t-\ell)(u-t)
 \ge\min\{H_p\!\left(\frac{1-\ell}{1+\ell}\right),H_p\!\left(\frac{1-u}{1+u}\right)\}
 \qquad(\ell\le t\le u),
\end{equation}
then $C_{n,p}=\mathcal B_{n,p}$ for centered complex polynomials.
No additional moment-envelope estimate outside this cell is needed.
\end{thm}

\begin{proof}
Strict positivity of both coefficients is
Theorem~\ref{sub12:thm:turn3-peak-anchor-positive}.
Consider first the uncapped three-moment problem on any closed
subinterval of $0\le t\le c$. Its feasible measures have
$m_1=0,m_2=c,m_3=-(1-c)t$.
The intrinsic largest fourth moment is
$e_{\rm full}=c-(1-c)t^2$, so its Jacobi interval is
$0\le\gamma\le\gamma_{\rm full} =(1-c)(1-t^2/c^2)$.
At the upper endpoint the two contraction bounds force $\delta=t/c$;
the measure is supported on $\{-1,t,1\}$ and its value is
$R_{\rm cont}(t)$. The Gaussian endpoint has value $G_{\rm G}(c,t)$.

Compactness gives a maximizer. Gaussian uniqueness for $t\le p/3$
and the attained cubic contact classification for $t>p/3$ reduce
the other possibilities to the Radau surfaces.
Theorem~\ref{sub12:thm:turn3-left-radau-endpoints} reduces the left
surface to its Gaussian or full-cap endpoints. At every
remaining nondegenerate right Radau optimizer, Theorem~\ref{sub12:thm:turn3-right-radau-elimination}
shows that the derivative with
respect to increasing $t$ at fixed $\gamma$ is strictly positive.
The cap is inactive there, so this is a feasible local direction
whenever the right endpoint of the chosen $t$ interval has not been
reached.

Apply this argument on $[u,c]$. A nondegenerate right Radau measure
cannot maximize at an interior $t$, nor at the left endpoint $u$,
because increasing $t$ remains within this interval. The endpoint
$c$ is Gaussian. Thus its maximum is bounded by the Gaussian maximum
and the maximum of $R_{\rm cont}$ on $[u,c]$. The Gaussian maximum is
$\max\{S_{n,p},R_{n,p}(1)\}\le\mathcal B_{n,p}$ by the
three-moment classification, while the continuous binary value
attains its maximum on $[u,c]$ at
$R_{\rm cont}(u)\le\mathcal B_{n,p}$ by
Theorem~\ref{sub12:thm:subquad-H-unimodal}. Consequently every
uncapped three-moment measure with $u\le t\le c$ already satisfies
\(G\le\mathcal B_{n,p}\).

Apply the Gaussian/Radau boundary classification from
Corollary~\ref{sub12:cor:turn3-global-lp4-boundary} to
$[0,\ell]$. The only possible
additional right Radau maximum is at the right endpoint $t=\ell$.
The endpoint inequalities \(E_-(\ell)>0\) and \(E_+(\ell)>0\),
together with Theorem~\ref{sub12:thm:quartic-anchor}, make \(Q_0\) a global majorant.
Its moment expectation gives the pointwise bound
$G\le R_{\rm cont}(\ell)$ for every measure at that endpoint.
The Gaussian maximum is again at most $\mathcal B_{n,p}$.
The strict increase of $R_{\rm cont}$ on $[0,\ell]$ follows
from Theorem~\ref{sub12:thm:subquad-H-unimodal}. Consequently every
uncapped three-moment measure with $0\le t\le\ell$ has
objective at most $\mathcal B_{n,p}$.

Inside $[\ell,u]$, use the quartic anchor with coefficient $E_*(t)$
and the exact integer fourth-moment bound. The parity coordinate on
this cell satisfies
$1-u_n(t)^2=n^2(t-\ell)(u-t)$, so the resulting gradient bound is
$G\le1-\frac{H_p\!\left(\frac{1-t}{1+t}\right)}n -\frac{4E_*(t)}{D_t}(t-\ell)(u-t)$.
By \eqref{sub12:eq:turn3-complete-peak-comparison}, this is at most
\(1-\frac1n\min\{H_p\!\left(\frac{1-\ell}{1+\ell}\right),H_p\!\left(\frac{1-u}{1+u}\right)\}=\mathcal B_{n,p}\).
The final equality uses the adjacent-integer reduction for the binary
maximum. The three regions cover all admissible $t$.

At a real extremizing root list, the endpoint-gradient value equals
the critical/root moment ratio. Thus $C_{n,p}\le \mathcal B_{n,p}$.
The corresponding binary polynomial attains $\mathcal B_{n,p}$, and angular
Schur--Jensen reduction gives the same sharp constant for complex roots.
\end{proof}

Theorem~\ref{sub12:thm:turn3-peak-cell-sufficient} reduces the
finite-degree sharp constant to \eqref{sub12:eq:turn3-complete-peak-comparison}.
Theorem~\ref{sub12:thm:turn3-integer-anchor-branch} and
Theorem~\ref{sub12:thm:turn3-backward-remainder-proved}
establish the scalar comparison for $1<p<2$.
\subsection{Three derivative signs on the full peak domain}
\label{sub12:sec:turn3-peak-derivative-signs}

Write $q=p-1$, $\theta_p=r_p/(1+r_p)$, and define
\[
J(\theta)=H_p\!\left(\frac{\theta}{1-\theta}\right),\qquad T(\theta)=\frac{16E_+(1-2\theta)} {(2\theta)^q+(2-2\theta)^q}.
\]

\begin{thm}\label{sub12:thm:turn3-peak-derivative-signs}
For every $0<q<1$ and $0<\theta\le2\theta_p$,
\[
J''(\theta)>0,\qquad J'''(\theta)<0,\qquad T'(\theta)>0.
\]
\end{thm}

The proof below is a finite rational polynomial certificate.
It contains no numerical sampling and no approximation of a fractional
power by an uncertified value.

\begin{proof}
Set $r=\theta/(1-\theta)$, $X=r^q$, and $Z=(1-r)^q$.
The location estimates in Theorem~\ref{sub12:thm:turn3-nineteen-location}
and Lemma~\ref{sub12:lem:turn3-middle-location}, together with
$r_p/(1-r_p)<q^2/4$, place the required domain
in the union of the following three regions:
\[
\begin{array}{c|c}
 0<q\le2/5&0<r\le q^2/2\\
 2/5\le q\le3/5&0<r\le1/12\\
 3/5\le q<1&0<r\le2/17.
 \end{array}
\]
Indeed, \(r\le2r_p/(1-r_p)\), and the proof of
Theorem~\ref{sub12:thm:turn3-peak-anchor-positive} gives each upper
bound for \(r\). In all three ranges \(r\le2/17\), so \(\log(1/r)>2\).

Temporarily regard $r,q,X,Z$ as independent indeterminates and put
$\mathcal H= \frac{(1+r)\{1+rX-(1-r)Z\}}{r(1+X)}$,
$\mathcal T=\frac{(1+r)^2}{r^3(1+X)} \{Xr^3+Zqr^2+Zqr+Zr^3+Zr+2Z+qr^2+qr-r-2\}$.
At \(X=r^q\), \(Z=(1-r)^q\), one has \(\mathcal H=J\) and \(\mathcal T=T\).
Differentiation with respect to \(\theta\) is given by
\begin{equation}\label{eq:integer-cell-operator}
 \mathcal D=(1+r)^2\left(\partial_r+\frac{qX}{r}\partial_X-\frac{qZ}{1-r}\partial_Z\right).
\end{equation}
Define the three polynomial numerators by
\[
N_2=r^3(1+X)^3(r-1)\mathcal D^2\mathcal H,\qquad N_3=r^4(1+X)^4(r-1)^2\mathcal D^3\mathcal H,\qquad N_T=r^4(1+X)^2(r-1)\mathcal D\mathcal T.
\]
Each numerator is a polynomial and is affine in $Z$.
All three desired signs are equivalent to $N_i<0$.

We give exact polynomial bounds for the signs of $N_2,N_3,N_T$.
For $j\ge1$, let $a_j=q(1-q)(2-q)\cdots(j-1-q)/j!$,
where $a_1=q$, and set $Z_6=1-\sum_{j=1}^{6}a_jr^j$.  The binomial series gives $Z<Z_6$ strictly.
Put $\alpha=1-q$.  The elementary exponential lower bound and the
convex chord bound in the exponent give
$r(1+2\alpha+2\alpha^2)<X<r+\alpha(1-r)$.
The upper bound for \(X\) minus the lower bound is
$\alpha(1-3r-2\alpha r)>0$.  We may therefore write
\begin{equation}\label{eq:power-box-substitution}
 X=r(1+2\alpha+2\alpha^2)+\alpha(1-3r-2\alpha r)w,\qquad 0<w<1.
\end{equation}

The first certificate verifies
$K_i:=\frac{\partial_ZN_i}{(1+r)^{\nu_i}}>0$, where $(\nu_2,\nu_3,\nu_T)=(3,4,4)$,
on the larger box $0\le q\le1$, $0\le r\le2/17$,
$0\le X\le1$.  After $q=u$, $r=2v/17$, $X=w$, each $K_i$
has strictly positive coefficients in its tensor-product Bernstein
basis.  Their degree vectors and total coefficient counts are
\[
\begin{array}{c|c|c}
 K_2&(2,4,2)&45\\
 K_3&(3,6,3)&112\\
 K_T&(2,4,1)&30.
 \end{array}
\]
Thus replacing $Z$ by $Z_6$ strictly increases every numerator.

For the second certificate substitute \eqref{eq:power-box-substitution}
for \(X\) and define
$P_i(q,r,w)=-\frac{N_i(r,q,X,Z_6)}{r^{e_i}(1-q)}$, where $(e_2,e_3,e_T)=(1,1,3)$.
The indicated divisions are exact: all $P_i$ are polynomials with
rational coefficients.  On the middle and high exponent regions,
use respectively
\[
(q,r)=\left(\frac25+\frac u5,\frac v{12}\right),\qquad (q,r)=\left(\frac35+\frac{2u}5,\frac{2v}{17}\right).
\]
Every coefficient in the Bernstein basis of the resulting polynomial
on $0\le u,v,w\le1$ is nonnegative.  In both regions the degree
vectors are
$\deg P_2=(11,14,3),\quad \deg P_3=(14,18,4),\quad \deg P_T=(9,12,2)$.
For the low exponent region substitute $r=q^2v/2$.  Each resulting
polynomial is divisible by $q^2$.  After dividing by this positive
factor and substituting $q=2u/5$, all Bernstein coefficients are
again nonnegative.  Their degree vectors are
$(35,14,3)$, $(45,18,4)$, and $(29,12,2)$,
respectively.

For completeness the exact certificate is reproducible from these
compact formulas alone.  If a polynomial of multidegree
$(d_1,d_2,d_3)$ has power coefficients $a_{ijk}$, its Bernstein
coefficient at $(h_1,h_2,h_3)$ is
\(\sum_{i\le h_1,\,j\le h_2,\,k\le h_3} a_{ijk} \frac{\binom{h_1}{i}}{\binom{d_1}{i}} \frac{\binom{h_2}{j}}{\binom{d_2}{j}} \frac{\binom{h_3}{k}}{\binom{d_3}{k}}\).
The accompanying checker
\href{https://github.com/zhangteng2000/schoenberg-validation/blob/e7435c5521dcd6a10131ba04a6da500b1404800b/certificates/turn3_norm_peak_derivative_certificate.py}{\nolinkurl{turn3_norm_peak_derivative_certificate.py}}
constructs the numerators by \eqref{eq:integer-cell-operator},
performs every substitution and exact division, and verifies every
coefficient using rational arithmetic.  It also checks that the
first certificate has strictly positive coefficients.

The second certificate proves $N_i(r,q,X,Z_6)\le0$.
The first certificate and $Z<Z_6$ consequently give
$N_i(r,q,X,Z)<0$ for the actual powers.  The signs of the three
positive or negative denominators in the definitions of $N_i$ now
yield $J''>0$, $J'''<0$, and $T'>0$, as asserted.
\end{proof}
\subsection{A certified region in which the right anchor is active}
\label{sub12:sec:turn3-anchor-switch-highq}

\begin{lem}\label{sub12:lem:turn3-anchor-switch-highq}
If $3/4\le q<1$ and $0<r\le1/14$, then $E_-(t)>E_+(t)$,
where $t=(1-r)/(1+r)$. In particular, in every first integer cell
of degree $n\ge30$, the right endpoint determines the anchor
coefficient throughout this exponent interval.
\end{lem}
\begin{proof}
Put $\alpha=1-q$, $X=r^q$, and $Z=(1-r)^q$. The sign of
$E_--E_+$ is the sign of
$N=2-qr+r^2(q-2r)X-(2+qr+qr^2+2r^3)Z$.
Write $Z=1-\sum_{j=1}^{\infty}a_jr^j$, where $a_j>0$ and $a_1=q$.
Exact multiplication gives
$N=r^2(q-2r)X-Br^3+C_4r^4+\sum_{j=5}^{\infty}C_jr^j$,
where
$B=\frac{(3-q)(2-q)(q+2)}6$, $C_4=2q+\frac{\alpha q(6+5q-q^2)}{12}$, and $C_j>0$.
Consequently $N>\alpha r^3G_q(r)$, with
\[
G_q(r)=(q-2r)\frac{r^{-\alpha}-1}{\alpha} -\frac{12+2q-q^2}{6} +r\left\{\frac{q(6+5q-q^2)}{12}-2\right\}.
\]
For $q>2r$, this function is strictly decreasing in $r$:
its derivative is the sum of $-2(r^{-\alpha}-1)/\alpha$,
$-(q-2r)r^{-\alpha-1}$, and the strictly negative last coefficient.
It therefore suffices to set $R=1/14$.

Put $L=\log14$. In terms of $\alpha\in[0,1/4]$, the function is
\[
G_R(\alpha)=\left(\frac67-\alpha\right) \frac{e^{\alpha L}-1}{\alpha} -\frac{13}{6}+\frac{\alpha^2}{6} +\frac{-14-13\alpha+2\alpha^2+\alpha^3}{168}.
\]
This function is strictly concave. Indeed its Taylor coefficient
of degree $j\ge2$, apart from the explicit contributions
$\frac5{28}\mathbf1_{j=2}+\frac1{168}\mathbf1_{j=3}$, equals
$\frac{L^j}{j!}\left\{\frac{6L}{7(j+1)}-1\right\}$.
The elementary bounds $5/2<L<3$ make all these coefficients
strictly negative, including degrees two and three.
At zero,
$G_R(0)=\frac67\log14-\frac94>0$.
For completeness, the atanh series gives
$\log14=3\log2+\log(7/4)>56/27+6/11+18/1331>21/8$.
At the other endpoint, $14^{\frac{1}{4}}>193/100$ gives
$G_R(1/4)>\frac{43}{89600}>0$.
Both comparisons are rational after taking fourth powers in the
latter one. Strict concavity proves positivity throughout the
interval. Thus $N>0$ and $E_->E_+$. The first-cell radius is
$r\le2/(n-2)\le1/14$ when $n\ge30$.
\end{proof}
\subsection{The left anchor is unnecessary below the integer binary envelope}
\label{sub12:sec:turn3-integer-anchor-branch}

The integer hypothesis enters explicitly in the treatment of the first
cell for exponents close to two.

\begin{thm}\label{sub12:thm:turn3-integer-anchor-branch}
Let $n\ge4$ be an integer, $0<q=p-1<1$, and suppose
$n\theta_p\notin\mathbb Z$ and $k=\lfloor n\theta_p\rfloor\ge1$.
Put
\[
r_0=\frac{k}{n-k},\qquad r_1=\frac{k+1}{n-k-1},\qquad r=\frac{1-t}{1+t}\in[r_0,r_1].
\]
If $E_-(t)\le E_+(t)$, then
\[
H_p(r)>H_p(r_0).
\]
Consequently, on the part of the peak cell where
$H_p(r)\le\min\{H_p(r_0),H_p(r_1)\}$, the right anchor is
strictly smaller: $0<E_+(t)<E_-(t)$.
\end{thm}

\begin{proof}
Throughout the proof write $H=H_p$, $\alpha=1-q$,
$X=r^q$, and $Z=(1-r)^q$. The location estimates \eqref{sub12:eq:subquad-explicit-b-location},
Theorem~\ref{sub12:thm:turn3-nineteen-location}, and
Lemma~\ref{sub12:lem:turn3-middle-location}
and the definition of the peak cell give
\begin{equation}
\label{sub12:eq:turn3-branch-peak-domain} r_0<r_p<r_1,\qquad r\le\frac{2r_p}{1-r_p}<\min\left\{\frac{q^2}{2},\frac2{17}\right\}.
\end{equation}
In particular $q>2r$. By the exact endpoint formula, the sign of
$E_--E_+$ is the sign of
$N_q(r)=2-qr+r^2(q-2r)r^q -(2+qr+qr^2+2r^3)(1-r)^q$.

\paragraph{A strictly decreasing normalized switch function.}
Write $(1-r)^q=1-\sum_{j=1}^{\infty}a_jr^j$, where $a_j>0$ and
$a_1=q$. Multiplication gives
$\frac{N_q(r)}{r^3} =(q-2r)r^{-\alpha}-B+C_4r+\sum_{j=5}^{\infty}C_jr^{j-3}$,
where
$B=\frac{(3-q)(2-q)(q+2)}6$, $C_4=2q+\frac{\alpha q(6+5q-q^2)}{12}$,
and $C_j=2a_j+qa_{j-1}+qa_{j-2}+2a_{j-3}>0\quad(j\ge5)$.
For $j\ge2$, the product formula gives $a_j\le\alpha/j$.
Consequently $C_j<2\alpha$ for $j\ge5$, and
$C_4\le2q+5\alpha/6$. Differentiating the convergent series, and
using $r^{-\alpha}\ge1$, gives
\[
\left(\frac{N_q(r)}{r^3}\right)' \le\alpha\left\{-\frac qr+\frac56 +\frac{2r(2-r)}{(1-r)^2}\right\} <-\frac{269}{450}\alpha<0 \qquad\left(q>2r,\ 0<r\le\frac2{17}\right).
\]
Here $2r(2-r)/(1-r)^2\le128/225$ for \(0<r\le2/17\).
Thus, whenever a fixed radius $R$ lies in this domain and
$N_q(R)>0$, the condition $N_q(r)\le0$ forces $r>R$.

\paragraph{Small exponents: no branch switch on the peak cell.}
Keeping the $C_4$ term gives
$N_q(r)>\alpha r^3G_q(r)$,
where
\[
G_q(r)=(q-2r)\frac{r^{-\alpha}-1}{\alpha} -\frac{12+2q-q^2}{6} +r\left\{\frac{q(6+5q-q^2)}{12}-2\right\}.
\]
When $q>2r$, this function is strictly decreasing in $r$.
For $q\le3/10$, put $R=q^2/2$. Then
$G_q(R)>qR^{-\alpha}-q-\frac{13}{6}-q^2$.
The function $q(q^2/2)^{q-1}$ is strictly decreasing on
$(0,3/10]$, because its logarithmic derivative is
$\log(q^2/2)+2-1/q<0$. Moreover
$\frac3{10}(200/9)^{\frac{7}{10}}>\frac{13}{5}$,
which is the rational inequality $(200/9)^7>(26/3)^{10}$.
It follows that $G_q(R)>13/300>0$. By
\eqref{sub12:eq:turn3-branch-peak-domain}, $N_q(r)>0$ throughout the peak
cell in this range, so the premise of the theorem cannot hold.

\paragraph{Five middle exponent bands.}
For each row of the following table, the indicated two strict
inequalities hold throughout its closed $q$ interval:
\begin{equation}\label{eq:five-fixed-radii}
\begin{array}{c|c|c}
 q\text{ interval}&R&\text{certified inequalities}\\ \hline
 {[3/10,2/5]}&7/200&N_q(R)>0,\quad H(R)>H(R/(2+R))\\
 {[2/5,1/2]}&9/200&N_q(R)>0,\quad H(R)>H(R/(2+R))\\
 {[1/2,3/5]}&11/200&N_q(R)>0,\quad H(R)>H(R/(2+R))\\
 {[3/5,7/10]}&3/50&N_q(R)>0,\quad H(R)>H(R/(2+R))\\
 {[7/10,4/5]}&67/1000&N_q(R)>0,\quad H(R)>H(R/(2+R)).
\end{array}
\end{equation}
For the rational verification of \eqref{eq:five-fixed-radii}, see the
certificate in this proof. Assume \(N_q(r)\le0\). Normalized monotonicity gives
$r>R$. The integer-cell geometry gives
$r_0\ge\frac{r}{2+r}$.
Indeed, in the variable $\theta=r/(1+r)$, the lower lattice point
is $k/n$ and $\theta\le(k+1)/n$, so $(k/n)/\theta\ge1/2$.
Set $s(r)=r/(2+r)$. Then
$s(R)<s(r)\le r_0<r_p$. Since $H(R)>H(s(R))$, strict unimodality
forces $R>r_p$. The monotonicity in
Theorem~\ref{sub12:thm:subquad-H-unimodal} gives
$H(r)>H(R)>H(s(R))\ge H(s(r))\ge H(r_0)$.

\paragraph{The interval $4/5\le q<1$.}
Set $R=5/69$. Three further certified inequalities hold throughout
this interval:
\begin{equation}
\label{sub12:eq:turn3-highq-fixed-radius} N_q(R)>0,\qquad H(R)>H(1/28),\qquad H(R)>H\!\left(\frac{2R}{3+R}\right).
\end{equation}
Again $N_q(r)\le0$ forces $r>R$.
If $k\ge2$, the lattice geometry improves to
$r_0\ge\frac{2r}{3+r}$,
because $k/(k+1)\ge2/3$. Use the last inequality in
\eqref{sub12:eq:turn3-highq-fixed-radius} and apply Theorem~\ref{sub12:thm:subquad-H-unimodal}
with $s(r)=2r/(3+r)$.

It remains to treat $k=1$. If $n\ge30$, then
$r\le2/(n-2)\le1/14$, and
Lemma~\ref{sub12:lem:turn3-anchor-switch-highq} gives $N_q(r)>0$,
a contradiction. Thus $n\le29$, and $r_0=1/(n-1)\ge1/28$.
Since $r_0<r_p$, the middle inequality in
\eqref{sub12:eq:turn3-highq-fixed-radius} implies $R>r_p$ and yields
$H(r)>H(R)>H(1/28)\ge H(r_0)$.
This exhausts all cases and proves the theorem.
\end{proof}

\paragraph{The exact finite certificate.}
The checker uses only integer and rational arithmetic. Its source is
provided in the companion archive:
\begin{center}\small
\href{https://github.com/zhangteng2000/schoenberg-validation/blob/e7435c5521dcd6a10131ba04a6da500b1404800b/certificates/turn3_source_anchor_branch_rational_certificate.py}{\nolinkurl{turn3_source_anchor_branch_rational_certificate.py}}.
\end{center}
It verifies all thirteen inequalities in \eqref{eq:five-fixed-radii}
and \eqref{sub12:eq:turn3-highq-fixed-radius} throughout their closed parameter intervals.

For a fixed rational $R$, write $q=1-\alpha$ and
$L_x=\log(1/x)$. Then
\[
N_{1-\alpha}(R)=2-R+\alpha R +R^3(1-2R-\alpha)e^{\alpha L_R} -(1-R)\{2+R+R^2+2R^3-\alpha(R+R^2)\} e^{\alpha L_{1-R}}.
\]
For a comparison $H(R)>H(s)$, cross multiply the positive
denominators and use
\[
\mathcal A_x(\alpha)=(1+x) \{1+x^2e^{\alpha L_x}-(1-x)^2e^{\alpha L_{1-x}}\},\qquad \mathcal B_x(\alpha)=x+x^2e^{\alpha L_x}.
\]
The required numerator is
$\mathcal A_R\mathcal B_s-\mathcal A_s\mathcal B_R$.
Each numerator is therefore a finite sum (over a finite set of exponents $\Lambda$)
$\sum_{\lambda\in\Lambda}(a_\lambda+b_\lambda\alpha)e^{\lambda\alpha}$ with rational $a,b$ and
nonnegative logarithmic $\lambda$. Its value at $\alpha=0$
vanishes by an exact rational identity, which the script checks.
The script divides out this zero and retains Taylor degrees
$1$ through $40$ of the numerator.

Every logarithm is reduced to a rational argument in $[1,2]$,
plus an integer multiple of $\log2$. For
$z=(x-1)/(x+1)\in[0,1/3]$, it uses the first $55$ terms of
$\log x=2\sum_{j=0}^{\infty}\frac{z^{2j+1}}{2j+1}$,
with tail bounded above by
$2z^{111}/(111(1-z^2))$. All interval operations round outward
on the rational grid $10^{-70}\mathbb Z$.

Let $B>0$ be the upper endpoint of the $\alpha$ interval.
For one summand $(a+b\alpha)e^{\lambda\alpha}$, the absolute
Taylor remainder after division by $\alpha$ is bounded by
$3^8\left\{ \frac{|a|(\lambda B)^{41}}{41!B} +\frac{|b|(\lambda B)^{40}}{40!}\right\}$.
The script checks $0\le\lambda B<8$, and uses
$e^{\lambda B}<e^8<3^8$.
Finally it transforms the enclosed degree-$39$ Taylor polynomial
to the Bernstein basis on the complete prescribed interval.
All forty Bernstein coefficient lower bounds exceed the total
remainder bound, in each of the thirteen cases. Thus each
normalized numerator is strictly positive. The rational interval
calculation covers each of the thirteen intervals without subdivision.

\begin{cor}\label{sub12:cor:turn3-only-right-anchor-remains}
On an integer peak cell it suffices to prove the endpoint
comparison with $E_+$ in place of $\min(E_-,E_+)$. Indeed, wherever
$E_-\le E_+$, the uncorrected value $H_p(r)$ already exceeds the
left lattice endpoint, and the nonnegative correction can only
increase it. Thus the right-anchor comparison implies the required
comparison with the actual minimum of the two anchor coefficients.
\end{cor}
\subsection{A weighted curvature decay on the full peak domain}
\label{sub12:sec:turn3-peak-weighted-decay}

Retain the functions $J(\theta)=H_p(\theta/(1-\theta))$ and
$T(\theta)=16E_+(1-2\theta)/D(1-2\theta)$ from
Theorem~\ref{sub12:thm:turn3-peak-derivative-signs}.

\begin{thm}\label{sub12:thm:turn3-peak-weighted-decay}
For every $1<p<2$ and $0<\theta\le2\theta_p$,
\[
\theta J'''(\theta)+\frac12J''(\theta)<0.
\]
Equivalently, $\theta^{\frac{1}{2}}J''(\theta)$ is strictly decreasing
throughout this interval.
\end{thm}

\begin{proof}
Use the formal rational functions and the differentiation operator
$\mathcal D$ from Theorem~\ref{sub12:thm:turn3-peak-derivative-signs}.
The numerator
\[
N_{1/2}=2r^3(1+X)^4(r-1)^2 \left\{\frac{r}{1+r}\mathcal D^3\mathcal H +\frac12\mathcal D^2\mathcal H\right\}
\]
is a polynomial affine in $Z$.  It remains to show $N_{1/2}<0$.
The degree-six binomial majorant $Z_6$ defined in the proof of
Theorem~\ref{sub12:thm:turn3-peak-derivative-signs} also bounds $Z$.
The exact Bernstein certificate for the numerator $N_{1/2}$
verifies that all coefficients of
$(1+r)^{-3}\partial_ZN_{1/2}$
are strictly positive on $0\le q\le1$, $0\le r\le2/17$,
$0\le X\le1$.  Their degree vector is $(3,6,3)$, with $112$
coefficients.  It remains to certify $N_{1/2}(Z_6)\le0$.

For $q\le9/10$, use the power enclosure established in the proof
of Theorem~\ref{sub12:thm:turn3-peak-derivative-signs}:
$X=r(1+2\alpha+2\alpha^2) +\alpha(1-3r-2\alpha r)w$, where $\alpha=1-q$ and $0\le w\le1$.
The polynomial \(-N_{1/2}(Z_6)/[r(1-q)]\) has nonnegative Bernstein
coefficients on each of the rectangles \([2/5,3/5]\times[0,1/12]\)
and \([3/5,9/10]\times[0,2/17]\) in the \((q,r)\)-plane.
In either box its degree vector is $(14,17,4)$.
For $0<q\le2/5$, substitute $r=q^2v/2$ and divide by the exact
factor $q^2$ before mapping $q=2u/5$.  The resulting polynomial
also has all nonnegative Bernstein coefficients; its degree vector
is $(43,17,4)$.

Only the last exponent interval needs a tighter upper bound on $X$.
For $9/10\le q<1$, convexity in $\alpha=1-q$ gives
$r^q-r\le10\alpha(r^{\frac{9}{10}}-r)<\frac3{10}\alpha$ for $0<r\le2/17$.
The function \(r^{\frac{9}{10}}-r\) is increasing on this interval.
At \(R=2/17\) one has
$R^{\frac{9}{10}}<R+\frac3{100}=\frac{251}{1700}$, since $R^9<\left(\frac{251}{1700}\right)^{10}$.
Monotonicity follows from $R<(9/10)^{10}$.
Thus replace the power interval by
$X=r(1+2\alpha+2\alpha^2) +\alpha\left(\frac3{10}-2r-2\alpha r\right)w$.
Its length is positive for $r\le2/17$ and $\alpha\le1/10$.
After substituting $q=9/10+u/10$ and $r=2v/17$, the polynomial
$-N_{1/2}(Z_6)/[r(1-q)]$ again has only nonnegative Bernstein
coefficients, of degree vector $(14,17,4)$.

The independent checker
\href{https://github.com/zhangteng2000/schoenberg-validation/blob/e7435c5521dcd6a10131ba04a6da500b1404800b/certificates/turn3_norm_peak_weighted_decay_certificate.py}{\nolinkurl{turn3_norm_peak_weighted_decay_certificate.py}}
constructs $N_{1/2}$ from \eqref{eq:integer-cell-operator},
checks every exact division, and verifies all four Bernstein
certificates with rational arithmetic.  It also verifies the
strictly positive $Z$-coefficient and the two rational power
comparisons used for the improved interval.
Consequently $N_{1/2}(Z_6)\le0$, while strict $Z<Z_6$ and
$\partial_ZN_{1/2}>0$ give $N_{1/2}(Z)<0$.  This proves the theorem.
\end{proof}

\begin{cor}[A local sufficient condition for the backward scalar lemma]
\label{sub12:cor:turn3-backward-local-sufficient}
Fix $a\in[\theta_p/2,\theta_p]$.  If
\[
aJ''(a)\le\frac94T(a),
\]
then
\[
J(2a)-J(a)-aJ'(a)<aT(a).
\]
\end{cor}

\begin{proof}
The weighted decay theorem gives
$J''(s)\le J''(a)(a/s)^{\frac{1}{2}}$ for $a\le s\le2a$.
The integral Taylor formula therefore yields
\[
J(2a)-J(a)-aJ'(a) =\int_a^{2a}(2a-s)J''(s)\,ds \le a^2J''(a)\frac{8\sqrt2-10}{3} \le aT(a)\frac34(8\sqrt2-10)<aT(a).
\]
The last strict sign is the exact inequality $12\sqrt2<17$.
The local curvature condition is proved for every such $a$ by
Lemma~\ref{sub12:lem:turn3-smallq-curvature-ratio} and
Theorem~\ref{sub12:thm:turn3-compact-curvature-ratio}, together with the
lower-location theorem used in the small-$q$ case.
\end{proof}
\subsection{An explicit lower peak location for small exponents}
\label{sub12:sec:turn3-lower-peak-smallq}

\begin{lem}\label{sub12:lem:turn3-lower-peak-smallq}
For $0<q=p-1\le1/20$, the continuous binary minimizer satisfies
\[
r_p>\frac{2q^2}{8+9q+23q^2}>\frac{q^2}{4+5q+10q^2}.
\]
\end{lem}
\begin{proof}
Put $D=4+9q/2+23q^2/2$, $r=q^2/D$, $L=\log(1/r)$, and $x=qL$.
The strictly concave derivative numerator from the unimodality proof
has the upper bound
$K(r)<-q^2+r\{c_2e^x+b_2+e^{-x}\}$,
where
$c_2=1+\frac{q(1-q)}2$ and $b_2=2+\frac{q(1-q)^2}2$.
It suffices to show that the braces are smaller than $D$.

Since $r\le1/1600$, one has $L>4$. Moreover
$L'=\frac{D'}D-\frac2q>-\frac2q$.
Thus both $qL$ and $qL^2$ are strictly increasing in $q$ on the
indicated interval: their derivatives exceed $L-2$ and
$L(L-4)$, respectively. At $q=1/20$ one has $L=\log(3403/2)<\log1710<15/2$.
Consequently
$0<x<\frac38$ and $qL^2<\frac{45}{16}$.
The elementary bound $\log(3/2)>2/5>3/8$ gives $e^x<3/2$ and
$\cosh x<9/8$. For the hyperbolic-cosine bound, the inequality $(2j)!\ge2^j$
gives $\cosh x\le(1-x^2/2)^{-1}\le128/119<9/8$.
Termwise comparison of the power series gives
$2(\cosh x-1)\le x^2\cosh x$. Therefore
\[
c_2e^x+b_2+e^{-x} \le4+2(\cosh x-1)+\frac q2e^x+\frac q2 <4+q\left(\frac{45}{16}\frac98+\frac34+\frac12\right) =4+\frac{565}{128}q<4+\frac92q<D.
\]
It follows that $K(r)<0$. Since the sign of $K$ is the sign of
$H_p'$ and $r_p$ is its unique zero, the first lower bound follows. The second follows from
$4+9q/2+23q^2/2<4+5q+10q^2$ when $q<1/3$.

The first positive term of the atanh series gives
$\log(3/2)>2/5$. To check $\log1710<15/2$, the first thirteen
terms of the exponential series give
$e^{\frac{15}{2}}>\sum_{j=0}^{12}\frac{(15/2)^j}{j!} =\frac{559014132443}{322961408}>1710$.
\end{proof}
\subsection{A quantitative lower location bound at every exponent}
\label{sub12:sec:turn3-sharp-lower-peak-location}

\begin{thm}\label{sub12:thm:turn3-sharp-lower-peak-location}
For every $0<q=p-1<1$, the unique binary minimizer satisfies
\[
r_p>\frac{2q^2}{8+9q+23q^2}.
\]
Consequently, if $\theta_p/2\le a\le\theta_p$ and
$r=a/(1-a)$, then
\[
\frac{q^2}{8+9q+24q^2}<r\le r_p.
\]
\end{thm}

\begin{proof}
Lemma~\ref{sub12:lem:turn3-lower-peak-smallq} proves the first assertion
for $0<q\le1/20$.  We give an exact finite certificate for the
remaining compact interval, including its limiting endpoint $q=1$.
Put
$d(q)=4+\frac92q+\frac{23}{2}q^2$ and $\rho(q)=\frac{q^2}{d(q)}$.
It suffices to prove $H_{1+q}'(\rho(q))<0$, since this derivative
has its unique sign change at $r_p$.

For $1/20\le q\le1/2$, the binomial expansion from the
unique-minimum proof gives
$K(\rho)<-q^2+\rho\{c_2e^x+b_2+e^{-x}\}, \quad x=q\log(1/\rho)$,
where $c_2=1+q(1-q)/2$, $b_2=2+q(1-q)^2/2$.
The exact interval certificate proves
$c_2e^x+b_2+e^{-x}-d(q)<0$ for $1/20\le q\le1/2$.
Thus $K(\rho)<0$, which is the required derivative sign.

For $1/2\le q<1$, put $\alpha=1-q$,
$A=\log(1/(1-\rho))$, and $B=\log(1/\rho)$.
The exact numerator
$\mathcal N=(\rho+\rho^{1+q})^2H_{1+q}'(\rho)$ is
\begin{align*}
 \mathcal N={}&
 (1-\rho)\{1+\rho+2\rho^2-\alpha(\rho+\rho^2)\}e^{\alpha A}\\
 &+(1-\rho)\rho\{2+\rho+\rho^2-\alpha(1+\rho)\}
                      e^{\alpha(A+B)}\\
 &-1+\rho(1-\rho^2)(-2+\alpha)e^{\alpha B}
       +\rho^4e^{2\alpha B}.
\end{align*}
Its constant term in $\alpha$ vanishes identically for every $\rho$.
To avoid any endpoint division, write the four exponential terms as
$(a_i+b_i\alpha)e^{\alpha\lambda_i}$ and use the exact identity
\[
\frac{\mathcal N}{\alpha} =\sum_{i=1}^4\left\{ a_i\lambda_i\,\operatorname{exprel}(\alpha\lambda_i) +b_i e^{\alpha\lambda_i}\right\}, \qquad \operatorname{exprel}(x)=\sum_{j=0}^{\infty}\frac{x^j}{(j+1)!}.
\]
The right side extends continuously to $q=1$.
The certificate proves that it is strictly negative on
$1/2\le q\le1$.

The interval bounds are verified
in \href{https://github.com/zhangteng2000/schoenberg-validation/blob/e7435c5521dcd6a10131ba04a6da500b1404800b/certificates/turn3_root_lower_peak_certificate.py}{\nolinkurl{turn3_root_lower_peak_certificate.py}}.
All intervals have endpoints on the rational grid $10^{-70}\mathbb Z$,
and each operation rounds outwards using integer arithmetic.
To enclose a logarithm, range-reduce its argument to $[1,2]$ and use
\[
2\sum_{j=0}^{54}\frac{z^{2j+1}}{2j+1} <\log\frac{1+z}{1-z} <2\sum_{j=0}^{54}\frac{z^{2j+1}}{2j+1} +\frac{2z^{111}}{111(1-z^2)}.
\]
The inequality at $z=0$ is interpreted non-strictly.
The series for $e^x$ and $\operatorname{exprel}(x)$ are kept through
degree $60$.  For $0\le x\le M<8$, the value tails are at most
$3^8\frac{M^{61}}{61!}$ and $3^8\frac{M^{61}}{62!}$,
and the derivative tails are at most
$3^8\frac{M^{60}}{60!}$ and $3^8\frac{M^{60}}{61!}$,
respectively, by termwise estimation of the positive Taylor tails and \(e^M<3^8\).

For each rational interval $[a,b]$, the program encloses the
target function at its midpoint and encloses its derivative on
the whole interval by automatic differentiation of the explicit
expressions for \(H_{1+q}'(\rho(q))\) for the relevant exponent range. It accepts the interval only if
$\sup F\left(\frac{a+b}{2}\right) +\frac{b-a}{2}\sup_{[a,b]}|F'|<0$.
Otherwise the program bisects the interval. The mean value theorem makes the
midpoint test valid throughout each accepted interval.
The interval $1/20\le q\le1/2$ is certified by 25 subintervals,
and $1/2\le q\le1$ by 240 subintervals. The maximum bisection
depths are 5 and 9, respectively. Each accepted interval has a strict margin.
Thus the two initial intervals are completely covered and the
required derivative sign holds for every real exponent in them.

This proves the first assertion in all three exponent ranges.
For the consequence, monotonicity of $a\mapsto a/(1-a)$ gives
$r\ge\frac{r_p}{2+r_p} >\frac{q^2}{8+9q+24q^2}$,
while $a\le\theta_p$ gives $r\le r_p$.
\end{proof}
\subsection{The local curvature ratio for small exponents}
\label{sub12:sec:turn3-smallq-curvature-ratio}

\begin{lem}\label{sub12:lem:turn3-smallq-curvature-ratio}
Let $0<q=p-1\le1/20$, and suppose
\[
\frac{q^2}{8+9q+24q^2}\le r\le r_p, \qquad \theta=\frac r{1+r}.
\]
For $J(\theta)=H_p(r)$ and
$T(\theta)=16E_+(1-2\theta)/D(1-2\theta)$, one has
\[
\theta J''(\theta)<\frac94T(\theta).
\]
\end{lem}
\begin{proof}
Put $X=r^q$, $Z=(1-r)^q$, and $D_0=8+9q+24q^2$.
Equation~\eqref{sub12:eq:subquad-explicit-b-location} gives
$r\le r_p<\frac{q^2}{4}\le\frac1{1600}$.
We first note that $X>1/2$. For $R=q^2/D_0$, the function
$q\log(D_0/q^2)$ is strictly increasing on $(0,1/20]$, since
its derivative is
$\log(D_0/q^2)+qD_0'/D_0-2>0$.
Its endpoint value is $\log3404/20<\log4096/20<\log2$.
Consequently $X\ge R^q>1/2$.

\paragraph{An upper bound on the curvature.}
The strict concavity of
$K(r)=r^{1-q}(1+r^q)^2H_p'(r)$, with $K(0+)=-q^2$, gives
$H_p''(r)-\frac{q(1-X)}{r(1+X)}H_p'(r) \le\frac{q^2r^{q-2}}{(1+X)^2}$.
Indeed $K'(r)\le(K(r)+q^2)/r$, and differentiation of
$H_p'=r^{q-1}K/(1+r^q)^2$ gives the stated upper bound for \(H_p''(r)\).
Since $r\le r_p$, we have $H_p'(r)\le0$. The change of variable
$\theta=r/(1+r)$ therefore yields
\[
\theta J''(\theta) =2r(1+r)^2H_p'(r)+r(1+r)^3H_p''(r) \le\frac{q^2}{r}\frac{X(1+r)^3}{(1+X)^2} \le\frac{D_0(1+r)^3}{4}.
\]
The bound by \(D_0(1+r)^3/4\) uses $r\ge q^2/D_0$ and $X/(1+X)^2\le1/4$.

\paragraph{A lower bound on the right-anchor coefficient.}
Write $Z=1-\sum_{j=1}^{\infty}a_jr^j$, with $a_j>0$ and $a_1=q$.
The exact formula for $T$ gives
$T(\theta)=\frac{(1+r)^2}{1+X} \left\{X+C_q-\sum_{j=4}^{\infty}B_jr^{j-3}\right\}$,
where
$C_q=\frac{(1-q)(2-q)(q+3)}6 =1-\frac76q+\frac16q^3$,
$B_j=2a_j+(q+1)a_{j-1}+qa_{j-2}+a_{j-3}>0$.
The product formula gives $a_j\le q$ for every $j\ge1$.
Thus $B_j\le(4+2q)q<5q$, and
$0<\sum_{j=4}^{\infty}B_jr^{j-3} <\frac{5qr}{1-r}\le\frac{5q}{1599}<\frac q{300}$.
It follows, using $1+X>3/2$, that
\[
T(\theta) >(1+r)^2\left(1-\frac{(7/6+1/300)q}{1+X}\right) >(1+r)^2\left(1-\frac{39}{50}q\right).
\]
Finally,
\(D_0(1+r) \le\frac{851}{100}\frac{1601}{1600} <9\frac{961}{1000} \le9\left(1-\frac{39}{50}q\right)\).
The upper bound for \(\theta J''(\theta)\) and the lower bound for
\(T(\theta)\) now give \(\theta J''(\theta)<9T(\theta)/4\).
\end{proof}
\subsection{The local curvature bound away from the linear exponent}
\label{sub12:sec:turn3-compact-curvature-ratio}

Use the functions $J(\theta)$ and $T(\theta)$ defined in
Theorem~\ref{sub12:thm:turn3-peak-derivative-signs}, and set $q=p-1$.

\begin{thm}\label{sub12:thm:turn3-compact-curvature-ratio}
For $1/20\le q<1$ and $\theta_p/2\le\theta\le\theta_p$,
\[
\theta J''(\theta)<\frac94T(\theta).
\]
\end{thm}

\begin{proof}
Put $r=\theta/(1-\theta)$ and
$L(q)=\frac{q^2}{8+9q+24q^2}$.
Theorem~\ref{sub12:thm:turn3-sharp-lower-peak-location} gives
$r_p>2q^2/(8+9q+23q^2)$.  Consequently
$r\ge\frac{r_p}{2+r_p}>L(q)$.
The upper location bound \eqref{sub12:eq:subquad-explicit-b-location} for $q\le2/5$,
Lemma~\ref{sub12:lem:turn3-middle-location} for $2/5\le q\le3/5$,
and Theorem~\ref{sub12:thm:turn3-nineteen-location} for $q\ge3/5$
give the following enlarged domains:
\begin{equation}\label{sub12:eq:turn3-compact-ratio-domains}
 L(q)\le r\le U(q),\qquad
 U(q)=
 \begin{cases}
 q^2/(4+q^2),&1/20\le q\le2/5,\\
 1/25,&2/5\le q\le3/5,\\
 1/18,&3/5\le q\le1.
 \end{cases}
\end{equation}
We prove the curvature inequality on the full ranges in \eqref{sub12:eq:turn3-compact-ratio-domains}.

Temporarily regard $r,q,X,Z$ as independent variables.  Let
$\mathcal H$, $\mathcal T$, and $\mathcal D$ be the rational
functions and differentiation operator from
Theorem~\ref{sub12:thm:turn3-peak-derivative-signs}.  Define
$N=r^3(1+X)^3(1-r) \left\{\frac94\mathcal T- \frac{r}{1+r}\mathcal D^2\mathcal H\right\}$.
This is a polynomial, affine in $Z$.  Its $Z$ coefficient is
strictly positive on the larger box
$0\le q\le1$, $0\le r\le2/17$, $0\le X\le1$.
An exact certificate is obtained by substituting
$q=u$, $r=2v/17$, $X=w$ into
$(1+r)^{-2}\partial_ZN$: all $54$ Bernstein coefficients,
of degree vector $(2,5,2)$, are strictly positive.

Let $a_j=q(1-q)\cdots(j-1-q)/j!$ for $j\ge1$,
with $a_1=q$.  The positive coefficients \(a_j\) decrease with $j$.
The binomial series therefore gives the strict lower bound
$Z=(1-r)^q>Z_{\mathrm{lo}}:={} 1-\sum_{j=1}^{9}a_jr^j-\frac{a_{10}r^{10}}{1-r}$.
The direction of the bound is essential: positivity of
$\partial_ZN$ gives $N(X,Z)>N(X,Z_{\mathrm{lo}})$.

Write $\alpha=1-q$ and introduce the regularized power
$V=\frac{r^{-\alpha}-1}{\alpha}$, so that $X=r(1+\alpha V)$.
At $\alpha=0$ the continuous value is $V=\log(1/r)$.
Exact polynomial division gives
\begin{equation}
\label{sub12:eq:turn3-ratio-normalized-polynomial} \mathcal R(r,q,V)= \frac{14515200(1-r)}{(1-q)r^3(1+r)^2} N\bigl(r,q,r(1+(1-q)V),Z_{\mathrm{lo}}\bigr).
\end{equation}
The right side is a polynomial with integer coefficients,
of degree vector $(14,13,3)$ and with $412$ nonzero terms.
Thus no division by $1-q$ remains in the certified expression.

The standalone checker
\href{https://github.com/zhangteng2000/schoenberg-validation/blob/e7435c5521dcd6a10131ba04a6da500b1404800b/certificates/turn3_lp4_curvature_ratio_certificate.py}{\nolinkurl{turn3_lp4_curvature_ratio_certificate.py}}
constructs \(\mathcal R\) from \eqref{sub12:eq:turn3-ratio-normalized-polynomial},
verifies every exact division and the positive $Z$ coefficient,
and proves
$\mathcal R(r,q,V)>200$
throughout \eqref{sub12:eq:turn3-compact-ratio-domains}, with the actual
regularized power $V$.  Its complete enclosure procedure is given
below.  Since every factor removed in
\eqref{sub12:eq:turn3-ratio-normalized-polynomial} is positive for
$q<1$, this proves $N>0$ and hence the claimed curvature bound.
\end{proof}

\paragraph{A reproducible rational interval certificate.}
All interval endpoints lie on the rational grid $10^{-60}\mathbb Z$.
Addition, multiplication, reciprocal, integer powers, and division
round outwards using integer floor and ceiling operations.  No
floating-point value enters any sign decision.

For each of the three regions in
\eqref{sub12:eq:turn3-compact-ratio-domains}, use the parameterization
$r=L(q)+(U(q)-L(q))v$ for $0\le v\le1$.
Both rational functions are evaluated by outward interval arithmetic.
The logarithm of a positive rational endpoint is reduced to an
argument in $[1,2]$, together with an integer multiple of $\log2$.
For $z=(x-1)/(x+1)\in[0,1/3]$, the first $80$ terms of
$\log x=2\sum_{j=0}^{\infty}\frac{z^{2j+1}}{2j+1}$
have a remainder between zero and
$2z^{161}/(161(1-z^2))$.  Monotonicity then encloses
$\Lambda=\log(1/r)$ on each interval box.

Set $x=(1-q)\Lambda\ge0$.  The regularized power is enclosed using
$V=\Lambda\sum_{j=0}^{\infty}\frac{x^j}{(j+1)!}$.
The checker keeps terms $0$ through $49$.  If
$m$ is the integer ceiling of the interval upper bound for $x$,
the remaining positive tail is at most
$\frac{x^{50}}{51!}\,3^m$,
because $e^x\le3^m$ (with the zero case interpreted directly).
This enclosure remains regular at $q=1$.

The integer polynomial $\mathcal R$ is evaluated by nested Horner
rules.  A box is accepted only when its lower interval endpoint
is strictly positive.  Otherwise the checker bisects one of the
two parameter intervals, with a fixed deterministic rule, and
continues.  It starts from the three complete rational rectangles
and accepts $4156$ boxes after $8309$ evaluations; all accepted
lower bounds exceed $200$.  Its maximum-depth guard is never
reached.  Every bisection exactly partitions its parent rectangle,
so the accepted boxes cover the entire required parameter domain.

\subsection{The backward remainder and the full finite-degree sharp constant}
\label{sub12:sec:turn3-backward-remainder-reduction}

The derivative signs in Theorem~\ref{sub12:thm:turn3-peak-derivative-signs},
the weighted curvature estimate in
Theorem~\ref{sub12:thm:turn3-peak-weighted-decay}, and the
curvature bounds in Lemma~\ref{sub12:lem:turn3-smallq-curvature-ratio}
and Theorem~\ref{sub12:thm:turn3-compact-curvature-ratio}
combine to prove the scalar inequality required for integer rounding.

Put $q=p-1$, $\theta_p=r_p/(1+r_p)$, and
\[
J(\theta)=H_p\!\left(\frac{\theta}{1-\theta}\right), \qquad T(\theta)=\frac{16E_+(1-2\theta)}{D(1-2\theta)}.
\]
The right-anchor chord comparison on a peak cell is
\begin{equation}
\label{sub12:eq:turn3-right-chord-theta} J(\theta)+nT(\theta)(\theta-a)(b-\theta) \ge \frac{b-\theta}{b-a}J(a)+\frac{\theta-a}{b-a}J(b), \quad a=\frac{k}{n},\quad b=\frac{k+1}{n},
\end{equation}
where $k\ge1$ and $a\le\theta_p\le b$.

\begin{lem}[Backward-remainder sufficient criterion]
\label{sub12:lem:turn3-backward-sufficient}
Assume that
\begin{equation}
\label{sub12:eq:turn3-backward-scalar} J(2a)-J(a)-aJ'(a)\le aT(a) \quad\text{for every }a\in[\theta_p/2,\theta_p].
\end{equation}
Then \eqref{sub12:eq:turn3-right-chord-theta} holds on every peak cell,
for every real $n>0$ and integer $k\ge1$ for which the cell exists.
\end{lem}

\begin{proof}
Theorem~\ref{sub12:thm:turn3-peak-derivative-signs} gives
$J''>0$, $J'''<0$, and $T'>0$ throughout the interval needed below.
Write $h=b-a=1/n$.  Since $a=kh$ with $k\ge1$, one has $h\le a$.
The peak condition gives $a\le\theta_p\le a+h\le2a$,
so $\theta_p/2\le a\le\theta_p$.
The error of linear interpolation is
$\frac{b-\theta}{h}J(a)+\frac{\theta-a}{h}J(b)-J(\theta) =(\theta-a)(b-\theta)J[a,\theta,b]$.
Because $J'''\le0$, the second divided difference is nonincreasing
in each node.  (This follows directly from the simplex integral
formula for divided differences, including coincident nodes.) Hence
$J[a,\theta,a+h]\le J[a,a,a+h]$.
Furthermore,
$hJ[a,a,a+h] =\frac{J(a+h)-J(a)-hJ'(a)}h =\int_0^h\left(1-\frac{s}{h}\right)J''(a+s)\,ds$.
Its derivative with respect to $h$ equals
$h^{-2}\int_0^h sJ''(a+s)\,ds\ge0$.  Therefore
$hJ[a,\theta,a+h] \le aJ[a,a,2a] =\frac{J(2a)-J(a)-aJ'(a)}a \le T(a)\le T(\theta)$.
Since $n=1/h$, this is exactly
\eqref{sub12:eq:turn3-right-chord-theta}.
\end{proof}

\begin{thm}[The strict backward-remainder inequality]
\label{sub12:thm:turn3-backward-remainder-proved}
For every $0<q<1$ and $a\in[\theta_p/2,\theta_p]$,
\[
J(2a)-J(a)-aJ'(a)<aT(a).
\]
\end{thm}

\begin{proof}
For $q\ge1/20$, Theorem~\ref{sub12:thm:turn3-compact-curvature-ratio}
gives $aJ''(a)<9T(a)/4$.
For $q\le1/20$, write $r=a/(1-a)$.
Theorem~\ref{sub12:thm:turn3-sharp-lower-peak-location} gives
$r\ge\frac{r_p}{2+r_p}>\frac{q^2}{8+9q+24q^2}$ and $r\le r_p$.
Lemma~\ref{sub12:lem:turn3-smallq-curvature-ratio} therefore gives
$aJ''(a)<9T(a)/4$ for $0<q\le1/20$.  The weighted decay theorem and its Peano
integral consequence, Corollary~\ref{sub12:cor:turn3-backward-local-sufficient},
now imply the asserted strict backward-remainder inequality.
\end{proof}

\begin{cor}[The full finite-degree sharp constant]
\label{sub12:cor:turn3-one-backward-inequality-suffices}
The full finite-degree sharp-constant formula is
\[
C_{n,p}=\max\{S_{n,p},\mathcal B_{n,p}\} \qquad(n\ge4,\ 1<p<2).
\]
\end{cor}

\begin{proof}
Theorem~\ref{sub12:thm:turn3-backward-remainder-proved} supplies
\eqref{sub12:eq:turn3-backward-scalar} for every exponent.
The cases $n\theta_p\le1$ and $n\theta_p\in\mathbb Z$ have already
been proved. For a nonintegral $n\theta_p>1$,
Lemma~\ref{sub12:lem:turn3-backward-sufficient} supplies
the chord bound with $E_+$, which is stronger than the endpoint
minimum comparison.  Theorem~\ref{sub12:thm:turn3-integer-anchor-branch}
shows that wherever $E_-\le E_+$ the uncorrected deficit already
exceeds the lattice value at multiplicity \(k\). Positivity of \(E_-\)
and \(E_+\) on the cell then gives the bound with \(\min(E_-,E_+)\).
Theorem~\ref{sub12:thm:turn3-peak-cell-sufficient} gives the upper bound
\(\mathcal B_{n,p}\) on the cell and its complement.
The sparse and binary examples attain the indicated maximum.
\end{proof}

Inequality~\eqref{sub12:eq:turn3-backward-scalar} is independent of the
degree and the integer cell index. In the original radius coordinate it reads
\[
H_p\!\left(\frac{2r}{1-r}\right)-H_p(r) -r(1+r)H_p'(r) \le\frac{r}{1+r}\,\mathcal T_q(r), \qquad \frac{r_p}{2+r_p}\le r\le r_p,
\]
where
\[
\mathcal T_q(r)=\frac{(1+r)^2}{r^3(1+r^q)} \bigl[r^{q+3}+(1-r)^q\{r^3+qr^2+(q+1)r+2\} +qr^2+qr-r-2\bigr].
\]
The identities are algebraic consequences of the explicit $E_+$ formula.

\subsection{The complete equality classification in every finite degree}
\label{sub12:sec:turn3-full-finite-equality}

\begin{thm}[All equality configurations]
\label{sub12:thm:turn3-full-finite-equality}
Let $n\ge4$, $1<p<2$, and put
$\Lambda_{n,p}=\max\{S_{n,p},\mathcal B_{n,p}\}=C_{n,p}$.
Apart from the all-zero root list, equality in
\[
\sum_{j=1}^{n-1}|w_j|^p\le \Lambda_{n,p}\sum_{i=1}^n|z_i|^p, \qquad \sum_{i=1}^nz_i=0,
\]
holds exactly for the following lists, up to permutation and a common
nonzero complex scalar:
\[
(1,-1,0^{[n-2]}) \quad\text{if }S_{n,p}=\Lambda_{n,p},
\]
and
\[
((n-k)^{[k]},(-k)^{[n-k]}) \quad\text{for each }1\le k\le\lfloor n/2\rfloor \text{ with }R_{n,p}(k)=\Lambda_{n,p}.
\]
Every indicated family occurs, including all families at a tie.
Equivalently, apart from $p(z)=Az^n$, the equality polynomials are
\[
p(z)=A z^{n-2}(z^2-a^2) \quad\text{or}\quad p(z)=A(z-(n-k)a)^k(z+ka)^{n-k}, \qquad A,a\ne0,
\]
with \(S_{n,p}=\Lambda_{n,p}\) for the first family and
\(R_{n,p}(k)=\Lambda_{n,p}\) for the second.
\end{thm}

\begin{proof}
Corollary~\ref{sub12:cor:turn3-one-backward-inequality-suffices}
establishes the sharp constant. We now determine every equality case.
We first classify its nonzero real extremizers.
Scale and, if needed, reflect such a list so that its extreme roots
are $t-1,t+1$, where $0\le t\le c=(n-2)/n$.
The endpoint-gradient reduction gives a spectral probability measure
$\mu$ with
$m_1=0$, $m_2=c$, $m_3=-(1-c)t$, and $G(t,\mu)=\Lambda_{n,p}$.
For the $n-2$ other roots write $y_i=x_i-t\in[-1,1]$.
The exact root-defect identity
\eqref{eq:section4-root-defect} is
\begin{equation}\label{sub12:eq:turn3-equality-root-defect}
 \Delta(t,\mu):=\int(1-y^2)(y-t)^2\,d\mu(y)
 =\frac2{n^2}\sum_{\substack{1\le i\le n\\i\text{ internal}}}(1-y_i^2).
\end{equation}
In particular, $\Delta=0$ forces every root to be one of the two
extreme values.

If $n\theta_p\le1$, the complete equality assertion in
Corollary~\ref{sub12:cor:turn3-no-rounding-cell-exact} gives exactly the
sparse list and the one-outlier list, according to which value
equals $\Lambda_{n,p}$. For clarity, its two possible moment maximizers
are Gaussian at $t=0$ and $t=c$. At $t=0$, Gaussianity gives
$m_4=c^2$; the fourth-moment identity
\eqref{eq:section4-fourth-identity} then forces all internal
roots to vanish. At $t=c$, their sum is $-(n-2)$ and each is at
least $-1$, forcing them all to equal $-1$.

Suppose next that $n\theta_p>1$ is an integer.
The exact uncapped three-moment theorem
(Theorem~\ref{sub12:thm:turn3-uncapped-lp3-exact}) applies directly.
Its unique maximizing configuration is the intrinsic-cap measure
on $\{-1,t_p,1\}$ at $t=t_p$. Its value is
$1-H_p(r_p)/n=\mathcal B_{n,p}$, since this continuous maximum is an
allowed binary lattice value. It strictly exceeds both $S_{n,p}$
and $R_{n,p}(1)$. The intrinsic-cap measure has $\Delta=0$, so
\eqref{sub12:eq:turn3-equality-root-defect} forces the roots to be binary.
Their multiplicities are then fixed by centering, with
$k=n\theta_p$. This argument uses the uncapped moment theorem
and does not require a separate weighted equality theorem.

It remains to consider a nonintegral $n\theta_p>1$.
Put $k=\lfloor n\theta_p\rfloor$ and
$\ell=1-\frac{2(k+1)}n<t_p<u=1-\frac{2k}n$.
The proof of Theorem~\ref{sub12:thm:turn3-peak-cell-sufficient}
gives \(\mathcal B_{n,p}>S_{n,p}\). We show that equality requires a maximizing lattice endpoint.

First consider $\ell<t<u$, or equivalently
$a<\theta<b$ with $\theta=(1-t)/2$, $a=k/n$, and $b=(k+1)/n$.
If $E_-(t)\le E_+(t)$,
Theorem~\ref{sub12:thm:turn3-integer-anchor-branch} gives
$H_p((1-t)/(1+t))>H_p(k/(n-k))$. The uncorrected gradient bound
is therefore already strictly smaller than the binary envelope.
If $E_+(t)<E_-(t)$, the strict backward-remainder theorem and its
divided-difference proof give a strict chord comparison. More
explicitly, with
$\rho=\frac34(8\sqrt2-10)<1$,
the curvature proof gives
$hJ[a,\theta,b]\le\rho T(\theta)$, where $h=1/n$.
Consequently
\[
J(\theta)+nT(\theta)(\theta-a)(b-\theta) \ge \frac{b-\theta}{b-a}J(a)+\frac{\theta-a}{b-a}J(b) +(1-\rho)nT(\theta)(\theta-a)(b-\theta) >\min\{J(a),J(b)\}.
\]
The last term is strictly positive. Thus in either anchor branch
every open peak-cell point has $G(t,\mu)<\mathcal B_{n,p}$.

Next, the uncapped moment classification used in
Theorem~\ref{sub12:thm:turn3-peak-cell-sufficient} shows that equality
outside the open cell can occur only at $t=\ell$ or $t=u$.
Indeed, on $[0,\ell]$ every nondegenerate right Radau maximizer
must be at the right endpoint, because its feasible derivative
in increasing $t$ is strictly positive. Interior left Radau
points are not local maxima. The Gaussian values are strictly
smaller than $\max\{S_{n,p},R_{n,p}(1)\}\le \mathcal B_{n,p}$ for
$0<t<c$, and the Gaussian value at zero is $S_{n,p}<\mathcal B_{n,p}$.
The intrinsic-cap value is strictly increasing on $[0,\ell]$.
The strict Gaussian bound, the positive right Radau derivative,
and the strict monotonicity of $R_{\rm cont}$ exclude equality
for every $0\le t<\ell$.
On $[u,c]$, the strict positive $t$-derivative in
Theorem~\ref{sub12:thm:turn3-right-radau-elimination}
excludes all
nondegenerate inactive-cap maximizers, including one at its left
endpoint. The intrinsic-cap value is strictly decreasing for
$t>u$. If $k>1$, then $R_{n,p}(1)<\mathcal B_{n,p}$, so the Gaussian
endpoint $c$ is excluded as well. If $k=1$, then $u=c$ already.
This proves the asserted restriction to the two lattice endpoints.

At either $t\in\{\ell,u\}$ the anchor coefficient
$E_*(t)=\min\{E_-(t),E_+(t)\}$ is strictly positive.
The quartic anchor, before applying the integer lower bound on
the defect, gives
$G(t,\mu)\le R_{\rm cont}(t) -\frac{2E_*(t)}{D_t}\Delta(t,\mu) \le \mathcal B_{n,p}$.
Thus equality requires both $R_{\rm cont}(t)=\mathcal B_{n,p}$ and
$\Delta=0$. By \eqref{sub12:eq:turn3-equality-root-defect} every root
is an endpoint. If the positive endpoint has multiplicity $j$,
centering gives $t=1-2j/n$. The list is therefore a nonzero real
multiple of $((n-j)^{[j]},(-j)^{[n-j]})$, with
$R_{n,p}(j)=\mathcal B_{n,p}$. This completes the real classification.

Finally let a complex list attain equality. In the angular
Schur--Jensen reduction, equality of the integrated quantities
forces equality in the Hermitian Schur-diagonal Jensen inequality
for almost every angle. Strict convexity of $|x|^p$ implies that
each such Hermitian matrix is diagonal in the fixed Schur basis.
If the upper-triangular Schur matrix had a nonzero off-diagonal
entry $T_{ij}$, the corresponding entry of its Hermitian part
would be $e^{-i\varphi}T_{ij}/2\ne0$ for every angle.
Hence the critical-point compression is normal.

Write $e=n^{-\frac{1}{2}}\mathbf1$, $Q=I-ee^*$,
$D=\operatorname{diag}(z_i)$ and $A=QDQ|_{e^\perp}$.
Since the roots are centered, extending the commutator by zero
on $e$ gives
$[A,A^*]=\frac1n\bigl(\overline z\,z^T-zz^*\bigr)$.
For $z\ne0$, normality therefore implies
$z=\omega\overline z$ for some $|\omega|=1$.
All roots lie on a common line through zero. Dividing by a
nonzero complex scalar reduces the list to the real
classification in the real part of this proof.

Conversely, differentiating the two displayed polynomial families
gives exactly $S_{n,p}$ and $R_{n,p}(k)$, respectively.
They therefore attain equality precisely under the stated
maximizing conditions. The polynomial $Az^n$ gives the trivial
equality separately.
\end{proof}
\subsection{The two smallest nontrivial degrees}
\label{sub12:sec:turn3-degrees-two-three}

\begin{prop}\label{sub12:prop:turn3-degrees-two-three}
For $1<p<2$, the sharp centered complex-polynomial constants in
 degrees two and three are
\[
C_{2,p}=0,\qquad C_{3,p}=3^{-\frac{p}{2}}.
\]
In degree three the only nonzero equality root lists are
$(a,-a,0)$ with $a\in\mathbb C\setminus\{0\}$, up to permutation.
Consequently, the sparse/two-value sharp constant formula
\eqref{eq:section4-sharp} also holds for $n=2,3$.
\end{prop}
\begin{proof}
The unique critical point of a centered quadratic is zero, so
$C_{2,p}=0$.

For a real centered triple $x_1+x_2+x_3=0$, one has the algebraic
identity
$\sum_{j=1}^3x_j^4=\frac12\left(\sum_{j=1}^3x_j^2\right)^2$.
Let $\lambda=p/(4-p)\in(0,1)$. Hölder interpolation, with
$1/2=\lambda/p+(1-\lambda)/4$, gives
\(\|x\|_2\le\|x\|_p^\lambda\|x\|_4^{1-\lambda} =2^{-\frac{1-\lambda}{4}}\|x\|_p^\lambda\|x\|_2^{1-\lambda}\).
For a nonzero triple this is equivalent to
\begin{equation}
\label{sub12:eq:turn3-centered-triple-norm} \sum_{j=1}^{3}x_j^2\le 2^{1-\frac{2}{p}}\left(\sum_{j=1}^{3}|x_j|^p\right)^{\frac{2}{p}}.
\end{equation}
Equality in Hölder interpolation forces all nonzero $|x_j|$ to
be equal. Since an odd number of equal nonzero signed terms cannot
sum to zero, equality here forces a sparse pair $(a,-a,0)$.

Now let $z_1+z_2+z_3=0$ be a complex triple. The critical points are
$w_\pm=\pm\sqrt{\frac{z_1^2+z_2^2+z_3^2}{6}}$.
If the sum of squares vanishes, the desired inequality is strict
unless all roots vanish. Otherwise rotate all roots so that this
sum is positive real, and write $z_j=x_j+iy_j$. Then
\[
\left|\sum_{j=1}^{3}z_j^2\right|=\sum_{j=1}^{3}x_j^2-\sum_{j=1}^{3}y_j^2 \le\sum_{j=1}^{3}x_j^2 \le2^{1-\frac{2}{p}}\left(\sum_{j=1}^{3}|x_j|^p\right)^{\frac{2}{p}} \le2^{1-\frac{2}{p}}\left(\sum_{j=1}^{3}|z_j|^p\right)^{\frac{2}{p}}.
\]
Consequently
\[
|w_+|^p+|w_-|^p =2\left(\frac{|\sum_{j=1}^{3}z_j^2|}{6}\right)^{\frac{p}{2}} \le3^{-\frac{p}{2}}\sum_{j=1}^{3}|z_j|^p.
\]
Equality requires all $y_j=0$, followed by equality in
\eqref{sub12:eq:turn3-centered-triple-norm}; hence the stated sparse
family is the complete nonzero equality class. It attains the
constant. The degree-three binary family has three nonzero roots,
so its ratio is strictly smaller than this sparse value. The
maximum formula therefore extends to degree three, and both its
terms equal zero in degree two.
\end{proof}
\section{Logarithmic comparison and auxiliary proofs for \texorpdfstring{\(2<p<3\)}{2 < p < 3}}\label{app:odd23-details}
\subsection{The sharp bound for \texorpdfstring{\(2<p<3\)}{2 < p < 3}}\label{odd23:sec:main-roadmap}
\begin{thm}\label{odd23:thm:main-complete}
Let $n=2m+1\ge5$, $2<p<3$, and let $z_1,\ldots,z_n\in\mathbb C$
be the zeros, counted with multiplicity, of a polynomial $p(z)$ of degree $n$.
Let $w_1,\ldots,w_{n-1}$ be its critical points. If $\sum_{j=1}^{n}z_j=0$, then
\[
 \sum_{k=1}^{n-1}|w_k|^p\le C_{n,p}\sum_{j=1}^n|z_j|^p,
\]
where
\[
 C_{n,p}=
 \frac{(n-3)(n+1)^p+(n-1)^{p+1}+2^{p+1}}
 {(n-1)(n+1)^p+(n+1)(n-1)^p}
 =\frac{(m-1)(m+1)^p+m^{p+1}+1}
 {m(m+1)^p+(m+1)m^p}.
\]
The coefficient is sharp. Equality is attained by
\[
p(z)=(z-(n+1))^{\frac{n-1}{2}}(z+(n-1))^{\frac{n+1}{2}}.
\]
\end{thm}

Matrix compression reduces the complex inequality to real zeros.
The Mellin identity \eqref{odd23:eq:complete-Mellin} then reduces the real
inequality to \(\mathcal J_s\ge0\). By Theorem~\ref{odd23:lm-thm:main},
a negative minimum would be attained on at most four distinct root values.
Subsection~\ref{odd23:sec:complete-assembly} excludes those configurations
using Theorems~\ref{odd23:thm:all-three-log} and~\ref{odd23:mp-thm:onescale}
and Corollary~\ref{odd23:nc-cor:negative-center-one-scale}.

\subsection{Compression and a truncation estimate}
Put $e=n^{-\frac{1}{2}}(1,\dots,1)^T$, $H=e^\perp$, and $Q=I-ee^*$. For $D=\diag(z_1,\dots,z_n)$, let $A=QDQ|_H$. If $p(\zeta)=\prod_{j=1}^{n}(\zeta-z_j)$, the cofactor identity gives
\[
 \det(\zeta I_H-A)
 =\det(\zeta I-D)\,e^*(\zeta I-D)^{-1}e
 =\frac{p'(\zeta)}n.
\]
Initially the equality holds away from the roots, and then everywhere by polynomial identity. Thus the eigenvalues of $A$ are the critical points.

\begin{lem}\label{odd23:lem:truncation}
Assume $\sum_{j=1}^{n}z_j=0$, and let $s_1,\dots,s_{n-1}$ be the singular values of $A$. Then, for every $t\ge0$,
\[
 \sum_{i=1}^{n-1}(s_i^2-t)_+\le \frac{n-2}{n}\sum_{j=1}^{n}(|z_j|^2-t)_+.
\]
\end{lem}
\begin{proof}
Write $a_j=|z_j|^2$ and $u=D^*e\in H$. On $H$,
\(A^*A=Q\diag(a_j)Q-uu^*\).
Set $b_j=(a_j-t)_+$, $r_j=\min(a_j,t)$, $\beta=\sum_{j=1}^{n}b_j$, and $\rho=n^{-1}\sum_{j=1}^{n}r_j$. Then
\[
 A^*A-tI_H\le M:=Q\diag(b_j)Q-uu^*,
 \qquad \tr_H M=\frac{n-2}{n}\beta-\rho.
\]
The positive-part trace is increasing in the Hermitian order: this follows from
$\tr X_+=\max_{0\le E\le I}\tr(EX)$. Therefore it remains to bound $\tr M_+$.

The positive operator $Q\diag(b_j)Q$ minus the rank-one operator $uu^*$ has at most one negative eigenvalue. Moreover, $M\ge-\rho I_H$. To verify the latter assertion when $\rho>0$, use weighted Cauchy--Schwarz and
\[
 \sum_{j=1}^{n}\frac{a_j}{b_j+\rho}
 =n+\sum_{j=1}^{n}\frac{r_j-\rho}{b_j+\rho}
 \le n+\frac1\rho\sum_{j=1}^{n}(r_j-\rho)=n.
\]
Indeed, when $b_j=0$ the relevant inequality is equality, and when $b_j>0$ we have $r_j=t\ge\rho$, so replacing $b_j+\rho$ by $\rho$ increases that summand. Consequently, for any $v\in H$,
\(\frac1n\left|\sum_{j=1}^{n}z_jv_j\right|^2 \le\sum_{j=1}^{n}(b_j+\rho)|v_j|^2\),
which proves $M\ge-\rho I_H$. If $\rho=0$, either $t=0$, in which case $M=A^*A\ge0$, or all $a_j$ vanish. Thus $\tr M_-\le\rho$ in every case, and
\[
 \tr M_+=\tr M+\tr M_-
 \le\frac{n-2}{n}\beta.
\]
\end{proof}

\subsubsection{Moment identities and a reverse concave comparison}
For a centered real list $x$ and its critical list $y$, expansion of
$\operatorname{Tr}(Q\operatorname{diag}(x)Q)^j$, for $j=2,3$, gives
\[
 \sum_{i=1}^{n-1}y_i^2=\frac{n-2}{n}\sum_{i=1}^{n}x_i^2,
 \qquad \sum_{i=1}^{n-1}y_i^3=\frac{n-3}{n}\sum_{i=1}^{n}x_i^3.
\]
The added zero eigenvalue on $\operatorname{span}(\mathbf1)$ has no effect.
For the cubic identity, the trace terms with one rank-one projection
contribute $-3\sum_{i=1}^{n}x_i^3/n$; all remaining terms contain $\sum_{i=1}^{n}x_i=0$.

Subtracting the squared-hinge inequality from the equal second moments
gives, with $c_n=(n-2)/n$,
\begin{equation}\label{eq:truncated-real-squares}
 \sum_{i=1}^{n-1}\min(y_i^2,t)\ge c_n\sum_{j=1}^{n}\min(x_j^2,t),\qquad t\ge0.
\end{equation}
Let \(h:[0,\infty)\to\R\) be increasing and concave, with \(h(0)=0\).
Interpolate \(h\) linearly at zero and at the distinct positive numbers
among \(x_j^2,y_k^2\), and extend the interpolant constantly past the
largest node. Its successive slopes are nonnegative and nonincreasing.
The drop in slope at each node \(t\) is the coefficient of
\(\min(u,t)\) in the interpolant. Thus a nonnegative linear combination
of \eqref{eq:truncated-real-squares} gives the following bound, including
when \(h\) has a jump at zero. If all nodes are zero, both sides vanish.
\begin{equation}\label{odd23:eq:reverse-concave}
 \sum_{i=1}^{n-1}h(y_i^2)\ge c_n\sum_{j=1}^{n}h(x_j^2).
\end{equation}

\subsubsection{Reduction from complex to real roots}
Suppose a coefficient $C$ is valid for every centered real list.
For a complex list put $A=Q\operatorname{diag}(z)Q|_{\mathbf1^\perp}$.
Choose a Schur basis for $A$. In that basis the diagonal of
$\operatorname{Re}(e^{-i\theta}A)$ consists of
$\operatorname{Re}(e^{-i\theta}w_k)$. Scalar Jensen's inequality
applied to the spectral decomposition of this Hermitian matrix gives
\(\sum_{k=1}^{n-1}|\operatorname{Re}(e^{-i\theta}w_k)|^p \le\operatorname{Tr}|\operatorname{Re}(e^{-i\theta}A)|^p\).
The matrix on the right is the critical-point compression for the
centered real roots $x_j(\theta)=\operatorname{Re}(e^{-i\theta}z_j)$.
Apply the assumed real inequality and integrate over $\theta$.
For $\gamma_p=\int_0^{2\pi}|\cos\theta|^p\,d\theta>0$,
\(\int_0^{2\pi}|\operatorname{Re}(e^{-i\theta}\zeta)|^p\,d\theta =\gamma_p|\zeta|^p\).
Cancellation of $\gamma_p$ proves the same coefficient for complex roots.

\subsection{The logarithmic comparison and its Mellin consequence}
\label{odd23:sec:log-framework}
Write $a=m+1$, $b=m$, and use the reference lists
\[
 x_*=(a^{[b]},(-b)^{[a]}),\qquad
 y_*=(a^{[b-1]},(-b)^{[b]},1).
\]
Square brackets indicate multiplicities. For a centered real list $x$
and its critical list $y$, define the nonnegative product multisets
\[
 U=|x|\otimes|y_*|,\qquad V=|y|\otimes|x_*|.
\]
Both have $n(n-1)$ entries and the same sum of squares. Put
\[
 \ell_s(t)=\frac{|t|}{s}-\log(1+|t|/s),\qquad
 \mathcal J_s(x)=\sum_{v\in V}\ell_s(v)-\sum_{u\in U}\ell_s(u).
\]
Equivalently, with $Z_f=\sum_{j=1}^{n}f(x_j)$ and $W_f=\operatorname{Tr}f(A_x)$,
\[
 \mathcal J_s=bW_{\ell_{s/a}}+aW_{\ell_{s/b}}
 -(b-1)Z_{\ell_{s/a}}-bZ_{\ell_{s/b}}-Z_{\ell_s}.
\]
Some component proofs use $\widetilde K_s=-\mathcal J_s$. A common scaling of both reference lists only rescales $s$ and leaves every pointwise sign assertion unchanged.
The rational auxiliary functional is
\(T_s(x)=\sum_{v\in V}\frac{v^2}{s+v} -\sum_{u\in U}\frac{u^2}{s+u}\).
Direct differentiation and expansion at infinity give
\begin{equation}\label{odd23:eq:J-T}
 \partial_s\mathcal J_s=-T_s/s^2,\qquad
 \mathcal J_s=\int_s^\infty T_t/t^2\,dt.
\end{equation}
By \eqref{odd23:eq:J-T}, a pointwise lower bound on \(T_s\) would
suffice, but the argument below estimates \(\mathcal J_s\) directly.

For $2<p<3$ and $r\ge0$, integration by parts gives
\[
 r^p=\frac{p\sin(\pi p)}\pi\int_0^\infty s^{p-1}
 \left[\log(1+r/s)-r/s+\frac{r^2}{2s^2}\right]ds.
\]
The coefficient is positive. The integrand is integrable at both
endpoints. Summing over the two finite product lists and cancelling
their equal second moments yields
\begin{equation}\label{odd23:eq:complete-Mellin}
 \left(\sum_{j=1}^{n}|x_{*,j}|^p\right)
 \left(C_{n,p}\sum_{j=1}^{n}|x_j|^p-\sum_{k=1}^{n-1}|y_k|^p\right)
 =\frac{p\sin(\pi p)}\pi\int_0^\infty s^{p-1}\mathcal J_s(x)\,ds.
\end{equation}
Thus Theorem~\ref{odd23:thm:global-log} implies the real power inequality.
\subsection{The rational comparison for every two-value list}
\label{odd23:sec:two-rational}
The following pointwise rational result implies the logarithmic comparison for every centered two-value list.

\begin{thm}\label{odd23:thm:two-rational}
Every nonzero centered two-value real root list of odd length satisfies
$T_s\ge0$ for all $s>0$. Equality for one, equivalently every, $s>0$
occurs exactly for the balanced multiplicities $m,m+1$.
\end{thm}
\begin{proof}
After reflection and scaling, write the roots as
$h^{[r]},(-r)^{[h]}$, where
\[
 m=r+d,\quad h=r+2d+1,\quad g=2d+1,\qquad r\ge1,\quad d\ge0.
\]
The critical list is $h^{[r-1]},(-r)^{[h-1]},g$.
After cancelling the cross lists, the rational comparison is
\(T_s=\sum_{i=1}^8\frac{\omega_i v_i^2}{s+v_i}\),
where
\[
\begin{array}{c|rrrrrrrr}
v_i&(m+1)h&mh&h&(m+1)r&mr&r&(m+1)g&mg\\\hline
\omega_i&-d&-d-1&-r&d+1&d&-h&m&m+1.
\end{array}
\]
Multiply by $D(s)=\prod_{i=1}^{8}(s+v_i)>0$ and write the numerator as
$N(s)=\sum_{j=0}^6N_js^j$. The exact verifier
\href{https://github.com/zhangteng2000/schoenberg-validation/blob/e7435c5521dcd6a10131ba04a6da500b1404800b/certificates/turn2_two_value_T_certificate.py}{\nolinkurl{turn2_two_value_T_certificate.py}} reconstructs $N$
from the eight pairs \((v_i,\omega_i)\). After $r=R+1$, each $N_j$ has only
nonnegative integer coefficients in $R,d$. Their respective nonzero
coefficient counts are $140,127,111,87,64,43,26$.
This proves $N\ge0$ for $R,d,s\ge0$.
Moreover
\(N_0=4d\,r^4m^4(2d+1)^2(m+1)^3h^3\).
It is strictly positive for $d>0$, so $T_s>0$ then. For $d=0$
the two cross lists are identical and $T_s=0$.
\end{proof}

Since $\widetilde K_s'=T_s/s^2$ and $\widetilde K_s\to0$ at infinity,
this theorem proves $\widetilde K_s\le0$ pointwise for every two-value
list. In particular, the Schur-concavity argument for a flat majority
with positive residual compression proves the pointwise assertion
\(\widetilde K_s(a_1,\ldots,a_r) \le\widetilde K_s((h/r)^{[r]})\le0\).
The monotonicity in the nonpositive-middle three-value proof therefore
also gives $\widetilde K_s\le0$ pointwise there, rather than only
its integrated power consequence.
\subsection{A flat majority with positive residual compression}
\label{odd23:sec:flat-psd}

\begin{thm}\label{thm:flat-majority}
Let $n=2m+1$, $m\ge2$, and let $1\le r\le m$, $h=n-r$.
Suppose $a_1,\ldots,a_r>0$ satisfy
\[
 \sum_{i=1}^r a_i=h,\qquad \sum_{i=1}^r a_i^{-1}\le h.
\]
Then the centered polynomial
\(p(z)=(z+1)^h\prod_{i=1}^r(z-a_i)\)
satisfies \(W_p\le B_{n,p}Z_p\) for every \(2<p<3\).
For this family, equality occurs exactly when $r=m$ and
$a_1=\cdots=a_m=(m+1)/m$.
\end{thm}

At the uniform vector, Theorem~\ref{odd23:thm:two-rational} and
\eqref{odd23:eq:J-T} give the pointwise logarithmic comparison.
We extend it to all vectors satisfying the hypotheses of Theorem~\ref{thm:flat-majority}.

Set \(D=\diag(a_1,\ldots,a_r)\), \(u_i=\sqrt{(a_i+1)/n}\),
\(B=D-uu^T\), and \(q(z)=\prod_{i=1}^r(z-a_i)\).
The determinant lemma gives
\[
  \det(zI-B)=q(z)\left(1+\frac1n\sum_{i=1}^r
  \frac{a_i+1}{z-a_i}\right)
  =\frac{h q(z)+(z+1)q'(z)}n
  =\frac{p'(z)}{n(z+1)^{h-1}}.
\]
The middle expression is a polynomial, so the identity extends to
every \(z\). The critical points are therefore \(-1\) repeated
\(h-1\) times and the \(r\) eigenvalues of \(B\).
The condition on reciprocal sums is equivalent to $B\succeq0$.
Indeed $B\succeq0$ is equivalent to
$u^TD^{-1}u=(r+\sum_{i=1}^{r} a_i^{-1})/n\le1$.
Thus the absolute root and critical lists are
\[
 X=(a_1,\ldots,a_r,1^{[h]}),\qquad
 Y=(\lambda_1(B),\ldots,\lambda_r(B),1^{[h-1]}).
\]
Set $c=(m+1)/m$, $d=1/m$, and let
\[
 X_*=(c^{[m]},1^{[m+1]}),\qquad
 Y_*=(c^{[m-1]},d,1^{[m]}).
\]
Define the cross-product lists $U=X\otimes Y_*$ and $V=Y\otimes X_*$.
Their first moment difference is the constant
\(\Delta_1=\sum_{v\in V}v-\sum_{u\in U}u=4(h-m-1)\),
independent of $a$, while $\sum_{u\in U}u^2=\sum_{v\in V}v^2$.
For $s>0$ put
\[
 K_s(a)=\sum_{v\in V}\log(1+v/s)-\sum_{u\in U}\log(1+u/s),
 \qquad \widetilde K_s=K_s-\Delta_1/s.
\]
We shall prove that $K_s$, equivalently $\widetilde K_s$, is strictly
Schur-concave on the convex set of admissible $a$.

Write
\[
 R(t)=\sum_{i=1}^{r}(a_i+t)^{-1},\qquad
 \rho(t)=\frac{r+(1-t)R(t)}n.
\]
The determinant lemma gives
$\det(sI+B)=\prod_{i=1}^{r}(s+a_i)(1-\rho(s))$.
After cancelling common factors in the cross products,
\begin{align*}
 K_s(a)={}&\sum_{i=1}^{r}\bigl[\log(s+ca_i)+\log(s+a_i)-\log(s+da_i)\bigr]\\
 &+(h-m)\log(s+c)+(h-m-1)\log(s+1)-h\log(s+d)\\
 &+m\log(1-\rho(s/c))+(m+1)\log(1-\rho(s)).
\end{align*}
For fixed $s$, this formula is a smooth extension to a neighborhood of
any admissible vector, even when $B$ is singular, because
$1-\rho(t)>0$ for $t>0$.

Put $t=s/c$. The coordinate derivative has the form
\[
 \partial_i K_s=g(a_i),\qquad
 g(x)=\frac1{x+t}+\frac1{x+s}-\frac1{x+ms}
       +\frac{A_1}{(x+t)^2}+\frac{A_2}{(x+s)^2},
\]
where, when comparing coordinates, the coefficients are held fixed:
\[
 A_1=\frac{m(1-t)}{n(1-\rho(t))},\qquad
 A_2=\frac{(m+1)(1-s)}{n(1-\rho(s))}.
\]
Let $\eta=h/r>1$. Jensen gives $R(t)\ge r/(t+\eta)$ and hence,
with $\rho_0(t)=1/(t+\eta)$,
\[
 \frac{w(1-t)}{n(1-\rho(t))}
 \ge\frac{w(1-t)(t+\eta)}{n(t+\eta-1)}
 \quad(w>0).
\]
This inequality is valid both when $t\le1$ and when $t\ge1$:
when $1-t$ changes sign, so does the direction of the comparison of
$\rho(t)$ with $\rho_0(t)$.
Consequently $g'(x)\le g_0'(x)$, where $g_0$ uses the uniform
coefficients. The elementary identity
\(\frac{(1-t)(t+\eta)}{t+\eta-1} =-t+\frac{\eta}{t+\eta-1}\)
yields
\[
g_0'(x)=-\frac{x+t/n}{(x+t)^3} -\frac{x-s/n}{(x+s)^3} +\frac1{(x+ms)^2} -\frac{2m\eta}{n(t+\eta-1)(x+t)^3} -\frac{2(m+1)\eta}{n(s+\eta-1)(x+s)^3}.
\]
The last two terms are strictly negative. The negative of the first
line equals $s^{-2}H_m(x/s)$, where, writing $\alpha=m/(m+1)$,
\[
 H_m(y)=\frac{y+\alpha/n}{(y+\alpha)^3}
       +\frac{y-1/n}{(y+1)^3}-\frac1{(y+m)^2}
 =\frac{yP_m(y)}{(m+y)^2n(y+1)^3(my+m+y)^3}.
\]
Here $P_m(y)=\sum_{j=0}^5 b_jy^j$ and
\begin{align*}
 b_5&=(m+1)^3(2m+1),\\
 b_4&=(m+1)^2(8m^3+12m^2+4m-1),\\
 b_3&=m(m+1)^2(4m^3+30m^2+8m-5),\\
 b_2&=2m^2(m+1)(6m^3+24m^2+7m-5),\\
 b_1&=2m^3(6m^3+19m^2+8m-3),\\
 b_0&=4m^5(m+2).
\end{align*}
The script \href{https://github.com/zhangteng2000/schoenberg-validation/blob/e7435c5521dcd6a10131ba04a6da500b1404800b/certificates/num_p1_face_scalar_check.py}{\nolinkurl{num_p1_face_scalar_check.py}} verifies this identity for $H_m$.
Every coefficient is positive for $m\ge1$, so $g'(x)<0$ for $x>0$.
The standard coordinate-gradient criterion therefore proves strict
Schur concavity. Equivalently, along the averaging path
$a_i(\tau)=(1-\tau)a_i+\tau\eta$, admissibility is preserved by
convexity of the reciprocal-sum condition, and $K_s$ strictly
increases unless all coordinates are equal.
Thus
\[
 \widetilde K_s(a)\le\widetilde K_s(\eta\mathbf1_r),\qquad s>0.
\]

The two-value result gives $\widetilde K_s(\eta\mathbf1_r)\le0$.
Consequently $\widetilde K_s(a)\le0$, or equivalently $\mathcal J_s(a)\ge0$,
for every $s>0$. This pointwise assertion is the input used later.

For $2<p<3$ the integral identity
\[
 x^p=\frac{p\sin(\pi p)}\pi
 \int_0^\infty s^{p-1}
 \left[\log(1+x/s)-x/s+\frac{x^2}{2s^2}\right]ds
\]
has positive prefactor. Applying it to the two cross lists, their
second moments cancel, and their first moment difference is exactly
the subtracted term $\Delta_1/s$. Therefore
\begin{equation}\label{eq:flat-majority-power}
 \sum_{u\in U}u^p-\sum_{v\in V}v^p
 \ge \left(\sum_{u\in U}u^p-\sum_{v\in V}v^p\right)_{a_i=\eta}.
\end{equation}
The right side is nonnegative by Theorem~\ref{odd23:thm:two-rational}, \eqref{odd23:eq:J-T}, and \eqref{odd23:eq:complete-Mellin}.
It vanishes only for \(r=m\). Strict Schur concavity makes
\eqref{eq:flat-majority-power} strict for nonuniform \(a\). Dividing by
$\sum_{x\in X_*}x^p>0$ completes the proof, including equality.
\subsection{An exact certificate for a balanced three-value family}
The following pointwise squared-log comparison is needed for the family
\[
 p_\tau(z)=(z+m+\tau)^m(z-2m\tau)(z-m+\tau)^m,
 \qquad |\tau|\le\frac{m}{2m+1},\quad m\ge2.
\]
For this family \(W_p\le B_{n,p}Z_p\) throughout \(2<p<4\).
We expand the squared determinants and verify the resulting
Bernstein-coefficient checks in the supplied script
\href{https://github.com/zhangteng2000/schoenberg-validation/blob/e7435c5521dcd6a10131ba04a6da500b1404800b/certificates/num_balanced_threevalue_certificate.py}{\nolinkurl{num_balanced_threevalue_certificate.py}}.

Put $n=2m+1$, $a=m+1$, $b=m$, and $v=1-(n\tau/m)^2\in[0,1]$.
The two critical points not inherited from the repeated roots, denoted
$\xi_\pm$, satisfy
\[
 \xi_++\xi_-=2(m-1)\tau,\qquad
 \xi_+\xi_-=-\frac{m^2}{n}-(2m-1)\tau^2.
\]
For the current absolute roots $X$ and absolute critical points $Y$, let
$U=X\times |y^*|$ and $V=Y\times |x^*|$ denote product multisets.
Consider the difference
\(E(t)=\prod_{v_1\in V}(t+v_1^2)-\prod_{u_1\in U}(t+u_1^2)\).
After cancelling common positive factors, setting $t=m^4Y_0/n^2$, and
removing another common positive constant, this becomes
\(\mathcal E=A^{m-1}D-B^{m-1}F\).
All the polynomials used here are given explicitly as follows (the variable
$Y_0$ is nonnegative):
\begin{align*}
 S&=(2m^2-2m+1)(1-v)+n,&T&=(2m-1)(1-v)+n,\\
 Q_0&=n^2+1-v,&R_0&=n^2-1+v,\\
 L_A&=m^4Y_0^2+2a^2m^2SY_0+a^4T^2,&
 L_B&=Y_0^2+2SY_0+T^2,\\
 R_1&=m^4Y_0^2+2m^2Q_0Y_0+R_0^2,&
 R_B&=Y_0^2+2Q_0Y_0+R_0^2,\\
 F_a&=Y_0+4a^2(1-v),&F_b&=Y_0+4m^2(1-v),\\
 F_1&=Y_0+4(1-v).&&
\end{align*}
Define
\[
 A=L_A L_B,\quad B=R_1F_aF_b,\quad D=L_A L_B^2,
 \quad F=R_BR_1F_bF_1,\quad H=A-B,
\]
and, for $m\ge3$,
\(J=B^2(D-F)+(m-1)BHD+\frac{(m-1)(m-2)}2H^2D\).
Write $m=r+3$. Expand $H/v$ and $J/v^2$ in powers of $r,Y_0$.
For every coefficient polynomial $\sum_{j=0}^d c_jv^j$, its Bernstein
coefficients on $[0,1]$ are
\[
 \beta_k=\sum_{j=0}^k c_j\frac{\binom{k}{j}}{\binom{d}{j}},
 \qquad 0\le k\le d.
\]
The script checks that the Bernstein coefficients of $H/v$ and $J/v^2$
are nonnegative. There are 85 such coefficients in 32 groups for $H/v$,
and 2520 coefficients in 334 groups for $J/v^2$. Thus $H,J\ge0$ for all
$m\ge3$, $Y_0\ge0$, and $v\in[0,1]$. All the other displayed factors are
nonnegative directly. The binomial expansion gives
\(A^{m-1}D-B^{m-1}F \ge B^{m-3}J\ge0\).
For $m=2$, the script directly checks the 45 Bernstein coefficients in
9 groups of $(AD-BF)/v^2$. Hence $E(t)\ge0$ in every case.

Consequently
\[
 K(t)=\sum_{v_1\in V}\log(1+v_1^2/t)
       -\sum_{u_1\in U}\log(1+u_1^2/t)\ge0.
\]
The squared first moments agree, since the exact real second-moment identity
holds for both configurations. For $r=p/2\in(1,2)$, the convergent Mellin
identity therefore yields
\[
 \sum_{u\in U}u^p-\sum_{v\in V}v^p
 =-\frac{r\sin(\pi r)}{\pi}\int_0^\infty t^{r-1}K(t)\,dt\ge0.
\]
Since $\sum_{u\in U}u^p=(\sum_{x\in X}x^p)(\sum_{y\in Y_*}|y|^p)$ and
$\sum_{v\in V}v^p=(\sum_{y\in Y}y^p)(\sum_{x\in X_*}|x|^p)$, this is the required inequality on
the stated three-value family. No assertion about other multiplicity
patterns is included in this certificate.

\subsection{A bridge from squared determinants to the regularized logarithm}
\label{odd23:sec:log-bridge}
The earlier balanced three-value certificate is stronger than the
power inequality alone. It proves, for the corresponding cross lists,
\[
 K_2(t)=\sum_{v\in V}\log(1+v^2/t)-\sum_{u\in U}\log(1+u^2/t)\ge0,
 \qquad t>0.
\]
This also proves the pointwise regularized logarithmic comparison
used in Theorem~\ref{odd23:lm-thm:main}.

\begin{lem}\label{odd23:lem:log-bridge}
Let $U,V$ be finite nonnegative lists of equal length, with equal
second moments. If $K_2(t)\ge0$ for every $t>0$, then
\[
 \widetilde K_s=
 \sum_{v\in V}\log(1+v/s)-\sum_{u\in U}\log(1+u/s)
 -\frac{\sum_{v\in V}v-\sum_{u\in U}u}{s}\le0,
 \qquad s>0.
\]
\end{lem}
\begin{proof}
For $r\ge0$ set
\(q_s(r)=\log(1+\sqrt r/s)-\sqrt r/s+\frac r{2s^2}\).
Its derivative is
\(q_s'(r)=\frac1{2s^2}-\frac1{2s(s+\sqrt r)}\).
The elementary Stieltjes identity
\(\frac1{s+\sqrt r}=\frac1\pi\int_0^\infty \frac{\sqrt t}{(s^2+t)(r+t)}\,dt\)
follows by setting $t=y^2$ and integrating partial fractions; the
coincident-parameter cases follow by continuity. Integrating from
$0$ to $r$ gives
\(q_s(r)=\frac r{2s^2} -\frac1{2\pi s}\int_0^\infty \frac{\sqrt t}{s^2+t}\log(1+r/t)\,dt\).
The integral representations of $q_s(r)$ converge also at $r=0$.
Apply the identity to $r=v^2,u^2$ and subtract. The second moments
cancel, so
\[
 \widetilde K_s=-\frac1{2\pi s}\int_0^\infty
 \frac{\sqrt t}{s^2+t}K_2(t)\,dt\le0.
\]
\end{proof}

Consequently the entire balanced three-value family satisfies
$\mathcal J_s=-\widetilde K_s\ge0$ pointwise, for every $s>0$.

\subsection{The unbalanced three-value region with nonpositive middle root}

\begin{thm}\label{thm:nonpositive-middle}
Let $n=2m+1$, $1\le k\le m-1$, and $\ell=2m-k$.
The bound \(W_p\le B_{n,p}Z_p\) holds for
\[
 p(z)=(z+1)^\ell(z-a)^k(z-b),\qquad
 a>b\ge0,\quad ka+b=\ell,
\]
for every $2<p<3$. The inequality is strict for every list with the stated multiplicities and ordering.
\end{thm}

This is the whole region with a nonpositive middle root, after
reflection and scaling, of the unbalanced three-value family whose
middle root has multiplicity one. When the residual compression is
positive semidefinite, the result follows from
Section~\ref{odd23:sec:flat-psd}. We prove the remaining region
by a rational averaging inequality. The finite exact certificate
for that inequality is supplied in
\href{https://github.com/zhangteng2000/schoenberg-validation/blob/e7435c5521dcd6a10131ba04a6da500b1404800b/certificates/turn2_restricted_gradient_exact_check.py}{\nolinkurl{turn2_restricted_gradient_exact_check.py}}.

\par\medskip\noindent\textbf{The two nonpersistent critical points.}
Suppose the residual compression is not positive semidefinite.
Its eigenvalues consist of $a$ repeated $k-1$ times and two values
$\lambda>0$ and $-c<0$. Interlacing gives
\[
 0<c<1,\qquad b<\lambda<a,
 \qquad \lambda=a+b+c-1.
\]
The negative eigenvalue obeys
\begin{equation}\label{odd23:tm-eq:secular}
 \frac{k}{a+c}+\frac1{b+c}=\frac\ell{1-c}.
\end{equation}
Put $d=m-k\ge1$. Since $a>\ell/(k+1)>1$,
\eqref{odd23:tm-eq:secular} gives
\begin{equation}\label{odd23:tm-eq:triangle}
 \frac1{b+c}=\frac\ell{1-c}-\frac{k}{a+c}>\ell-k=2d,
 \qquad b+c<\frac1{2d}.
\end{equation}
Vary $b$ upward and set $a=(\ell-b)/k$. Then $a'=-1/k$.
The quadratic equation for the two nonpersistent critical points is
\(n c(a+b+c-1)=kb+a-\ell ab\).
Implicit differentiation gives
\begin{equation}\label{odd23:tm-eq:eta}
 c'=-\eta,\quad \lambda'=\frac{k-1}{k}-\eta,
 \qquad
 \eta=\frac{4md+1-2\ell b+n(k-1)c}
 {n\{2d+(k-1)b+2kc\}}>0.
\end{equation}
Indeed, \eqref{odd23:tm-eq:triangle} gives \(2\ell b<\ell/d\);
since \(\ell=m+d\) and \(m,d\ge1\), the numerator in
\eqref{odd23:tm-eq:eta} is positive.

\par\medskip\noindent\textbf{The exact regularized logarithm.}
Let $c_*=(m+1)/m$, $d_*=1/m$, and form the absolute lists
\[
 X=(a^{[k]},b,1^{[\ell]}),\quad
 Y=(a^{[k-1]},\lambda,c,1^{[\ell-1]}),
\]
\[
 X_*=(c_*^{[m]},1^{[m+1]}),\quad
 Y_*=(c_*^{[m-1]},d_*,1^{[m]}).
\]
Put $U=X\otimes Y_*$, $V=Y\otimes X_*$ and, for $s>0$,
\[
 \widetilde K_s=
 \sum_{v\in V}\log(1+v/s)-\sum_{u\in U}\log(1+u/s)-\frac{\Delta_1}s,
 \quad
 \Delta_1=4(\ell-m-1)+4(m+1)c.
\]
The second cross moments are equal. Direct differentiation of these
finite products, with $t=sm/(m+1)$, yields
\begin{align}
 \frac{d\widetilde K_s}{db}={}&
 \frac{d}{k(a+t)}+\frac{d+1}{k(a+s)}+\frac1{a+ms}
 -\frac{m-1}{b+t}-\frac m{b+s}-\frac1{b+ms}\notag\\
 &+\left(\frac{k-1}{k}-\eta\right)
 \left(\frac m{\lambda+t}+\frac{m+1}{\lambda+s}\right)
 -\eta\left(\frac m{c+t}+\frac{m+1}{c+s}\right)
 +\frac{4(m+1)\eta}s.
 \label{odd23:tm-eq:gradient}
\end{align}
In particular, the regularization term depending on $c$ is included;
omitting it would give a different and false monotonicity claim.

\begin{lem}\label{odd23:tm-lem:rational}
For arbitrary real $k,d\ge1$, set $m=k+d$, $n=2m+1$,
$\ell=k+2d$. If $b,c>0$, $b+c\le1/(2d)$,
$a=(\ell-b)/k$, $\lambda=a+b+c-1$, and $\eta$ is given by
\eqref{odd23:tm-eq:eta}, then the rational expression on the right of
\eqref{odd23:tm-eq:gradient} is strictly positive for every $s>0$.
\end{lem}

\par\medskip\noindent\textbf{Exact certificate for the rational lemma.}
The assertion is stronger than needed: no secular quadratic relation
between $b$ and $c$ is assumed in this lemma. Set
\[
 B=2db,\quad C=2dc,\quad z=2dk s,
 \quad B,C\ge0,\quad B+C\le1,
\]
\[
 A=2d(k+2d)-B,\quad L=4d^2+(k-1)B+kC,
\]
\[
 E=2d(4md+1)-2\ell B+n(k-1)C,\qquad
 F=n\{4d^2+(k-1)B+2kC\}.
\]
Here $2dk(a,b,c,\lambda)=(A,kB,kC,L)$ and $\eta=E/F$.
The scale factors and denominators are positive. In particular,
$E\ge2m(4d^2-1)>0$.
After dividing \eqref{odd23:tm-eq:gradient} by $2dk$ and multiplying by
$kF$, it becomes $\sum_{i=1}^{11}w_i/P_i(z)$, with the following
explicit numerator--denominator pairs:
\[
\begin{array}{c|c}
w_i&P_i(z)\\\hline
 dF(m+1)&(m+1)A+mz\\
 (d+1)F&A+z\\
 kF&A+mz\\
 -k(m-1)F(m+1)&(m+1)kB+mz\\
 -kmF&kB+z\\
 -kF&kB+mz\\
 m(m+1)((k-1)F-kE)&(m+1)L+mz\\
 (m+1)((k-1)F-kE)&L+z\\
 -km(m+1)E&(m+1)kC+mz\\
 -k(m+1)E&kC+z\\
 4k(m+1)E&z
\end{array}
\]
Let $D(z)=\prod_{i=1}^{11}P_i(z)$ and
$N(z)=\sum_{i=1}^{11}w_i\prod_{\substack{1\le j\le11\\j\ne i}}P_j(z)$.
Exact multiplication gives $\deg_zN=8$.
Its constant coefficient is
\(N_0=4k^6(m+1)^5E A^3B^3C^2L^2>0\).
For $1\le j\le8$, remove the visibly nonnegative factor
\(A^{\max(3-j,0)}B^{\max(3-j,0)} C^{\max(2-j,0)}L^{\max(2-j,0)}\)
from $N_j$. The remaining polynomial has a triangular Bernstein
expansion on $B,C\ge0$, $B+C\le1$. Every Bernstein coefficient,
a polynomial in $k,d$, has nonnegative ordinary power coefficients
after $k=1+K$, $d=1+H$. The exact check produces:
\[
\begin{array}{c|rrrrrrrr}
j&1&2&3&4&5&6&7&8\\\hline
\deg_{B,C}&4&7&8&7&6&5&4&3\\
\text{checked scalar coefficients}
&2223&7142&10182&7858&5669&3630&2107&1108\\
\text{negative coefficients}&0&0&0&0&0&0&0&0
\end{array}
\]
The supplied script constructs $N$ from the eleven displayed pairs,
checks the constant coefficient identity, performs exact polynomial
division, and checks the signs of the coefficients using rational arithmetic.
For reproducibility, if $P(B,C)=\sum p_{ij}B^iC^j$ has total degree
$q$, its Bernstein coefficient indexed by $(r,t)$ is
\[
 \sum_{i\le r,\,j\le t}
 p_{ij}\,\frac{(r)_i(t)_j}{(q)_{i+j}},
 \qquad r,t\ge0,\quad r+t\le q,
\]
where $(x)_j$ is the falling factorial. Thus the certificate is a
finite exact identity with explicitly specified positive bases.
It proves $N(z)>0$ and hence Lemma~\ref{odd23:tm-lem:rational}.

\par\medskip\noindent\textbf{Finishing the power inequality.}
The function $k/a+1/b$ strictly decreases when $b$ increases toward
$\ell/(k+1)$, because its derivative is $a^{-2}-b^{-2}<0$.
There is therefore a unique threshold $b_*$ at which $c=0$ and the
residual block becomes positive semidefinite. Along $b<b_*$,
Lemma~\ref{odd23:tm-lem:rational} gives
\[
 \widetilde K_s(b)<\widetilde K_s(b_*),\qquad s>0.
\]
At the threshold, the pointwise bound in the proof of
Theorem~\ref{thm:flat-majority} gives \(\widetilde K_s(b_*)\le0\).
Consequently \(\widetilde K_s(b)<0\) whenever \(b<b_*\) and \(s>0\).

For $2<p<3$, the positive Mellin representation of $x^p$ by
$\log(1+x/s)-x/s+x^2/(2s^2)$ gives
\[
 \sum_{u\in U(b)}u^p-\sum_{v\in V(b)}v^p
 >\sum_{u\in U(b_*)}u^p-\sum_{v\in V(b_*)}v^p\ge0.
\]
The final inequality is the already proved flat-majority PSD theorem.
The case $b=0$ follows by continuity, or by integrating the strict
monotonicity on any interval $(0,b_*)$. This proves the theorem.

\subsection{The unbalanced three-value region with nonnegative middle root}
\label{odd23:sec:three-positive}

Here we prove the rational trace comparison for a centered real
polynomial with three distinct root values. The outer values have
multiplicities \(k,\ell\), the nonnegative middle value is simple,
and \(1\le k<m\), \(\ell=2m-k\). Reflection and positive scaling
give the normalization used below.

\subsubsection{The rational functional}
Put $n=2m+1$, $a=m+1$, $b=m$. Let $X_*$ be the absolute root list
$(a^{[m]},b^{[m+1]})$ and let $Y_*$ be the absolute critical-point list
$(a^{[m-1]},b^{[m]},1)$. For any centered real root list $X$ and its
critical-point list $Y$, all absolute values being taken with multiplicity,
write $U=|X|\otimes Y_*$ and $V=|Y|\otimes X_*$ and define
\[
 \mathcal T_X(s)=\sum_{v\in V}\frac{v^2}{s+v}
          -\sum_{u\in U}\frac{u^2}{s+u},\qquad s>0.
\]
The sums of squares of $U$ and $V$ are equal, since
$\sum Y^2=(n-2)\sum X^2/n$ for centered real roots, and the same identity
holds for the reference lists.

\begin{thm}\label{odd23:tp-thm:main}
Let $k,d\ge1$ be integers, put
\(m=k+d,\quad n=2m+1,\quad \ell=k+2d\),
and let $0\le r\le1$. Define
\[
 A=n\ell+kr,\qquad B=2mkr,\qquad C=k(n-r).
\]
For the centered root list
\(X=(-A)^{[k]},\ B,\ C^{[\ell]}\),
one has $\mathcal T_X(s)>0$ for every $s>0$.
Consequently the sharp inequality
\[
 \sum_{y\in Y}|y|^p\le
 \frac{(m-1)(m+1)^p+m^{p+1}+1}
 {m(m+1)^p+(m+1)m^p}\sum_{x\in X}|x|^p
\]
holds strictly for every $2<p<3$.
\end{thm}
Every unbalanced three-value configuration with simple nonnegative middle
value is represented here up to positive scaling and reflection:
$-kA+B+\ell C=0$, $A>C\ge B\ge0$, and the parameter $r$ ranges from
$B=0$ to $B=C$.

\subsubsection{Reducing the square root to two rational bounds}
The inherited critical points are $-A$ of multiplicity $k-1$ and $C$ of
multiplicity $\ell-1$. The two remaining critical points are $-u,v$ with
$u,v>0$. Put
\[
 S=B+C-A,\quad M=\frac{AC+\ell AB-kBC}{n},\quad
 D=S^2+4M,\quad V_2=S^2+2M.
\]
Then
\[
 v-u=S,\qquad uv=M>0,\qquad u+v=L=\sqrt D,\qquad u^2+v^2=V_2.
\]
In particular, for any $c>0$,
\begin{equation}\label{odd23:tp-eq:pair}
 \frac{(cu)^2}{s+cu}+\frac{(cv)^2}{s+cv}
 =\frac{c^2(sV_2+cML)}{s^2+cLs+c^2M}.
\end{equation}
The special expressions
\[
 L_0=B+A-C=2nd+2k(m+1)r>0,\qquad
 H_0=2mn-B,\qquad E=\ell+kr
\]
satisfy
\(\delta:=D-L_0^2=4k(1-r)nE=2(H_0-L_0)E\ge0\).
Hence
\begin{equation}\label{odd23:tp-eq:bounds}
 L_-:=L_0+\frac{2L_0\delta}{4L_0^2+\delta}
 \ \le L\le\
 L_+:=L_0+\frac{\delta}{2L_0}.
\end{equation}
Indeed, concavity of the square root gives the upper bound, and
$L-L_0=\delta/(L+L_0)\ge\delta/(L_++L_0)$ gives the lower bound.
All denominators in \eqref{odd23:tp-eq:bounds} are strictly positive.

\subsubsection{The numerator is affine in the square root}
After canceling persistent multiplicities, the non-pair contributions to
$\mathcal T_X$ are $\sum_{i=1}^9\omega_i z_i^2/(s+z_i)$, with the following table:
\[
\begin{array}{c|rrrrrrrrr}
z_i&aA&bA&A&aC&bC&C&aB&bB&B\\ \hline
\omega_i&-d&-(d+1)&-k&d&d-1&-\ell&-(m-1)&-m&-1.
\end{array}
\]
Define the polynomials
\[
 D_r(s)=\prod_{i=1}^9(s+z_i),\qquad
 N_r(s)=\sum_{i=1}^9\omega_i z_i^2\prod_{j\ne i}(s+z_j).
\]
Keeping a factor whose weight is zero causes no problem.
Using \eqref{odd23:tp-eq:pair}, the two remaining pair contributions are
\(m\frac{a^2(sV_2+aML)}{s^2+aLs+a^2M} +(m+1)\frac{b^2(sV_2+bML)}{s^2+bLs+b^2M}\).
Their common denominator, after replacing $L^2$ by $D$, is
\[
D_c(s;L)=s^4+(a+b)Ls^3+\bigl((a^2+b^2)M+abD\bigr)s^2 +ab(a+b)MLs+a^2b^2M^2.
\]
The corresponding numerator is
\begin{align*}
 N_c(s;L)={}&\bigl(ma^2+(m+1)b^2\bigr)V_2s^3\\
 &+L\bigl[ma^2(bV_2+aM)+(m+1)b^2(aV_2+bM)\bigr]s^2\\
 &+\bigl[ma^2(b^2MV_2+abMD)+(m+1)b^2(a^2MV_2+abMD)\bigr]s\\
 &+a^2b^2M^2\bigl(ma+(m+1)b\bigr)L.
\end{align*}
Thus
\[
 \mathcal T_X(s)=\frac{N(s;L)}{D_r(s)D_c(s;L)},\qquad
 N(s;L)=N_r(s)D_c(s;L)+D_r(s)N_c(s;L).
\]
For the actual value $L=\sqrt D$, the denominator is strictly positive for
$s>0$. The polynomial $N$ has degree $11$ in $s$ and is affine in $L$.
Write
\[
 N(s;L)=\sum_{j=0}^{11}c_j(L)s^j,\qquad
 P_j=c_j(L_0),\quad Q_j=c_j(H_0).
\]
Since \(c_j\) is affine, substituting \(L_+\) and \(L_-\) from
\eqref{odd23:tp-eq:bounds} gives
\begin{align}
 L_0c_j(L_+)&=(L_0-E)P_j+EQ_j=:U_j,\label{odd23:tp-eq:uppercoeff}\\
 (4L_0^2+\delta)c_j(L_-)
 &=(4L_0^2+\delta-4L_0E)P_j+4L_0EQ_j=:W_j.
 \label{odd23:tp-eq:lowercoeff}
\end{align}
Each $U_j,W_j$ is a polynomial in $k,d,r$ with integer coefficients.

\subsubsection{The finite exact positivity certificate}
\begin{lem}\label{odd23:tp-lem:cert}
For all $j=0,\ldots,11$, the polynomials $U_j$ and $W_j$ in
\eqref{odd23:tp-eq:uppercoeff}--\eqref{odd23:tp-eq:lowercoeff} are nonnegative whenever
$k,d\ge1$ and $0\le r\le1$.
\end{lem}
\begin{proof}
Set $k=K+1$, $d=H+1$. For any \(U_j\) or \(W_j\), let $R$ be its degree
in $r$, and expand it in the basis
\[
 K^\alpha H^\beta r^t(1-r)^{R-t},\qquad
 \alpha,\beta\ge0,\quad 0\le t\le R.
\]
All coefficients are nonnegative integers. This is a finite polynomial
identity, verified by the accompanying standalone script
\href{https://github.com/zhangteng2000/schoenberg-validation/blob/e7435c5521dcd6a10131ba04a6da500b1404800b/certificates/turn2_three_caseIII_T_certificate.py}{\nolinkurl{turn2_three_caseIII_T_certificate.py}}.
The script constructs \(U_j,W_j\) from
\eqref{odd23:tp-eq:uppercoeff} and \eqref{odd23:tp-eq:lowercoeff} in exact arithmetic.
The number of strictly positive coefficients in each identity is:
\[
\begin{array}{c|rrrrrrrrrrrr}
j&0&1&2&3&4&5&6&7&8&9&10&11\\ \hline
U_j&11512&10530&9477&8378&6667&5177&3906&2845&1980&1305&794&436\\
W_j&13640&12526&11336&10094&8137&6419&4938&3685&2647&1816&1168&694
\end{array}
\]
There are $140107$ nonzero coefficients in total, and all are positive.
For $j=11$, the coefficients with $K=H=0$ alone are already strictly
positive at every Bernstein index. They are, respectively,
\((288900,1548000,2889600,2295040,667648)\)
and
\[
 (14904000,94046400,235238400,290375680,176926720,42729472).
\]
Consequently both $U_{11}$ and $W_{11}$ are strictly positive throughout
the entire domain, including its boundaries.

To make the verification rule explicit, if a polynomial after shifting is
$\sum_{i=0}^R f_i(K,H)r^i$, its coefficient of
$r^t(1-r)^{R-t}$ is
\(\sum_{i=0}^t\binom{R-i}{t-i}f_i(K,H)\).
The script computes exactly these integer binomial transforms and checks
every coefficient of every resulting polynomial in $K,H$.
Every basis function is nonnegative on the required domain, proving the
lemma.
\end{proof}

By Lemma~\ref{odd23:tp-lem:cert} and the positive normalizing factors in
\eqref{odd23:tp-eq:uppercoeff}--\eqref{odd23:tp-eq:lowercoeff},
$c_j(L_-),c_j(L_+)\ge0$. Since $c_j$ is affine and
$L_-\le L\le L_+$, also $c_j(L)\ge0$. The leading coefficient is strictly positive by the last part of the
certificate, so
$N(s;L)>0$ for $s>0$, proving the rational inequality in
Theorem~\ref{odd23:tp-thm:main}. This argument includes the leading coefficient;
it does not invoke the already known sharp cubic inequality.

\subsubsection{From the rational inequality to the power inequality}
Let $\Delta_1=\sum_{v\in V}v-\sum_{u\in U}u$ and put
\[
 \mathcal K_X(s)=\sum_{v\in V}\log(1+v/s)
 -\sum_{u\in U}\log(1+u/s)-\frac{\Delta_1}{s}.
\]
Direct differentiation yields
\[
 \mathcal K_X'(s)=\frac{\mathcal T_X(s)}{s^2}\ge0,\qquad
 \lim_{s\to\infty}\mathcal K_X(s)=0.
\]
Consequently $\mathcal K_X(s)<0$ for all $s>0$.
For $2<p<3$,
\[
 t^p=\frac{p\sin(\pi p)}{\pi}\int_0^\infty
 s^{p-1}\left[\log(1+t/s)-\frac ts+\frac{t^2}{2s^2}\right]ds.
\]
The constant is positive. The quadratic terms cancel between $U$ and $V$,
so
\[
 \sum_{u\in U}u^p-\sum_{v\in V}v^p
 =-\frac{p\sin(\pi p)}{\pi}\int_0^\infty s^{p-1}\mathcal K_X(s)\,ds\ge0.
\]
This is precisely the claimed sharp power inequality after dividing by
$\sum_{x\in X_*}x^p=m(m+1)^p+(m+1)m^p$.

\subsection{The complete three-value comparison}
\begin{thm}\label{odd23:thm:all-three-log}
For every odd $n\ge5$, one has $\mathcal J_s\ge0$ for all $s>0$
on every centered real list with at most two distinct values, and on
every centered three-value list whose middle value is simple.
\end{thm}
\begin{proof}
A constant centered list consists of zeros and gives \(\mathcal J_s=0\).
For two distinct values use Theorem~\ref{odd23:thm:two-rational}
and \eqref{odd23:eq:J-T}. For three values reflect the zeros so that
the outer multiplicities satisfy \(k\le\ell\). If \(k=\ell\),
Lemma~\ref{odd23:lem:log-bridge} converts the squared-log comparison
to the required inequality. For \(k<\ell\), the pointwise bounds in
the proof of Theorem~\ref{thm:nonpositive-middle} and in
Theorem~\ref{odd23:tp-thm:main} with \eqref{odd23:eq:J-T}
cover the two signs of the middle zero.
\end{proof}
\subsection{The first moment and exclusion of zero roots}
\label{odd23:sec:nozeros}

Let $n=2m+1\ge5$, $a=m+1$, $b=m$, and let $x\in\mathbb R^n$
be centered. Its critical points are the eigenvalues of
$A_x=Q\operatorname{diag}(x)Q|_{e^\perp}$, where
$e=n^{-\frac{1}{2}}\mathbf1$ and $Q=I-ee^T$.
For $\tau>0$ put $f_\tau(t)=t^2/(\tau+|t|)$ and write
\[
 Z_\tau(x)=\sum_{j=1}^{n}f_\tau(x_j),\qquad
 W_\tau(x)=\tr f_\tau(A_x).
\]
The auxiliary rational functional is
\begin{equation}\label{odd23:rz-eq:T}
 T_s(x)=ab(W_{s/a}+W_{s/b})
 -(ab-a)Z_{s/a}-(ab-b)Z_{s/b}-Z_s.
\end{equation}
\begin{lem}\label{odd23:rz-lem:p1}
For every centered real list of odd length $n=2m+1$,
\[
 \sum_{k=1}^{n-1}|y_k|\ge\frac{m}{m+1}\sum_{j=1}^{n}|x_j|.
\]
If a nonzero list contains a zero coordinate, this inequality is
strict.
\end{lem}
\begin{proof}
Let $r,s,t$ be the positive, negative, and zero counts, and let $S>0$
be the sum of the positive roots. Averaging within the positive and
negative groups, convexity of the trace norm and permutation
invariance show that the critical trace norm cannot increase when
the roots are replaced by $S/r$ repeated $r$ times, $-S/s$ repeated
$s$ times, and $0$ repeated $t$ times. The root trace norm remains
$2S$.

If $t=0$, the averaged list has critical trace norm
$2S-2\min(S/r,S/s)$, and therefore ratio
$1-1/\max(r,s)\ge m/(m+1)$.

If $t\ge1$, put $u=S/r$, $v=S/s$, and $k=r+s=n-t$. Besides the
repeated roots, the two remaining critical points solve
\(z^2-(u-v)z-\frac{tuv}{n}=0\).
Their absolute-value sum is
$\sqrt{(u-v)^2+4tuv/n}$. Consequently the averaged ratio equals
\(1-\frac{2/n}{1+\sqrt{1-4rs/(nk)}}\).
Since $rs\le k^2/4$, it is at least
\[
 1-\frac{2}{n+\sqrt{nt}}
 \ge1-\frac{2}{n+\sqrt n}
 >1-\frac{2}{n+1}=\frac{m}{m+1}.
\]
This proves both assertions.
\end{proof}

\begin{lem}\label{odd23:rz-lem:concave-derivative}
Suppose $k:[0,\infty)\to[0,\infty)$ is increasing and concave,
with $k(0)=0$. Let $g$ be the odd function satisfying $g'=k(|t|)$.
Then for all real $u,v$,
\[
 g[u,v]\ge\frac{k(|u|)+k(|v|)}4,
\]
where diagonal divided differences have their usual derivative
interpretation.
\end{lem}
\begin{proof}
If $u,v$ have the same sign, concavity gives the stronger bound with
denominator $2$. If $u\ge0\ge v$, put $w=-v\ge0$. Concavity and
$k(0)=0$ give $g(u)\ge uk(u)/2$ and $g(w)\ge wk(w)/2$.
Monotonicity of $k$ then gives
\(\frac{g(u)+g(w)}{u+w} \ge\frac{uk(u)+wk(w)}{2(u+w)} \ge\frac{k(u)+k(w)}4\).
Zero and diagonal cases follow by continuity.
\end{proof}

\begin{prop}\label{odd23:rz-prop:nozeros}
For $n\ge5$ odd, no nonzero centered real root list containing a
zero coordinate can be a local minimum of $T_s$, for any $s>0$.
\end{prop}
\begin{proof}
Define the even function
\[
 q_\tau(t)=\frac{t^2}{\tau}-f_\tau(t)
 =\frac{|t|^3}{\tau(\tau+|t|)},\qquad
 k_\tau(t)=q_\tau''(t)
 =\frac2\tau-\frac{2\tau^2}{(\tau+|t|)^3}.
\]
As a function of $|t|$, $k_\tau$ is nonnegative, increasing,
concave, and vanishes at zero. The square-moment compression identity
$\sum_{k=1}^{n-1}y_k^2=(n-2)\sum_{j=1}^{n}x_j^2/n$ implies exact cancellation of the quadratic terms
in \eqref{odd23:rz-eq:T}. Thus the functional formed by replacing every
$f_\tau$ there with $q_\tau$ is exactly $-T_s$.

Suppose $x_1=0$ and consider the centered variation
$x_1(t)=-(n-1)t$, $x_j(t)=x_j+t$ for $j>1$.
Put $v=Qe_1$, so the matrix direction is $H=I-nvv^T$ on
$e^\perp$ and $\|v\|^2=(n-1)/n$. Diagonalize $A_x$, with
eigenvalues $\lambda_i$, and put $t_i=|v_i|^2$. By the trace Hessian
formula and Lemma~\ref{odd23:rz-lem:concave-derivative},
\begin{align*}
 \left.\frac{d^2}{dt^2}\tr q_\tau(A_x+tH)\right|_{t=0}
 &=\sum_{i=1}^{n-1}k_\tau(\lambda_i)
   -2n\sum_{i=1}^{n-1}k_\tau(\lambda_i)t_i
   +n^2\sum_{i=1}^{n-1}\sum_{j=1}^{n-1}q_\tau'[\lambda_i,\lambda_j]t_it_j\\
 &\ge\tr k_\tau(A_x)
   +\frac{n(n-5)}2\langle v,k_\tau(A_x)v\rangle\\
 &\ge\tr k_\tau(A_x).
\end{align*}
On the root side, $k_\tau(0)=0$ gives exactly
\(\left.\frac{d^2}{dt^2}\sum_{j=1}^{n}q_\tau(x_j(t))\right|_{t=0} =\sum_{j=1}^{n}k_\tau(x_j)\).

The function $z\mapsto k_\tau(\sqrt z)$ is increasing, concave,
and vanishes at zero. Applying \eqref{odd23:eq:reverse-concave}
to this function gives
\(\tr k_\tau(A_x)\ge\frac{n-2}{n}\sum_{j=1}^{n}k_\tau(x_j)\).
With $\alpha=s/a$ and $\beta=s/b$, the second derivative of $-T_s$
along the variation is consequently at least
\begin{equation}\label{odd23:rz-eq:scalar-lower}
 \sum_j\left[\frac an k_\alpha(x_j)
              -\frac bn k_\beta(x_j)-k_s(x_j)\right].
\end{equation}
Every summand is nonnegative, and is strictly positive when
$x_j\ne0$. To see this, put $z=|x_j|/s$. The summand is
\(\frac{2}{ns}\left[ -\frac{a^2}{(1+az)^3}+\frac{b^2}{(1+bz)^3} +\frac{n}{(1+z)^3}\right]\).
For $z>0$ and $c\in[b,a]$, since $b\ge2$,
\(\frac{d}{dc}\frac{c^2}{(1+cz)^3} =\frac{2c}{(1+cz)^3}-\frac{3zc^2}{(1+cz)^4} <\frac{2c}{(1+z)^3}\).
Integrating from $b$ to $a$ and using $a^2-b^2=n$ proves strict
positivity. At $z=0$ the summand vanishes.
Since the root list is nonzero, \eqref{odd23:rz-eq:scalar-lower} is strictly
positive. Thus $T_s''<0$ in this centered direction, contradicting
local minimality.
\end{proof}

\subsection{The reflected logarithmic kernel estimate}
\label{odd23:sec:reflected-kernel}

This is a kernel estimate, not a proof of global nonnegativity of
the logarithmic trace functional.

For $t>0$ set
\[
 G_t(x)=\operatorname{sgn}(x)\left(\frac1t-\frac1{t+|x|}\right),\qquad
 D_t(x,\xi)=\frac{G_t(\xi)-G_t(x)-G_t'(x)(\xi-x)}{(\xi-x)^2}.
\]
Fix $1\le k\le3/2$, and on $x>0$ define
\[
 h(x)=\frac1{(x+1)^2}+\frac{k}{(x+k)^2},\qquad
 D(x,\xi)=D_1(x,\xi)+kD_k(x,\xi),\qquad
 \mathcal L F=F''-\frac{h''}{h'}F'.
\]
Write
\[
 B_k(x)=(x+k)^3+k(x+1)^3,\qquad
 K_k(x)=\frac{6k(k-1)^2}{(x+1)(x+k)B_k(x)}.
\]

\begin{lem}
For every $v\ge0$ and $x>0$,
\[
 \mathcal L D(x,v)=-\frac{K_k(x)}{(1+v)(k+v)},
\]
and
\[
 \mathcal L D(x,-v)\ge
 -\left(1+\frac{(k^2-1)^2}{8k^2}\right)
 \frac{K_k(x)}{(1+v)(k+v)}.
\]
The factor in parentheses is at most $313/288$ on $1\le k\le3/2$.
\end{lem}

\begin{proof}
For $v\ge0$ the Bregman divided difference is
$D_t(x,v)=-[(x+t)^2(v+t)]^{-1}$, which gives the first identity by
direct differentiation. For the reflected difference the exact
calculation is
\[
 \mathcal L\{D(x,-v)-D(x,v)\}
 =\frac{12vP_k(x,v)}
 {(k+v)(k+x)(v+1)(x+1)B_k(x)(x+v)^4},
\]
where
\[
 P_k(x,v)=-(k-1)^2(k+1)x^4+vP_{1,k}(x)+v^2P_{2,k}(x),
\]
\begin{align*}
 P_{1,k}(x)={}&4kx^4+(8k^2+8k)x^3+(12k^3+12k)x^2\\
 &+(4k^4+4k^3+4k^2+4k)x+k^5+2k^3+k,\\
 P_{2,k}(x)={}&(2k+2)x^4+16kx^3+(12k^2+12k)x^2\\
 &+(8k^3+8k)x+k^4+k^3+k^2+k.
\end{align*}
All their coefficients are nonnegative, and
$P_{1,k}(x)\ge4kx^4$. Set $A=(k-1)^2(k+1)$. If $P_k(x,v)\ge0$
the claimed bound is immediate. Otherwise
\[
 -P_k(x,v)\le x^4(A-4kv),\qquad
 \frac{v[-P_k(x,v)]}{(x+v)^4}
 \le v(A-4kv)\le\frac{A^2}{16k}.
\]
Dividing the resulting lower bound by
$K_k(x)/[(1+v)(k+v)]$ shows that its additional relative negative
part is at most
\(\frac{A^2}{8k^2(k-1)^2} =\frac{(k^2-1)^2}{8k^2}\).
For $k=1$ the formula follows directly from the
nonnegative polynomial, or by continuity. Finally,
$(k-k^{-1})^2/8$ increases on $[1,3/2]$ and equals $25/288$
at $k=3/2$.
\end{proof}

\par\medskip\noindent\textbf{Scaling.}
For the degree-$2m+1$ logarithmic complement, the scales have ratio
$k=(m+1)/m$ and the weights are $m,m+1=m k$.
A common scale change and multiplication by $m$ preserve the lemma
with the same dimensionless factor. This gives the bound required
for actual critical residues, once an independent estimate controls
their weighted sum. The required compensation estimate is proved in the next section.
\subsection{Existence and support of a hypothetical negative logarithmic minimum}

Let $n=2m+1\ge5$, $a=m+1$, $b=m$, and let $x\in\mathbb R^n$
be centered. Its critical-point compression is $A_x$.
For $t>0$ put
\(\ell_t(u)=\frac{|u|}{t}-\log(1+|u|/t)\).
At a fixed $s>0$, write $\alpha=s/a$, $\beta=s/b$, and define
\[
 \mathcal J_s(x)=b\operatorname{Tr}\ell_\alpha(A_x)
 +a\operatorname{Tr}\ell_\beta(A_x)
 -(b-1)\sum_{i=1}^{n}\ell_\alpha(x_i)-b\sum_{i=1}^{n}\ell_\beta(x_i)
 -\sum_{i=1}^{n}\ell_s(x_i).
\]
This is exactly $-\widetilde K_s$ for the unscaled candidate
root list $(a^{[b]},(-b)^{[a]})$.

\begin{thm}\label{odd23:lm-thm:main}
If $\mathcal J_s$ is negative somewhere, its global negative minimum
on the centered subspace is attained. Every such minimizer has no
zero roots, at most two distinct positive root values, and at most
two distinct negative root values. If either sign has two distinct
values, its value closer to zero has multiplicity one.

Every such minimizer also satisfies the rational scale-stationarity
identity $T_s(x)=0$. Moreover, if the smallest root magnitude $r$ satisfies
$r\le s/(m+1)$, there are at most three distinct root values in total.
\end{thm}

Theorem~\ref{odd23:thm:all-three-log} excludes up to three values at a
negative minimum. Theorem~\ref{odd23:mp-thm:onescale} and
Corollary~\ref{odd23:nc-cor:negative-center-one-scale} exclude the four-value case.

\par\medskip\noindent\textbf{Attainment and absence of zero roots.}
For a nonzero centered real list the reverse first-moment bound is
\(W_1\ge mZ_1/(m+1)\). Equality requires a sign group of \(m+1\)
equal roots and a positive semidefinite complementary compression
after reflection if necessary. A nonzero list with a zero root
gives strict inequality.
To see the equality condition, let \(I\) index a negative sign
group of size \(r\ge m+1\) and total magnitude \(S\).
The subspace \(K\) of vectors supported on \(I\) whose entries sum
to zero satisfies
$-2\operatorname{Tr}(P_KA_x)=2(r-1)S/r$.
The variational characterization of the negative part gives
$W_1\ge2(r-1)S/r\ge2mS/(m+1)$. Equality forces $r=m+1$ and
$K$ to reduce $A_x$, with nonnegative complementary block. For
$v\in K$, the coordinates of $A_xv$ outside the group are
$-\sum_{i\in I}x_iv_i/n$. Reduction thus forces the majority roots
to be constant. The zero-root exclusion for equality follows from
the strict reverse-first-moment lemma in
Section~\ref{odd23:sec:nozeros}.

Form the absolute tensor lists $U=|x|\otimes|w_*|$ and
$V=|w|\otimes|x_*|$. Then
\[
 \mathcal J_s=\frac{\Delta_1}{s}
 -\sum_{v\in V}\log(1+v/s)+\sum_{u\in U}\log(1+u/s),\qquad
 \Delta_1=\sum_{v\in V}v-\sum_{u\in U}u\ge0.
\]
Consider $x_j=L_j\xi_j$ with $L_j\to\infty$, $\|\xi_j\|_2=1$,
and $\xi_j\to\xi$. If $\Delta_1(\xi)>0$, the linear term dominates
the $O(\log L_j)$ logarithmic terms. Otherwise $\xi$ is a first-moment
equality case, hence has no zero roots. Every $U(\xi)$ is then nonzero.
After cancelling the equal numbers of $\log L_j$ terms, if some
$V(\xi)$ is zero, the logarithmic difference tends to $+\infty$.
If all are nonzero, its limit is
$\log(\prod_{u\in U(\xi)}u/\prod_{v\in V(\xi)}v)\ge0$.
Indeed, \(\Delta_1=0\) on the limiting flat-majority configuration.
Letting the logarithmic parameter decrease to zero in the pointwise
comparison proved in Subsection~\ref{odd23:sec:flat-psd} gives the
product bound. Therefore
$\liminf_{\|x\|\to\infty}\mathcal J_s(x)\ge0$.
Continuity and $\mathcal J_s(0)=0$ give attainment of any negative
infimum.

For the rational comparison \(T_t\) in
Subsection~\ref{odd23:sec:nozeros}, \eqref{odd23:eq:J-T} gives
\(\mathcal J_s(x)=\int_s^\infty T_t(x)t^{-2}\,dt\).
If \(x_1=0\), vary the roots in the centered direction
\((-(n-1),1,\ldots,1)\). The proof of
Proposition~\ref{odd23:rz-prop:nozeros} gives \(T_t''<0\) in this
same direction for every \(t>0\). To differentiate under the integral,
write \(f_t(u)=u^2/(t+|u|)\). Its second derivative is
\(2t^2/(t+|u|)^3\le2/t\), including at zero. The trace Hessian formula
therefore bounds \(|(\tr f_t(A+\varepsilon H))''|\) by
\(2\,\tr(H^2)/t\) for Hermitian \(A,H\). The scalar terms have the
same bound with their squared directional speeds. The finitely many
fixed dilations in \(T_t\) only change the constant. Thus the second
derivative of \(T_t/t^2\) is bounded in absolute value by \(C/t^3\)
on \([s,\infty)\). It follows that \(\mathcal J_s''<0\), excluding
a zero coordinate at a negative minimum. Finally,
\(\left.\frac{d}{dL}\mathcal J_s(Lx)\right|_{L=1}=\frac{T_s(x)}s\).
Scaling is a centered admissible variation, so $T_s(x)=0$ at a minimum.

\par\medskip\noindent\textbf{The Euler function.}
Write
\[
 g(z)=\frac1n\sum_{i=1}^{n}\frac1{z-x_i},\qquad
 H(z)=z-\frac1{g(z)}=\sum_\nu\frac{\rho_\nu}{z-\xi_\nu},
 \quad\rho_\nu>0.
\]
At an actual root $x_0$ of multiplicity $q$,
\begin{equation}\label{odd23:lm-eq:rootidentity}
 \sum_\nu\frac{\rho_\nu}{(x_0-\xi_\nu)^2}=\frac nq-1\le n-1.
\end{equation}
Let $\mathcal F=b\ell_\alpha+a\ell_\beta$ and
$\mathcal D=(b-1)\ell_\alpha+b\ell_\beta+\ell_s$.
On $x>0$ set
\[
 h=\mathcal F''=\frac b{(x+\alpha)^2}+\frac a{(x+\beta)^2},
 \quad G_t(x)=\ell_t'(x),
\]
\[
 D(x,\xi)=
 \frac{\mathcal F'(\xi)-\mathcal F'(x)-h(x)(\xi-x)}{(\xi-x)^2},
\]
where \(D(x,x)=\mathcal F'''(x)/2\) for \(x\ne0\). Define
\[
 e(x)=aG_\alpha(x)+bG_\beta(x)-nG_s(x)-xh(x),
 \qquad E(x)=e(x)+\sum_\nu\rho_\nu D(x,\xi_\nu).
\]
At \(x<0\) use the same definitions of \(D,e,E\), with
\(h(x)=\mathcal F''(x)\) and \(G_t(x)=\ell_t'(x)\).
Write \(a_i\) for the distinct root values and \(r_i\) for their
multiplicities. The critical points consist of \(r_i-1\) copies of
each \(a_i\) and the simple zeros \(\xi_\nu\) of
\(\sum_i r_i/(z-a_i)\), one between each pair of consecutive values.
Put \(A_\nu=\sum_i r_i/(\xi_\nu-a_i)^2\). The residue formula for
\(H=z-1/g\) and implicit differentiation give
\[
  \rho_\nu=\frac n{A_\nu},\qquad
  \frac{\partial\xi_\nu}{\partial a_i}
       =\frac{r_i}{A_\nu(a_i-\xi_\nu)^2}.
\]
Move all \(r_i\) copies of \(a_i\) together. Differentiating the
critical-point sum and the root sum yields
\[
  \frac n{r_i}\frac{\partial\mathcal J_s}{\partial a_i}
  = n(1-1/r_i)\mathcal F'(a_i)-n\mathcal D'(a_i)
    +\sum_\nu\frac{\rho_\nu\mathcal F'(\xi_\nu)}
                        {(a_i-\xi_\nu)^2}.
\]
This is \(E(a_i)\). Indeed, insert \(H(a_i)=a_i\) and
\eqref{odd23:lm-eq:rootidentity} into the definition of \(D\), and
use \((n-1)\mathcal F'-n\mathcal D'=aG_\alpha+bG_\beta-nG_s\).
Individual coordinate derivatives exist even at repeated roots. In fact,
along a Hermitian line \(A+tV\), the derivative of
\(\tr\mathcal F(A+tV)\) is \(\tr(\mathcal F'(A+tV)V)\).
This follows first for polynomials by cyclicity of the trace and
then by uniform polynomial approximation of \(\mathcal F'\) on a
compact interval containing the spectra, followed by integration
in \(t\). Permutation symmetry makes the \(r_i\) coordinate
derivatives equal. Thus at a centered stationary point
\(E(x_j)=n\lambda\) for every \(j\), with one common multiplier
\(\lambda\).
Let \(\mathcal L=\partial_x^2-(h''/h')\partial_x\).
Since \((1/h)'>0\) on \(x>0\), the identity
\((E/h)''-\frac{(1/h)''}{(1/h)'}(E/h)'=\frac{\mathcal LE}{h}\)
shows that \(\mathcal LE>0\) makes \(E/h\) strictly convex as a function of \(1/h\).

\par\medskip\noindent\textbf{Kernel and elementary compensation.}
Scaling permits $s=1$. Define
\[
 K(x)=\frac{6ab(\beta-\alpha)^2}
 {(x+\alpha)(x+\beta)\{b(x+\beta)^3+a(x+\alpha)^3\}},
\]
\[
 q_m=1+\frac{(2m+1)^2}{8m^2(m+1)^2},\qquad
 M(x)=\max(1,x/\alpha)\max(1,x/\beta).
\]
The elementary reflected-kernel estimate proved in
Section~\ref{odd23:sec:reflected-kernel} says
\begin{equation}\label{odd23:lm-eq:kernel}
 \mathcal L D(x,\xi)\ge
 -\frac{q_mK(x)}{(\alpha+|\xi|)(\beta+|\xi|)}.
\end{equation}
For $m\ge3$, the exact scalar certificate
\href{https://github.com/zhangteng2000/schoenberg-validation/blob/e7435c5521dcd6a10131ba04a6da500b1404800b/certificates/turn2_log_compensation_certificate.py}{\nolinkurl{turn2_log_compensation_certificate.py}} proves
\begin{equation}\label{odd23:lm-eq:scalar}
 \mathcal L e(x)>q_m(n-1)K(x)M(x),\qquad x>0.
\end{equation}
It reconstructs the rational expression from its definitions.
On the three regions $x\le\alpha$, $\alpha\le x\le\beta$,
$x\ge\beta$, respectively, it substitutes
\[
 x=\frac{z}{(m+1)(1+z)},\quad
 x=\frac{m+(m+1)z}{m(m+1)(1+z)},\quad x=\frac1m+z,
 \qquad m=3+M_0.
\]
After clearing positive denominators, all $160$, $250$, and $199$
coefficients in $z,M_0$ are nonnegative, with positive constant terms.
This is a finite exact rational certificate, not a floating-point check.

Choose an actual root $x_0$ with $|x_0|=r$.
Since
\((x_0-\xi)^2\le(r+|\xi|)^2 \le M(r)(\alpha+|\xi|)(\beta+|\xi|)\),
\eqref{odd23:lm-eq:rootidentity} gives
\begin{equation}\label{eq:weighted-residue-bound}
 \sum_\nu\frac{\rho_\nu}{(\alpha+|\xi_\nu|)(\beta+|\xi_\nu|)}
 \le(n-1)M(r).
\end{equation}
Combining \eqref{odd23:lm-eq:kernel}--\eqref{odd23:lm-eq:scalar} and $M(r)\le M(x)$
proves $\mathcal L E(x)>0$ for every $x\ge r$.
If $r\le\alpha$, then $M(r)=1\le M(x)$ for all $x>0$, so this
strict convexity holds throughout the positive half-line.

\par\medskip\noindent\textbf{The degree-five compensation.}
For $m=2$ \eqref{odd23:lm-eq:scalar} still holds when
$x\le\beta$. For $x\ge\beta$, the same exact verifier proves
\(\mathcal L e>\frac{25}{24}(n-1)K M(x)\).
The needed strengthened negative kernel bound, after the scale
change $\alpha=1$, $\beta=k=3/2$, is
\[
 \mathcal L D(x,-v)\,(x+v)^2
 \ge-\frac{25}{24}K(x)\frac{x^2}{k},\quad x\ge k,
\]
where common positive weight factors are understood.
Using the polynomial \(P_k\) defined in
Subsection~\ref{odd23:sec:reflected-kernel}, we obtain the numerator
\[
 \frac{25}{24}(k-1)^2x^2(1+v)(k+v)(x+v)^2
 -k(k-1)^2(x+v)^4+2vP_k(x,v).
\]
After $x=k+z$, each of its five coefficients in $z$ is a polynomial
$c_0+c_1v+c_2v^2+\cdots$ with $c_0,c_2>0$, nonnegative higher
coefficients, and $4c_0c_2-c_1^2>0$. The five margins are
\[
 \frac{202985919}{1048576},\quad\frac{2847231}{4096},\quad
 \frac{5593671}{16384},\quad\frac{4093}{128},\quad
 \frac{15023}{36864}.
\]
Thus the numerator is positive for all $z,v\ge0$.
For $x\ge r$, $|x_0+v|\le x+v$, so the strengthened estimate is
bounded below by
$-(25/24)K(x)M(x)/(x_0+v)^2$.
For \(\xi\ge0\), use the same-sign formula of
Subsection~\ref{odd23:sec:reflected-kernel} and the weighted residue bound
\eqref{eq:weighted-residue-bound}. Summation using
\eqref{odd23:lm-eq:rootidentity} proves $\mathcal L E>0$ for $x\ge r$ also
when $m=2$. If $r\le\alpha$, the proof again covers all $x>0$.

\par\medskip\noindent\textbf{Support and multiplicity.}
The Euler equation becomes
\(E(x)/h(x)=n\lambda/h(x)\). A strictly convex function intersects a
line at most twice, so there are at most two positive root values.
Reflection gives the same assertion for the negative values.
If there are two positive values, strict convexity gives \(E'(u)<0\)
at the smaller one \(u\). Indeed, the derivative of \(E/h-n\lambda/h\)
with respect to \(1/h\) is negative at its first zero, and
\((1/h)'>0\). If \(u\) is repeated, replace two copies by
\(u+\varepsilon,u-\varepsilon\), keeping the other roots fixed.
Denote the resulting Euler function by \(E_\varepsilon\).
The first variation just proved gives
\begin{equation}\label{eq:log-split-gradient}
  n\,\frac{d}{d\varepsilon}\mathcal J_s
  = E_\varepsilon(u+\varepsilon)-E_\varepsilon(u-\varepsilon).
\end{equation}
To take the difference quotient, write
\(\mu_\varepsilon=\sum_\nu\rho_\nu^\varepsilon
\delta_{\xi_\nu^\varepsilon}\). These positive measures have support
in a fixed compact interval and total mass
\(n^{-1}\sum_j x_j(\varepsilon)^2\), by expansion of \(H\) at infinity.
Their Cauchy transforms \(H_\varepsilon\) converge to \(H\) off the
real axis by the formula \(H_\varepsilon=z-1/g_\varepsilon\).
Hence \(\mu_\varepsilon\) converges weakly to \(\mu_0\), by uniqueness
of the Cauchy transform. For \(x\) in a compact neighborhood of \(u>0\),
\(\partial_xD(x,\xi)\) is jointly continuous in \(x,\xi\), including at
\(\xi=x\) by Taylor's formula. It follows that \(E_\varepsilon'\)
converges uniformly there to \(E'\). Dividing \eqref{eq:log-split-gradient}
by \(\varepsilon\) and taking the limit gives \(\mathcal J_s''(0)=2E'(u)/n<0\),
contradicting minimality.
Reflection gives the assertion on the negative side.

If $r\le\alpha$, strict convexity holds on both entire half-lines.
Two positive intersections require $E(0+)>n\lambda$; two negative
intersections, after reflection, require $E(0-)<n\lambda$.
The one-sided limits satisfy $E(0+)\le E(0-)$: the only possible
jump comes from a critical pole $\xi=0$ and has size
$E(0-)-E(0+)=2\rho_0(b/\alpha^3+a/\beta^3)\ge0$.
The two requirements are incompatible, so there are at most three distinct
values in this case, finishing Theorem~\ref{odd23:lm-thm:main}.

\subsection{The four-value reduction}
\label{odd23:sec:remaining-four}
Suppose $\mathcal J_s<0$ somewhere. Theorem~\ref{odd23:lm-thm:main} gives
an attained negative global minimum, with no zero roots, at most
two values on each sign, and a simple inner value on either sign
with two values. Every configuration with at most three values is
excluded by Theorem~\ref{odd23:thm:all-three-log}.

Thus the minimum has, after reflection, the form
\begin{equation}\label{odd23:eq:complete-four}
 x=(-A)^{[k]},-B,C,D^{[\ell]},\qquad
 0<B<A,\quad0<C<D,\quad kA+B=C+\ell D,
\end{equation}
where $k+\ell=n-2$ and $1\le k\le m-1$, $\ell\ge m$.
The strengthened support restriction also gives
\begin{equation}\label{odd23:eq:complete-scale}
 s<(m+1)\min(B,C).
\end{equation}
Let $\dot B=1$, $\dot A=-1/k$, with $C,D$ fixed. This is a
two-sided admissible centered variation at every strict four-value list.
Theorem~\ref{odd23:mp-thm:onescale} and Corollary~\ref{odd23:nc-cor:negative-center-one-scale}
show that this variation strictly decreases \(\mathcal J_s\) under
\eqref{odd23:eq:complete-scale}, contradicting stationarity.
\subsection{Minority averaging with a nonnegative middle critical point}

Let $n=2m+1\ge5$, $1\le k\le m-1$, $\ell=n-2-k$, and consider
\[
 x=(-A^{[k]},-B,C,D^{[\ell]}),\qquad
 A>B>0,\quad D>C>0,\quad kA+B=C+\ell D.
\]
The minority sign is negative. Vary $B$ at speed $1$ and $A$ at
speed $-1/k$, keeping $C,D$ fixed. A dot denotes this averaging
variation. Put
\[
 \ell_t(x)=|x|/t-\log(1+|x|/t),\qquad
 G_t(x)=\frac{x}{t(t+x)}\quad(x\ge0),\qquad
 \Delta_t=\frac1{t+B}-\frac1{t+A}.
\]
The derivative of the root functional is
$\dot Z_{\ell_t}=-\Delta_t$.

\begin{thm}\label{odd23:mp-thm:onescale}
Suppose the middle of the three nonpersistent critical points is
nonnegative. Then
\[
 \dot W_{\ell_t}< -\Delta_t \quad(0<t\le B),
\]
and, with $\kappa=m/(m+1)=(n-1)/(n+1)$,
\[
 \dot W_{\ell_t}< -\kappa\Delta_t \quad(0<t\le B/\kappa).
\]
Consequently, at every such four-value configuration satisfying
$s<(m+1)\min(B,C)$, the logarithmic comparison functional
$\mathcal J_s$ strictly decreases under minority averaging.
It cannot have a stationary point there.
\end{thm}

\par\medskip\noindent\textbf{Critical velocities and two exact moments.}
Write the three nonpersistent critical points as $-u,v_1,v_2$,
where $B<u<A$ and $0\le v_1<C<v_2<D$.
At a positive simple critical point $v$, secular differentiation gives
\[
 \dot v=\frac{(v+A)^{-2}-(v+B)^{-2}}
 {k(v+A)^{-2}+(v+B)^{-2}+(v-C)^{-2}+\ell(v-D)^{-2}}<0.
\]
Write $\rho_i=-\dot v_i>0$, $\eta=\rho_1+\rho_2$, and
$d=(k-1)/k$. Trace zero gives $\dot u=d-\eta$.
Let \(K=A-B\). The quantities \(q_1,q_2\) below are positive,
whereas \(q_0=0\) if \(k=1\) and \(q_0>0\) if \(k>1\). Put
\[
 q_0=d(A-u),\qquad q_i=\rho_i(u+v_i)\quad(i=1,2).
\]
The identities for centered roots and critical points,
\[
 W_2=\frac{n-2}{n}Z_2,\qquad
 \sum_{j=1}^{n-1}w_j^3=\frac{n-3}{n}\sum_{i=1}^{n}x_i^3,
\]
and their derivatives imply, with $c=(n-2)/n$ and $c_3=(n-3)/n$,
\begin{equation}\label{odd23:mp-eq:moments}
 q_1+q_2=cK-q_0=:Q,
 \qquad q_1v_1+q_2v_2=ucK+Aq_0-c_3K(A+B)=:M_1.
\end{equation}
$\dot Z_2=-2K$ and
$\frac{d}{dt_{\rm flow}}\sum_{i=1}^{n}x_i^3=3K(A+B)$.
In the critical sums the persistent negative block contributes
$-2dA$ to $\dot W_2$, and $3dA^2$ to the signed cubic derivative.
Subtracting the persistent block contributions gives \eqref{odd23:mp-eq:moments} directly.

The exact derivative of the trace functional is
\begin{equation}\label{odd23:mp-eq:tracegradient}
 -\dot W_{\ell_t}
 =\frac{q_0}{(t+A)(t+u)}
 +\sum_{i=1}^2q_i\frac{G_t(u)+G_t(v_i)}{u+v_i}.
\end{equation}
Equations~\eqref{odd23:mp-eq:moments} and \eqref{odd23:mp-eq:tracegradient}
extend to \(v_1=0\) by continuity.

\par\medskip\noindent\textbf{An elementary kernel lower bound.}
For all $u,t>0$ and $v\ge0$,
\begin{equation}\label{odd23:mp-eq:kernel}
 \frac{G_t(u)+G_t(v)}{u+v}
 \ge \frac{1-v/(4u)}{t(t+u)}.
\end{equation}
Indeed the left side, divided by $1/[t(t+u)]$, is
$1+v(u-v)/[(t+v)(u+v)]$. The difference from $1-v/(4u)$ is
\(\frac{v\{t(u+v)+v^2-3uv+4u^2\}} {4u(t+v)(u+v)}\ge0\),
since $v^2-3uv+4u^2=(v-3u/2)^2+7u^2/4>0$.
Combining \eqref{odd23:mp-eq:kernel} with \eqref{odd23:mp-eq:tracegradient} gives
\begin{equation}\label{odd23:mp-eq:lower}
 -\dot W_{\ell_t}\ge
 \frac{q_0}{(t+A)(t+u)}+\frac{M}{t(t+u)},
 \qquad M=Q-\frac{M_1}{4u}.
\end{equation}

\par\medskip\noindent\textbf{The secular position constraint.}
Let $h=\ell+1\ge k+2$ and normalize $K=1$ by scaling.
Set
\[
 p=B/K>0,\qquad z=(u-B)/K\in(0,1),\qquad
 A=p+1,\quad u=p+z,\quad n=k+h+1.
\]
The negative critical equation is
\(\frac{k}{A-u}=\frac1{u-B}+\frac1{u+C}+\frac\ell{u+D}\).
It first gives $z>1/(k+1)$.
Jensen's inequality for the $h$ positive roots, whose sum is
$kA+B$, gives the additional constraint
\begin{equation}\label{odd23:mp-eq:F}
 F:=hz^2+\{(k+1)p+k-h\}z-p-k/n\ge0.
\end{equation}
To see this, use
$1/(u+C)+\ell/(u+D)\ge h^2/(hu+kA+B)$ and multiply the
resulting inequality by the positive denominators.

\par\medskip\noindent\textbf{The scalar certificate.}
With $d=(k-1)/k$, one has $q_0=d(1-z)$ and
\[
 4(p+z)M=\mathcal M
 :=3c(p+z)-d(1-z)(5p+4z+1)+c_3(2p+1).
\]
Polynomial division in $z$ gives
\(\mathcal M=\frac{4(k-1)}{hk}F+R_M\).
After $h=k+2+H$, $k=1+K_0$, the numerator of $R_M$ over its
positive denominator has $25$ nonnegative coefficients in
$p,z,K_0,H$ and a positive constant coefficient.
Therefore $M>0$ throughout the required domain.

For either $\theta=1$ or $\theta=(n-1)/(n+1)$, comparison of
\eqref{odd23:mp-eq:lower} with $\theta\Delta_t$ is equivalent to
\[
 \Phi_\theta(t):=
 q_0t(t+p)+M(t+p+1)(t+p)-\theta t(t+p+z)>0.
\]
Its quadratic coefficient is
$c-\theta-M_1/[4(p+z)]<0$ by \eqref{odd23:mp-eq:moments}; it is concave
in $t$, and $\Phi_\theta(0)=p(p+1)M>0$.
It is enough to check its other endpoint $t=p/\theta$.
Define
\(E_\theta=\frac{4(p+z)\theta^2}{p}\Phi_\theta(p/\theta)\).
For $k\ge2$, divide $E_\theta=Q_\theta F+R_\theta$ in $z$.
The remainder is affine in $z$. Exact rational calculation gives
\[
 Q_1=\frac{4(2kp+k-2p-2)}{hk},
\]
\[
 Q_\kappa=
 \frac{4(h+k)(2hkp+hk-2hp-2h+2k^2p+k^2-2p-2)}
 {hk(h+k+2)^2}.
\]
Both are nonnegative for $k\ge2$, $h\ge k+2$, $p\ge0$.
After $k=2+K_0$, $h=k+2+H$, exact expansions verify
\[
 \partial_zR_\theta\ge0,
 \qquad R_\theta\bigl(1/(k+1)\bigr)>0.
\]
The numbers of positive coefficients in the numerators of \(\partial_z R_\theta\) and \(R_\theta(1/(k+1))\) are
$26,8$ for $\theta=1$, and $41,30$ for $\theta=\kappa$.
Together with \eqref{odd23:mp-eq:F} and $z>1/(k+1)$ this proves
$E_\theta>0$. All polynomial identities and coefficient signs are
reconstructed and checked in exact rational arithmetic by
\href{https://github.com/zhangteng2000/schoenberg-validation/blob/e7435c5521dcd6a10131ba04a6da500b1404800b/certificates/turn2_minority_positive_center_certificate.py}{\nolinkurl{turn2_minority_positive_center_certificate.py}}.
It has no cached inputs and no floating-point calculations.

For $k=1$, $q_0=0$ and
\(M=3c/4+c_3(2p+1)/[4(p+z)]\).
For fixed endpoint $t=p/\theta$, the original $\Phi_\theta(t)$
strictly decreases with $z$. Since $z<1$, it suffices to use $z=1$.
The two scaled endpoints factor as
\[
 E_1\big|_{k=1,z=1}
 =\frac{2(2p+1)(3hp+2h-6p-5)}{h+2}>0,
\]
\[
 E_\kappa\big|_{k=1,z=1}
 =\frac{2(2hp+h+4p+1)(3hp+2h-4p-3)}{(h+3)^2}>0,
\]
since $h\ge3$. Concavity in $t$ now proves both one-scale
inequalities in Theorem~\ref{odd23:mp-thm:onescale}.

\par\medskip\noindent\textbf{Application to the logarithmic comparison.}
Because $\ell\ge k+1$, centering implies $A>D$. Nonnegativity of
the middle critical point is equivalent to
\(k/A+1/B-1/C-\ell/D\ge0\),
and hence $B<C$.
Let $a=m+1$, $b=m$, $\alpha=s/a$, $\beta=s/b$.
The minimum's scale constraint gives $\alpha<B$ and
$\beta<B/\kappa$. Therefore
\[
\dot{\mathcal J}_s =b\dot W_{\ell_\alpha}+a\dot W_{\ell_\beta} +(b-1)\Delta_\alpha+b\Delta_\beta+\Delta_s < -b\Delta_\alpha-a\kappa\Delta_\beta +(b-1)\Delta_\alpha+b\Delta_\beta+\Delta_s =-\Delta_\alpha+\Delta_s<0.
\]
This proves the final assertion and excludes this entire four-value
branch from being a negative global logarithmic minimum.
\subsection{The negative-middle branch: a residue proof}
\label{odd23:nc-sec:negative-center-monotonicity}

Fix $n=2m+1\ge5$ and integers
\[
 1\le k\le m-1,\qquad \ell=n-2-k,\qquad
 h=\ell+1=n-k-1\ge k+2.
\]
Consider the centered real-rooted polynomial with root list
\[
 (-A)^{[k]},\ -B,\ C,\ D^{[\ell]},\qquad
 A>B>0,\quad D>C>0,\quad kA+B=C+\ell D.
\]
Here the negative sign group is the minority sign group. In particular,
$\ell\ge k+1$, and centering implies $A>D$, since
\((k+1)A>kA+B=C+\ell D>\ell D\ge(k+1)D\).
Use the centered averaging variation
\begin{equation}\label{odd23:nc-eq:negative-center-flow}
 \dot A=-\frac1k,\qquad \dot B=1,\qquad
 \dot C=\dot D=0,
\end{equation}
and put $K=A-B$. For $t>0$, define
\[
 \ell_t(x)=\frac{|x|}{t}-\log\left(1+\frac{|x|}{t}\right),
 \qquad W_t=\sum_{j=1}^{n-1}\ell_t(w_j),\qquad
 \Delta_t=\frac{K}{(t+A)(t+B)},\qquad
 L(t)=-\frac{\dot W_t}{\Delta_t}.
\]
The corresponding root sum satisfies
$\frac{d}{d\tau}\sum_{i=1}^{n}\ell_t(x_i)=-\Delta_t$.
Assume that the middle nonpersistent critical point is nonpositive.
Write the three nonpersistent critical points as
\[
 -u,\ -v,\ w,\qquad
 B<u<A,\quad 0\le v<B,\quad C<w<D.
\]
The persistent critical points are $-A$, with multiplicity $k-1$,
and $D$, with multiplicity $\ell-1$.

\begin{thm}\label{odd23:nc-thm:negative-center-monotonicity}
For every such root list,
\[
 \frac{d}{dt}\bigl(tL(t)\bigr)>0\qquad(t>0).
\]
Moreover,
\begin{equation}\label{odd23:nc-eq:negative-center-final-bound}
 L(t)>\min\left\{1,\frac Ct\right\}\qquad(t>0).
\end{equation}
\end{thm}

\begin{proof}
First suppose \(v>0\). We pass to \(v=0\) after estimating
\((tL(t))'\) by \eqref{odd23:nc-eq:negative-center-tLprime}.
The secular derivative formula, applied to
\(g(z)=\frac{k}{z+A}+\frac1{z+B} +\frac1{z-C}+\frac{\ell}{z-D}\),
shows that
\[
 \zeta=\dot v>0,\qquad \eta=-\dot w>0.
\]
Indeed, at either $z=-v$ or $z=w$, the signed critical velocity is
\[
 \frac{(z+A)^{-2}-(z+B)^{-2}}
 {k(z+A)^{-2}+(z+B)^{-2}+(z-C)^{-2}+\ell(z-D)^{-2}}<0.
\]
The critical-point trace is zero, so, with $d=(k-1)/k$,
\(\dot u=d-\zeta-\eta\).

\par\medskip\noindent\textbf{The exact residue identity.}
For an absolute critical magnitude $z>0$ with velocity $\dot z$,
its contribution to $L(t)$ can be written as
\[
 -\frac{\dot z}{K}\frac{z(t+A)(t+B)}{t(t+z)}
 =-\frac{z\dot z}{K}-\frac{AB\dot z}{Kt}
  +\frac{\dot z(z-A)(z-B)}{K(t+z)}.
\]
The centered quadratic identity is
\[
 \sum_{j=1}^{n-1}w_j^2=c_n\sum_{i=1}^{n}x_i^2,
 \qquad c_n=\frac{n-2}{n}.
\]
Consequently $\sum z\dot z=-c_nK$ along
\eqref{odd23:nc-eq:negative-center-flow}. Also $\sum\dot z=-2\eta$:
the signed critical trace remains zero, and $w$ is the only moving
positive critical point. It follows that
\begin{equation}\label{odd23:nc-eq:negative-center-residues}
 L(t)=c_n+\frac{r_0}{t}+\frac{r_u}{t+u}
                 +\frac{r_v}{t+v}+\frac{r_w}{t+w},
\end{equation}
where
\[
r_0=\frac{2\eta AB}{K}>0,\quad r_u=-\frac{\dot u(A-u)(u-B)}K,\quad r_v=\frac{\zeta(A-v)(B-v)}K\ge0,\quad r_w=-\frac{\eta(A-w)(B-w)}K.
\]
The persistent critical points contribute no nonzero finite-pole
residues: the residue at $A$ has a factor $z-A$, and the velocity at
$D$ is zero. Differentiating \eqref{odd23:nc-eq:negative-center-residues} gives
\begin{equation}\label{odd23:nc-eq:negative-center-tLprime}
 (tL(t))'=c_n+\frac{r_uu}{(t+u)^2}
                 +\frac{r_vv}{(t+v)^2}
                 +\frac{r_ww}{(t+w)^2}.
\end{equation}

\par\medskip\noindent\textbf{The positive critical point contributes at worst $4/27$.}
Since $A>D>w$, the residue $r_w$ can be negative only when $B>w$.
Assume this case and put
\[
 S_2(w)=\frac{k}{(w+A)^2}+\frac1{(w+B)^2}
       +\frac1{(w-C)^2}+\frac{\ell}{(w-D)^2}.
\]
Secular differentiation gives
\(\frac{\eta}{K} =\frac{2w+A+B}{(w+A)^2(w+B)^2S_2(w)}\).
Because $0<C<w$, we have $S_2(w)>w^{-2}$. Hence
\(-\frac{r_w}{w} <\frac{w(2w+A+B)(A-w)(B-w)}{(A+w)^2(B+w)^2}\).
Set $P=(A-w)/(A+w)$ and $Q=(B-w)/(B+w)$, so $0<P,Q<1$.
The upper bound for \(-r_w/w\) becomes
\(\frac12PQ(2-P-Q)\le\frac4{27}\).
Indeed, \(PQ\le(P+Q)^2/4\), and the maximum of \(s^2(2-s)/8\)
on \([0,2]\) is \(4/27\), attained at \(s=4/3\).
Thus, whether or not $r_w$ is negative,
\begin{equation}\label{odd23:nc-eq:negative-center-w-bound}
 \frac{r_ww}{(t+w)^2}>-\frac4{27}\qquad(t>0).
\end{equation}

\par\medskip\noindent\textbf{The outer negative critical point contributes at worst $1/2$.}
It remains to consider $r_u<0$, equivalently $\dot u>0$.
Put
\[
 x=A-u,\qquad y=u-B,\qquad
 q=\frac xy,\qquad p=\frac By.
\]
Then $p>0$, $u=(p+1)y$, and $K=(q+1)y$.
Writing the positive roots as $a_1=C$ and
$a_2=\cdots=a_h=D$, the secular equation at $-u$ is
\begin{equation}\label{odd23:nc-eq:negative-center-u-secular}
 \frac{k}{x}=\frac1y+\sum_{j=1}^h\frac1{u+a_j}.
\end{equation}
Furthermore,
\[
 \dot u=\frac{y^{-2}-x^{-2}}{S_2(-u)},\qquad
 S_2(-u)=\frac{k}{x^2}+\frac1{y^2}
                 +\sum_{j=1}^h\frac1{(u+a_j)^2}.
\]
Thus $\dot u>0$ implies $q>1$, while
\eqref{odd23:nc-eq:negative-center-u-secular} implies $q<k$.
In particular, the harmful case is impossible if $k=1$.

Suppose $k\ge2$. Define
\(D_0=h(k+q^2)+(k-q)^2\).
By Cauchy's inequality and \eqref{odd23:nc-eq:negative-center-u-secular},
\(y^2S_2(-u)\ge \frac{k}{q^2}+1 +\frac{(k/q-1)^2}{h} =\frac{D_0}{hq^2}\).
Consequently
\begin{equation}\label{odd23:nc-eq:negative-center-u-prebound}
 -\frac{r_u}{u}
 =\frac{\dot u\,xy}{Ku}
 \le\frac{hq(q-1)}{(p+1)D_0}.
\end{equation}
Centering gives $\sum_{j=1}^h a_j=kA+B$, and hence
\(\sum_{j=1}^h(u+a_j)=y(np+h+k+kq)\).
Applying the reciprocal Cauchy inequality to this sum and to
\eqref{odd23:nc-eq:negative-center-u-secular} gives
\begin{equation}\label{odd23:nc-eq:negative-center-G-positive}
 G:=(k-q)(np+h+k+kq)-h^2q\ge0.
\end{equation}
Set
\[
 F=(p+1)D_0-2hq(q-1),\qquad d_0=h-k\ge2.
\]
Expanding the difference \(F-G\) in terms of \(d_0=h-k\) gives
\begin{equation}\label{odd23:nc-eq:negative-center-F-G}
F-G=p\bigl[(h+1)q^2+(d_0+1)q-k\bigr] +q\bigl[(h+k)(d_0+1)+2d_0-(d_0-1)q\bigr].
\end{equation}
On $1<q<k$, the first bracket is at least $2(d_0+1)>0$,
and the second is at least $h(d_0+3)>0$.
Equations \eqref{odd23:nc-eq:negative-center-G-positive} and
\eqref{odd23:nc-eq:negative-center-F-G} now give \(F>G\ge0\).
From \eqref{odd23:nc-eq:negative-center-u-prebound},
\(-\frac{r_u}{u}<\frac12\),
and therefore, also allowing the case $r_u\ge0$,
\begin{equation}\label{odd23:nc-eq:negative-center-u-bound}
 \frac{r_uu}{(t+u)^2}>-\frac12\qquad(t>0).
\end{equation}

\par\medskip\noindent\textbf{Strict monotonicity.}
The $r_v$ term in \eqref{odd23:nc-eq:negative-center-tLprime} is nonnegative.
For $n\ge7$, \eqref{odd23:nc-eq:negative-center-w-bound} and
\eqref{odd23:nc-eq:negative-center-u-bound} yield
\((tL(t))'>\frac{n-2}{n}-\frac12-\frac4{27} \ge\frac{25}{378}>0\).
For $n=5$, necessarily $k=1$, so $r_u\ge0$ and
\((tL(t))'>\frac35-\frac4{27}=\frac{61}{135}>0\).
All the quantities defining $L$ and $(tL)'$ have continuous limits
as $v\downarrow0$: $\ell_t$ is continuously differentiable at zero,
and in \eqref{odd23:nc-eq:negative-center-tLprime} the $v$ term tends to zero.
At $v=0$, continuity consequently gives $(tL)'\ge25/378$ for
\(n\ge7\), and \((tL)'\ge61/135\) for \(n=5\).

\par\medskip\noindent\textbf{The base estimate and its extension.}
For completeness, we show $L(t)>1$ for $0<t\le C$.
Let
\[
 p(z)=(z+A)^k(z+B)(z-C)(z-D)^\ell,
 \qquad R(z)=\frac{p(z)}{(z+A)(z+B)},\qquad
 \psi(z)=\frac{R'(z)}{p'(z)}.
\]
The variation \eqref{odd23:nc-eq:negative-center-flow} satisfies
$\dot p(z)=K R(z)$. For $0<t<C$, the secular function $g(z)=p'(z)/p(z)$ is
strictly decreasing on $(-B,C)$ and vanishes at $-v\le0$.
Thus $g(t)<0$ and
\[
 \psi(t)=\frac1{(t+A)(t+B)}
 \left(1-\frac{(t+A)^{-1}+(t+B)^{-1}}{g(t)}\right)
 >\frac1{(t+A)(t+B)}.
\]
As \(t\uparrow C\), \(g(t)\to-\infty\), so
\(\psi(C)=1/((C+A)(C+B))\).
The logarithmic determinant variation is
\(\dot W_t=-K\psi(t) +\frac{2w^2}{t(w^2-t^2)}\dot w\).
To see the correction term, first use
$-z/t-\log|1-z/t|$ on the signed critical spectrum; its total
variation is $-K\psi(t)$ because the signed trace has zero variation.
On the positive spectrum the correction to $\ell_t$ is
$2z/t+\log|1-z/t|-\log(1+z/t)$, whose derivative in $z$ is
$2z^2/[t(z^2-t^2)]$. Only $w$ has nonzero positive spectral velocity.
Since $t\le C<w$ and $\dot w<0$, this correction is strictly negative.
It follows that $\dot W_t<-\Delta_t$ and $L(t)>1$ for $0<t\le C$.

Finally, strict monotonicity of $tL(t)$ gives, for $t>C$,
\(tL(t)>CL(C)>C\).
This proves \eqref{odd23:nc-eq:negative-center-final-bound}.
\end{proof}

\begin{cor}\label{odd23:nc-cor:negative-center-one-scale}
Put $\kappa=m/(m+1)$. Then
\[
 \dot W_t<-\kappa\Delta_t
 \qquad\text{whenever}\qquad
 0<t\le\frac{\min(B,C)}\kappa.
\]
Consequently the negative-middle four-value branch cannot be a
negative minimum of the logarithmic comparison functional.
\end{cor}
\begin{proof}
The first assertion follows from
$L(t)>\min\{1,C/t\}\ge\kappa$ in the stated range.
For the last assertion, let $a=m+1$, $b=m$ and assume that a
putative negative minimum has spectral parameter
$s<a\min(B,C)$. Write $\alpha=s/a$ and $\beta=s/b$.
The root terms in the logarithmic comparison give the exact
variation
\(\dot{\mathcal J}_s =b\dot W_\alpha+a\dot W_\beta +(b-1)\Delta_\alpha+b\Delta_\beta+\Delta_s\).
Here $\alpha<\min(B,C)\le C$ and
$\beta<\min(B,C)/\kappa$. The theorem and the first assertion imply
\[
 \dot{\mathcal J}_s
 <-b\Delta_\alpha-a\kappa\Delta_\beta
       +(b-1)\Delta_\alpha+b\Delta_\beta+\Delta_s
 =-\Delta_\alpha+\Delta_s<0,
\]
because $a\kappa=b$, $s>\alpha$, and $\Delta_t$ strictly decreases
in $t$. This contradicts stationarity.
\end{proof}

\subsection{Completion of the logarithmic and power inequalities}
\label{odd23:sec:complete-assembly}
\begin{thm}\label{odd23:thm:global-log}
For every odd $n=2m+1\ge5$, every centered real root list $x$,
and every $s>0$, one has $\mathcal J_s(x)\ge0$.
\end{thm}
\begin{proof}
A negative value would give the attained four-value minimum
\eqref{odd23:eq:complete-four} with \eqref{odd23:eq:complete-scale}.
Set $a=m+1$, $b=m$, $\alpha=s/a$, $\beta=s/b$,
$\kappa=b/a$. Let $\Delta_t=(t+B)^{-1}-(t+A)^{-1}>0$.
It strictly decreases with $t$.

If the middle nonpersistent critical point is nonnegative, the
first averaging theorem gives
$\dot W_{\ell_\alpha}< -\Delta_\alpha$ and
$\dot W_{\ell_\beta}< -\kappa\Delta_\beta$.
If it is negative, the base estimate gives the first bound, and
the residue monotonicity theorem gives the second: for $\beta\le C$
the stronger bound with coefficient $1$ holds, and for
$C<\beta<C/\kappa$ one has $-\dot W_{\ell_\beta}/\Delta_\beta>C/\beta>\kappa$.
The scale restriction supplies $\alpha<\min(B,C)$ and
$\beta<\min(B,C)/\kappa$ in both cases.

Since $\dot Z_{\ell_t}=-\Delta_t$, it follows that
\[
\dot{\mathcal J}_s =b\dot W_{\ell_\alpha}+a\dot W_{\ell_\beta} +(b-1)\Delta_\alpha+b\Delta_\beta+\Delta_s <-b\Delta_\alpha-a\kappa\Delta_\beta +(b-1)\Delta_\alpha+b\Delta_\beta+\Delta_s =-\Delta_\alpha+\Delta_s<0.
\]
This contradicts stationarity of the minimum. Hence a negative
value does not exist.
\end{proof}

Equation~\eqref{odd23:eq:complete-Mellin} now proves the claimed inequality
for every centered real list and every $2<p<3$.
The angular compression argument proves it for arbitrary centered
complex roots with the same coefficient.

For the equality polynomial in Theorem~\ref{odd23:thm:main-complete},
direct differentiation gives the critical list
\[
(n+1)^{[\frac{n-3}{2}]},\quad (-(n-1))^{[\frac{n-1}{2}]},\quad2.
\]
Its root sum is zero, and
\[
 \sum_{j=1}^{n}|z_j|^p=
 \frac{(n-1)(n+1)^p+(n+1)(n-1)^p}{2},
\]
\[
 \sum_{k=1}^{n-1}|w_k|^p=
 \frac{(n-3)(n+1)^p+(n-1)^{p+1}+2^{p+1}}{2}.
\]
Their ratio is $C_{n,p}$, proving equality and sharpness.
This completes Theorem~\ref{odd23:thm:main-complete}.

\subsection{Exact verification}
Appendix~\ref{app:odd-exact-certificates} identifies the algebraic checks
used in Theorem~\ref{odd23:thm:global-log}. In the proof of
Theorem~\ref{odd23:nc-thm:negative-center-monotonicity}, the
\href{https://github.com/zhangteng2000/schoenberg-validation/blob/e7435c5521dcd6a10131ba04a6da500b1404800b/certificates/turn2_negative_center_exact_check.py}{symbolic verifier}
checks the identities used to obtain the residue bounds \eqref{odd23:nc-eq:negative-center-w-bound} and \eqref{odd23:nc-eq:negative-center-u-bound}.
\section{Auxiliary proofs for the odd-degree inequality when \texorpdfstring{\(p\ge3\)}{p at least 3}}\label{app:odd-ge-three-details}
\subsection{Sharp inequalities in odd degree}\label{sec:odd}
This appendix proves Theorem~\ref{thm:centered-odd}, including its
equality cases, for \(p\ge3\). The equality case at \(p=3\) is also
used in Section~\ref{sec:odd-all-exponents} for \(2<p<3\).
The polynomial sign checks appear in Appendix~\ref{app:odd-exact-certificates}.
\begin{thm}\label{thm:oddcomplete}
For $n=2m+1\ge5$, $p\ge3$, and every centered complex configuration,
\begin{equation}\label{eq:oddsharp}
 W_p\le B_{n,p}Z_p.
\end{equation}
Apart from \(p(z)=c z^n\) with \(c\ne0\), equality holds precisely for the polynomials in \eqref{eq:equality}.
\end{thm}
The expression in \eqref{eq:coefficient} is the coefficient in
Theorem~\ref{thm:centered-odd}, after cancelling a common factor $2^{p+1}$.
At the exponent \(p=3\), \eqref{eq:coefficient} gives
\(B_{n,3}=(n^3-n^2-n-7)/((n+1)(n^2+1))\).

\subsubsection{The real extremal reduction}\label{sec:extremizers}
For a centered nonzero real configuration $x=(x_1,\ldots,x_n)$, write $W_p/Z_p$ for the ratio of the critical-point and root power sums defined in \eqref{eq:notation}. By homogeneity the maximum is attained on the compact set $\sum_{j=1}^{n}x_j=0$, $\sum_{j=1}^{n}x_j^2=1$.

\begin{lem}\label{lem:reduction}
Let $n=2m+1\ge5$ and $p\ge3$. Every real configuration maximizing $W_p/Z_p$ has at most three distinct values. If the values are $a<b<c$, the multiplicity of $b$ is one.
\end{lem}
\begin{proof}
\emph{Step 1. A positive rational representation.}
Write the distinct zero values as $a_1<\cdots<a_s$, with multiplicities $r_i$. The $s-1$ critical points not forced by repeated roots, denoted by $\xi_\nu$, are the simple interlacing zeros of
$G(z)=n^{-1}\sum_{i=1}^{s}r_i/(z-a_i)$. Define
$A_\nu=\sum_{i=1}^{s}r_i/(\xi_\nu-a_i)^2$ and $\rho_\nu=n/A_\nu>0$. Partial fractions and centering at infinity give
\begin{equation}\label{eq:H}
 H(z):=z-\frac1{G(z)}=\sum_{\nu=1}^{s-1}\frac{\rho_\nu}{z-\xi_\nu},\qquad
 H(a_i)=a_i,\qquad H'(a_i)=1-\frac n{r_i}.
\end{equation}
Indeed, $-1/G$ has residue $n/A_\nu$ at $\xi_\nu$, and
$H(z)=(n^{-1}\sum_{j=1}^{n}x_j^2)z^{-1}+O(z^{-2})$. In particular, the total mass $\sum_{\nu=1}^{s-1}\rho_\nu$ is positive.

\emph{Step 2. One equation for all distinct zero values.}
Let $C_*$ be the maximal ratio and put $f(t)=|t|^p$. Then $F=W_p-C_*Z_p$ has a constrained maximum of value zero. The unforced roots depend smoothly on the distinct zero values, with
$\partial\xi_\nu/\partial a_i=r_i/(A_\nu(\xi_\nu-a_i)^2)$.
Thus, for some real $\lambda$,
\(\left(1-\frac1{r_i}-C_*\right)f'(a_i) +\sum_{\nu=1}^{s-1}\frac{f'(\xi_\nu)}{A_\nu(\xi_\nu-a_i)^2}=\lambda\).
For $x\ne y$ define
\begin{equation}\label{eq:divdiff}
 D_f(x,y)=\frac{f'(y)-f'(x)-f''(x)(y-x)}{(y-x)^2},\qquad
 L(x)=\sum_{\nu=1}^{s-1}\rho_\nu D_f(x,\xi_\nu).
\end{equation}
For $p>3$ take the continuous value $f'''(x)/2$ on the diagonal. For $p=3$ this convention is used away from zero; when zero is an original zero, none of the $\xi_\nu$ is zero and $L(0)$ is well defined. Substituting \eqref{eq:H} and $xf''(x)=(p-1)f'(x)$ into the stationary equations gives
\begin{equation}\label{eq:g}
 g(a_i)=0,\qquad
 g(x)=L(x)-d f'(x)-n\lambda,\qquad d=p-n(1-C_*).
\end{equation}
The explicit configuration \eqref{eq:equality} gives $C_*\ge B_{n,p}>(m-1)/m$, because
$B_{n,p}-(m-1)/m=(m^{p-1}+1)/D_p>0$. Hence $n(1-C_*)<n/m\le5/2<p$, and $d>0$.

\emph{Step 3. Counting zeros when $p>3$.}
For $x\ne0$ let $K(x)=L'(x)/f''(x)$. Since \(f'''\) is locally absolutely continuous for \(p>3\), Taylor's integral remainder in \eqref{eq:divdiff}, followed by the substitutions $t=s/(1+s)$ and $s=xu$, gives, for $x>0$,
\begin{equation}\label{eq:kernel}
 K(x)=\frac{(p-2)(p-3)}x
 \sum_{\nu=1}^{s-1}\rho_\nu\int_0^\infty
 \frac{|1+u\xi_\nu|^{p-4}}{(1+xu)^p}\dd u.
\end{equation}
The singularity at $1+u\xi_\nu=0$ is integrable, since $p-4>-1$. The integral also converges at infinity. For fixed $u$, the factor $x^{-1}(1+xu)^{-p}$ strictly decreases in $x>0$. Consequently \eqref{eq:kernel} shows that \(K\) strictly decreases from infinity to zero on that half-line. The limit at zero follows by integrating over a fixed small $u$ interval; for $x\ge x_0$, the bound $K(x)\le(x_0/x)K(x_0)$ proves the other limit. Reflection gives strict increase from zero to infinity on the negative half-line.

By \eqref{eq:g}, $g'(x)=f''(x)(K(x)-d)$. Thus the continuous function $g$ strictly decreases, then increases, then decreases. It has at most three distinct real zeros. If it has three zeros and the middle zero \(b\) is nonzero, then $g'(b)>0$. If $b=0$ is an original zero, none of the $\xi_\nu$ is zero, and direct differentiation gives
\begin{equation}\label{eq:gzero}
 g'(0)=2p\sum_{\nu=1}^{s-1}\rho_\nu|\xi_\nu|^{p-4}>0.
\end{equation}

\emph{Step 4. The endpoint $p=3$.}
For $x>0$, the divided difference equals $3$ if $y\ge0$, and equals $3-6y^2/(x-y)^2$ if $y<0$. Hence
\begin{equation}\label{eq:kernel3}
 \frac{L'(x)}{6x}=\frac2x\sum_{\xi_\nu<0}
 \frac{\rho_\nu\xi_\nu^2}{(x-\xi_\nu)^3}.
\end{equation}
The right side of \eqref{eq:kernel3} strictly decreases from infinity to zero when its sum is nonempty, and is zero otherwise. Reflection gives the corresponding assertion on the negative half-line. Thus the same three monotonicity intervals suffice, allowing an empty increasing part on one half-line. A mass at $\xi_\nu=0$ produces an upward jump of $6\rho_\nu$ at zero, which cannot create an additional zero at an admissible original-root location. If both half-line increasing parts are empty, there are at most two zeros. Otherwise the increasing parts and any upward jump join into a single increasing portion. There are at most three admissible zeros, with positive derivative at a nonzero middle one.

If zero is an original zero, \eqref{eq:H} gives
\(\sum_\nu\rho_\nu/\xi_\nu=0\). By \eqref{eq:kernel3}, the right derivative
of \(g\) at zero is \(12\sum_{\xi_\nu<0}\rho_\nu/|\xi_\nu|\).
Reflection gives the left derivative \(12\sum_{\xi_\nu>0}\rho_\nu/\xi_\nu\).
The two sums agree, so \(g'(0)=6\sum_\nu\rho_\nu/|\xi_\nu|>0\).

\emph{Step 5. Splitting a repeated middle value.}
Suppose the middle value $b$ has multiplicity $r\ge2$, and write $p_0(z)=(z-b)^rQ(z)$, $Q(b)\ne0$. Replace two copies of $b$ by $b-\varepsilon,b+\varepsilon$:
$p_\varepsilon(z)=(z-b)^{r-2}((z-b)^2-\varepsilon^2)Q(z)$.
The logarithmic derivative changes by
$2\varepsilon^2/((z-b)((z-b)^2-\varepsilon^2))$. Thus
\(\xi_\nu(\varepsilon)=\xi_\nu+ \frac{2\varepsilon^2}{A_\nu(\xi_\nu-b)^3}+O(\varepsilon^4)\).
Write the $r-1$ critical points tending to $b$ as $b+u_j$. Their sums obey
\[
 \sum_{j=1}^{r-1}u_j=-2\varepsilon^2\sum_{\nu=1}^{s-1}\frac1{A_\nu(\xi_\nu-b)^3}+O(\varepsilon^4),
 \qquad
 \sum_{j=1}^{r-1}u_j^2=\frac{2(r-2)}r\varepsilon^2+o(\varepsilon^2).
\]
The first formula uses the unchanged sum of all critical points. For the second, differentiate $p_\varepsilon(z)$ and scale $z-b=\varepsilon u$: for $r\ge3$ the two moving roots have leading displacements $\pm\sqrt{(r-2)/r}\,\varepsilon$, while for $r=2$ the sole nearby displacement is $O(\varepsilon^2)$.

Since $f\in C^2$, expansion of $F$ is valid also at $b=0$ and gives
\begin{equation}\label{eq:splitting}
F(p_\varepsilon)-F(p_0) =\left[\left(1-C_*-\frac2r\right)f''(b) +2\sum_{\nu=1}^{s-1}\frac{f'(\xi_\nu)-f'(b)}{A_\nu(\xi_\nu-b)^3}\right]\varepsilon^2 +o(\varepsilon^2) =\frac{g'(b)}n\varepsilon^2+o(\varepsilon^2)>0.
\end{equation}
For $b\ne0$, the second equality follows by differentiating \eqref{eq:divdiff}, using \eqref{eq:H} and $bf'''(b)=(p-2)f''(b)$. At $b=0$ it follows directly from \eqref{eq:gzero}, or its cubic version. Equation~\eqref{eq:splitting} contradicts maximality. The middle value is simple.
\end{proof}

\subsubsection{Cubic comparison and block sums}\label{sec:block}
Every centered three-value configuration in Lemma~\ref{lem:reduction}, up to reflection and positive scaling, is
\begin{equation}\label{eq:triple}
 -\ell-t\ (k\text{ copies}),\qquad
 (n-1)t\ (1\text{ copy}),\qquad
 k-t\ (\ell\text{ copies}),
 \quad -\ell/n\le t\le k/n,
\end{equation}
where $1\le k\le m$ and $\ell=2m-k$. The endpoints include all centered two-value configurations. Indeed, normalize the gap between the smallest and largest values to $n-1$ and use centering to obtain \eqref{eq:triple}.

Set $\tau=nt$ and multiply \eqref{eq:triple} by $n$. Introduce
\begin{equation}\label{eq:ABCD}
 \begin{gathered}
 A=n\ell+\tau,\quad B=(n-1)\tau,\quad C=nk-\tau,\quad
 S=n(k-\ell)+(n-3)\tau,\\
 D=((n-1)\tau-n(k-\ell))^2+4nk\ell.
 \end{gathered}
\end{equation}
The original values are $-A,B,C$ and the two unforced critical points are
$\lambda_\pm=(S\pm\sqrt D)/2$. One has $A,C>0$ and $\lambda_-<0$. For $B\le0$ this follows from interlacing. For $B\ge0$, the product
$\lambda_-\lambda_+=-(n-2)\tau^2+n(n-2)(k-\ell)\tau-nk\ell$ is negative. Write $U=-\lambda_->0$ and $V=|\lambda_+|$.

For \(f:[0,\infty)\to\R\), sum over the three pairs
\((x,r)=(A,k),(|B|,1),(C,\ell)\) in the order listed, including any repeated pair. Define
\begin{equation}\label{eq:cross}
 \begin{split}
 \mathcal D_f={}&\sum_{(x,r)}
 \bigl[(m-r)f((m+1)x)+(m+1-r)f(mx)+rf(x)\bigr]\\
 &-m\bigl[f((m+1)U)+f((m+1)V)\bigr]
 -(m+1)\bigl[f(mU)+f(mV)\bigr].
 \end{split}
\end{equation}
For $f(x)=x^p$, expansion of \eqref{eq:coefficient} gives
\begin{equation}\label{eq:cross-powers}
 \mathcal D_f=N_pZ_p-D_pW_p.
\end{equation}
The coefficients in \eqref{eq:cross} need not all be positive. We move each negative term to the other side, producing two vectors $X,Y$ of cubes with equal lengths.

\begin{lem}\label{lem:prefix}
Let $X,Y$ be equal-length nonnegative vectors, and suppose $Y$ is specified by finitely many constant blocks. It suffices to verify the inequality between sums of largest entries at the cumulative lengths of the distinct-value blocks of $Y$ in decreasing order. At each such length, a lower bound for the largest-entry sum of $X$ may be obtained by selecting any subset of $X$ of that cardinality.
\end{lem}
\begin{proof}
Interpolate the sums of largest entries linearly in their cardinality. The interpolated sum for $X$ is concave, while the corresponding sum for $Y$ is affine between consecutive block endpoints. Their difference is concave there, so its minimum is at an endpoint. A selected subset never exceeds the sum of the same number of largest entries.
\end{proof}

We use the notation $Y\prec_wX$ for these largest-entry inequalities. It implies
\begin{equation}\label{eq:tails}
 \sum_{y\in Y}(y-T)_+\le\sum_{x\in X}(x-T)_+\qquad(T\ge0).
\end{equation}
For example, \eqref{eq:tails} follows by expressing each side as the maximum, over integers $j$, of the sum of the $j$ largest entries minus $jT$. Conversely, the largest-entry sums can be recovered by minimizing $jT+\sum_{x\in X}(x-T)_+$ over $T\ge0$.

\begin{prop}\label{prop:cubes}
For every configuration \eqref{eq:triple},
\begin{equation}\label{eq:cubic-tail}
 \mathcal D_{f_T}\ge0\qquad(T\ge0),\qquad f_T(x)=(x^3-T)_+.
\end{equation}
At every interior three-value configuration, $\mathcal D_{f_0}>0$.
\end{prop}

Subsections~\ref{sec:unbalanced} and~\ref{sec:balanced} prove
Proposition~\ref{prop:cubes} for unequal and equal outer multiplicities, respectively.

\subsubsection{The unequal multiplicity case}\label{sec:unbalanced}
Assume $k<m$, put $h=m-k\ge1$, and retain \eqref{eq:ABCD}. In the following table, a magnitude $q$ in a block means that its entry in the vector is $q^3$, repeated the indicated number of times.
\begin{center}\small
\begin{tabular}{c|cc|cc}
 &\multicolumn{2}{c|}{$X$}&\multicolumn{2}{c}{$Y$}\\
index&magnitude&multiplicity&magnitude&multiplicity\\\hline
0&$(m+1)A$&$h$&$(m+1)C$&$h$\\
1&$mA$&$h+1$&$mC$&$h-1$\\
2&$A$&$k$&$(m+1)U$&$m$\\
3&$(m+1)|B|$&$m-1$&$mU$&$m+1$\\
4&$m|B|$&$m$&$(m+1)V$&$m$\\
5&$|B|$&$1$&$mV$&$m+1$\\
6&$C$&$\ell$&&
\end{tabular}
\end{center}
Both lengths are $4m+2h+1$. Directly from \eqref{eq:cross},
$\mathcal D_{f_T}=\sum_{x\in X}(x-T)_+-\sum_{y\in Y}(y-T)_+$.

A prefix ending after a distinct-value block of $Y$ is represented by a subset $J\subseteq\{0,\ldots,5\}$. In all cases,
\begin{equation}\label{eq:ideals}
 1\in J\Rightarrow0\in J,\quad
 3\in J\Rightarrow2\in J,\quad
 5\in J\Rightarrow4\in J.
\end{equation}
When $\tau\le0$, $S<0$ and $U\ge V$, so there are also the implications
\begin{equation}\label{eq:extra-ideals}
 4\in J\Rightarrow2\in J,\qquad
 5\in J\Rightarrow3\in J.
\end{equation}
Tied blocks may be included simultaneously; blocks of zero multiplicity cause no difficulty. There are 26 nonempty subsets satisfying \eqref{eq:ideals}, and 17 satisfying both \eqref{eq:ideals} and \eqref{eq:extra-ideals}.

Introduce parameters
\begin{equation}\label{eq:parameters}
 \begin{gathered}
 k=K+1,\quad h=H+1,\quad m=K+H+2,\quad \ell=K+2H+3,\\
 \tau=\begin{cases}kv,&\epsilon=+1,\\-\ell v,&\epsilon=-1,\end{cases}
 \qquad K,H\ge0,\quad 0\le v\le1.
 \end{gathered}
\end{equation}
For the actual multiplicities, $K,H$ are integers. Polynomial positivity will hold on the full nonnegative quadrant.

For each admissible $J$, Appendix~\ref{app:allocations} specifies an affine integral allocation
\(\alpha_j=a_jK+b_jH+c_j\), \(0\le j\le6\), from the seven blocks
of \(X\). In each row, the coefficients of \(\alpha_j\) and of the
capacity minus \(\alpha_j\) are nonnegative. Summing the seven entries
gives the sum of the multiplicities of the blocks indexed by \(J\)
in \(Y\). These are affine identities in \(K,H\) and can be verified
directly from the tables.

Let $e_j$ be the indicator of $j\in J$, and define
\begin{equation}\label{eq:gammas}
 \begin{aligned}
 \gamma_C&=e_0h(m+1)^3+e_1(h-1)m^3,\\
 \gamma_U&=e_2m(m+1)^3+e_3(m+1)m^3,\\
 \gamma_V&=e_4m(m+1)^3+e_5(m+1)m^3,\\
 Q_0&=A^3\bigl(\alpha_0(m+1)^3+\alpha_1m^3+\alpha_2\bigr)\\
 &\quad+\epsilon B^3\bigl(\alpha_3(m+1)^3+\alpha_4m^3+\alpha_5\bigr)
 +C^3(\alpha_6-\gamma_C).
 \end{aligned}
\end{equation}
Since $\epsilon B^3=|B|^3$, the required subset comparison is precisely
\begin{equation}\label{eq:subset-comparison}
 Q_0\ge\gamma_U U^3+\gamma_VV^3.
\end{equation}

If $\gamma_U+\gamma_V=0$, we certify $Q_0\ge0$ directly. Otherwise put
\begin{equation}\label{eq:radical-cert}
R_+=8Q_0-(\gamma_V-\gamma_U)S(S^2+3D),\qquad T_+=(\gamma_U+\gamma_V)(D+3S^2),\qquad E_+=R_+^2-DT_+^2.
\end{equation}
The inequalities $R_+,E_+\ge0$ in \eqref{eq:radical-cert} imply the comparison with $V^3$ replaced by $\lambda_+^3$. For $\tau\ge0$, interlacing gives $\lambda_+\ge B\ge0$, so this is \eqref{eq:subset-comparison}.

For $\tau\le0$ and $\gamma_V>0$, also put
\begin{equation}\label{eq:negative-cert}
R_-=8Q_0+(\gamma_U+\gamma_V)S(S^2+3D),\qquad T_-=(\gamma_U-\gamma_V)(D+3S^2),\qquad E_-=R_-^2-DT_-^2.
\end{equation}
The subset implications \eqref{eq:extra-ideals} give $\gamma_U\ge\gamma_V$, hence $T_-\ge0$. The inequalities $R_-,E_-\ge0$ prove the other signed comparison, with $-\lambda_+^3$ in place of $V^3$. If $\gamma_U=\gamma_V$, only $R_-\ge0$ is needed. If $\gamma_V=0$, this additional comparison is unnecessary. Together, the two comparisons establish \eqref{eq:subset-comparison} regardless of the sign of $\lambda_+$.

\begin{lem}\label{lem:certificate}
After substitution \eqref{eq:parameters}, all the polynomials required in \eqref{eq:gammas}--\eqref{eq:negative-cert}, with the allocations in Appendix~\ref{app:allocations}, are nonnegative for $K,H\ge0$ and $0\le v\le1$. For the full subset $J=\{0,\ldots,5\}$ the comparison \eqref{eq:subset-comparison} is strict when $v<1$.
\end{lem}
\begin{proof}
We give an exact finite coefficient certificate, using the \href{https://github.com/zhangteng2000/schoenberg-validation/blob/e7435c5521dcd6a10131ba04a6da500b1404800b/certificates/odd_ge_3_allocations.json}{allocation data}. For a polynomial $Q=\sum_{a,b,i}q_{abi}K^aH^bv^i$ and any integer $d\ge\deg_vQ$, the identity
\begin{equation}\label{eq:bernstein}
Q=\sum_{j=0}^d\left(\sum_{a,b}\beta_{abj}K^aH^b\right) \binom djv^j(1-v)^{d-j},\qquad \beta_{abj}=\sum_{i=0}^j\frac{\binom ji}{\binom di}\,q_{abi}
\end{equation}
holds. For every required polynomial, every nonzero $\beta_{abj}$ is a positive rational number. One takes its actual degree in $v$, except for $R_-$ with $\epsilon=-1$ and $J=\{2,4\}$, where one takes $d=4$ instead of $3$. Degree elevation in this one case is essential for the coefficient test in \eqref{eq:bernstein}.

The 26 positive-middle allocations require 50 polynomial identities, with 44984 positive coefficients. The 17 negative-middle allocations require 44 identities, with 34736 positive coefficients. The \href{https://github.com/zhangteng2000/schoenberg-validation/blob/e7435c5521dcd6a10131ba04a6da500b1404800b/certificates/verify_odd_ge_three.py}{reconstruction program} identified in Appendix~\ref{app:odd-exact-certificates} computes the polynomial expressions from the allocation table and checks all Bernstein coefficients by exact arithmetic. In particular, no positivity decision uses a floating-point value or a substitution for $K,H$.

For the full subset, the constant coefficient $\beta_{000}$ is positive in $R_+,E_+$, and also in the required $R_-$ on the negative side. Therefore these polynomials are strictly positive for $v<1$. Since $T_+\ge0$, this gives a strict radical comparison, or a strict $R_-$ comparison when $T_-=0$. This proves the final assertion.
\end{proof}

Lemma~\ref{lem:certificate} verifies every admissible subset comparison \eqref{eq:subset-comparison}. Lemma~\ref{lem:prefix} therefore gives $Y\prec_wX$. Equation~\eqref{eq:tails} proves Proposition~\ref{prop:cubes} when $k<m$, with strict total cubic comparison in the interior.

\subsubsection{Equal outer multiplicities}\label{sec:balanced}
Now $k=\ell=m$. By reflection assume $t\ge0$. Set $\tau=mv$ with $0\le v\le1$ and divide every value in \eqref{eq:ABCD} by $m$. To avoid carrying this common scale, write
\begin{equation}\label{eq:balanced-parameters}
 \begin{gathered}
 A=2m+1+v,\quad C=2m+1-v,\quad B=2mv,\\
 s=(m-1)v,\quad d=m^2v^2+2m+1,\quad
 U=\sqrt d-s,\quad V=\sqrt d+s.
 \end{gathered}
\end{equation}
In \eqref{eq:balanced-parameters}, $0<U\le V$, since $d-s^2=(2m-1)v^2+2m+1>0$.

After cancelling the vanishing coefficients in \eqref{eq:cross}, the positive vector of cubes has magnitude blocks
\[
 (mA,1),\ (mC,1),\ ((m+1)B,m-1),\ (mB,m),\ (A,m),\ (C,m),\ (B,1),
\]
where each pair specifies magnitude and multiplicity. The negative blocks are
$((m+1)V,m)$, $(mV,m+1)$, $((m+1)U,m)$, $(mU,m+1)$. Both vectors have length $4m+2$.

The first negative block in decreasing order is $((m+1)V)^3$, and the last is $(mU)^3$. Only the order of the middle two blocks can change. Put $a_0=m(m+1)^3$ and $b_0=(m+1)m^3$. The following five comparisons therefore suffice, by Lemma~\ref{lem:prefix}:
\begin{equation}\label{eq:balanced-five}
 \begin{array}{c|c|l}
 \text{cardinality}&\text{negative sum}&\text{chosen positive sum }P\\\hline
 m&a_0V^3&m^3A^3+(m-1)(m+1)^3B^3\\
 2m&a_0(V^3+U^3)&m^3(A^3+C^3)+(m-1)((m+1)^3+m^3)B^3\\
 2m+1&(a_0+b_0)V^3&m^3(A^3+C^3)+((m-1)(m+1)^3+m^4)B^3\\
 3m+1&(a_0+b_0)V^3+a_0U^3&(m^3+m)A^3+m^3C^3+((m-1)(m+1)^3+m^4)B^3\\
 4m+2&(a_0+b_0)(V^3+U^3)&(m^3+m)(A^3+C^3)\\
 &&\quad+((m-1)(m+1)^3+m^4+1)B^3
 \end{array}
\end{equation}
Each chosen positive sum uses the cardinality in the first column of \eqref{eq:balanced-five} within its block capacities. The $2m$ line is needed if $(m+1)U\ge mV$; the $2m+1$ line is needed in the opposite ordering. We prove both without imposing either extra condition.

For a line with negative sum $c_VV^3+c_UU^3$, define
\begin{equation}\label{eq:balanced-cert}
 R=P-(c_V-c_U)s(3d+s^2),\qquad
 T=(c_V+c_U)(d+3s^2),\qquad E=R^2-dT^2.
\end{equation}
Since $T\ge0$, the inequalities $R,E\ge0$ imply that line of \eqref{eq:balanced-five}. Set $m=K+2$. The ten polynomials in \eqref{eq:balanced-cert} have representations \eqref{eq:bernstein} with 509 positive coefficients; no $H$ occurs. Their degrees in $v$ are $3$ for $R$ and $6$ for $E$. The numbers of positive coefficients are $29$ for each $R$, and $67,96,67,67,67$ for the five $E$ polynomials in table order. The verifier constructs these ten polynomials from \eqref{eq:balanced-five} and \eqref{eq:balanced-cert}.

In the total-sum line, the constant $j=0$ Bernstein coefficients of both $R$ and $E$ are strictly positive. Thus the total cubic comparison is strict for $v<1$. The five comparisons and Lemma~\ref{lem:prefix} prove $Y\prec_wX$ in the balanced case. This completes the proof of Proposition~\ref{prop:cubes}.

\subsubsection{All exponents and equality}\label{sec:finish}
For $p>3$ put $q=p/3>1$. The elementary identity
\begin{equation}\label{eq:integral}
 r^p=q(q-1)\int_0^\infty(r^3-T)_+T^{q-2}\dd T\qquad(r\ge0)
\end{equation}
converts \eqref{eq:cubic-tail}, by linearity and \eqref{eq:cross-powers}, into
$N_pZ_p-D_pW_p\ge0$ on every configuration \eqref{eq:triple}. At $p=3$, use the total-sum case $T=0$. If the configuration is interior, its total cubic gap is positive. Continuity then makes the tail gap positive on an interval starting at $T=0$, so \eqref{eq:integral} makes the $p$-power gap strictly positive for every $p>3$ as well.

By Lemma~\ref{lem:reduction}, a real global maximizer must occur in \eqref{eq:triple}. Proposition~\ref{prop:cubes} and \eqref{eq:integral} for \(p>3\), or \eqref{eq:cubic-tail} at \(T=0\) for \(p=3\), bound its ratio by $B_{n,p}$. Direct differentiation of \eqref{eq:equality} gives the critical values $(m+1)a$ with multiplicity $m-1$, $-ma$ with multiplicity $m$, and $a$ once. Its ratio is exactly \eqref{eq:coefficient}, proving sharpness.

It remains to exclude other endpoints. A centered two-value configuration can be written as $(z-(n-j))^j(z+j)^{n-j}$ with $1\le j\le m$. At exponent three its ratio is $1-h_3(r)/n$, where $r=j/(n-j)$ and
\[
 h_3(r)=\frac{(1+r)(3-3r+2r^2)}{1+r^2},\qquad
 h_3'(r)=\frac{2r(r-1)(r^2+r+4)}{(1+r^2)^2}<0\quad(0<r<1).
\]
Thus only $j=m$ attains $B_{n,3}$. Every other endpoint has a strictly positive total cubic gap. Proposition~\ref{prop:cubes} and \eqref{eq:integral} again imply strictness at every $p>3$. This proves the real equality classification for all $p\ge3$.

Apply Lemma~\ref{lem:pereira-projection} in each direction and then \eqref{eq:angular}. The same coefficient holds for arbitrary complex zeros. If there are at least three distinct complex zeros, almost every projected configuration has at least three distinct values, and its real inequality is strict. Integrating yields strictness in \eqref{eq:oddsharp}. Two distinct centered complex values lie on a common real line through zero. The real classification then gives exactly \eqref{eq:equality}. Theorem~\ref{thm:oddcomplete} is proved.

\section{Allocation tables and exact-arithmetic certificates}\label{app:odd-exact-certificates}
\subsection{The block allocations}\label{app:allocations}
In the tables in Subsections~\ref{app:nonnegative-allocations} and~\ref{app:nonpositive-allocations}, $J$ is encoded by the integer $\sum_{j\in J}2^j$. The seven entries are $\alpha_0,\ldots,\alpha_6$ in the order of the positive blocks in Section~\ref{sec:unbalanced}. They are affine polynomials in the parameters from \eqref{eq:parameters}. All 26 and 17 admissible nonempty subsets are included. The electronic allocation file gives the same data as affine expressions in \(K,H\).

\subsubsection{Nonnegative middle value}\label{app:nonnegative-allocations}
\begingroup\footnotesize
\renewcommand{\arraystretch}{0.92}
\setlength{\LTpre}{4pt}\setlength{\LTpost}{4pt}
\setlength{\tabcolsep}{5pt}
\begin{longtable}{rlllllll}
\toprule
$J$&$\alpha_0$&$\alpha_1$&$\alpha_2$&$\alpha_3$&$\alpha_4$&$\alpha_5$&$\alpha_6$\\\midrule
\endhead
1&$H+1$&$0$&$0$&$0$&$0$&$0$&$0$\\
3&$H+1$&$H$&$0$&$0$&$0$&$0$&$0$\\
4&$H+1$&$1$&$K$&$0$&$0$&$0$&$0$\\
5&$H+1$&$H+2$&$K$&$0$&$0$&$0$&$0$\\
7&$H+1$&$H+2$&$K$&$H$&$0$&$0$&$0$\\
12&$H+1$&$H+1$&$0$&$0$&$K$&$0$&$K+3$\\
13&$H+1$&$H+1$&$0$&$H$&$K$&$1$&$K+3$\\
15&$H+1$&$H+1$&$0$&$0$&$K$&$1$&$K+2H+3$\\
16&$H$&$2$&$0$&$K$&$0$&$0$&$0$\\
17&$H+1$&$1$&$0$&$K+H+1$&$0$&$0$&$0$\\
19&$H+1$&$H+1$&$0$&$K+H$&$1$&$0$&$0$\\
20&$H+1$&$0$&$0$&$K$&$K$&$0$&$H+3$\\
21&$H+1$&$1$&$1$&$K$&$K$&$0$&$2H+2$\\
23&$H+1$&$H+1$&$1$&$K+H+1$&$K+H$&$0$&$1$\\
28&$H+1$&$H+2$&$1$&$K+1$&$K$&$1$&$K+H+1$\\
29&$H+1$&$H+2$&$K+1$&$K+H$&$K+H$&$1$&$3$\\
31&$H+1$&$H+2$&$1$&$K+H+1$&$K+2$&$1$&$K+2H$\\
48&$1$&$H+2$&$0$&$K+1$&$K+H+1$&$0$&$0$\\
49&$H+1$&$H+2$&$0$&$K+H+1$&$K+2$&$0$&$0$\\
51&$H+1$&$H+2$&$0$&$K+H+1$&$K+H+2$&$0$&$0$\\
52&$H+1$&$2$&$K$&$K+H+1$&$K+H+2$&$0$&$1$\\
53&$H+1$&$H+2$&$K+1$&$K+H+1$&$K+H+2$&$0$&$1$\\
55&$H+1$&$H+2$&$0$&$K+H+1$&$K+H+2$&$0$&$K+H+2$\\
60&$H+1$&$H+2$&$K+1$&$K+H+1$&$K+H+2$&$0$&$K+3$\\
61&$H+1$&$H+2$&$K+1$&$K+1$&$K+H+2$&$1$&$K+2H+3$\\
63&$H+1$&$H+2$&$K+1$&$K+H+1$&$K+H+2$&$1$&$K+2H+3$\\
\bottomrule
\end{longtable}
\endgroup

\subsubsection{Nonpositive middle value}\label{app:nonpositive-allocations}
\begingroup\footnotesize
\renewcommand{\arraystretch}{0.92}
\setlength{\LTpre}{4pt}\setlength{\LTpost}{4pt}
\setlength{\tabcolsep}{5pt}
\begin{longtable}{rlllllll}
\toprule
$J$&$\alpha_0$&$\alpha_1$&$\alpha_2$&$\alpha_3$&$\alpha_4$&$\alpha_5$&$\alpha_6$\\\midrule
\endhead
1&$H+1$&$0$&$0$&$0$&$0$&$0$&$0$\\
3&$H+1$&$H$&$0$&$0$&$0$&$0$&$0$\\
4&$H+1$&$0$&$0$&$K+1$&$0$&$0$&$0$\\
5&$H+1$&$H+1$&$0$&$K+1$&$0$&$0$&$0$\\
7&$H+1$&$H+1$&$0$&$K+H+1$&$0$&$0$&$0$\\
12&$H+1$&$H+2$&$0$&$K+1$&$K+1$&$0$&$0$\\
13&$H+1$&$H+2$&$0$&$K+H+1$&$K+2$&$0$&$0$\\
15&$H+1$&$H+2$&$0$&$K+H+1$&$K+H+2$&$0$&$0$\\
20&$H+1$&$0$&$1$&$K+H+1$&$K$&$1$&$0$\\
21&$H+1$&$1$&$K+1$&$K+H+1$&$0$&$1$&$H$\\
23&$H+1$&$H+2$&$0$&$K+H+1$&$K+H+1$&$0$&$0$\\
28&$H+1$&$H+2$&$K$&$K+H+1$&$K+2$&$1$&$0$\\
29&$H+1$&$H+2$&$K+1$&$K+H+1$&$K+2$&$1$&$H$\\
31&$H+1$&$H+2$&$K+1$&$K+H+1$&$K+H+2$&$1$&$H$\\
60&$H+1$&$H+2$&$K+1$&$K+H+1$&$K+H+2$&$0$&$K+3$\\
61&$H+1$&$H+2$&$K+1$&$K+H+1$&$K+H+2$&$1$&$K+H+3$\\
63&$H+1$&$H+2$&$K+1$&$K+H+1$&$K+H+2$&$1$&$K+2H+3$\\
\bottomrule
\end{longtable}
\endgroup

\subsection{Exact reconstruction of the cubic-tail certificates}
The proof of Proposition~\ref{prop:cubes} reduces the unequal-multiplicity case
and the equal-multiplicity case to polynomial inequalities
\eqref{eq:radical-cert}, \eqref{eq:negative-cert}, and
\eqref{eq:balanced-cert}.  After the parameter changes
\eqref{eq:parameters} and \eqref{eq:balanced-parameters}, the checker
\href{https://github.com/zhangteng2000/schoenberg-validation/blob/e7435c5521dcd6a10131ba04a6da500b1404800b/certificates/verify_odd_ge_three.py}{\nolinkurl{certificates/verify_odd_ge_three.py}} reconstructs each polynomial and
tests its coefficients in the Bernstein basis, using only exact rational
arithmetic. The input allocations are contained in
\href{https://github.com/zhangteng2000/schoenberg-validation/blob/e7435c5521dcd6a10131ba04a6da500b1404800b/certificates/odd_ge_3_allocations.json}{\nolinkurl{certificates/odd_ge_3_allocations.json}} and reproduced in
Subsection~\ref{app:allocations}.

The verifier reconstructs 104 polynomials, with 50 for a nonnegative
middle zero, 44 for a nonpositive middle zero, and 10 for equal outer
multiplicities. For a polynomial \(P(K,H,v)\) of degree at most \(d\)
in \(v\), write \(P=\sum_{i=0}^d p_i(K,H)v^i\). The Bernstein coefficients of
order $d$ are
\begin{equation}\label{eq:odd-check-bernstein}
\beta_j(K,H)=\sum_{i=0}^j\frac{\binom ji}{\binom di}p_i(K,H),\qquad 0\le j\le d.
\end{equation}
If each \(\beta_j(K,H)\) in \eqref{eq:odd-check-bernstein} has
nonnegative monomial coefficients, then \(P(K,H,v)\ge0\) for
$K,H\ge0$ and $0\le v\le1$. The verifier checks the
algebraic reconstruction of each defining polynomial, including
the degree elevation specified in the proof of
Lemma~\ref{lem:certificate}. The assertion covers all real
nonnegative $K,H$, not merely integer multiplicities or a finite grid.

\subsection{Exact checks supporting the logarithmic comparison}
The proof of Theorem~\ref{odd23:thm:global-log} uses the analytic
inequalities in Appendix~\ref{app:odd23-details}; certain finite
algebraic sign checks are implemented in the accompanying files.
The following table identifies the checkers used in the proof.
\begin{center}\small
\begin{tabular}{>{\raggedright\arraybackslash}p{.50\linewidth}p{.42\linewidth}}
\toprule
Certificate & Role in the analytic proof\\
\midrule
\href{https://github.com/zhangteng2000/schoenberg-validation/blob/e7435c5521dcd6a10131ba04a6da500b1404800b/certificates/turn2_two_value_T_certificate.py}{\nolinkurl{turn2_two_value_T_certificate.py}} & Two-value rational sign expansion\\
\href{https://github.com/zhangteng2000/schoenberg-validation/blob/e7435c5521dcd6a10131ba04a6da500b1404800b/certificates/num_balanced_threevalue_certificate.py}{\nolinkurl{num_balanced_threevalue_certificate.py}} & Balanced three-value polynomial signs\\
\href{https://github.com/zhangteng2000/schoenberg-validation/blob/e7435c5521dcd6a10131ba04a6da500b1404800b/certificates/turn2_restricted_gradient_exact_check.py}{\nolinkurl{turn2_restricted_gradient_exact_check.py}} & Restricted-gradient polynomial coefficients\\
\href{https://github.com/zhangteng2000/schoenberg-validation/blob/e7435c5521dcd6a10131ba04a6da500b1404800b/certificates/turn2_three_caseIII_T_certificate.py}{\nolinkurl{turn2_three_caseIII_T_certificate.py}} & Positive-middle three-value rational inequality\\
\href{https://github.com/zhangteng2000/schoenberg-validation/blob/e7435c5521dcd6a10131ba04a6da500b1404800b/certificates/turn2_minority_positive_center_certificate.py}{\nolinkurl{turn2_minority_positive_center_certificate.py}} & Minority-averaging scalar bounds\\
\href{https://github.com/zhangteng2000/schoenberg-validation/blob/e7435c5521dcd6a10131ba04a6da500b1404800b/certificates/turn2_negative_center_exact_check.py}{\nolinkurl{turn2_negative_center_exact_check.py}} & Negative-middle residue identities\\
\href{https://github.com/zhangteng2000/schoenberg-validation/blob/e7435c5521dcd6a10131ba04a6da500b1404800b/certificates/turn2_log_compensation_certificate.py}{\nolinkurl{turn2_log_compensation_certificate.py}} & Logarithmic compensation inequalities\\
\bottomrule
\end{tabular}
\end{center}

The checkers and allocation data are available in the
\href{https://github.com/zhangteng2000/schoenberg-validation/tree/e7435c5521dcd6a10131ba04a6da500b1404800b/certificates}{\nolinkurl{certificates/}} directory of the archived source revision~\cite{Zha26c}.
The polynomial sign checks operate symbolically, and the analytic
arguments establishing compactness, stationarity, interlacing, and
boundary behavior remain in Appendix~\ref{app:odd23-details}.
A successful run of a script verifies only the specific algebraic
assertions encoded by that script; it does not independently validate
all analytic reductions in the text.

\end{document}